\documentclass[12pt]{book}
 \usepackage{amsmath,amssymb,amsbsy,amstext,amscd,amsfonts}
 \usepackage{exscale,calrsfs,times,multind,theorem}
 \usepackage{color}
 \usepackage{article}
 \usepackage{pictex,epsfig}
 \makeatletter
 \renewcommand{\@evenhead}{\textsc{\thepage}\hfill {\em\scriptsize\rightmark}
        \hfill \today \,}
\renewcommand{\@oddhead}{ {\em\scriptsize\today \hfill Martingale Hardy spaces and summability of the one dimensional Vilenkin-Fourier series} \hfill \textsc{\thepage}\,}
\renewcommand{\@evenfoot}{\scriptsize\em
\hfill  \hfill,}
\renewcommand{\@oddfoot}{\scriptsize\em
\hfill  G.Tephnadze \hfill,}

 \input article.def

\begin{document}
\def \R {{\mathbb {R}}}

\def \R {{\mathbb {R}}}

\newpage

\thispagestyle{empty}

\begin{center}
{\huge \textbf{Ph.D. Thesis}} \vspace{2cm}

{\Large \textbf{Martingale Hardy spaces}} \vspace{0.5cm}

{\Large \textbf{\ and summability of the one dimensional}}

\vspace{0.5cm}

{\Large \textbf{\ Vilenkin-Fourier series}}
\end{center}

\vspace{3cm}

\begin{center}
\textbf{George Tephnadze}
\end{center}

\vspace{2.0cm}

\begin{center}
{\ Department of Engineering Sciences and Mathematics\\[0pt]
Lule\aa\ University of Technology\\[0pt]
SE-971~87~Lule\aa ,~Sweden\\[0pt]
and \\[0pt]
Department of Mathematics\\[0pt]
Faculty of Exact and Natural Sciences \\[0pt]
Ivane Javakhishvili Tbilisi State University\\[0pt]
Chavchavadze str. 1, Tbilisi 0128\\[0pt]
{giorgitephnadze@gmail.com}}\\[0pt]
\vspace{3cm} \textbf{8 October, 2015}
\end{center}

\newpage

\thispagestyle{empty}

\begin{center}
{\scriptsize George Tephnadze}
\end{center}

\newpage

\tableofcontents

\newpage

\thispagestyle{empty}


\newpage

\begin{center}
{\Large \textbf{Abstract}}
\end{center}

\vspace{0.5cm}

The classical theory of Fourier series deals with decomposition of a
function into sinusoidal waves. Unlike these continuous waves the Vilenkin
(Walsh) functions are rectangular waves. Such waves have already been used
frequently in the theory of signal transmission, multiplexing, filtering,
image enhancement, codic theory, digital signal processing and pattern
recognition. The development of the theory of Vilenkin-Fourier series has
been strongly influenced by the classical theory of trigonometric series.
Because of this it is inevitable to compare results of Vilenkin series to
those on trigonometric series. There are many similarities between these
theories, but there exist differences also. Much of these can be explained
by modern abstract harmonic analysis, which studies orthonormal systems from the
point of view of the structure of a topological group.

In this PhD thesis we discuss, develop and apply this fascinating theory
connected to modern harmonic analysis. In particular we make new estimations
of Vilenkin-Fourier coefficients and prove some new results concerning
boundedness of maximal operators of partial sums. Moreover, we derive
necessary and sufficient conditions for the modulus of continuity so that
norm convergence of the partial sums is valid and develop new methods to
prove Hardy type inequalities for the partial sums with respect to the
Vilenkin systems. We also do the similar investigation for the Fej\'er means.
Furthermore, we investigate some N\"orlund means but only in the case  when their coefficients are monotone. Some well-know examples of N\"orlund means are
Fej\'er means, Ces\`aro means and N\"orlund logarithmic means. In addition, we
consider Riesz logarithmic means, which are not example of N\"orlund
means. It is also proved that these results are the best possible in a
special sense. As applications both some well-known and new results are
pointed out.

This PhD is written as a monograph consisting of four Chapters:
Preliminaries, Fourier coefficients and partial sums of Vilenkin-Fourier
series on martingale Hardy spaces, Vilenkin-Fej\'er means on
martingale Hardy spaces, Vilenkin-N\"orlund means on martingale Hardy
spaces. It is based on 15 papers with the candidate as author or coauthor,
but also some new results are presented for the first time.

In Chapter 1 we first present some definitions and notations, which are
crucial for our further investigations. After that we also define some
summabilitity methods and remind about some classical facts and results. We investigate some well-known results and prove new estimates for the kernels of these summabilitity methods, which are very important to prove our main results. Moreover, we define martingale Hardy spaces and construct martingales, which help us to prove sharpness of our main results in the later chapters.

Chapter 2 is devoted to present and prove some new and known results about Vilenkin-Fourier coefficients and partial sums of martingales in  Hardy spaces.
First, we show that Fourier coefficients of martingales are not uniformly bounded when $0<p<1$. By applying these results we prove some known Hardy and Paley type inequalities with a new method. After that we investigate partial sums with respect to Vilenkin systems and prove boundedness of maximal operators of partial sums. Moreover, we find necessary and sufficient conditions for the modulus of continuity for which norm convergence of partial sums hold and we present a new proof of a Hardy type inequality for it.

In Chapter 3 we investigate some analogous problems concerning the partial sums of  Fej\'er means. First we consider some weighted maximal operators of Fej\'er means and prove some boundedness results for them. After that we apply these results to find necessary and sufficient conditions for the modulus of continuity for which norm convergence of Fej\'er means hold. Finally, we prove some new Hardy type inequalities for Fej\'er means, which is a main part of this PhD thesis. We also prove sharpness of all our main results in this Chapter.

In Chapter 4 we consider boundedness of maximal operators of N\"orlund means. After that we prove some strong convergence theorems for these summablility methods.
Since Fej\'er means, Ces\`aro means are well-know examples of N\"orlund means some well-known and new results are pointed out. We also investigate Riesz and N\"orlund logarithmic means simultaneously at the end of this chapter.

\newpage

\begin{center}
{\Large \textbf{Preface}}
\end{center}

\vspace{0.10cm}

This PhD thesis is written as a monograph based on the following
publications: \vspace{0.10cm} \newline

[{1] G. Tephnadze, }The maximal operators of logarithmic means of
one-dimensional Vilenkin-Fourier series, Acta Math. Acad. Paedagog. Nyh\'azi.
27 (2011), {no. 2, }245-256.

[{2] G. Tephnadze, }Fej\'er means of Vilenkin-Fourier series. Studia Sci.
Math. Hungar. 49 (2012), no. 1, 79-90.

[{3] G. Tephnadze, }A note on the Fourier coefficients and partial sums of
Vilenkin-Fourier series. Acta Math. Acad. Paedagog. Nyh\'azi. 28 (2012), no.
2, 167-176.

[{4] G. Tephnadze, }On the maximal operators of Vilenkin-Fej\'er means.
Turkish J. Math. 37 (2013), no. 2, 308-318.

[{5] G. Tephnadze, }On the maximal operators of Vilenkin-Fej\'er means on
Hardy spaces. Math. Inequal. Appl. 16 (2013), no. 1, 301-312.

[{6] G. Tephnadze, On the Vilenkin-Fourier coefficients, Georgian Math. J.
20 (2013), no. 1, 169-177.}

[{7] G. Tephnadze, }On the partial sums of Vilenkin-Fourier series. J.
Contemp. Math. Anal. 49 (2014), no. 1, 23-32.

[{8] G. Tephnadze, }A note on the norm convergence by Vilenkin-Fej\'er means.
Georgian Math. J. 21 (2014), no. 4, 511--517.

[{9] I. Blahota and G. Tephnadze, }Strong convergence theorem for Vilenkin-Fej\'e%
r means. Publ. Math. Debrecen 85 (2014), no. 1-2, 181-196.

[10] {L. E. Persson and G. Tephnadze, A note on Vilenkin-Fej\'er means, on the
martingale Hardy space }$H_{p},$ Bulletin of TICMI 18 (2014), no. 1, 55-64.

[{11] L. E. Persson, G. Tephnadze and P. Wall, Maximal operators of Vilenkin-N\"o%
rlund means, }J. Fourier Anal. Appl. 21 (2015), no. 1, 76-94.

[{12] L. E. Persson, G. Tephnadze and P. Wall,} Some new $\left(
H_{p},L_{p}\right) $ type inequalities of maximal operators of Vilenkin-N\"o%
rlund means with non-decreasing coefficients, J. Math. Inequal. (to appear).

[13] L. E. Persson and G. Tephnadze,{\textit{\ }}A sharp boundedness result
concerning some maximal operators of Vilenkin-Fej\'er means, Mediterr. J.
Math. (to appear).

[14] I. Blahota, L. E. Persson and G. Tephnadze, On the N\"orlund means of
Vilenkin-Fourier series, Czechoslovak Math. J. (to appear).

[{15] I. Blahota and G. Tephnadze, }\textit{\ }A note on maximal operators of
Vilenkin-N\"orlund means, Acta Math. Acad. Paedagog. Nyh\'azi. (to appear).

Remark: Also some new results which can not be found in these papers appear
in this PhD thesis for the first time.

\newpage

\begin{center}
{\Large \textbf{Acknowledgement}}
\end{center}

\vspace{1.0cm}

It is a pleasure to express my warmest thanks to my supervisors Professors
Lars-Erik Persson, Ushangi Goginava and Peter Wall for the attention to my
work, their valuable remarks and suggestions and for their constant support
and help.

I am also grateful to Doctors Yulia Koroleva, John Fabricius and Professors
Niklas Grip, Maria Alessandra Ragusa and Natasha Samko for helping me with
several practical things. Moreover, I appreciate very much the warm and
friendly atmosphere at the Department of Mathematics at Lule{\aa }
University of Technology, which helps me to do research more effectively.

I thank my Georgian colleagues from Tbilisi State University for their
interest in my research and for many fruitful discussions. In particular, I thank Professor Amiran Gogatishvili, for his support and his valuable remarks. Moreover, I am very grateful to Georgian Mathematical Union and its President Professor Roland Duduchava for all support they have given me.

I also want to pronounce that the agreement about scientific collaboration
and PhD education between Tbilisi State University and Lule{\aa } University of Technology has been very important. Especially, I express my deepest gratitude to Professor Ramaz Bochorishvili at Tbilisi State University and Professor Elisabet Kassfeldt at Lule{\aa } University of Technology for their creative and hard work to realize this important agreement.

I thank Swedish Institute for economic support, which was very important to
be able to travel and be accommodated in Sweden during this PhD program.
Moreover, I thank Shota Rustaveli National Science Foundation for financial support.

Finally, I thank my family for their love, understanding, patience and long lasting constant support.

\newpage

\thispagestyle{empty}

\begin{center}
{\scriptsize George Tephnadze}
\end{center}

\newpage

\section{\textbf{Preliminaries}\protect\bigskip}

\vspace{0.5cm}

\subsection{Basic notations}

Denote by $%
\mathbb{N}
_{+}$ the set of the positive integers, $%
\mathbb{N}
:=%
\mathbb{N}
_{+}\cup \{0\}.$ Let $m:=(m_{0,}$ $m_{1},\ldots )$ be a sequence of positive
integers not less than 2. Denote by
\begin{equation*}
Z_{m_{k}}:=\{0,1,\ldots ,m_{k}-1\}
\end{equation*}%
the additive group of integers modulo $m_{k}$.

Define the group $G_{m}$ as the complete direct product of the groups $%
Z_{m_{i}}$ with the product of the discrete topologies of $Z_{m_{i}}$.

The direct product $\mu $ of the measures
\begin{equation*}
\mu _{k}\left( \{j\}\right) :=1/m_{k}\text{ \ \ \ }(j\in Z_{m_{k}})
\end{equation*}%
is the Haar measure on $G_{m_{k\text{ }}}$with $\mu \left( G_{m}\right) =1.$

In this paper we discuss bounded Vilenkin groups,\textbf{\ }i.e. the case
when $\sup_{n}m_{n}<\infty .$

The elements of $G_{m}$ are represented by sequences
\begin{equation*}
x:=\left( x_{0},x_{1},\ldots ,x_{j},\ldots \right) \ \left( x_{j}\in
Z_{m_{j}}\right) .
\end{equation*}

It is easy to give a base for the neighborhoods of $G_{m}:$

\begin{eqnarray*}
\text{ }I_{0}\left( x\right) &:&=G_{m},\text{ \ } \\
I_{n}(x) &:&=\{y\in G_{m}\mid y_{0}=x_{0},\ldots
,y_{n-1}=x_{n-1}\}\,\,\left( x\in G_{m},\text{ }n\in
\mathbb{N}
\right) .
\end{eqnarray*}%
Let

\begin{equation*}
e_{n}:=\left( 0,\ldots ,0,x_{n}=1,0,\ldots \right) \in G_{m}\qquad \left(
n\in
\mathbb{N}
\right) .
\end{equation*}

If we define $I_{n}:=I_{n}\left( 0\right) ,$\ for \ $n\in \mathbb{N}$ and $\
\overline{I_{n}}:=G_{m}$ $\backslash $ $I_{n},$ then%
\begin{equation}
\overline{I_{N}}=\overset{N-1}{\underset{s=0}{\bigcup }}I_{s}\backslash
I_{s+1}=\left( \overset{N-2}{\underset{k=0}{\bigcup }}\overset{N-1}{\underset%
{l=k+1}{\bigcup }}I_{N}^{k,l}\right) \bigcup \left( \underset{k=1}{%
\bigcup\limits^{N-1}}I_{N}^{k,N}\right) ,  \label{1.1}
\end{equation}%
where
\begin{equation*}
I_{N}^{k,l}:=\left\{
\begin{array}{l}
\text{ }I_{N}(0,\ldots ,0,x_{k}\neq 0,0,...,0,x_{l}\neq 0,x_{l+1\text{ }%
},\ldots ,x_{N-1\text{ }},\ldots ), \\
\text{for }k<l<N, \\
\text{ }I_{N}(0,\ldots ,0,x_{k}\neq 0,x_{k+1}=0,\ldots ,x_{N-1\text{ }}=0,%
\text{ }x_{N\text{ }},\ldots ), \\
\text{for }l=N.%
\end{array}%
\text{ }\right.
\end{equation*}

The norm (or quasi-norm when $0<p<1$) of the space $L_{p}(G_{m})$ $\left(
0<p<\infty \right) $ is defined by
\begin{equation*}
\left\Vert f\right\Vert _{p}:=\left( \int_{G_{m}}\left\vert f\right\vert
^{p}d\mu \right) ^{1/p}.
\end{equation*}%
\qquad

The space $weak-L_{p}\left( G_{m}\right) $ consists of all measurable
functions $f,$ for which
\begin{equation*}
\left\Vert f\right\Vert _{weak-L_{p}}:=\underset{\lambda >0}{\sup }\lambda
\mu \left( f>\lambda \right) ^{1/p}<+\infty .
\end{equation*}

The norm of the space of continuous functions $C(G_{m})$ \ is defined by
\qquad \qquad \thinspace\
\begin{equation*}
\left\Vert f\right\Vert _{C}:=\underset{x\in G_{m}}{\sup }\left\vert
f(x)\right\vert <c<\infty .
\end{equation*}

The best approximation of $f\in L_{p}(G_{m})$ ($1\leq p\in \infty $) is
defined as%
\begin{equation*}
E_{n}\left( f,L_{p}\right) :=\inf_{\psi \in \emph{P}_{n}}\left\Vert f-\psi
\right\Vert _{p},
\end{equation*}%
where $\emph{P}_{n}$ is set of all Vilenkin polynomials of order less than $%
n\in \mathbb{N}$.

The modulus of continuity of $f\in L_{p}\left( G_{m}\right) $ and $\ f\in
C\left( G_{m}\right) $ are defined by

\begin{equation*}
\omega _{p}\left( \frac{1}{M_{n}},f\right) :=\sup\limits_{h\in
I_{n}}\left\Vert f\left( \cdot -h\right) -f\left( \cdot \right) \right\Vert
_{p}
\end{equation*}

and
\begin{equation*}
\omega _{C}\left( \frac{1}{M_{n}},f\right) :=\underset{h\in I_{n}}{\sup }%
\left\Vert f\left( \cdot -h\right) -f\left( \cdot \right) \right\Vert _{C},
\end{equation*}%
respectively.

If we define the so-called generalized number system based on $m$ in the
following way :
\begin{equation*}
M_{0}:=1,\ M_{k+1}:=m_{k}M_{k}\,\,\,\ \ (k\in
\mathbb{N}
),
\end{equation*}%
then every $n\in
\mathbb{N}
$ can be uniquely expressed as
\begin{equation*}
n=\sum_{j=0}^{\infty }n_{j}M_{j},
\end{equation*}%
where $n_{j}\in Z_{m_{j}}$ $(j\in
\mathbb{N}
_{+})$ and only a finite number of $n_{j}^{^{\prime }}$s differ from zero.

\bigskip

Let
\begin{equation*}
\left\langle n\right\rangle :=\min \{j\in \mathbb{N}:n_{j}\neq 0\}\text{ \
and \ }\left\vert n\right\vert :=\max \{j\in \mathbb{N}:n_{j}\neq 0\},
\end{equation*}%
that is $M_{\left\vert n\right\vert }\leq n\leq M_{\left\vert n\right\vert
+1}.$ Set
\begin{equation*}
d\left( n\right) :=\left\vert n\right\vert -\left\langle n\right\rangle ,%
\text{ \ for \ all \ \ }n\in \mathbb{N}.
\end{equation*}

For the natural number $n=\sum_{j=1}^{\infty }n_{j}M_{j},$ we define
functions $v$ and $v^{\ast }$ by (for details see Lukomskii \cite{luko})

\begin{equation*}
v\left( n\right) :=\sum_{j=1}^{\infty }\left\vert \delta _{j+1}-\delta
_{j}\right\vert +\delta _{0},\text{ \ \ }v^{\ast }\left( n\right)
:=\sum_{j=1}^{\infty }\delta _{j}^{\ast },
\end{equation*}%
where
\begin{equation*}
\delta _{j}=sign\left( n_{j}\right) =sign\left( \ominus n_{j}\right) \text{
\ and \ \ \ }\delta _{j}^{\ast }=\left\vert \ominus n_{j}-1\right\vert
\delta _{j}
\end{equation*}%
and $\ominus $\ is the inverse operation for
\begin{equation*}
a_{k}\oplus b_{k}:=(a_{k}+b_{k})\text{mod}m_{k}.
\end{equation*}

Next, we introduce on $G_{m}$ an orthonormal systems, which are called
Vilenkin systems.

At first, we define the complex-valued function $r_{k}\left( x\right)
:G_{m}\rightarrow
\mathbb{C}
,$ the generalized Rademacher functions, by%
\begin{equation*}
r_{k}\left( x\right) :=\exp \left( 2\pi ix_{k}/m_{k}\right) ,\text{ }\left(
i^{2}=-1,x\in G_{m},\text{ }k\in
\mathbb{N}
\right) .
\end{equation*}

Now, define Vilenkin systems$\,\,\,\psi :=(\psi _{n}:n\in
\mathbb{N}
)$ on $G_{m}$ as:
\begin{equation*}
\psi _{n}(x):=\prod\limits_{k=0}^{\infty }r_{k}^{n_{k}}\left( x\right)
,\,\,\ \ \,\left( n\in
\mathbb{N}
\right) .
\end{equation*}

The Vilenkin systems are orthonormal and complete in $L_{2}\left( G_{m}\right)
$ (see e.g. Vilenkin\textit{\ }\cite{Vi}).

Specifically, we call this system the Walsh-Paley system when $m\equiv 2.$

Next, we introduce some analogues of the usual definitions in
Fourier-analysis. If $f\in L_{1}\left( G_{m}\right) $ we can define the
Fourier coefficients, the partial sums of the Fourier series, the Dirichlet
kernels with respect to Vilenkin systems in the usual manner:

\begin{equation*}
\widehat{f}\left( n\right) :=\int_{G_{m}}f\overline{\psi }_{n}d\mu ,\,\ \ \
\ \ \ \ \ \ \ \,\left( n\in
\mathbb{N}
\right) ,
\end{equation*}%
\begin{equation*}
S_{n}f:=\sum_{k=0}^{n-1}\widehat{f}\left( k\right) \psi _{k},\text{ \ \ \ }%
\left( n\in
\mathbb{N}
_{+}\right) ,
\end{equation*}%
\begin{equation*}
D_{n}:=\sum_{k=0}^{n-1}\psi _{k\text{ }},\ \ \ \ \ \ \ \ \ \ \ \left( n\in
\mathbb{N}
_{+}\right) ,
\end{equation*}%
respectively.

The $n$-th Lebesgue constant is defined in the following way:

\begin{equation*}
L_{n}:=\left\Vert D_{n}\right\Vert _{1}.
\end{equation*}

\subsection{Definition and examples of N\"orlund means and its maximal
operators}

Let $\{q_{k}:k\in \mathbb{N}\}$ be a sequence of nonnegative numbers. The $n$%
-th N\"orlund means for the Fourier series of $f$ \ is defined by

\begin{equation}
t_{n}f:=\frac{1}{Q_{n}}\overset{n}{\underset{k=1}{\sum }}q_{n-k}S_{k}f,
\label{1.2}
\end{equation}%
where \
\begin{equation*}
Q_{n}:=\sum_{k=0}^{n-1}q_{k}.
\end{equation*}

We always assume that $q_{0}>0$ and $\ $%
\begin{equation*}
\lim_{n\rightarrow \infty }Q_{n}=\infty .\
\end{equation*}

In this case it is well-known that the summability method generated by $%
\{q_{k}:k\geq 0\}$ is regular if and only if
\begin{equation}
\underset{n\rightarrow \infty }{\lim }\frac{q_{n-1}}{Q_{n}}=0.  \label{112}
\end{equation}

Concerning this fact and related basic results, we refer to \cite{moo}.

The next remark is due to Persson,
Tephnadze and Wall \cite{ptw2}:

\begin{remark}
\label{theorem1n}a) Let the sequence $\{q_{k}:k\in \mathbb{N}\}$ be
non-increasing. Then the summability method generated by $\{q_{k}:k\in
\mathbb{N}\}$ is regular.

b) Let the sequence $\{q_{k}:k\in \mathbb{N}\}$ be non-decreasing. Then the
summability method generated by $\{q_{k}:k\in \mathbb{N}\}$ is not always
regular.
\end{remark}

{\bf Proof}:
Let the sequence $\{q_{k}:k\in \mathbb{N}\}$ be non-increasing. Then
\begin{equation*}
\frac{q_{n-1}}{Q_{n}}\leq \frac{q_{n-1}}{nq_{n-1}}=\frac{1}{n}\rightarrow 0,%
\text{ when }n\rightarrow \infty .
\end{equation*}

According to (\ref{112}) we conclude that in this case the summability
method is regular.

Now, we prove part b) of theorem and construct N\"orlund mean, with
non-decreasing coefficients $\{q_{k}:k\in \mathbb{N}\},$ which is not
regular.

Let $\{q_{k}=2^{k}:k\in \mathbb{N}\}.$ Then
\begin{equation*}
Q_{n}=\overset{n-1}{\underset{k=0}{\sum }}2^{k}=2^{n}-1\leq 2^{n}
\end{equation*}%
and%
\begin{equation*}
\frac{q_{n-1}}{Q_{n}}=\frac{2^{n-1}}{2^{n}-1}\geq \frac{2^{n-1}}{2^{n}}=%
\frac{1}{2}\nrightarrow 0,\text{ when }n\rightarrow \infty .
\end{equation*}

By using again (\ref{112}) we obtain that when the sequence $\{q_{k}:k\in
\mathbb{N}\}$ is non-decreasing, then the summability method is not always
regular.

The proof is complete.
\QED

\bigskip

Let $t_{n}$ be N\"orlund means with monotone and bounded sequence $%
\{q_{k}:k\in \mathbb{N}\},$ such that
\begin{equation*}
q:=\lim_{n\rightarrow \infty }q_{n}>c>0.
\end{equation*}

Then, if the sequence $\{q_{k}:k\in \mathbb{N}\}$ is non-decreasing, we get
that
\begin{equation*}
nq_{0}\leq Q_{n}\leq nq.
\end{equation*}

In the case when the sequence $\{q_{k}:k\in \mathbb{N}\}$ is non-increasing,
then
\begin{equation}
nq\leq Q_{n}\leq nq_{0}.  \label{monotone0}
\end{equation}

In both cases we can conclude that
\begin{equation}
\frac{q_{n-1}}{Q_{n}}=O\left( \frac{1}{n}\right) ,\text{ \ \ when \ \ \ }%
n\rightarrow \infty .  \label{monotone1}
\end{equation}

In the special case when $\{q_{k}=1:k\in \mathbb{N}\},$ we get Fej\'er means
\begin{equation*}
\sigma _{n}f:=\frac{1}{n}\sum_{k=1}^{n}S_{k}f\,.
\end{equation*}%
$\,$

The $\left( C,\alpha \right) $-means (Ces\`aro means) of the Vilenkin-Fourier
series are defined by
\begin{equation*}
\sigma _{n}^{\alpha }f:=\frac{1}{A_{n}^{\alpha }}\overset{n}{\underset{k=1}{%
\sum }}A_{n-k}^{\alpha -1}S_{k}f,
\end{equation*}%
where \
\begin{equation*}
A_{0}^{\alpha }:=0,\text{ \ \ }A_{n}^{\alpha }:=\frac{\left( \alpha
+1\right) ...\left( \alpha +n\right) }{n!},~~\alpha \neq -1,-2,...
\end{equation*}

It is well-known that (see e.g. Zygmund \cite{13})%
\begin{equation}
A_{n}^{\alpha }=\overset{n}{\underset{k=0}{\sum }}A_{n-k}^{\alpha -1},
\label{node0}
\end{equation}%
\begin{equation}
A_{n}^{\alpha }-A_{n-1}^{\alpha }=A_{n}^{\alpha -1},\ \ \ \ \ \
A_{n}^{\alpha }\backsim n^{\alpha }.  \label{node01}
\end{equation}

In the literature, there is the notion of Riesz means ($\left( R,\alpha
\right) $-means) of the Fourier series. Let $\beta _{n}^{\alpha }$ denote
the N\"orlund mean, where
\begin{equation*}
\left\{ q_{0}=1,\text{ }q_{k}=k^{\alpha -1}:k\in \mathbb{N}_{+}\right\} ,
\end{equation*}%
that is%
\begin{equation*}
\beta _{n}^{\alpha }f:=\frac{1}{Q_{n}}\overset{n}{\underset{k=1}{\sum }}%
\left( n-k\right) ^{\alpha -1}S_{k}f,\text{ \ \ }0<\alpha <1.
\end{equation*}

It is obvious that

\begin{equation}
\frac{\left\vert q_{n}-q_{n+1}\right\vert }{n^{\alpha -2}}=O\left( 1\right)
,\ \ \text{when \ }n\rightarrow \infty .  \label{node1}
\end{equation}

and
\begin{equation}
\frac{q_{0}}{Q_{n}}=O\left( \frac{1}{n^{\alpha }}\right) ,\text{ \ when \ }%
n\rightarrow \infty .  \label{node2}
\end{equation}

The $n$-th N\"orlund logarithmic mean $L_{n}$ and the Riesz logarithmic mean $%
R_{n}$ are defined by
\begin{equation*}
L_{n}f:=\frac{1}{l_{n}}\sum_{k=1}^{n-1}\frac{S_{k}f}{n-k},\text{ \ \ }%
R_{n}f:=\frac{1}{l_{n}}\sum_{k=1}^{n-1}\frac{S_{k}f}{k},
\end{equation*}%
\ respectively, where%
\begin{equation*}
l_{n}:=\sum_{k=1}^{n-1}\frac{1}{k}.
\end{equation*}

Up to now we have considered N\"orlund mean in the case when the sequence $\
\{q_{k}:k\in \mathbb{N}\}$ is bounded but now we consider N\"orlund
summabilities with unbounded sequence $\{q_{k}:k\in \mathbb{N}\}.$

Let $\alpha \in
\mathbb{R}
_{+},$ $\beta \in
\mathbb{N}
_{+}$ and
\begin{equation*}
\log ^{(\beta )}x:=\overset{\beta \text{ times}}{\overbrace{\log ...\log }}x.
\end{equation*}

If we define the sequence $\{q_{k}:k\in \mathbb{N}\}$ by
\begin{equation*}
\left\{ q_{0}=0\text{ and }q_{k}=\log ^{\left( \beta \right) }k^{\alpha
}:k\in
\mathbb{N}
_{+}\right\} ,
\end{equation*}%
then we get the class of N\"orlund means with non-decreasing coefficients:
\begin{equation*}
\kappa _{n}^{\alpha ,\beta }f:=\frac{1}{Q_{n}}\sum_{k=1}^{n}\log ^{\left(
\beta \right) }\left( n-k\right) ^{\alpha }S_{k}f.
\end{equation*}%
First we note that $\kappa _{n}^{\alpha ,\beta }$ are
well-defined for every $n\in
\mathbb{N}
_{+}$, if we rewrite them as:
\begin{equation*}
\kappa _{n}^{\alpha ,\beta }f=\sum_{k=1}^{n}\frac{\log ^{\left( \beta
\right) }\left( n-k\right) ^{\alpha }}{Q_{n}}S_{k}f.
\end{equation*}

It is obvious that

\begin{equation*}
\frac{n}{2}\log ^{\left( \beta \right) }\frac{n^{\alpha }}{2^{\alpha }}\leq
Q_{n}\leq n\log ^{\left( \beta \right) }n^{\alpha }.
\end{equation*}

It follows that%
\begin{equation}
\frac{q_{n-1}}{Q_{n}}\leq \frac{c\log ^{\left( \beta \right) }\left(
n-1\right) ^{\alpha }}{n\log ^{\left( \beta \right) }n^{\alpha }}=O\left(
\frac{1}{n}\right) \rightarrow 0,\text{ when \ }n\rightarrow \infty .
\label{node00}
\end{equation}

For the function $\ f$ \ we consider the following maximal operators:%
\begin{equation*}
S^{\ast }f:=\sup_{n\in
\mathbb{N}
}\left\vert S_{n}f\right\vert ,\text{ \ \ \ \ }t^{\ast }f:=\sup_{n\in
\mathbb{N}
}\left\vert t_{n}f\right\vert ,\text{ \ \ \ \ \ }\sigma ^{\ast
}f:=\sup_{n\in
\mathbb{N}
}\left\vert \sigma _{n}f\right\vert ,
\end{equation*}%
\begin{equation*}
\text{\ }\sigma ^{\alpha ,\ast }f:=\sup_{n\in
\mathbb{N}
}\left\vert \sigma _{n}^{\alpha }f\right\vert ,\text{ \ \ }R^{\ast
}f:=\sup_{n\in
\mathbb{N}
}\left\vert R_{n}f\right\vert ,\text{ \ \ }L^{\ast }f:=\sup_{n\in
\mathbb{N}
}\left\vert L_{n}f\right\vert ,
\end{equation*}

\begin{equation*}
\kappa ^{\alpha ,\beta ,\ast }f:=\sup_{n\in
\mathbb{N}
}\left\vert \kappa _{n}^{\alpha ,\beta }f\right\vert ,\text{ \ \ }\beta
^{\alpha ,\ast }f:=\sup_{n\in
\mathbb{N}
}\left\vert \beta _{n}^{\alpha }f\right\vert .
\end{equation*}

We also consider the following weighted maximal operators:%
\begin{equation*}
\overset{\sim }{t}_{p,\alpha }^{\ast }f:=\sup_{n\in \mathbb{N}}\frac{%
\left\vert t_{n}f\right\vert }{\left( n+1\right) ^{1/p-1-\alpha }},\ \ \ \ \
\ \ \ \ \ \ \ (0<p<1/\left( 1+\alpha \right) ,\text{ }0<\alpha \leq 1),
\end{equation*}%
\begin{equation*}
\widetilde{t}_{\alpha }^{\ast }f:=\sup_{n\in \mathbb{N}}\frac{\left\vert
t_{n}f\right\vert }{\log ^{1+\alpha }\left( n+1\right) },\text{ \ \ }\
\left( 0<\alpha \leq 1\right) ,
\end{equation*}%
\begin{equation*}
\overset{\sim }{\sigma }_{p}^{\alpha ,\ast }f:=\sup_{n\in \mathbb{N}}\frac{%
\left\vert \sigma _{n}^{\alpha }f\right\vert }{\left( n+1\right)
^{1/p-1-\alpha }},\ \ \ \ \ \ \ \ \ \ \ \text{ }(0<p<1/\left( 1+\alpha
\right) ,\text{ }0<\alpha \leq 1),
\end{equation*}%
\begin{equation*}
\overset{\sim }{\beta }_{p}^{\alpha ,\ast }f:=\sup_{n\in \mathbb{N}}\frac{%
\left\vert \beta _{n}^{\alpha }f\right\vert }{\left( n+1\right)
^{1/p-1-\alpha }},\text{ \ \ \ \ \ \ \ \ \ \ \ }(0<p<1/\left( 1+\alpha
\right) ,\text{ }0<\alpha \leq 1),
\end{equation*}%
\begin{equation*}
\widetilde{\sigma }^{\alpha ,\ast }f:=\sup_{n\in \mathbb{N}}\frac{\left\vert
\sigma _{n}^{\alpha }f\right\vert }{\log ^{1+\alpha }\left( n+1\right) },%
\text{ \ \ }\widetilde{\beta }^{\alpha ,\ast }f:=\sup_{n\in \mathbb{N}}\frac{%
\left\vert \beta _{n}^{\alpha }f\right\vert }{\log ^{1+\alpha }\left(
n+1\right) },\text{ }\left( 0<\alpha <1\right) ,
\end{equation*}%
\
\begin{equation*}
\widetilde{\sigma }_{p}^{\ast }:=\sup_{n\in \mathbb{N}}\frac{\left\vert
\sigma _{n}\right\vert }{\left( n+1\right) ^{1/p-2}},\text{ \ \ \ }%
\widetilde{\kappa }_{p}^{\alpha ,\beta ,\ast }f:=\sup_{n\in \mathbb{N}}\frac{%
\left\vert \kappa _{n}^{\alpha ,\beta }f\right\vert }{\left( n+1\right)
^{1/p-2}},\text{ \ }\left( 0<p<1/2\right) ,
\end{equation*}%
and%
\begin{equation*}
\widetilde{\sigma }^{\ast }f:=\sup_{n\in \mathbb{N}}\frac{\left\vert \sigma
_{n}f\right\vert }{\log ^{2}\left( n+1\right) },\text{ \ \ \ \ }\widetilde{%
\kappa }^{\alpha ,\beta ,\ast }f:=\sup_{n\in \mathbb{N}}\frac{\left\vert
\kappa _{n}^{\alpha ,\beta }f\right\vert }{\log ^{2}\left( n+1\right) }.
\end{equation*}%
\

\subsection{Dirichlet kernels and Lebesgue constants with respect to
Vilenkin systems}

\bigskip It is easy to see that%
\begin{equation*}
S_{n}f\left( x\right) =\int_{G_{m}}f\left( t\right) \sum_{k=0}^{n-1}\psi
_{k}\left( x-t\right) d\mu \left( t\right)
\end{equation*}%
\begin{equation*}
=\int_{G_{m}}f\left( t\right) D_{n}\left( x-t\right) d\mu \left( t\right)
=\left( f\ast D_{n}\right) \left( x\right) .
\end{equation*}

The next well-known identities with respect to Dirichlet kernels (see Lemmas %
\ref{dn1} and \ref{dn2.1}, Corollaries \ref{dn2.2} and \ref{dn2.3}) will be
used many time in the proofs of our main results. The first equality can be
found in Vilenkin \cite{Vi} and the second identity in G\`{a}t and\textit{\ }%
Goginava \cite{Ga2}:

\begin{lemma}
\label{dn1}Let $n\in
\mathbb{N}
$, $1\leq s_{n}\leq m_{n}-1.$ Then%
\begin{equation}
D_{j+M_{n}}=D_{M_{n}}+\psi _{M_{n}}D_{j}=D_{M_{n}}+r_{n}D_{j},\text{ \ \ }%
\,j\leq \left( m_{n}-1\right) M_{n}  \label{dn21}
\end{equation}%
and
\begin{equation}
D_{s_{n}M_{n}-j}=D_{s_{n}M_{n}}-\psi _{s_{n}M_{n}-1}\overline{D_{j}},\text{ }%
\,j<M_{n}.  \label{dn22}
\end{equation}
\end{lemma}

\begin{lemma}
\label{dn2.1}Let $n\in \mathbb{N}$ and $1\leq s_{n}\leq m_{n}-1.$ Then \
\begin{equation}
D_{s_{n}M_{n}}=D_{M_{n}}\sum_{k=0}^{s_{n}-1}\psi
_{kM_{n}}=D_{M_{n}}\sum_{k=0}^{s_{n}-1}r_{n}^{k}  \label{9dn}
\end{equation}%
and%
\begin{equation}
D_{n}=\psi _{n}\left( \sum_{j=0}^{\infty
}D_{M_{j}}\sum_{k=m_{j}-n_{j}}^{m_{j}-1}r_{j}^{k}\right) ,  \label{2dn}
\end{equation}%
for $n=\sum_{i=0}^{\infty }n_{i}M_{i}.$
\end{lemma}

\begin{corollary}
\label{dn2.2}Let $n\in \mathbb{N}.$ Then%
\begin{equation*}
D_{M_{n+1}}=\prod_{k=0}^{n}\left( \sum_{s=0}^{m_{k}-1}r_{k}^{s}\right) .
\end{equation*}
\end{corollary}

\begin{corollary}
\label{dn2.3}Let $n\in \mathbb{N}.$ Then%
\begin{equation*}
D_{M_{n}}\left( x\right) \text{\thinspace }=\left\{
\begin{array}{ll}
M_{n}, & x\in I_{n}, \\
0, & x\notin I_{n}.%
\end{array}%
\right.
\end{equation*}
\end{corollary}

We also need the following estimate (see Tephnadze \cite{tep6}):

\begin{lemma}
\label{dn2.6}\bigskip Let $x\in I_{s}\backslash I_{s+1},$ $s=0,...,N-1.$ Then%
\begin{equation*}
\int_{I_{N}}\left\vert D_{n}\left( x-t\right) \right\vert d\mu \left(
t\right) \leq \frac{cM_{s}}{M_{N}},
\end{equation*}%
where $c$ is an absolute constant.
\end{lemma}

{\bf Proof}:
By combining (\ref{2dn}) in Lemma \ref{dn2.1} and Corollary \ref{dn2.3} we
have that
\begin{equation*}
\left\vert D_{n}\left( x\right) \right\vert \leq \underset{j=0}{\overset{l}{%
\sum }}n_{j}D_{M_{j}}\left( x\right) =\underset{j=0}{\overset{l}{\sum }}%
n_{j}M_{j}\leq cM_{l}.
\end{equation*}

Since $t\in I_{N}$ and $x\in I_{s}\backslash I_{s+1},$ $s=0,...,N-1,$ we
obtain that $x-t\in I_{s}\backslash I_{s+1}$. By using the estimate above we get
that
\begin{equation*}
\left\vert D_{n}\left( x-t\right) \right\vert \leq cM_{s}
\end{equation*}%
and%
\begin{equation*}
\int_{I_{N}}\left\vert D_{n}\left( x-t\right) \right\vert d\mu \left(
t\right) \leq \,\frac{cM_{s}}{M_{N}}.
\end{equation*}

The proof is complete.
\QED

\bigskip To study the Dirichlet kernels we need an estimate of some sums of
Rademacher functions. This Lemma can be found in Persson and Tephnadze \cite%
{pt2}.

\begin{lemma}
\label{rad4}Let $n\in \mathbb{N}$, and $x_{n}=1.$ Then
\begin{equation*}
\left\vert \sum_{u=0}^{s_{n}-1}r_{n}^{u}\left( x\right) \right\vert \geq 1,%
\text{ for some }1\leq s_{n}\leq m_{n}-1
\end{equation*}%
and%
\begin{equation*}
\left\vert \sum_{u=1}^{s_{n}-1}r_{n}^{u}\left( x\right) \right\vert \geq 1,%
\text{ \ for some }2\leq s_{n}\leq m_{n}-1.
\end{equation*}
\end{lemma}

{\bf Proof}:
Let $x_{n}=1.$ Then we readily get that
\begin{equation*}
\left\vert \sum_{u=0}^{s_{n}-1}r_{n}^{u}\left( x\right) \right\vert
=\left\vert \frac{r_{n}^{s_{n}}\left( x\right) -1}{r_{n}\left( x\right) -1}%
\right\vert
\end{equation*}%
\begin{equation*}
=\frac{\sin \left( \pi s_{n}x_{n}/m_{n}\right) }{\sin \left( \pi
x_{n}/m_{n}\right) }=\frac{\sin \left( \pi s_{n}/m_{n}\right) }{\sin \left(
\pi /m_{n}\right) }\geq 1.
\end{equation*}

Analogously, we can prove that%
\begin{equation*}
\left\vert \sum_{u=1}^{s_{n}-1}r_{n}^{u}\left( x\right) \right\vert
=\left\vert r_{n}\left( x\right) \frac{r_{n}^{s_{n}-1}\left( x\right) -1}{%
r_{n}\left( x\right) -1}\right\vert =\frac{\sin \left( \pi \left(
s_{n}-1\right) /m_{n}\right) }{\sin \left( \pi /m_{n}\right) }\geq 1.
\end{equation*}

The proof is complete.
\QED

The next Lemma can be found in Tephnadze \cite{tep6}:

\begin{lemma}
\label{dn2.6.2}\bigskip Let $n\in \mathbb{N}$, $\left\vert n\right\vert \neq
\left\langle n\right\rangle $ and $x\in I_{\left\langle n\right\rangle
}\backslash I_{\left\langle n\right\rangle +1}.$ Then
\begin{equation*}
\left\vert D_{n}\right\vert =\left\vert D_{n-M_{\left\vert n\right\vert
}}\right\vert \geq M_{\left\langle n\right\rangle }.
\end{equation*}
\end{lemma}

{\bf Proof}:
\bigskip Let $x\in I_{\left\langle n\right\rangle }\backslash I_{\left\langle
n\right\rangle +1}.$\textbf{\ }Since
\begin{equation*}
n=n_{\left\langle n\right\rangle }M_{\left\langle n\right\rangle
}+\sum_{j=\left\langle n\right\rangle }^{\left\vert n\right\vert
-1}n_{j}M_{j}+n_{\left\vert n\right\vert }M_{\left\vert n\right\vert }
\end{equation*}%
and%
\begin{equation*}
n-M_{\left\vert n\right\vert }=n_{\left\langle n\right\rangle
}M_{\left\langle n\right\rangle }+\sum_{j=\left\langle n\right\rangle
}^{\left\vert n\right\vert -1}n_{j}M_{j}+\left( n_{\left\vert n\right\vert
}-1\right) M_{\left\vert n\right\vert },  \label{n1}
\end{equation*}%
if we apply Lemma \ref{rad4}, Corollary \ref{dn2.3} and (\ref{2dn}) in Lemma %
\ref{dn2.1} we can conclude that
\begin{equation*}
\left\vert D_{n-M_{\left\vert n\right\vert }}\right\vert
\end{equation*}%
\begin{equation*}
\geq \left\vert \psi _{\left\langle n\right\rangle }D_{M_{\left\langle
n\right\rangle }}\sum_{s=m_{\left\langle n\right\rangle }-n_{\left\langle
n\right\rangle }}^{m_{\left\langle n\right\rangle }-1}r_{\left\langle
n\right\rangle }^{s}\right\vert -\left\vert \psi _{\left\langle
n\right\rangle }\sum_{j=\left\langle n\right\rangle +1}^{\left\vert
n\right\vert }D_{M_{j}}\sum_{s=m_{j}-n_{j}}^{m_{j}-1}r_{j}^{s}\right\vert
\end{equation*}%
\begin{equation*}
=\left\vert D_{M_{\left\langle n\right\rangle }}\sum_{s=m_{\left\langle
n\right\rangle }-n_{\left\langle n\right\rangle }}^{m_{\left\langle
n\right\rangle }-1}r_{\left\langle n\right\rangle }^{s}\right\vert
=\left\vert D_{M_{\left\langle n\right\rangle }}r_{\left\langle
n\right\rangle }^{m_{\left\langle n\right\rangle }-\alpha _{\left\langle
n\right\rangle }}\sum_{s=0}^{n_{\left\langle n\right\rangle
}-1}r_{\left\langle n\right\rangle }^{s}\right\vert
\end{equation*}%
\begin{equation*}
=D_{M_{\left\langle n\right\rangle }}\left\vert \sum_{s=0}^{n_{\left\langle
n\right\rangle }-1}r_{\left\langle n\right\rangle }^{s}\right\vert \geq
D_{M_{\left\langle n\right\rangle }}\geq M_{\left\langle n\right\rangle }.
\end{equation*}%
Analogously we can show that
\begin{equation*}
\left\vert D_{n}\right\vert =D_{M_{\left\langle n\right\rangle }}\left\vert
\sum_{s=0}^{n_{\left\langle n\right\rangle }}r_{\left\langle n\right\rangle
}^{s}\right\vert \geq M_{\left\langle n\right\rangle }
\end{equation*}%
and the proof is complete.
\QED

The next Lemma can be found in Lukomskii \cite{luko}:

\begin{lemma}
\label{theorem1}Let $n=\sum_{i=1}^{\infty }n_{i}M_{i}$. Then
\begin{equation*}
\frac{1}{4\lambda }v\left( n\right) +\frac{1}{\lambda }v^{\ast }\left(
n\right) +\frac{1}{2\lambda }\leq L_{n}\leq \frac{3}{2}v\left( n\right)
+4v^{\ast }\left( n\right) -1,
\end{equation*}%
where $\lambda :=\sup_{n\in \mathbb{N}}m_{n}.$
\end{lemma}

The next result for the Walsh system can be found in the book \cite{sws} and for
bounded Vilenkin systems in the book \cite{AVD}.

\begin{corollary}
\label{dn5}\bigskip\ Let $q_{n_{k}}$\bigskip $%
=M_{2n_{k}}+M_{2n_{k}-2}+M_{2}+M_{0}$. Then%
\begin{equation*}
\frac{n_{k}}{2\lambda }\leq \left\Vert D_{q_{n_{k}}}\right\Vert _{1}\leq
\lambda n_{k},
\end{equation*}%
where $\lambda :=\sup_{n\in \mathbb{N}}m_{n}.$
\end{corollary}

{\bf Proof}:
The proof readily follows by just using Theorem \ref{theorem1} and the following identity
$
v\left( q_{n_{k}}\right) =2n_{k}.
$
Thus, we leave out the details.
\QED

\subsection{Fej\'er Kernels with respect to Vilenkin systems}

It is obvious that%
\begin{equation*}
\sigma _{n}f\left( x\right) =\frac{1}{n}\overset{n-1}{\underset{k=0}{\sum }}%
\left( D_{k}\ast f\right) \left( x\right)
\end{equation*}%
\begin{equation*}
=\left( f\ast K_{n}\right) \left( x\right) =\int_{G_{m}}f\left( t\right)
K_{n}\left( x-t\right) d\mu \left( t\right) .
\end{equation*}

We frequently use the following well-known result, which was proved in G\'{a}%
t\textit{\ }\cite{gat}:

\begin{lemma}
\label{lemma2} Let $n>t,$ $t,n\in \mathbb{N}.$ Then%
\begin{equation*}
K_{M_{n}}\left( x\right) =\left\{
\begin{array}{ll}
\frac{M_{t}}{1-r_{t}\left( x\right) }, & x\in I_{t}\backslash I_{t+1},\text{
}x-x_{t}e_{t}\in I_{n}, \\
\frac{M_{n}+1}{2}, & x\in I_{n}, \\
0, & \text{otherwise. }%
\end{array}%
\right.
\end{equation*}%

The proof of the next lemma can easily be done by using last lemma (c.f. also the book \cite{AVD} and Tephnadze \cite%
{tep2, tep3}:
\end{lemma}

\begin{lemma}
\label{lemma222} Let $n\in \mathbb{N}$ and $x\in I_{N}^{k,l},$ where $k<l.$
Then%
\begin{equation}
K_{M_{n}}\left( x\right) =0,\,\,\text{if \thinspace \thinspace }n>l.
\label{star1}
\end{equation}%
and%
\begin{equation}
\left\vert K_{M_{n}}\left( x\right) \right\vert \leq cM_{k}.  \label{star2}
\end{equation}%
Moreover,%
\begin{equation}
\int_{G_{m}}\left\vert K_{M_{n}}\right\vert d\mu \leq c<\infty ,
\label{star3}
\end{equation}%
where $c$ is an absolute constant.
\end{lemma}

We also need the following useful result:

\begin{lemma}
\label{lemma6kn}Let $t,s_{n},n\in \mathbb{N},$ and $1\leq s_{n}\leq m_{n}-1$%
. Then\bigskip
\begin{equation}
s_{n}M_{n}K_{s_{n}M_{n}}=\sum_{l=0}^{s_{n}-1}\left(
\sum_{i=0}^{l-1}r_{n}^{i}\right) M_{n}D_{M_{n}}+\left(
\sum_{l=0}^{s_{n}-1}r_{n}^{l}\right) M_{n}K_{M_{n}}  \label{mag}
\end{equation}%
and%
\begin{equation}
\left\vert K_{s_{n}M_{n}}\left( x\right) \right\vert \geq \frac{M_{n}}{2\pi
s_{n}},\text{ for }x\in I_{n+1}\left( e_{n-1}+e_{n}\right) .  \label{100kn1}
\end{equation}%
Moreover, if $x\in I_{t}\backslash I_{t+1},$ \ $x-x_{t}e_{t}\notin I_{n}$
and $n>t,$ then
\begin{equation*}
K_{s_{n}M_{n}}(x)=0.
\end{equation*}
\end{lemma}

\begin{remark}
This result was proved by Blahota and Tephnadze \cite{bt1}, but here we will give a completely
different and simpler proof.
\end{remark}

{\bf Proof}:
We can write that%
\begin{equation}
s_{n}M_{n}K_{s_{n}M_{n}}=\sum_{l=0}^{s_{n}-1}\sum_{k=lM_{n}}^{\left(
l+1\right) M_{n}-1}D_{k}  \label{skn}
\end{equation}%
\begin{equation*}
=\sum_{l=0}^{s_{n}-1}\sum_{k=lM_{n}}^{\left( l+1\right)
M_{n}-1}D_{k}=\sum_{l=0}^{s_{n}-1}\sum_{k=0}^{M_{n}-1}D_{k+lM_{n}}.
\end{equation*}

Let $0\leq k<M_{n}$. Then, by using (\ref{9dn}) in Lemma \ref{dn1} we have
that\qquad\
\begin{equation*}
D_{k+lM_{n}}=\sum_{m=0}^{lM_{n}-1}\psi _{m}+\sum_{m=lM_{n}}^{lM_{n}+k-1}\psi
_{m}
\end{equation*}%
\begin{equation*}
=D_{lM_{n}}+\sum_{m=0}^{k-1}\psi
_{m+lM_{n}}=D_{lM_{n}}+r_{n}^{l}\sum_{m=0}^{k-1}\psi _{m}
\end{equation*}%
\begin{equation*}
=\left( \sum_{i=0}^{l-1}r_{n}^{i}\right) D_{M_{n}}+r_{n}^{l}D_{k}.
\end{equation*}

According to (\ref{skn}) we readily get that
\begin{equation*}
s_{n}M_{n}K_{s_{n}M_{n}}=\sum_{l=0}^{s_{n}-1}\sum_{k=0}^{M_{n}-1}D_{k+lM_{n}}
\end{equation*}%
\begin{equation*}
=\sum_{l=0}^{s_{n}-1}\sum_{k=0}^{M_{n}-1}\left( \left(
\sum_{i=0}^{l-1}r_{n}^{i}\right) D_{M_{n}}+r_{n}^{l}D_{k}\right)
\end{equation*}%
\begin{equation*}
=\sum_{l=0}^{s_{n}-1}\left( \sum_{i=0}^{l-1}r_{n}^{i}\right)
M_{n}D_{M_{n}}+\left( \sum_{l=0}^{s_{n}-1}r_{n}^{l}\right)
\sum_{k=0}^{M_{n}-1}D_{k}
\end{equation*}%
\begin{equation*}
=\sum_{l=0}^{s_{n}-1}\left( \sum_{i=0}^{l-1}r_{n}^{i}\right)
M_{n}D_{M_{n}}+\left( \sum_{l=0}^{s_{n}-1}r_{n}^{l}\right) M_{n}K_{M_{n}}.
\end{equation*}

The first part of the Lemma is proved.

Now, let $t,s_{n},n\in \mathbb{N},\ n>t,\ x\in I_{t}\backslash I_{t+1}$. If $%
x-x_{t}e_{t}\notin I_{n}$, then by combining Corollary \ref{dn2.3} and Lemma %
\ref{lemma2} we obtain that
\begin{equation}
D_{M_{n}}(x)=K_{M_{n}}(x)=0.  \label{star4}
\end{equation}

Let $x\in I_{n+1}\left( e_{n-1}+e_{n}\right) .$ By Lemma \ref{lemma2} we
have that%
\begin{equation*}
\left\vert K_{M_{n}}\left( x\right) \right\vert =\frac{M_{n-1}}{\left\vert
1-r_{n-1}\left( x\right) \right\vert }=\frac{M_{n-1}}{2\sin \pi /m_{n-1}}.
\end{equation*}

By combining Lemmas \ref{rad4} and \ref{lemma2} and the first part of Lemma %
\ref{lemma6kn} we immediately get that%
\begin{equation*}
\left\vert s_{n}M_{n}K_{s_{n}M_{n}}\left( x\right) \right\vert =\left\vert
\left( \sum_{l=0}^{s_{n}-1}r_{n}^{l}\left( x\right) \right)
M_{n}K_{M_{n}}\left( x\right) \right\vert
\end{equation*}%
\begin{equation*}
\frac{M_{n}M_{n-1}}{\left\vert 1-r_{n}\left( x\right) \right\vert }=\frac{%
M_{n}M_{n-1}}{2\sin \pi /m_{n-1}}\geq \frac{M_{n}M_{n-1}m_{n-1}}{2\pi }\geq
\frac{M_{n}^{2}}{2\pi }.
\end{equation*}

Now, let $t,s_{n},n\in \mathbb{N},\ n>t,\ x\in I_{t}\backslash I_{t+1}$. If $%
x-x_{t}e_{t}\notin I_{n}$, then, by using the first part of Lemma \ref%
{lemma6kn}, with identities (\ref{star4}) we immediately get that%
\begin{equation*}
K_{s_{n}M_{n}}(x)=0.
\end{equation*}

The proof is complete.\vspace{0pt}\qquad
\QED

\bigskip In the same paper Blahota and Tephnadze \cite{bt1} also proved the
following result:

\begin{lemma}
\label{lemma4}Let $n=\sum_{i=1}^{r}s_{n_{i}}M_{n_{i}}$, where $%
n_{1}>n_{2}>\dots >n_{r}\geq 0$ and $1\leq s_{n_{i}}<m_{n_{i}}$ \ for all $%
1\leq i\leq r$ as well as $n^{(k)}=n-\sum_{i=1}^{k}s_{n_{i}}M_{n_{i}}$,
where $0<k\leq r$. Then
\begin{equation*}
nK_{n}=\sum_{k=1}^{r}\left( \prod_{j=1}^{k-1}r_{n_{j}}^{s_{n_{j}}}\right)
s_{n_{k}}M_{n_{k}}K_{s_{n_{k}}M_{n_{k}}}+\sum_{k=1}^{r-1}\left(
\prod_{j=1}^{k-1}r_{n_{j}}^{s_{n_{j}}}\right) n^{(k)}D_{s_{n_{k}}M_{n_{k}}}.
\end{equation*}
\end{lemma}

{\bf Proof}:
Let $k,n\in \mathbb{N},\ 0\leq k<M_{n}$. If we use the identity (\ref{dn21})
in Lemma \ref{dn1} we readily get that
\begin{equation*}
nK_{n}=\sum_{k=1}^{n}D_{k}=\sum_{k=1}^{s_{n_{1}}M_{n_{1}}}D_{k}+%
\sum_{k=s_{n_{1}}M_{n_{1}}+1}^{n}D_{k}
\end{equation*}%
\begin{equation*}
=s_{n_{1}}M_{n_{1}}K_{s_{n_{1}}M_{n_{1}}}+%
\sum_{k=1}^{n^{(1)}}D_{k+s_{n_{1}}M_{n_{1}}}
\end{equation*}%
\begin{equation*}
=s_{n_{1}}M_{n_{1}}K_{s_{n_{1}}M_{n_{1}}}+\sum_{k=1}^{n^{(1)}}\left(
D_{s_{n_{1}}M_{n_{1}}}+r_{n_{1}}^{s_{n_{1}}}D_{k}\right)
\end{equation*}%
\begin{equation*}
=s_{n_{1}}M_{n_{1}}K_{s_{n_{1}}M_{n_{1}}}+n^{(1)}D_{s_{n_{1}}M_{n_{1}}}+r_{n_{1}}^{s_{n_{1}}}n^{(1)}K_{n^{(1)}}.
\end{equation*}%
If we calculate $n^{(1)}K_{n^{(1)}}$ in similar way, we get that
\begin{equation*}
n^{(1)}K_{n^{(1)}}=s_{n_{2}}M_{n_{2}}K_{s_{n_{2}}M_{n_{2}}}+n^{(2)}D_{s_{n_{2}}M_{n_{2}}}+r_{n_{2}}^{s_{n_{2}}}n^{(2)}K_{n^{(2)}},
\end{equation*}%
so
\begin{equation*}
nK_{n}=s_{n_{1}}M_{n_{1}}K_{s_{n_{1}}M_{n_{1}}}+r_{n_{1}}^{s_{n_{1}}}s_{n_{2}}M_{n_{2}}K_{s_{n_{2}}M_{n_{2}}}
\end{equation*}%
\begin{equation*}
+r_{n_{1}}^{s_{n_{1}}}r_{n_{2}}^{s_{n_{2}}}n^{(2)}K_{n^{(2)}}+n^{(1)}D_{s_{n_{1}}M_{n_{1}}}+r_{n_{1}}^{s_{n_{1}}}n^{(2)}D_{s_{n_{2}}M_{n_{2}}}.
\end{equation*}

By using this method successively with $n^{(2)}K_{n^{(2)}},\dots
,n^{(r-1)}K_{n^{(r-1)}}$, we obtain that
\begin{equation*}
nK_{n}=\sum_{k=1}^{r}\left( \prod_{j=1}^{k-1}r_{n_{j}}^{s_{n_{j}}}\right)
s_{n_{k}}M_{n_{k}}K_{s_{n_{k}}M_{n_{k}}}
\end{equation*}%
\begin{equation*}
+\left( \prod_{j=1}^{r}r_{n_{j}}^{s_{n_{j}}}\right)
n^{(r)}K_{n^{(r)}}+\sum_{k=1}^{r-1}\left(
\prod_{j=1}^{k-1}r_{n_{j}}^{s_{n_{j}}}\right) n^{(k)}D_{s_{n_{k}}M_{n_{k}}}.
\end{equation*}

Since $n^{(r)}=0$ we conclude that the proof is complete.\vspace{0pt}
\QED

\begin{corollary}
\bigskip \label{lemma7kn}Let $n\in \mathbb{N}.$ Then%
\begin{equation}
n\left\vert K_{n}\right\vert \leq c\sum_{l=\left\langle n\right\rangle
}^{\left\vert n\right\vert }M_{l}\left\vert K_{M_{l}}\right\vert \leq
c\sum_{l=0}^{\left\vert n\right\vert }M_{l}\left\vert K_{M_{l}}\right\vert
\label{fn5}
\end{equation}%
and%
\begin{equation}
\sup_{n}\int_{G_{m}}\left\vert K_{n}\right\vert d\mu \leq c<\infty ,
\label{fn4}
\end{equation}%
where $c$ is an absolute constant.
\end{corollary}

\begin{remark}
Corollary \ref{lemma7kn} is known (see the book \textit{\ }\cite{AVD}), but it is also a simple consequence of
Lemmas  \ref{lemma222} and \ref{lemma4}.
\end{remark}

The next lemma can be found as a part of a more general result by Persson and Tephnadze \cite{pt2}:

\begin{lemma}
\label{lemma8ccc}Let $n\in \mathbb{N}$, $\left\langle n\right\rangle \neq
\left\vert n\right\vert $ and \ $x\in I_{\left\langle n\right\rangle
+1}\left( e_{\left\langle n\right\rangle -1}+e_{\left\langle n\right\rangle
}\right) .$ Then
\begin{equation*}
\left\vert nK_{n}\right\vert =\left\vert \left( n-M_{\left\vert n\right\vert
}\right) K_{-M_{\left\vert n\right\vert }}\right\vert \geq \frac{%
M_{\left\langle n\right\rangle }^{2}}{2\pi \lambda },
\end{equation*}%
where $\lambda :=\sup m_{n}.$
\end{lemma}

{\bf Proof}:
Let $x\in $ $I_{_{\left\langle n\right\rangle +1}}^{\left\langle
n\right\rangle -1,\left\langle n\right\rangle }.$ \textbf{\ }Since
\begin{equation*}
n=n_{\left\langle n\right\rangle }M_{\left\langle n\right\rangle
}+\sum_{j=\left\langle n\right\rangle }^{\left\vert n\right\vert
-1}n_{j}M_{j}+n_{\left\vert n\right\vert }M_{\left\vert n\right\vert }
\end{equation*}%
and%
\begin{equation*}
n-M_{\left\vert n\right\vert }=n_{\left\langle n\right\rangle
}M_{\left\langle n\right\rangle }+\sum_{j=\left\langle n\right\rangle
}^{\left\vert n\right\vert -1}n_{j}M_{j}+\left( n_{\left\vert n\right\vert
}-1\right) M_{\left\vert n\right\vert },
\end{equation*}%
if we combine (\ref{9dn}), (\ref{100kn1}) and invoke Corollary \ref{dn2.3}, Lemmas \ref{lemma2} and \ref%
{lemma4} we obtain that
\begin{equation*}
n\left\vert K_{n}\right\vert =\left( n-M_{\left\vert n\right\vert }\right)
\left\vert K_{n-M_{\left\vert n\right\vert }}\right\vert
\end{equation*}%
\begin{equation*}
=\left\vert \left( \overset{\left\langle n\right\rangle -1}{\underset{j=1}{%
\prod }}\psi _{M_{n_{j}}}^{j}\right) s_{\left\langle n\right\rangle
}M_{\left\langle n\right\rangle }K_{s_{\left\langle n\right\rangle
}M_{\left\langle n\right\rangle }}\right\vert =\left\vert s_{\left\langle
n\right\rangle }M_{\left\langle n\right\rangle }K_{s_{\left\langle
n\right\rangle }M_{\left\langle n\right\rangle }}\right\vert
\end{equation*}%
\begin{equation*}
=\left\vert \sum_{l=0}^{s_{\left\langle n\right\rangle }-1}\left(
\sum_{i=0}^{l-1}r_{\left\langle n\right\rangle }^{i}\right) M_{\left\langle
n\right\rangle }D_{M_{\left\langle n\right\rangle }}+\left(
\sum_{l=0}^{s_{\left\langle n\right\rangle }-1}r_{\left\langle
n\right\rangle }^{l}\right) M_{\left\langle n\right\rangle
}K_{M_{\left\langle n\right\rangle }}\right\vert
\end{equation*}%
\begin{equation*}
=\left\vert \left( \sum_{l=0}^{s_{\left\langle n\right\rangle
}-1}r_{\left\langle n\right\rangle }^{l}\right) M_{\left\langle
n\right\rangle }K_{M_{\left\langle n\right\rangle }}\right\vert \geq
\left\vert M_{\left\langle n\right\rangle }K_{M_{\left\langle n\right\rangle
}}\right\vert \geq \frac{M_{\left\langle n\right\rangle }^{2}}{2\pi \lambda }%
.
\end{equation*}

The proof is complete.
\QED

The next result is proved in Blahota, G\'{a}t and Goginava\textit{\ }\cite%
{BGG} and \cite{BGG2}, (see also Tephnadze \cite{tep5}):

\begin{lemma}
\label{lemma3} Let $2<n\in \mathbb{N}_{+},$ $k\leq s<n$ and $%
q_{n}=M_{2n}+M_{2n-2}+...+M_{2}+M_{0}.$ Then%
\begin{equation*}
q_{n-1}\left\vert K_{q_{n-1}}(x)\right\vert \geq \frac{M_{2k}M_{2s}}{8},
\end{equation*}%
for%
\begin{equation*}
x\in I_{2n}\left( 0,...,x_{2k}\neq 0,0,...,0,x_{2s}\neq
0,x_{2s+1},...,x_{2n-1}\right) ,
\end{equation*}%
\begin{equation*}
k=0,1,...,n-3,\text{ \ \ }s=k+2,k+3,...,n-1.
\end{equation*}
\end{lemma}

The next lemma can be found in Blahota and Tephnadze \cite{bt1} (see also
Tephnadze \cite{tep2} and \cite{tep3}):

\begin{lemma}
\label{lemma5}\ Let $x\in I_{N}^{k,l},$ $k=0,\dots ,N-2,$ $l=k+1,\dots ,N-1.$
Then
\begin{equation*}
\int_{I_{N}}\left\vert K_{n}\left( x-t\right) \right\vert d\mu \left(
t\right) \leq \frac{cM_{l}M_{k}}{nM_{N}}.
\end{equation*}

Let $x\in I_{N}^{k,N},$ $k=0,\dots ,N-1.$ Then%
\begin{equation*}
\int_{I_{N}}\left\vert K_{n}\left( x-t\right) \right\vert d\mu \left(
t\right) \leq \frac{cM_{k}}{M_{N}},
\end{equation*}%
where $c$ is an absolute constant.
\end{lemma}

{\bf Proof}:
Let $x\in I_{N}^{k,l},$ for $0\leq k<l\leq N-1$ and $t\in I_{N}.$ Since $%
x-t\in $ $I_{N}^{k,l}$ and $n\geq M_{N},$ by combining Lemma \ref{lemma2}
and (\ref{fn5}) in Corollary \ref{lemma7kn} we obtain that%
\begin{equation*}
n\left\vert K_{n}\left( x\right) \right\vert
\end{equation*}%
\begin{equation*}
\leq c\underset{i=0}{\overset{l}{\sum }}M_{i}\int_{I_{N}}\left\vert
K_{M_{i}}\left( x-t\right) \right\vert d\mu \left( t\right) \leq c\overset{l}%
{\underset{i=0}{\sum }}M_{i}M_{k}\leq cM_{k}M_{l}
\end{equation*}%
and%
\begin{equation}
\int_{I_{N}}\left\vert K_{n}\left( x-t\right) \right\vert d\mu \left(
t\right) \leq \frac{cM_{k}M_{l}}{nM_{N}}.  \label{888}
\end{equation}

Let $x\in I_{N}^{k,N}$. Then by applying Lemma \ref{lemma2} and (\ref{fn5})
in Corollary \ref{lemma7kn}, we have that%
\begin{equation}
\int_{I_{N}}n\left\vert K_{n}\left( x-t\right) \right\vert d\mu \left(
t\right) \leq \underset{i=0}{\overset{\left\vert n\right\vert }{\sum }}%
M_{i}\int_{I_{N}}\left\vert K_{M_{i}}\left( x-t\right) \right\vert d\mu
\left( t\right) .  \label{999}
\end{equation}

Let%
\begin{equation*}
\left\{
\begin{array}{l}
x=\left( 0,\ldots ,0,x_{k}\neq 0,\ldots
,x_{N-1}=0,x_{N},x_{N+1},x_{q},...,x_{\left\vert n\right\vert -1},\ldots
\right) , \\
t=\left( 0,\ldots ,0,x_{N},\ldots ,x_{q-1},t_{q}\neq x_{q},t_{q+1},\ldots
,t_{\left\vert n\right\vert -1},\ldots \right) ,\,\,q=N,\ldots ,\left\vert
n\right\vert -1.\text{ }%
\end{array}%
\right.
\end{equation*}

By using Lemmas \ref{lemma2} and \ref{lemma222} in (\ref{999}) it is easy to
see that%
\begin{equation}
\int_{I_{N}}\left\vert K_{n}\left( x-t\right) \right\vert d\mu \left(
t\right)  \label{11110}
\end{equation}%
\begin{equation*}
\leq \frac{c}{n}\underset{i=0}{\overset{q-1}{\sum }}M_{i}\int_{I_{N}}M_{k}d%
\mu \left( t\right) \leq \frac{cM_{k}M_{q}}{nM_{N}}\leq \frac{cM_{k}}{M_{N}}.
\end{equation*}

Let%
\begin{equation*}
\left\{
\begin{array}{l}
\text{ }x=\left( 0,\ldots ,0,x_{m}\neq 0,0,\ldots
,0,x_{N},x_{N+1},x_{q},\ldots ,x_{\left\vert n\right\vert -1},\ldots \right)
,\text{ } \\
t=\left( 0,,\ldots ,x_{N}=0,\ldots ,x_{_{\left\vert n\right\vert -1}},\ldots
\right) .%
\end{array}%
\right.
\end{equation*}

If we apply again Lemmas \ref{lemma2} and \ref{lemma222} in (\ref{999}) we
obtain that%
\begin{equation}
\int_{I_{N}}\left\vert K_{n}\left( x-t\right) \right\vert d\mu \left(
t\right) \leq \frac{c}{n}\overset{\left\vert n\right\vert -1}{\underset{i=0}{%
\sum }}M_{i}\int_{I_{N}}M_{k}d\mu \left( t\right) \leq \frac{cM_{k}}{M_{N}}.
\label{11111}
\end{equation}

By combining (\ref{888}), (\ref{11110}) and (\ref{11111}) we can complete
the proof.
\QED

\bigskip Also the next lemma is due to Tephnadze \cite{tep2, tep3}, but it
is also a simple consequence of Lemma \ref{lemma5}.

\begin{lemma}
\label{lemma5aa}\ Let $x\in I_{N}^{k,l},$ $k=0,\dots ,N-1,$ $l=k+1,\dots ,N.$
Then%
\begin{equation*}
\int_{I_{N}}\left\vert K_{n}\left( x-t\right) \right\vert d\mu \left(
t\right) \leq \frac{cM_{l}M_{k}}{M_{N}^{2}},\text{ \ for }n\geq M_{N},
\end{equation*}%
where $c$ is an absolute constant.
\end{lemma}

{\bf Proof}:
Since $n\geq M_{N}$ if we apply Lemma \ref{lemma5} we immediately get a proof
of this estimate.
\QED

\subsection{N\"orlund Kernels with respect to Vilenkin systems}

A representation
\begin{equation*}
t_{n}f\left( x\right) =\underset{G}{\int }f\left( t\right) F_{n}\left(
x-t\right) d\mu \left( t\right)
\end{equation*}%
plays a central role in the sequel, where

\begin{equation*}
F_{n}:=\frac{1}{Q_{n}}\overset{n}{\underset{k=1}{\sum }}q_{n-k}D_{k}
\end{equation*}%
is the so-called N\"orlund kernel.

In this section we study N\"orlund kernels with respect to Vilenkin
systems. The next results (Lemmas \ref{lemma0nn}-\ref{lemma5aaa}) are due to
Persson, Tephnadze and Wall \cite{ptw2}:

\begin{lemma}
\label{lemma0nn}Let $\{q_{k}:k\in \mathbb{N}\}$ be a sequence of
non-decreasing numbers, satisfying the condition
\begin{equation}
\frac{q_{n-1}}{Q_{n}}=O\left( \frac{1}{n}\right) ,\text{ \ \ when \ \ }\
n\rightarrow \infty .  \label{fn01}
\end{equation}%
\textit{Then}
\begin{equation*}
\left\vert F_{n}\right\vert \leq \frac{c}{n}\left\{ \sum_{j=0}^{\left\vert
n\right\vert }M_{j}\left\vert K_{M_{j}}\right\vert \right\} ,
\end{equation*}%
where $c$ is an absolute constant.
\end{lemma}

{\bf Proof}:
First, we invoke Abel transformation to obtain the following
identities%
\begin{equation}
Q_{n}:=\overset{n-1}{\underset{j=0}{\sum }}q_{j}=\overset{n}{\underset{j=1}{%
\sum }}q_{n-j}\cdot 1=\overset{n-1}{\underset{j=1}{\sum }}\left(
q_{n-j}-q_{n-j-1}\right) j+q_{0}n  \label{2b}
\end{equation}%
and%
\begin{equation}
F_{n}=\frac{1}{Q_{n}}\left( \overset{n-1}{\underset{j=1}{\sum }}\left(
q_{n-j}-q_{n-j-1}\right) jK_{j}+q_{0}nK_{n}\right) .  \label{2bb}
\end{equation}%
Let the sequence $\{q_{k}:k\in \mathbb{N}\}$ be non-decreasing. Then, by using (%
\ref{fn01}), we get that
\begin{equation*}
\frac{1}{Q_{n}}\left( \overset{n-1}{\underset{j=1}{\sum }}\left\vert
q_{n-j}-q_{n-j-1}\right\vert +q_{0}\right)
\end{equation*}%
\begin{equation*}
\leq \frac{1}{Q_{n}}\left( \overset{n-1}{\underset{j=1}{\sum }}\left(
q_{n-j}-q_{n-j-1}\right) +q_{0}\right) \leq \frac{q_{n-1}}{Q_{n}}\leq \frac{c%
}{n}.
\end{equation*}

Under condition (\ref{fn01}) if we apply (\ref{fn5}) in Corollary \ref%
{lemma7kn} and use the equalities (\ref{2b}) and (\ref{2bb}) we immediately get
that%
\begin{equation*}
\left\vert F_{n}\right\vert \leq \left( \frac{1}{Q_{n}}\left( \overset{n-1}{%
\underset{j=1}{\sum }}\left\vert q_{n-j}-q_{n-j-1}\right\vert +q_{0}\right)
\right) \sum_{i=0}^{\left\vert n\right\vert }M_{i}\left\vert
K_{M_{i}}\right\vert
\end{equation*}%
\begin{equation*}
=\left( \frac{1}{Q_{n}}\left( \overset{n-1}{\underset{j=1}{\sum }}\left(
q_{n-j}-q_{n-j-1}\right) +q_{0}\right) \right) \sum_{i=0}^{\left\vert
n\right\vert }M_{i}\left\vert K_{M_{i}}\right\vert
\end{equation*}%
\begin{equation*}
\leq \frac{q_{n-1}}{Q_{n}}\sum_{i=0}^{\left\vert n\right\vert
}M_{i}\left\vert K_{M_{i}}\right\vert \leq \frac{c}{n}\sum_{i=0}^{\left\vert
n\right\vert }M_{i}\left\vert K_{M_{i}}\right\vert .
\end{equation*}

The proof is complete by just combining the estimates above.
\QED

\begin{corollary}
\label{Corollary3nn} Let $\{q_{k}:k\in \mathbb{N}\}$ be a sequence of
non-decreasing numbers. \textit{Then}
\begin{equation*}
\sup_{n}\int_{G_{m}}\left\vert F_{n}\right\vert d\mu \leq c<\infty ,
\end{equation*}%
where $c$ is an absolute constant.
\end{corollary}

{\bf Proof}:
If we apply (\ref{fn4}) in Corollary \ref{lemma7kn} and invoke the
identities (\ref{2b}) and (\ref{2bb}) we readily get the proof.
So, we leave out the details.
\QED

\begin{lemma}
\label{lemma00nn} Let $n\geq M_{N}$ and $\{q_{k}:k\in \mathbb{N}\}$ be a
sequence of non-decreasing numbers. \textit{Then}
\begin{equation*}
\left\vert \frac{1}{Q_{n}}\overset{n}{\underset{j=M_{N}}{\sum }}%
q_{n-j}D_{j}\right\vert \leq \frac{c}{M_{N}}\left\{ \sum_{j=0}^{\left\vert
n\right\vert }M_{j}\left\vert K_{M_{j}}\right\vert \right\} ,
\end{equation*}%
where $c$ is an absolute constant.
\end{lemma}

{\bf Proof}:
Let $M_{N}\leq j\leq n.$ By using (\ref{fn5}) in Corollary \ref{lemma7kn} we
get that%
\begin{equation*}
\left\vert K_{j}\right\vert \leq \frac{1}{j}\sum_{l=0}^{\left\vert
j\right\vert }M_{l}\left\vert K_{M_{l}}\right\vert \leq \frac{1}{M_{N}}%
\sum_{l=0}^{\left\vert n\right\vert }M_{l}\left\vert K_{M_{l}}\right\vert .
\end{equation*}

Let the sequence $\{q_{k}:k\in \mathbb{N}\}$ be non-decreasing. Then%
\begin{equation*}
\overset{n-1}{\underset{j=M_{N}}{\sum }}\left\vert
q_{n-j}-q_{n-j-1}\right\vert j+q_{0}n\leq \overset{n-1}{\underset{j=0}{\sum }%
}\left\vert q_{n-j}-q_{n-j-1}\right\vert j+q_{0}n
\end{equation*}%
\begin{equation*}
=\overset{n-1}{\underset{j=1}{\sum }}\left( q_{n-j}-q_{n-j-1}\right)
j+q_{0}n=Q_{n}.
\end{equation*}

By using Abel transformation we can write that%
\begin{equation*}
\left\vert \frac{1}{Q_{n}}\overset{n}{\underset{j=M_{N}}{\sum }}%
q_{n-j}D_{j}\right\vert
\end{equation*}%
\begin{equation*}
=\left\vert \frac{1}{Q_{n}}\left( \overset{n-1}{\underset{j=M_{N}}{\sum }}%
\left( q_{n-j}-q_{n-j-1}\right) jK_{j}+q_{0}nK_{n}\right) \right\vert
\end{equation*}%
\begin{equation*}
\left( \frac{1}{Q_{n}}\left( \overset{n-1}{\underset{j=M_{N}}{\sum }}%
\left\vert q_{n-j}-q_{n-j-1}\right\vert j+q_{0}n\right) \right) \frac{1}{%
M_{N}}\sum_{i=0}^{\left\vert n\right\vert }M_{i}\left\vert
K_{M_{i}}\right\vert
\end{equation*}%
\begin{equation*}
\leq \frac{1}{M_{N}}\sum_{i=0}^{\left\vert n\right\vert }M_{i}\left\vert
K_{M_{i}}\right\vert .
\end{equation*}%
The proof is complete.
\QED

\begin{lemma}
\label{lemma5a}\ Let $x\in I_{N}^{k,l},$ $k=0,\dots ,N-2,$ $l=k+1,\dots ,N-1$
and $\{q_{k}:k\in \mathbb{N}\}$ be a sequence of non-decreasing numbers,
satisfying condition (\ref{fn01}). Then%
\begin{equation*}
\int_{I_{N}}\left\vert F_{n}\left( x-t\right) \right\vert d\mu \left(
t\right) \leq \frac{cM_{l}M_{k}}{nM_{N}}.
\end{equation*}%
Let $x\in I_{N}^{k,N},$ $k=0,\dots ,N-1.$ Then%
\begin{equation*}
\int_{I_{N}}\left\vert F_{n}\left( x-t\right) \right\vert d\mu \left(
t\right) \leq \frac{cM_{k}}{M_{N}}.
\end{equation*}%
Here $c$ is an absolute constant.
\end{lemma}

{\bf Proof}:
Let $x\in I_{N}^{k,l},$ for $0\leq k<l\leq N-1$ and $t\in I_{N}.$ First, we observe that $x-t\in $ $I_{N}^{k,l}.$ Next, we apply Lemma \ref%
{lemma0nn} and invoke (\ref{star1}) and (\ref{star2}) in Lemma \ref{lemma222}
to obtain that%
\begin{equation}
\int_{I_{N}}\left\vert F_{n}\left( x-t\right) \right\vert d\mu \left(
t\right)  \label{8881}
\end{equation}%
\begin{equation*}
\leq \frac{c}{n}\underset{i=0}{\overset{\left\vert n\right\vert }{\sum }}%
M_{i}\int_{I_{N}}\left\vert K_{M_{i}}\left( x-t\right) \right\vert d\mu
\left( t\right)
\end{equation*}%
\begin{equation*}
\leq \frac{c}{n}\int_{I_{N}}\overset{l}{\underset{i=0}{\sum }}M_{i}M_{k}d\mu
\left( t\right) \leq \frac{cM_{k}M_{l}}{nM_{N}}
\end{equation*}%
and the first estimate is proved.

Now, let $x\in I_{N}^{k,N}$. Since $x-t\in I_{N}^{k,N}$ for $t\in I_{N},$ by
combining Lemmas \ref{lemma2} and \ref{lemma0nn} with (\ref{star1}) and (\ref%
{star2}) in \ref{lemma222}, we have that%
\begin{equation}
\int_{I_{N}}\left\vert F_{n}\left( x-t\right) \right\vert d\mu \left(
t\right)  \label{111100}
\end{equation}%
\begin{equation*}
\leq \frac{c}{n}\underset{i=0}{\overset{\left\vert n\right\vert }{\sum }}%
M_{i}\int_{I_{N}}\left\vert K_{M_{i}}\left( x-t\right) \right\vert d\mu
\left( t\right)
\end{equation*}%
\begin{equation*}
\leq \frac{c}{n}\overset{\left\vert n\right\vert -1}{\underset{i=0}{\sum }}%
M_{i}\int_{I_{N}}M_{k}d\mu \left( t\right) \leq \frac{cM_{k}}{M_{N}}.
\end{equation*}

By combining (\ref{8881}) and (\ref{111100}) we complete the proof.
\QED

\begin{lemma}
\label{lemma5b}\ Let $n\geq M_{N},$ $\ x\in I_{N}^{k,l},$ $k=0,\dots ,N-1,$ $%
l=k+1,\dots ,N$ and $\{q_{k}:k\in \mathbb{N}\}$ be a sequence of
non-decreasing sequence, satisfying condition (\ref{fn01}). Then%
\begin{equation*}
\int_{I_{N}}\left\vert F_{n}\left( x-t\right) \right\vert d\mu \left(
t\right) \leq \frac{cM_{l}M_{k}}{M_{N}^{2}},
\end{equation*}%
where $c$ is an absolute constant.
\end{lemma}

{\bf Proof}:
Since $n\geq M_{N}$ if we apply Lemma \ref{lemma5a} we immediately get the  proof.
\QED

Next, we state analogical estimate, but now without any restriction like (%
\ref{fn01}):

\begin{lemma}
\label{lemma5aaa}\ Let $x\in I_{N}^{k,l},$ $k=0,\dots ,N-1,$ $l=k+1,\dots ,N$
and $\{q_{k}:k\in \mathbb{N}\}$ be a sequence of non-decreasing sequence.
Then%
\begin{equation*}
\int_{I_{N}}\left\vert \frac{1}{Q_{n}}\overset{n}{\underset{j=M_{N}}{\sum }}%
q_{n-j}D_{j}\left( x-t\right) \right\vert d\mu \left( t\right) \leq \frac{%
cM_{l}M_{k}}{M_{N}^{2}},
\end{equation*}%
where $c$ is an absolute constant.
\end{lemma}

{\bf Proof}:
Let $x\in I_{N}^{k,l},$ for $0\leq k<l\leq N-1$ and $t\in I_{N}.$ Since $%
x-t\in $ $I_{N}^{k,l}$ and $n\geq M_{N},$ if we combine Lemmas \ref{lemma2} and %
\ref{lemma00nn} and invoke (\ref{star1}) and (\ref{star2}) in Lemma \ref%
{lemma222} we readily obtain that%
\begin{equation}
\int_{I_{N}}\left\vert \frac{1}{Q_{n}}\overset{n}{\underset{j=M_{N}}{\sum }}%
q_{n-j}D_{j}\left( x-t\right) \right\vert d\mu \left( t\right)  \label{8881a}
\end{equation}%
\begin{equation*}
\leq \frac{c}{M_{N}}\underset{i=0}{\overset{\left\vert n\right\vert }{\sum }}%
M_{i}\int_{I_{N}}\left\vert K_{M_{i}}\left( x-t\right) \right\vert d\mu
\left( t\right)
\end{equation*}%
\begin{equation*}
\leq \frac{c}{M_{N}}\overset{l}{\underset{i=0}{\sum }}M_{i}%
\int_{I_{N}}M_{k}d\mu \left( t\right) \leq \frac{cM_{k}M_{l}}{M_{N}^{2}}
\end{equation*}%
and the first estimate is proved.

Now, let $x\in I_{N}^{k,N}$. Since $x-t\in I_{N}^{k,N}$ for $t\in I_{N},$ if
we apply again Lemmas \ref{lemma2}, \ref{lemma0nn} and (\ref{star1}) and (%
\ref{star2}) in Lemma \ref{lemma222} we can conclude that%
\begin{equation}
\int_{I_{N}}\left\vert \frac{1}{Q_{n}}\overset{n}{\underset{j=M_{N}}{\sum }}%
q_{n-j}D_{j}\left( x-t\right) \right\vert d\mu \left( t\right)
\label{111100a}
\end{equation}%
\begin{equation*}
\leq \frac{c}{M_{N}}\underset{i=0}{\overset{\left\vert n\right\vert }{\sum }}%
M_{i}\int_{I_{N}}\left\vert K_{M_{i}}\left( x-t\right) \right\vert d\mu
\left( t\right)
\end{equation*}%
\begin{equation*}
\leq \frac{c}{M_{N}}\overset{\left\vert n\right\vert -1}{\underset{i=0}{\sum
}}M_{i}\int_{I_{N}}M_{k}d\mu \left( t\right) \leq \frac{cM_{k}}{M_{N}}.
\end{equation*}

By combining (\ref{8881a}) and (\ref{111100a}) we complete the proof.
\QED

The next results (Lemmas \ref{lemma1n}-\ref{lemma3n}) are due to
Blahota, Persson and Tephnadze \cite{bpt1}:

\begin{lemma}
\label{lemma1n}Let $s_{n}M_{n}<r\leq \left( s_{n}+1\right) M_{n},$ where $%
1\leq s_{n}\leq m_{n}-1.$ Then
\begin{equation*}
Q_{r}F_{r}=Q_{r}D_{s_{n}M_{n}}-\psi _{s_{n}M_{n}-1}\overset{s_{n}M_{n}-2}{%
\underset{l=1}{\sum }}\left( q_{r-s_{n}M_{n}+l}-q_{r-s_{n}M_{n}+l+1}\right) l%
\overline{K_{l}}
\end{equation*}%
\begin{equation*}
-\psi _{s_{n}M_{n}-1}\left( s_{n}M_{n}-1\right) q_{r-1}\overline{%
K_{s_{n}M_{n}-1}}+\psi _{s_{n}M_{n}}Q_{r-s_{n}M_{n}}F_{r-s_{n}M_{n}}.
\end{equation*}
\end{lemma}

\begin{remark}
We note that Lemma \ref{lemma1n} is true for every N\"orlund mean, without
any restriction on the generative sequence $\{q_{k}:k\in \mathbb{N}\}.$
\end{remark}

{\bf Proof}:
Let $s_{n}M_{n}<r\leq \left( s_{n}+1\right) M_{n},$ where $1\leq s_{n}\leq
m_{n}-1.$ It is easy to see that
\begin{equation}
Q_{r}F_{r}=\overset{r}{\underset{k=1}{\sum }}q_{r-k}D_{k}=\overset{s_{n}M_{n}%
}{\underset{l=1}{\sum }}q_{r-l}D_{l}+\overset{r}{\underset{l=s_{n}M_{n}+1}{%
\sum }}q_{r-l}D_{l}:=I+II.  \label{nor1}
\end{equation}

We apply (\ref{dn22}) in Lemma \ref{dn1} and invoke Abel transformation to
obtain that
\begin{equation}
I=\overset{s_{n}M_{n}-1}{\underset{l=0}{\sum }}%
q_{r-s_{n}M_{n}+l}D_{s_{n}M_{n}-l}  \label{nor2}
\end{equation}%
\begin{equation*}
=\overset{s_{n}M_{n}-1}{\underset{l=1}{\sum }}%
q_{r-s_{n}M_{n}+l}D_{s_{n}M_{n}-l}+q_{r-s_{n}M_{n}}D_{s_{n}M_{n}}
\end{equation*}%
\begin{equation*}
=D_{s_{n}M_{n}}\overset{s_{n}M_{n}-1}{\underset{l=0}{\sum }}%
q_{r-s_{n}M_{n}+l}-\psi _{s_{n}M_{n}-1}\overset{s_{n}M_{n}-1}{\underset{l=1}{%
\sum }}q_{r-s_{n}M_{n}+l}\overline{D_{l}}
\end{equation*}%
\begin{equation}
=\left( Q_{r}-Q_{r-s_{n}M_{n}}\right) D_{s_{n}M_{n}}  \notag
\end{equation}%
\begin{equation*}
-\psi _{s_{n}M_{n}-1}\overset{s_{n}M_{n}-2}{\underset{l=1}{\sum }}\left(
q_{r-s_{n}M_{n}+l}-q_{r-s_{n}M_{n}+l+1}\right) l\overline{K_{l}}
\end{equation*}%
\begin{equation*}
-\psi _{s_{n}M_{n}-1}q_{r-1}\left( s_{n}M_{n}-1\right) \overline{%
K_{s_{n}M_{n}-1}}.
\end{equation*}

By using (\ref{dn21}) in Lemma \ref{dn1} we can rewrite $II$ as
\begin{equation}
II=\overset{r-s_{n}M_{n}}{\underset{l=1}{\sum }}%
q_{r-s_{n}M_{n}-l}D_{l+s_{n}M_{n}}  \label{nor3}
\end{equation}%
\begin{equation*}
=Q_{r-s_{n}M_{n}}D_{s_{n}M_{n}}+\psi
_{s_{n}M_{n}}Q_{r-s_{n}M_{n}}F_{r-s_{n}M_{n}}.
\end{equation*}

The proof is complete by just combining (\ref{nor1}-\ref{nor3}).
\QED

\begin{lemma}
\label{lemma0nn1}Let $\{q_{k}:k\in \mathbb{N}\}$ be a sequence of
non-increasing numbers satisfying the condition
\begin{equation}
\frac{1}{Q_{n}}=O\left( \frac{1}{n}\right) ,\text{ \ \ when \ \ }\
n\rightarrow \infty .  \label{fn0}
\end{equation}%
\textit{Then}
\begin{equation*}
\left\vert F_{n}\right\vert \leq \frac{c_{\alpha }}{n}\left\{
\sum_{j=0}^{\left\vert n\right\vert }M_{j}\left\vert K_{M_{j}}\right\vert
\right\} ,
\end{equation*}%
where $c$ is an absolute constant.
\end{lemma}

{\bf Proof}:
Let the sequence $\{q_{k}:k\in \mathbb{N}\}$ be non-increasing satisfying
condition (\ref{fn0}). Then
\begin{equation*}
\frac{1}{Q_{n}}\left( \overset{n-1}{\underset{j=1}{\sum }}\left\vert
q_{n-j}-q_{n-j-1}\right\vert +q_{0}\right)
\end{equation*}%
\begin{equation*}
\leq \frac{1}{Q_{n}}\left( \overset{n-1}{\underset{j=1}{\sum }}-\left(
q_{n-j}-q_{n-j-1}\right) +q_{0}\right)
\end{equation*}%
\begin{equation*}
\leq \frac{2q_{0}-q_{n-1}}{Q_{n}}\leq \frac{2q_{0}}{Q_{n}}\leq \frac{c}{n}.
\end{equation*}

If we apply (\ref{fn5}) in Corollary \ref{lemma7kn} and invoke equalities (%
\ref{2b}) and (\ref{2bb}) we immediately get that%
\begin{equation*}
\left\vert F_{n}\right\vert \leq \left( \frac{1}{Q_{n}}\left( \overset{n-1}{%
\underset{j=1}{\sum }}\left\vert q_{n-j}-q_{n-j-1}\right\vert +q_{0}\right)
\right) \sum_{i=0}^{\left\vert n\right\vert }M_{i}\left\vert
K_{M_{i}}\right\vert
\end{equation*}%
\begin{equation*}
=\left( \frac{1}{Q_{n}}\left( \overset{n-1}{\underset{j=1}{\sum }}-\left(
q_{n-j}-q_{n-j-1}\right) +q_{0}\right) \right) \sum_{i=0}^{\left\vert
n\right\vert }M_{i}\left\vert K_{M_{i}}\right\vert
\end{equation*}%
\begin{equation*}
\leq \frac{2q_{0}-q_{n-1}}{Q_{n}}\sum_{i=0}^{\left\vert n\right\vert
}M_{i}\left\vert K_{M_{i}}\right\vert \leq \frac{2q_{0}}{Q_{n}}%
\sum_{i=0}^{\left\vert n\right\vert }M_{i}\left\vert K_{M_{i}}\right\vert
\end{equation*}%
\begin{equation*}
\leq \frac{c}{n}\sum_{i=0}^{\left\vert n\right\vert }M_{i}\left\vert
K_{M_{i}}\right\vert .
\end{equation*}

The proof is complete by combining the estimates above.
\QED

\begin{corollary}
\label{corollary3n9} Let $\{q_{k}:k\in \mathbb{N}\}$ be a sequence of
non-increasing numbers satisfying condition (\ref{fn0}). Then%
\begin{equation*}
\sup_{n}\int_{G_{m}}\left\vert F_{n}\right\vert d\mu \leq c_{\alpha }<\infty
,
\end{equation*}%
where $c_{\alpha }$ is an absolute constant depending only on $\alpha .$
\end{corollary}

{\bf Proof}:
By applying Lemma \ref{lemma0nn1} we readily get the proof. So,
we leave out the details.
\QED

\begin{lemma}
\label{lemma2n}Let $\{q_{k}:k\in \mathbb{N}\}$ be a sequence of
non-increasing numbers, $0<\alpha <1,\ \ $and $\ $%
\begin{equation}
\frac{1}{Q_{n}}=O\left( \frac{1}{n^{\alpha }}\right) ,\text{ \ \ when }\
n\rightarrow \infty  \label{6a}
\end{equation}%
and
\begin{equation}
q_{n}-q_{n+1}=O\left( \frac{1}{n^{2-\alpha }}\right) ,\text{ \ when }\
n\rightarrow \infty .  \label{7a}
\end{equation}%
\textit{Then}
\begin{equation*}
\left\vert F_{n}\right\vert \leq \frac{c_{\alpha }}{n^{\alpha }}\left\{
\sum_{j=0}^{\left\vert n\right\vert }M_{j}^{\alpha }\left\vert
K_{M_{j}}\right\vert \right\} ,
\end{equation*}%
where $c_{\alpha }$ is an absolute constant depending only on $\alpha .$
\end{lemma}

{\bf Proof}:
According to the fact that $\{q_{k}:k\in \mathbb{N}\}$ is a sequence of non-negative and
non-increasing numbers we have two cases:
\begin{equation*}
\begin{array}{ll}
1. & \lim_{k\rightarrow \infty }q_{k}\geq c>0, \\
2. & \lim_{k\rightarrow \infty }q_{k}=0,%
\end{array}%
\end{equation*}

In the first case we obtain that (\ref{monotone0}) and (\ref{monotone1}) are
satisfied. Since the case
\begin{equation*}
\frac{n}{Q_{n}}=O\left( 1\right) ,\ \text{when }\ n\rightarrow \infty ,
\end{equation*}%
have already been considered in Lemma \ref{lemma0nn}, we can exclude it.

Hence, we may assume that
\begin{equation}
q_{n}=o(1),\text{ \ \ \ \ when \ }n\rightarrow \infty .  \label{102a}
\end{equation}%
By combining (\ref{7a}) and (\ref{102a}) we immediately get that
\begin{equation*}
q_{n}=\overset{\infty }{\underset{l=n}{\sum }}\left( q_{l}-q_{l+1}\right)
\leq \overset{\infty }{\underset{l=n}{\sum }}\frac{c}{l^{2-\alpha }}\leq
\frac{c}{n^{1-\alpha }}
\end{equation*}%
and%
\begin{equation*}
Q_{n}=\overset{n-1}{\underset{l=0}{\sum }}q_{l}\leq \overset{n}{\underset{l=1%
}{\sum }}\frac{c}{l^{1-\alpha }}\leq cn^{\alpha }.
\end{equation*}%
Let $M_{n}<k\leq M_{n+1}.$ It is easy to see that
\begin{equation}
Q_{k}\left\vert D_{sM_{n}}\right\vert \leq cM_{n}^{\alpha }\left\vert
D_{sM_{n}}\right\vert  \label{a1}
\end{equation}%
and
\begin{equation}
\left( sM_{n}-1\right) q_{k-1}\left\vert K_{sM_{n}-1}\right\vert  \label{a2}
\end{equation}%
\begin{equation*}
\leq ck^{\alpha -1}M_{n}\left\vert K_{sM_{n}-1}\right\vert \leq
cM_{n}^{\alpha }\left\vert K_{sM_{n}-1}\right\vert .
\end{equation*}

Let
\begin{equation*}
n=s_{n_{1}}M_{n_{1}}+s_{n_{2}}M_{n_{2}}+\dots +s_{n_{r}}M_{n_{r}},\ \
n_{1}>n_{2}>\dots >n_{r},
\end{equation*}%
and
\begin{equation*}
n^{\left( k\right) }=s_{n_{k+1}}M_{n_{k+1}}+\dots +s_{n_{r}}M_{n_{r}},\
1\leq s_{n_{l}}\leq m_{l}-1,\ l=1,\dots ,r.
\end{equation*}

By combining (\ref{a1}), (\ref{a2}) and Lemma \ref{lemma1n} we have that%
\begin{equation*}
\left\vert Q_{n}F_{n}\right\vert
\end{equation*}%
\begin{equation*}
\leq c_{\alpha }\left( M_{n_{1}}^{\alpha }\left\vert
D_{s_{n_{1}}M_{n_{1}}}\right\vert +\overset{s_{n_{1}}M_{n_{1}}-1}{\underset{%
l=1}{\sum }}\left\vert \left( n^{\left( 1\right) }+l\right) ^{\alpha
-2}\right\vert \left\vert lK_{l}\right\vert +M_{n_{1}}^{\alpha }\left\vert
K_{s_{n_{1}}M_{n_{1}}-1}\right\vert \right).
\end{equation*}%
\begin{equation*}
+c_{\alpha }\left\vert Q_{n^{\left( 1\right) }}F_{n^{\left( 1\right)
}}\right\vert
\end{equation*}

By repeating this process $r$ times we get that
\begin{equation*}
\left\vert Q_{n}F_{n}\right\vert
\end{equation*}%
\begin{equation*}
\leq c_{\alpha }\overset{r}{\underset{k=1}{\sum }}\left( M_{n_{k}}^{\alpha
}\left\vert D_{s_{n_{k}}M_{n_{k}}}\right\vert +\overset{s_{n_{k}}M_{n_{k}}-1}%
{\underset{l=1}{\sum }}\left( n^{\left( k\right) }+l\right) ^{\alpha
-2}\left\vert lK_{l}\right\vert +M_{n_{k}}^{\alpha }\left\vert
K_{s_{n_{k}}M_{n_{k}}-1}\right\vert \right)
\end{equation*}%
\begin{equation*}
:=I+II+III.
\end{equation*}

We\ combine Corollary \ref{dn2.3} and Lemma \ref{lemma2} and invoke (\ref%
{9dn}) in Lemma \ref{dn2.1} to obtain that%
\begin{equation*}
I\leq c_{\alpha }\overset{\left\vert n\right\vert }{\underset{k=1}{\sum }}%
M_{k}^{\alpha }\left\vert D_{s_{k}M_{k}}\right\vert
\end{equation*}%
\begin{equation*}
\leq c_{\alpha }\overset{\left\vert n\right\vert }{\underset{k=1}{\sum }}%
M_{k}^{\alpha }\left\vert D_{M_{k}}\right\vert \leq c_{\alpha }\overset{%
\left\vert n\right\vert }{\underset{k=1}{\sum }}M_{k}^{\alpha }\left\vert
K_{M_{k}}\right\vert
\end{equation*}

and
\begin{equation*}
III\leq c_{\alpha }\overset{r}{\underset{k=1}{\sum }}M_{n_{k}}^{\alpha
-1}\left\vert
M_{n_{k}}K_{s_{n_{k}}M_{n_{k}}}-M_{n_{k}}D_{s_{n_{k}}M_{n_{k}}}\right\vert
\end{equation*}%
\begin{equation*}
\leq c_{\alpha }\overset{r}{\underset{k=1}{\sum }}M_{n_{k}}^{\alpha
}\left\vert K_{s_{n_{k}}M_{n_{k}}}\right\vert +c_{\alpha }\overset{r}{%
\underset{k=1}{\sum }}M_{n_{k}}^{\alpha }\left\vert
D_{s_{n_{k}}M_{n_{k}}}\right\vert
\end{equation*}%
\begin{equation*}
\leq c_{\alpha }\overset{r}{\underset{k=1}{\sum }}M_{k}^{\alpha }\left\vert
K_{M_{k}}\right\vert .
\end{equation*}

Next, we can rewrite $II$ as
\begin{equation*}
II=c_{\alpha }\overset{r}{\underset{k=1}{\sum }}\overset{%
s_{n_{k+1}}M_{n_{k+1}}-1}{\underset{l=1}{\sum }}\left( n^{\left( k\right)
}+l\right) ^{\alpha -2}\left\vert lK_{l}\right\vert
\end{equation*}%
\begin{equation*}
+c_{\alpha }\overset{r}{\underset{k=1}{\sum }}\overset{s_{n_{k}}M_{n_{k}}-1}{%
\underset{l=s_{n_{k+1}}M_{n_{k+1}}}{\sum }}\left( n^{\left( k\right)
}+l\right) ^{\alpha -2}\left\vert lK_{l}\right\vert
\end{equation*}%
\begin{equation*}
:=II_{1}+II_{2}.
\end{equation*}

For $II_{1}$ we find that%
\begin{equation*}
II_{1}\leq c_{\alpha }\overset{r}{\underset{k=1}{\sum }}s_{n_{k+1}}^{\alpha
-2}M_{n_{k+1}}^{\alpha -2}\overset{s_{n_{k+1}}M_{n_{k+1}}-1}{\underset{l=1}{%
\sum }}\left\vert lK_{l}\right\vert
\end{equation*}%
\begin{equation*}
\leq c_{\alpha }\overset{n_{1}}{\underset{k=1}{\sum }}M_{k}^{\alpha -2}%
\overset{M_{k}-1}{\underset{l=1}{\sum }}\left\vert lK_{l}\right\vert
\end{equation*}%
\begin{equation*}
=c_{\alpha }\overset{n_{1}}{\underset{k=1}{\sum }}M_{k}^{\alpha -2}\overset{k%
}{\underset{i=1}{\sum }}\overset{M_{i}-1}{\underset{l=M_{i-1}}{\sum }}%
\left\vert lK_{l}\right\vert
\end{equation*}%
\begin{equation*}
\leq c_{\alpha }\overset{n_{1}}{\underset{k=1}{\sum }}M_{k}^{\alpha -2}%
\overset{k}{\underset{i=1}{\sum }}M_{i}\overset{i}{\underset{j=0}{\sum }}%
M_{j}\left\vert K_{M_{j}}\right\vert
\end{equation*}%
\begin{equation*}
\leq c_{\alpha }\overset{n_{1}}{\underset{k=0}{\sum }}M_{k}^{\alpha -1}%
\overset{k}{\underset{j=0}{\sum }}M_{j}|K_{M_{j}}|
\end{equation*}%
\begin{equation*}
=c_{\alpha }\overset{n_{1}}{\underset{j=0}{\sum }}M_{j}|K_{M_{j}}|\overset{%
n_{1}}{\underset{k=j}{\sum }}M_{k}^{\alpha -1}
\end{equation*}%
\begin{equation*}
\leq c_{\alpha }\overset{n_{1}}{\underset{j=0}{\sum }}M_{j}^{\alpha
}\left\vert K_{M_{j}}\right\vert .
\end{equation*}

Moreover,%
\begin{equation*}
II_{2}\leq c_{\alpha }\overset{r}{\underset{k=1}{\sum }}\overset{%
s_{n_{k}}M_{n_{k}}-1}{\underset{l=s_{n_{k+1}}M_{n_{k+1}}}{\sum }}l^{\alpha
-2}\left\vert lK_{l}\right\vert
\end{equation*}%
\begin{equation*}
\leq c_{\alpha }\overset{r}{\underset{k=1}{\sum }}\overset{n_{k}}{\underset{%
i=n_{k+1}+1}{\sum }}M_{i}^{\alpha -2}\overset{M_{i+1}-1}{\underset{l=M_{i}}{%
\sum }}\left\vert lK_{l}\right\vert
\end{equation*}%
\begin{equation*}
\leq c_{\alpha }\overset{r}{\underset{k=1}{\sum }}\overset{n_{k}}{\underset{%
i=n_{k+1}+1}{\sum }}M_{i}^{\alpha -2}M_{i}\overset{i}{\underset{j=0}{\sum }}%
M_{j}\left\vert K_{M_{j}}\right\vert
\end{equation*}%
\begin{equation*}
=c_{\alpha }\overset{r}{\underset{k=1}{\sum }}\overset{n_{k}}{\underset{%
i=n_{k+1}+1}{\sum }}M_{i}^{\alpha -1}\overset{i}{\underset{j=0}{\sum }}%
M_{j}\left\vert K_{M_{j}}\right\vert
\end{equation*}%
\begin{equation*}
\leq c_{\alpha }\overset{n_{1}}{\underset{i=1}{\sum }}M_{i}^{\alpha -1}%
\overset{i}{\underset{j=0}{\sum }}M_{j}\left\vert K_{M_{j}}\right\vert
\end{equation*}%
\begin{equation*}
\leq c_{\alpha }\overset{n_{1}}{\underset{j=0}{\sum }}M_{j}^{\alpha
}\left\vert K_{M_{j}}\right\vert .
\end{equation*}

The proof is complete by just combining the estimates above.
\QED

\begin{corollary}
\label{corollary3n}Let $0<\alpha \leq 1$ and $\{q_{k}:k\in \mathbb{N}\}$ be
a sequence of non-increasing numbers satisfying the conditions (\ref{6a}) and (%
\ref{7a}). Then%
\begin{equation*}
\sup_{n}\int_{G_{m}}\left\vert F_{n}\right\vert d\mu \leq c_{\alpha }<\infty
,
\end{equation*}%
where $c_{\alpha }$ is an absolute constant depending only on $\alpha .$
\end{corollary}

{\bf Proof}:
If we use Lemma \ref{lemma2n} we readily get the proof. Thus, we leave out the details.
\QED

\begin{lemma}
\label{lemma3n}Let $0<\alpha \leq 1$ and $\{q_{k}:k\in \mathbb{N}\}$ be a
sequence of non-increasing numbers satisfying the conditions (\ref{6a}) and (%
\ref{7a}). Then
\begin{equation*}
\int_{I_{N}}\left\vert F_{m}\left( x-t\right) \right\vert d\mu \left(
t\right) \leq \frac{c_{\alpha }M_{l}^{\alpha }M_{k}}{m^{\alpha }M_{N}},\text{
\ for }x\in I_{N}^{k,l},
\end{equation*}%
where \ \ $k=0,\dots ,N-2,$ $l=k+2,\dots ,N-1$. Moreover,
\begin{equation*}
\int_{I_{N}}\left\vert F_{m}\left( x-t\right) \right\vert d\mu \left(
t\right) \leq \frac{c_{\alpha }M_{k}}{M_{N}},\text{ \ \textit{for} \ }x\in
I_{N}^{k,N},\text{ \ }k=0,\dots ,N-1.
\end{equation*}%
Here $c_{\alpha }$ is an absolute constant depending only on $\alpha .$
\end{lemma}

{\bf Proof}:
Let $x\in I_{N}^{k,l},$ where \ \ $k=0,\dots ,N-2,$ $l=k+2,\dots ,N-1.$
Since $x-t\in I_{N}^{k,l},$ for $t\in I_{N}$ \ if we apply Lemma \ref%
{lemma2n} with (\ref{star1}) and (\ref{star2}) in Lemma \ref{lemma222} we
can conclude that%
\begin{equation}
\left\vert F_{m}\left( x-t\right) \right\vert \leq \frac{c_{\alpha }}{%
m^{\alpha }}\overset{l}{\underset{i=0}{\sum }}M_{i}^{\alpha }\left\vert
K_{M_{i}}\left( x-t\right) \right\vert  \label{7n}
\end{equation}%
\begin{equation*}
\leq \frac{c_{\alpha }}{m^{\alpha }}\overset{l}{\underset{i=0}{\sum }}%
M_{i}^{\alpha }M_{k}\leq \frac{c_{\alpha }M_{l}^{\alpha }M_{k}}{m^{\alpha }}.
\end{equation*}

Let $x\in I_{N}^{k,l},$ for some $0\leq k<l\leq N-1.$ Since $x-t\in
I_{N}^{k,l},$ for $t\in $ $I_{N}$ and $m\geq M_{N}$ from (\ref{7n}) we
readily obtain that
\begin{equation}
\int_{I_{N}}\left\vert F_{m}\left( x-t\right) \right\vert d\mu \left(
t\right) \leq \frac{c_{\alpha }M_{l}^{\alpha }M_{k}}{m^{\alpha }M_{N}}.
\label{8n}
\end{equation}

Let $x\in I_{N}^{k,N},$ $k=0,\dots ,N-1.$ Then, by applying Lemma \ref%
{lemma2n}, we have that%
\begin{equation}
\int_{I_{N}}\left\vert F_{m}\left( x-t\right) \right\vert d\mu \left(
t\right)  \label{star10}
\end{equation}%
\begin{equation*}
\leq \frac{c_{\alpha }}{m^{\alpha }}\underset{i=0}{\overset{\left\vert
m\right\vert }{\sum }}M_{i}^{\alpha }\int_{I_{N}}\left\vert K_{M_{i}}\left(
x-t\right) \right\vert d\mu \left( t\right) .
\end{equation*}

Let $x\in I_{N}^{k,N},$ $k=0,\dots ,N-1,$ $t\in $ $I_{N}$ \ and $x_{q}\neq
t_{q},$ where $N\leq q\leq \left\vert m\right\vert -1.$ By using Lemma \ref%
{lemma2} and estimate (\ref{star10}) we get that%
\begin{equation}
\int_{I_{N}}\left\vert F_{m}\left( x-t\right) \right\vert d\mu \left(
t\right)  \label{10n}
\end{equation}%
\begin{equation*}
\leq \frac{c_{\alpha }}{m^{\alpha }}\underset{i=0}{\overset{q-1}{\sum }}%
M_{i}^{\alpha }\int_{I_{N}}M_{k}d\mu \left( t\right) \leq \frac{c_{\alpha
}M_{k}M_{q}^{\alpha }}{m^{\alpha }M_{N}}\leq \frac{c_{\alpha }M_{k}}{M_{N}}.
\end{equation*}

Let $x\in I_{N}^{k,N},$ $k=0,\dots ,N-1,$ $t\in I_{N}$ and $%
x_{N}=t_{N},\dots ,x_{\left\vert m\right\vert -1}=t_{\left\vert m\right\vert
-1}.$ By again applying Lemma \ref{lemma2} and estimate (\ref{star10}) we
have that%
\begin{equation}
\int_{I_{N}}\left\vert F_{m}\left( x-t\right) \right\vert d\mu \left(
t\right) \leq \frac{c_{\alpha }}{m^{\alpha }}\overset{\left\vert
m\right\vert -1}{\underset{i=0}{\sum }}M_{i}^{\alpha }\int_{I_{N}}M_{k}d\mu
\left( t\right) \leq \frac{c_{\alpha }M_{k}}{M_{N}}.  \label{11n}
\end{equation}

The proof follows by combining (\ref{8n}), (\ref{10n}) and (\ref{11n}).
\QED

\begin{corollary}
\label{corollary3na}Let $0<\alpha <1,$ $m\geq M_{N}$ and $\{q_{k}:k\in
\mathbb{N}\}$ be a sequence of non-increasing numbers satisfying the conditions
(\ref{6a}) and (\ref{7a}). Then there exists an absolute constant $c_{\alpha
}$, depending only on $\alpha ,$ such that
\begin{equation*}
\int_{I_{N}}\left\vert F_{m}\left( x-t\right) \right\vert d\mu \left(
t\right) \leq \frac{c_{\alpha }M_{l}^{\alpha }M_{k}}{M_{N}^{1+\alpha }},%
\text{ \ for }x\in I_{N}^{k,l}
\end{equation*}%
where \ \ $k=0,\dots ,N-1,$ $l=k+2,\dots ,N$.
\end{corollary}

{\bf Proof}:
The proof readily follows Lemma \ref{lemma3n}
if we use additional condition $n\geq M_{N}$.
\QED

\begin{remark}
For some sequences $\{q_{k}:k\in \mathbb{N}\}$ of non-increasing numbers
conditions (\ref{6a}) and (\ref{7a}) can be true or false independently.
\end{remark}

\subsection{Introduction to the theory of Martingale Hardy spaces}

The $\sigma $-algebra generated by the intervals
\begin{equation*}
\left\{ I_{n}\left( x\right) :x\in G_{m}\right\}
\end{equation*}
will be denoted by $\digamma _{n}\left( n\in
\mathbb{N}
\right) .$

A sequence $f=\left( f^{\left( n\right) }:n\in \mathbb{N}\right) $ of
integrable functions $f^{\left( n\right) }$ is said to be a martingale with
respect to the $\sigma $-algebras $\digamma _{n}\left( n\in
\mathbb{N}
\right) $ if (for details see e.g. Weisz \cite{We1})

$1)$ \ $f_{n}$\ is $\digamma _{n}$\ measurable for all $n\in \mathbb{N},$

$2)$\ $S_{M_{n}}f_{m}=f_{n}$\ for all $n\leq m.$

The martingale $f=\left( f^{\left( n\right) },n\in \mathbb{N}\right) $ is
said to be $L_{p}$-bounded ($0<p\leq \infty $) if $f^{\left( n\right) }\in
L_{p}$ and

\begin{equation*}
\left\Vert f\right\Vert _{p}:=\underset{n\in \mathbb{N}}{\sup }\left\Vert
f_{n}\right\Vert _{p}<\infty .
\end{equation*}

If $f\in L_{1}\left( G_{m}\right) ,$ then it is easy to show that the
sequence $F=\left( E_{n}f:n\in \mathbb{N}\right) $ is a martingale. This
type of martingales is called regular. If $1\leq p\leq \infty $ and $f\in
L_{p}\left( G_{m}\right) $ then $f=\left( f^{\left( n\right) },n\in \mathbb{N%
}\right) $ is $L_{p}$-bounded and
\begin{equation*}
\underset{n\rightarrow \infty }{\lim }\left\Vert E_{n}f-f\right\Vert _{p}=0
\end{equation*}%
and consequently $\left\Vert F\right\Vert _{p}=\left\Vert f\right\Vert _{p}$
(see \cite{W}). The converse of the latest statement holds also if $%
1<p\leq \infty $ (see \cite{W}): for an arbitrary $L_{p}$-bounded martingale
$f=\left( f^{\left( n\right) },n\in \mathbb{N}\right) $ there exists a
function $f\in L_{p}\left( G_{m}\right) $ for which $f^{\left( n\right)
}=E_{n}f.$ If $p=1,$ then there exists a function $f\in L_{1}\left(
G_{m}\right)$ of the preceding type if and only if $f$ is uniformly
integrable (see \cite{W}), namely, if
\begin{equation*}
\underset{y\rightarrow \infty }{\lim }\underset{n\in \mathbb{N}}{\sup }%
\int_{\left\{ \left\vert f_{n}\right\vert >y\right\} }\left\vert f_{n}\left(
x\right) \right\vert d\mu \left( x\right) =0.
\end{equation*}

Thus the map $f\rightarrow f:=$ $\left( E_{n}f:n\in \mathbb{N}\right) $ is
isometric from $L_{p}$ onto the space of $L_{p}$-bounded martingales when $%
1<p\leq \infty .$ Consequently, these two spaces can be identified with each
other. Similarly, the space $L_{1}\left( G_{m}\right) $ can be identified
with the space of uniformly integrable martingales.

Analogously, the martingale $f=\left( f^{\left( n\right) },n\in \mathbb{N}%
\right) $ is said to be $weak-L_{p}$-bounded ($0<p\leq \infty $) if $%
f^{\left( n\right) }\in L_{p}$ and

\begin{equation*}
\left\Vert f\right\Vert _{weak-L_{p}}:=\underset{n\in \mathbb{N}}{\sup }%
\left\Vert f_{n}\right\Vert _{weak-L_{p}}<\infty .
\end{equation*}

The maximal function of a martingale $f$ \ is defined by
\begin{equation*}
f^{\ast }:=\sup_{n\in
\mathbb{N}
}\left\vert f^{(n)}\right\vert .
\end{equation*}

In the case $f\in L_{1}(G_{m}),$ the maximal functions are also given by
\begin{equation*}
f^{\ast }\left( x\right) :=\sup_{n\in \mathbb{N}}\frac{1}{\left\vert
I_{n}\left( x\right) \right\vert }\left\vert \int_{I_{n}\left( x\right)
}f\left( u\right) d\mu \left( u\right) \right\vert .
\end{equation*}

For $0<p<\infty $ \ the Hardy martingale spaces $H_{p}$ consist of all
martingales for which
\begin{equation*}
\left\Vert f\right\Vert _{H_{p}}:=\left\Vert f^{\ast }\right\Vert
_{p}<\infty .
\end{equation*}

If $f=\left( f^{\left( n\right) }:n\in
\mathbb{N}
\right) $ is a martingale, then the Vilenkin-Fourier coefficients must be
defined in a slightly different manner:
\begin{equation*}
\widehat{f}\left( i\right) :=\lim_{k\rightarrow \infty
}\int_{G_{m}}f^{\left( k\right) }\overline{\psi }_{i}d\mu .
\end{equation*}

The next Lemma can be found in \cite{We3} (see also book \cite{sws}):

\begin{lemma}
\label{lemma2.3.4}If $f\in L_{1},$ then the sequence $F:=\left(
S_{M_{n}}f:n\in \mathbb{N}\right) $ is a martingale and
\begin{equation*}
\left\Vert F\right\Vert _{H_{p}}\sim \left\Vert \sup_{n\in \mathbb{N}%
}\left\vert S_{M_{n}}f\right\vert \right\Vert _{p}.
\end{equation*}%
Moreover, if $F:=$ $\left( S_{M_{n}}f:n\in \mathbb{N}\right) $ is a regular
martingale generated by $f\in L_{1},$ then \qquad \qquad \qquad \qquad
\end{lemma}

\begin{equation*}
\widehat{F}\left( k\right) =\int_{G}f\left( x\right) \psi _{k}\left(
x\right) d\mu \left( x\right) =\widehat{f}\left( k\right) ,\text{ \qquad }%
k\in \mathbb{N}.
\end{equation*}

A bounded measurable function $a$ is a p-atom if there exist an interval $I$
such that \qquad
\begin{equation*}
\int_{I}ad\mu =0,\text{ \ \ }\left\Vert a\right\Vert _{\infty }\leq \mu
\left( I\right) ^{-1/p},\text{ \ \ supp}\left( a\right) \subset I.
\end{equation*}

Next, we note that the Hardy martingale spaces $H_{p}\left( G_{m}\right) $
for $0<p\leq 1$ have atomic characterizations (see e.g. Weisz \cite{We1,We3}%
):

\begin{lemma}
\label{lemma2.1} A martingale $f=\left( f^{\left( n\right) }:n\in \mathbb{N}%
\right) $ is in $H_{p}\left( 0<p\leq 1\right) $ if and only if there exist a
sequence $\left( a_{k},k\in \mathbb{N}\right) $ of p-atoms and a sequence $%
\left( \mu _{k}:k\in \mathbb{N}\right) $ of real numbers such that, for
every $n\in \mathbb{N},$%
\begin{equation}
\qquad \sum_{k=0}^{\infty }\mu _{k}S_{M_{n}}a_{k}=f^{\left( n\right) },\text{
\ \ a.e.,}  \label{condmart}
\end{equation}%
where
\begin{equation*}
\qquad \sum_{k=0}^{\infty }\left\vert \mu _{k}\right\vert ^{p}<\infty .
\end{equation*}%
Moreover,
\begin{equation*}
\left\Vert f\right\Vert _{H_{p}}\backsim \inf \left( \sum_{k=0}^{\infty
}\left\vert \mu _{k}\right\vert ^{p}\right) ^{1/p},
\end{equation*}%
where the infimum is taken over all decomposition of $f=\left( f^{\left(
n\right) }:n\in \mathbb{N}\right) $ of the form (\ref{condmart}).
\end{lemma}

By using atomic characterization it can be easily proved that the following
Lemmas hold (see e.g. Weisz \cite{We3}):

\begin{lemma}
\label{lemma2.2} Suppose that an operator $T$ is sub-linear and for some $%
0<p\leq 1$%
\begin{equation*}
\int\limits_{\overset{-}{I}}\left\vert Ta\right\vert ^{p}d\mu \leq
c_{p}<\infty
\end{equation*}%
\textit{for every }$p$\textit{-atom }$a$\textit{, where }$I$\textit{\
denotes the support of the atom. If }$T$\textit{\ is bounded from }$%
L_{\infty \text{ }}$\textit{\ to }$L_{\infty },$\ then\textit{\ }%
\begin{equation*}
\left\Vert Tf\right\Vert _{p}\leq c_{p}\left\Vert f\right\Vert _{H_{p}}.
\end{equation*}%
Moreover, if $p<1,$ then we have the weak (1,1) type estimate
\begin{equation*}
\lambda \mu \left\{ x\in G_{m}:\text{ }\left\vert Tf\left( x\right)
\right\vert >\lambda \right\} \leq \left\Vert f\right\Vert _{1}
\end{equation*}%
for all $f\in L_{1}.$
\end{lemma}

\begin{lemma}
\label{lemma2.3} Suppose that an operator $T$ is sub-linear and for some $%
0<p\leq 1$%
\begin{equation*}
\underset{\lambda >0}{\sup }\lambda ^{p}\mu \left\{ x\in \overset{-}{I}:%
\text{ }\left\vert Tf\right\vert >\lambda \right\} \leq c_{p}<+\infty
\end{equation*}%
\textit{for every }$p$\textit{-atom }$a$\textit{, where }$I$\textit{\ denote
the support of the atom. If }$T$\textit{\ is bounded from }$L_{\infty \text{
}}$\textit{\ to }$L_{\infty },$\textit{\ then}
\begin{equation*}
\left\Vert Tf\right\Vert _{weak-L_{p}}\leq c_{p}\left\Vert f\right\Vert
_{H_{p}}.
\end{equation*}%
Moreover, if $p<1,$ then
\begin{equation*}
\lambda \mu \left\{ x\in G_{m}:\text{ }\left\vert Tf\left( x\right)
\right\vert >\lambda \right\} \leq \left\Vert f\right\Vert _{1},
\end{equation*}%
for all $f\in L_{1}.$
\end{lemma}

\bigskip The concept of modulus of continuity in $H_{p}$ $\left( p>0\right) $
is defined by
\begin{equation*}
\omega _{H_{p}}\left( \frac{1}{M_{n}},f\right) :=\left\Vert
f-S_{M_{n}}f\right\Vert _{H_{p}}.
\end{equation*}

We need to understand the meaning of the expression $f-S_{M_{n}}f$ where $f$ is a martingale and $S_{M_{n}}f$ is function. So, we give an explanation in the following remark:

\begin{remark}
\label{lemma2.3.6} Let $0<p\leq 1.$ Since
\begin{equation*}
S_{M_{n}}f=f^{\left( n\right) },\text{ \ for }f=\left( f^{\left( n\right)
}:n\in
\mathbb{N}
\right) \in H_{p}
\end{equation*}%
and
\begin{equation*}
\left( S_{M_{k}}f^{\left( n\right) }:k\in
\mathbb{N}
\right) =\left( S_{M_{k}}S_{M_{n}}f,k\in
\mathbb{N}
\right)
\end{equation*}%
\begin{equation*}
=\left( S_{M_{0}}f,\ldots ,S_{M_{n-1}}f,S_{M_{n}}f,S_{M_{n}}f,\ldots \right)
\end{equation*}%
\begin{equation*}
=\left( f^{\left( 0\right) },\ldots ,f^{\left( n-1\right) },f^{\left(
n\right) },f^{\left( n\right) },\ldots \right)
\end{equation*}%
we obtain that%
\begin{equation*}
f-S_{M_{n}}f=\left( f^{\left( k\right) }-S_{M_{k}}f:k\in
\mathbb{N}
\right)
\end{equation*}%
is a martingale, for which%
\begin{equation}
\left( f-S_{M_{n}}f\right) ^{\left( k\right) }=\left\{
\begin{array}{ll}
0, & k=0,.\ldots ,n, \\
f^{\left( k\right) }-f^{\left( n\right) }, & k\geq n+1,%
\end{array}%
\right.  \label{g100}
\end{equation}
\end{remark}

Watari \cite{wat} showed that there are strong connections between
\begin{equation*}
\omega _{p}\left( \frac{1}{M_{n}},f\right) ,\text{ \ }E_{M_{n}}\left(
L_{p},f\right) \text{ \ \ and\ \ \ }\left\Vert f-S_{M_{n}}f\right\Vert
_{p},\ \ p\geq 1,\text{ \ }n\in \mathbb{N}.
\end{equation*}

In particular,%
\begin{equation}
\frac{1}{2}\omega _{p}\left( \frac{1}{M_{n}},f\right) \leq \left\Vert
f-S_{M_{n}}f\right\Vert _{p}\leq \omega _{p}\left( \frac{1}{M_{n}},f\right)
\label{eqvi}
\end{equation}%
and%
\begin{equation*}
\frac{1}{2}\left\Vert f-S_{M_{n}}f\right\Vert _{p}\leq E_{M_{n}}\left(
L_{p},f\right) \leq \left\Vert f-S_{M_{n}}f\right\Vert _{p}.
\end{equation*}

\begin{remark}
\label{lemma2.3.7}Since
\begin{equation*}
\left\Vert f\right\Vert _{H_{p}}\sim \left\Vert f\right\Vert _{p},
\end{equation*}%
when $p>1$, by applying (\ref{eqvi}), we obtain that%
\begin{equation*}
\omega _{H_{p}}\left( \frac{1}{M_{n}},f\right) \sim \omega _{p}\left( \frac{1%
}{M_{n}},f\right) .
\end{equation*}
\end{remark}

\vspace{1.0cm}

\subsection{{\protect\bigskip  Examples of $%
p$-atoms and $H_{p}$ maringales}}

\begin{equation*}
\end{equation*}

The next two Examples can be found in the papers \cite{BGG} and \cite{BGG2}
by Blahota, G\'{a}t and Goginava (see also :

\begin{example}
\label{example2.1}Let $0<p\leq 1$ and $\lambda =\sup_{n}m_{n}.$ Then the
function
\begin{equation*}
a_{k}:=\frac{M_{\alpha _{k}}^{1/p-1}}{\lambda }\left( D_{M_{\alpha
_{k}+1}}-D_{M_{_{\alpha _{k}}}}\right)
\end{equation*}%
is a p-atom. Moreover,
\begin{equation*}
\left\Vert a_{k}\right\Vert _{H_{p}}\leq 1.
\end{equation*}
\end{example}

{\bf Proof}:
Here we present the proof from Tephnadze\textit{\ }\cite{tep2}
and \cite{tep3}. Since
\begin{equation*}
\text{supp}(a_{k})=I_{\alpha _{k}},\text{ \ \ }\int_{I_{\alpha
_{k}}}a_{k}d\mu =0
\end{equation*}%
and
\begin{equation*}
\left\Vert a_{k}\right\Vert _{\infty }\leq \frac{M_{\alpha _{k}}^{1/p-1}}{%
\lambda }M_{\alpha _{k}+1}\leq M_{\alpha _{k}}^{1/p}=(\text{supp }%
a_{k})^{-1/p},
\end{equation*}%
we conclude that $a_{k}$ is a p-atom.

Moreover, by using the orthonormality of Vilenkin functions, we find that
\begin{equation*}
S_{M_{n}}\left( D_{M_{\alpha _{k}+1}}-D_{M_{\alpha _{k}}}\right) =\left\{
\begin{array}{ll}
0, & n=0,...,\alpha _{k}, \\
D_{M_{\alpha _{k}+1}}-D_{M_{\alpha _{k}}}, & n\geq \alpha _{k}+1,%
\end{array}%
\right.
\end{equation*}%
and
\begin{equation*}
\sup\limits_{n\in \mathbb{N}}\left\vert S_{M_{n}}\left( D_{M_{\alpha
_{k}+1}}\left( x\right) -D_{M_{\alpha _{k}}}\left( x\right) \right)
\right\vert
\end{equation*}%
\begin{equation*}
=\left\vert D_{M_{\alpha _{k}+1}}\left( x\right) -D_{M_{\alpha _{k}}}\left(
x\right) \right\vert ,\text{ \ for all \ }x\in G_{m}.
\end{equation*}

If we invoke Lemma \ref{lemma2.3.4} we obtain that%
\begin{equation*}
\left\Vert a_{k}\right\Vert _{H_{p}}=\frac{M_{\alpha _{k}}^{1/p-1}}{\lambda }%
\left\Vert D_{M_{\alpha _{k}+1}}-D_{M_{\alpha _{k}}}\right\Vert _{H_{p}}
\end{equation*}%

\begin{equation*}
=\frac{M_{\alpha _{k}}^{1/p-1}}{\lambda }\left\Vert D_{M_{\alpha
_{k}+1}}-D_{M_{\alpha _{k}}}\right\Vert _{p}
\end{equation*}%
\begin{equation*}
=\frac{M_{\alpha _{k}}^{1/p-1}}{\lambda }\left( \int_{I_{_{\alpha
_{k}}}\backslash \text{ }I_{_{\alpha _{k}+1}}}M_{\alpha _{k}}^{p}d\mu
+\int_{I_{_{\alpha _{k}+1}}}\left( M_{\alpha _{k}+1}-M_{\alpha _{k}}\right)
^{p}d\mu \right) ^{1/p}
\end{equation*}%
\begin{equation*}
=\frac{M_{\alpha _{k}}^{1/p-1}}{\lambda }\left( \frac{m_{\alpha _{k}}-1}{%
M_{\alpha _{k}+1}}M_{\alpha _{k}}^{p}+\frac{\left( m_{\alpha _{k}}-1\right)
^{p}}{M_{\alpha _{k}+1}}M_{\alpha _{k}}^{p}\right) ^{1/p}
\end{equation*}%
\begin{equation*}
\leq \frac{M_{\alpha _{k}}^{1/p-1}}{\lambda }\cdot M_{\alpha
_{k}}^{1-1/p}\leq 1.
\end{equation*}

The proof is complete.
\QED

\begin{example}
\label{example2.2}Let $0<p\leq 1$ and
\begin{equation*}
f_{k}=D_{M_{2n_{k}+1}}-D_{M_{_{2n_{k}}}}.
\end{equation*}%
Then%
\begin{equation}
\widehat{f}_{k}\left( i\right) =\left\{
\begin{array}{ll}
1, & i=M_{_{2n_{k}}},...,M_{_{2n_{k}+1}}-1, \\
0, & \text{otherwise,}%
\end{array}%
\right.  \label{13l}
\end{equation}%
and
\begin{equation}
S_{i}f_{k}=\left\{
\begin{array}{ll}
D_{i}-D_{M_{_{2n_{k}}}}, & i=M_{_{2n_{k}}}+1,...,M_{_{2n_{k}+1}}-1, \\
f_{k}, & i\geq M_{2n_{k}+1}, \\
0, & \text{otherwise. }%
\end{array}%
\right.  \label{14l}
\end{equation}%
Moreover,
\begin{equation}
\left\Vert f_{k}\right\Vert _{H_{p}}\leq \lambda M_{_{2n_{k}}}^{1-1/p},
\label{15l}
\end{equation}%
where $\lambda =\sup_{n}m_{n}.$
\end{example}

{\bf Proof}:
The proof follows by using Example \ref%
{example2.1} in the case when $\alpha _{k}=2n_{k}$. We leave out the details.
\QED

The next three examples of regular martingales will be used frequently and
it can be found in Persson and Tephnadze \cite{pt1}:

\begin{example}
\label{example002}Let $M_{k}\leq n<M_{k+1}$ and $S_{n}f$ be the n-th partial sum
with respect to Vilenkin systems, where $f\in H_{p}$ \ for some $0<p\leq 1.$
Then $S_{n}f\in L_{1}$ for every fixed $n\in \mathbb{N}$ and%
\begin{equation*}
\left\Vert S_{n}f\right\Vert _{H_{p}}\leq \left\Vert \sup_{0\leq l\leq
k}\left\vert S_{M_{l}}f\right\vert \right\Vert _{p}+\left\Vert
S_{n}f\right\Vert _{p}
\end{equation*}%
\begin{equation*}
\leq \left\Vert \widetilde{S}_{\#}^{\ast }f\right\Vert _{p}+\left\Vert
S_{n}f\right\Vert _{p}.
\end{equation*}
\end{example}

{\bf Proof}:
We consider the following martingale
\begin{equation*}
f_{\#}=\left( S_{M_{l}}S_{n},\text{ }l\in \mathbb{N}\right)
\end{equation*}%
\begin{equation*}
f_{\#}:=\left( S_{M_{k}}S_{n}f,\text{ }k\geq 1\right) =\left(
S_{M_{0}},...,S_{M_{k}}f,...,\text{ }S_{n}f,...,S_{n}f,..\right) .
\end{equation*}%
It immediately follows that%
\begin{equation*}
\left\Vert S_{n}f\right\Vert _{H_{p}}\leq \left\Vert \sup_{0\leq l\leq
k}\left\vert S_{M_{l}}f\right\vert \right\Vert _{p}+\left\Vert
S_{n}f\right\Vert _{p}\leq \left\Vert \widetilde{S}_{\#}^{\ast }f\right\Vert
_{p}+\left\Vert S_{n}f\right\Vert _{p}.
\end{equation*}%
The proof is complete.
\QED

\begin{example}
\label{example003}Let $M_{k}\leq n<M_{k+1}$ and $\sigma _{n}f$ be n-th Fej\'er
means with respect to Vilenkin systems, where $f\in H_{p}$ \ for some $0<p\leq
1.$ Then $\sigma _{n}f\in L_{1}$ for every fixed $n\in \mathbb{N}$ and%
\begin{equation*}
\left\Vert \sigma _{n}f\right\Vert _{H_{p}}\leq \left\Vert \sup_{0\leq l\leq
k}\left\vert \sigma _{M_{l}}f\right\vert \right\Vert _{p}+\left\Vert \sigma
_{n}f\right\Vert _{p}\leq \left\Vert \widetilde{\sigma }_{\#}^{\ast
}f\right\Vert _{p}+\left\Vert \sigma _{n}f\right\Vert _{p}.
\end{equation*}
\end{example}

{\bf Proof}:
We consider the following martingale
\begin{equation*}
f_{\#}=\left( S_{M_{k}}\sigma _{n}f,\text{ }k\geq 1\right)
\end{equation*}%
\begin{equation*}
=\left( \frac{M_{0}\sigma _{M_{0}}}{n},...,\frac{M_{k}\sigma _{M_{k}}f}{n}%
,...,\text{ }\sigma _{n}f,...,\sigma _{n}f,...\right) .
\end{equation*}
and if we follow analogous steps as in Example \ref{example002} we readily can
complete the proof. Thus, we leave out the details.
\QED

The next five Examples of martingales will be used many times to prove
sharpness of our main results (c.f. the papers \cite{tep1}, \cite{tep4},
\cite{tep6}, \cite{tep7}, \cite{tep9} by Tephnadze).

\begin{example}
\label{example2.5}Let $0<p\leq 1,$ $\lambda =\sup_{n}m_{n},$ $\left\{
\lambda _{k}:k\in \mathbb{N}\right\} $ be a sequence of real numbers $\mathbb{R%
}$, such that%
\begin{equation}
\sum_{k=0}^{\infty }\left\vert \lambda _{k}\right\vert ^{p}\leq c_{p}<\infty
  \label{2AA}
\end{equation}%
and $\left\{ a_{k}:k\in \mathbb{N}\right\} $ be a sequence of $p$-atoms,
defined by%
\begin{equation}
a_{k}:=\frac{M_{2\alpha _{k}}^{1/p-1}}{\lambda }\left( D_{M_{2\alpha
_{k}+1}}-D_{M_{_{2\alpha _{k}}}}\right) .  \label{3AA}
\end{equation}

Then $\,f=\left( f^{\left( n\right) }:\text{ }n\in \mathbb{N}\right) ,$
where \qquad
\begin{equation*}
f^{\left( n\right) }:=\sum_{\left\{ k:\text{ }2\alpha _{k}<n\right\}
}\lambda _{k}a_{k}
\end{equation*}%
is a martingale, $f\in H_{p}\ $and%
\begin{equation}
\widehat{f}(j)  \label{5AA}
\end{equation}%
\begin{equation*}
=\left\{
\begin{array}{ll}
\frac{\lambda _{k}M_{2\alpha _{k}}^{1/p-1}}{\lambda }, & j\in \left\{
M_{2\alpha _{k}},...,\text{ ~}M_{2\alpha _{k}+1}-1\right\} ,\text{ }k\in
\mathbb{N}_{+}, \\
0, & \text{\thinspace }j\notin \bigcup\limits_{k=1}^{\infty }\left\{
M_{2\alpha _{k}},...,\text{ ~}M_{2\alpha _{k}+1}-1\right\} .%
\end{array}%
\right.
\end{equation*}%
Moreover,
\begin{equation}
\omega _{H_{p}}\left( \frac{1}{M_{n}},f\right) =O\left( \sum_{\left\{ k;%
\text{ }2\alpha _{k}\geq n\right\} }^{\infty }\left\vert \lambda
_{k}\right\vert ^{p}\right) ^{1/p},\text{\ \ when \ \ }n\rightarrow \infty ,
\label{4AA}
\end{equation}%
for all $0<p\leq 1.$

Let $M_{2\alpha _{l-1}+1}\leq j\leq M_{2\alpha _{l}},$\ $l\in \mathbb{N}%
_{+}. $ Then%
\begin{equation}
S_{j}f=S_{M_{2\alpha _{l-1}+1}}=\sum_{\eta =0}^{l-1}\frac{\lambda _{\eta
}M_{2\alpha _{\eta }}^{1/p-1}}{\lambda }\left( D_{M_{_{2\alpha _{\eta
}+1}}}-D_{M_{_{2\alpha _{\eta }}}}\right) .  \label{5BB}
\end{equation}

Let $M_{2\alpha _{l}}\leq j<M_{2\alpha _{l}+1},$ $l\in \mathbb{N}_{+}.$ Then%
\begin{equation}
S_{j}f=S_{M_{2\alpha _{l}}}+\frac{\lambda _{l}M_{2\alpha _{l}}^{1/p-1}\psi
_{2\alpha _{l}}D_{j-M_{_{2\alpha _{l}}}}}{\lambda }  \label{5CC}
\end{equation}%
\begin{equation*}
=\sum_{\eta =0}^{l-1}\frac{\lambda _{\eta }M_{2\alpha _{\eta }}^{1/p-1}}{%
\lambda }\left( D_{M_{_{2\alpha _{\eta }+1}}}-D_{M_{_{2\alpha _{\eta
}}}}\right) +\frac{\lambda _{l}M_{2\alpha _{l}}^{1/p-1}\psi _{M_{_{2\alpha
_{l}}}}D_{j-M_{_{2\alpha _{l}}}}}{\lambda }.
\end{equation*}
\end{example}

{\bf Proof}:
Since%
\begin{equation*}
S_{M_{n}}a_{k}=\left\{
\begin{array}{ll}
a_{k}, & \alpha _{k}<n, \\
0, & \alpha _{k}\geq n,%
\end{array}%
\right.
\end{equation*}%
by applying Lemma \ref{lemma2.1} and (\ref{2AA}) we can conclude that $f\in
H_{p}.$ Furthermore, in view of (\ref{g100}) in Remark \ref{lemma2.3.6}
we immediately get also the estimate of $H_{p}$ modulus of continuity of the
martingale $f$.

Let \thinspace $j\notin \bigcup\limits_{k=1}^{\infty }\left\{ M_{2\alpha
_{k}},...,\text{ ~}M_{2\alpha _{k}+1}-1\right\} .$ By using orthonormality
of Vilenkin functions we find that%
\begin{equation*}
\widehat{a_{k}}(j)=\int_{G_{m}}\frac{M_{2\alpha _{k}}^{1/p-1}}{\lambda }%
\left( D_{M_{2\alpha _{k}+1}}-D_{M_{_{2\alpha _{k}}}}\right) \overline{\psi }%
_{i}d\mu =0
\end{equation*}%
and%
\begin{equation}
\widehat{f}(j)=\lim_{n\rightarrow \infty }\int_{G_{m}}\sum_{\left\{ k:\text{
}2\alpha _{k}<n\right\} }\lambda _{k}a_{k}\overline{\psi }_{i}d\mu
\label{fc1}
\end{equation}%
\begin{equation*}
=\lim_{n\rightarrow \infty }\sum_{\left\{ k:\text{ }2\alpha _{k}<n\right\}
}\lambda _{k}\int_{G_{m}}a_{k}\overline{\psi }_{i}d\mu =0.
\end{equation*}

Let $j\in \left\{ M_{2\alpha _{k}},...,\text{ ~}M_{2\alpha _{k}+1}-1\right\}
,$ for some $k\in \mathbb{N}$. Then, for every $l\in \mathbb{N}$ we can
conclude that
\begin{equation*}
\widehat{a_{l}}(j)=\frac{M_{2\alpha _{l}}^{1/p-1}}{\lambda }%
\int_{G_{m}}\left( D_{M_{2\alpha _{l}+1}}-D_{M_{_{2\alpha _{l}}}}\right)
\overline{\psi }_{i}d\mu =\frac{M_{2\alpha _{l}}^{1/p-1}}{\lambda }\delta
_{l,k},
\end{equation*}%
where $\delta _{l,k}$ is Kronecker symbol. It follows that
\begin{equation}
\widehat{f}(j)=\lim_{n\rightarrow \infty }\sum_{\left\{ l:\text{ }2\alpha
_{k}<n\right\} }\lambda _{l}\int_{G_{m}}a_{l}\overline{\psi }_{i}d\mu
\label{fc2}
\end{equation}%
\begin{equation*}
=\lim_{n\rightarrow \infty }\sum_{\left\{ l:\text{ }2\alpha _{l}<n\right\} }%
\frac{\lambda _{l}M_{2\alpha _{l}}^{1/p-1}}{\lambda }\delta _{l,k}=\frac{%
\lambda _{k}M_{2\alpha _{k}}^{1/p-1}}{\lambda }.
\end{equation*}

Let $M_{2\alpha _{l-1}+1}+1\leq $ $j\leq M_{2\alpha _{l}}.$ In the view of (%
\ref{5BB}) we can conclude that
\begin{equation}
S_{j}f=S_{M_{2\alpha _{l-1}+1}}f=\sum_{v=0}^{M_{2\alpha _{l-1}+1}-1}\widehat{%
f}(v)\psi _{v}  \label{20AA}
\end{equation}%
\begin{equation*}
=\sum_{\eta =0}^{l-1}\sum_{v=M_{2\alpha _{\eta }}}^{M_{2\alpha _{\eta }+1}-1}%
\frac{\lambda _{\eta }M_{2\alpha _{\eta }}^{1/p-1}}{\lambda }\psi _{v}
\end{equation*}%
\begin{equation*}
=\sum_{\eta =0}^{l-1}\frac{\lambda _{\eta }M_{2\alpha _{\eta }}^{1/p-1}}{%
\lambda }\left( D_{M_{_{2\alpha _{\eta }+1}}}-D_{M_{_{2\alpha _{\eta
}}}}\right) .
\end{equation*}

Let $M_{2\alpha _{l}}\leq j<M_{2\alpha _{l}+1}.$ We apply (\ref{5CC})
and invoke (\ref{dn21}) in Lemma \ref{dn1} to obtain that%
\begin{equation}
S_{j}f=S_{M_{2\alpha _{l}}}f+\sum_{v=M_{2\alpha _{l}}}^{j-1}\widehat{f}%
(v)\psi _{v}  \label{21AA}
\end{equation}%
\begin{equation*}
=S_{M_{2\alpha _{l}}}f+\frac{\lambda _{l}M_{2\alpha _{l}}^{1/p-1}}{\lambda }%
\sum_{v=M_{2\alpha _{l}}}^{j-1}\psi _{v}
\end{equation*}%
\begin{equation*}
=S_{M_{2\alpha _{l}}}f+\frac{\lambda _{l}M_{2\alpha _{l}}^{1/p-1}}{\lambda }%
\left( D_{j}-D_{M_{_{2\alpha _{l}}}}\right)
\end{equation*}%
\begin{equation*}
=S_{M_{2\alpha _{l}}}f+\frac{\lambda _{l}M_{2\alpha _{l}}^{1/p-1}\psi
_{M_{_{2\alpha _{l}}}}D_{j-M_{_{2\alpha _{l}}}}}{\lambda }.
\end{equation*}

If we take (\ref{20AA}) obtained by $j=M_{2\alpha _{l}}$ into
account we can rewrite (\ref{21AA}) as
\begin{equation*}
S_{j}f=\sum_{\eta =0}^{l-1}\frac{\lambda _{\eta }M_{2\alpha _{\eta }}^{1/p-1}%
}{\lambda }\left( D_{M_{_{2\alpha _{\eta }+1}}}-D_{M_{_{2\alpha _{\eta
}}}}\right)
\end{equation*}%
\begin{equation*}
+\frac{\lambda _{l}M_{2\alpha _{l}}^{1/p-1}}{\lambda }\left(
D_{j}-D_{M_{_{2\alpha _{l}}}}\right) .
\end{equation*}

By using (\ref{g100}) we get that%
\begin{equation}
\left( f-S_{M_{n}}f\right) ^{\left( k\right) }=\left\{
\begin{array}{ll}
0, & k=0,.\ldots ,n, \\
f^{\left( k\right) }-f^{\left( n\right) }, & k\geq n+1,%
\end{array}%
\right.  \label{g1001}
\end{equation}%
\begin{equation*}
=\left\{
\begin{array}{ll}
0, & k=0,\ldots ,n, \\
\sum_{l=n+1}^{k}\mu _{l}S_{M_{n}}a_{l}, & k\geq n+1.%
\end{array}%
\right.
\end{equation*}

According to Lemma \ref{lemma2.1} we conclude that (\ref{g1001}) is an
atomic decomposition of the martingale $f-S_{M_{n}}f\in H_{p}$ and
\begin{equation*}
\omega _{H_{p}}\left( \frac{1}{M_{n}},f\right) \backsim \left(
\sum_{k=n+1}^{\infty }\left\vert \mu _{k}\right\vert ^{p}\right)
^{1/p}<\infty .
\end{equation*}

The proof is complete.
\QED

Sometimes, we will use the martingale defined in the next example.

\begin{example}
\label{example2.6}Let $0<p\leq 1,$ $\lambda =\sup_{n\in \mathbb{N}}m_{n},$ $%
\left\{ \lambda _{k}:k\in \mathbb{N}\right\} $ be a sequence of real numbers $%
\mathbb{R}$, such that%
\begin{equation}
\sum_{k=0}^{\infty }\left\vert \lambda _{k}\right\vert ^{p}\leq c_{p}<\infty
.  \label{2aa}
\end{equation}%
and $\left\{ a_{k}:k\in \mathbb{N}\right\} $ be a sequence of $p$-atoms,
defined by%
\begin{equation*}
a_{k}:=\frac{M_{\left\vert \alpha _{k}\right\vert }^{1/p-1}}{\lambda }\left(
D_{M_{\left\vert \alpha _{k}\right\vert +1}}-D_{M_{_{\left\vert \alpha
_{k}\right\vert }}}\right) ,
\end{equation*}%
where $\left\vert \alpha _{k}\right\vert :=\max $ $\{j\in \mathbb{N}:$ $%
\left( \alpha _{k}\right) _{j}\neq 0\}$ and $\left( \alpha _{k}\right) _{j}$
denotes the $j$-th binary coefficient of $\alpha _{k}\in \mathbb{N}_{+}.$ Then $%
\,f=\left( f^{\left( n\right) }:\text{ }n\in \mathbb{N}\right) ,$ where
\qquad
\begin{equation*}
f^{\left( n\right) }:=\sum_{\left\{ k:\text{ }\left\vert \alpha
_{k}\right\vert <n\right\} }\lambda _{k}a_{k}.
\end{equation*}
is a martingale, $f\in H_{p}\ $and
\begin{equation}
\omega _{H_{p}}\left( \frac{1}{M_{n}},f\right) =O\left( \sum_{\left\{ k:%
\text{ }\left\vert \alpha _{k}\right\vert \geq n\right\} }^{\infty
}\left\vert \lambda _{k}\right\vert ^{p}\right) ^{1/p},\text{\ \ when \ \ }%
n\rightarrow \infty,  \label{2aa0}
\end{equation}%
for all $0<p\leq 1.$ Moreover%
\begin{equation}
\widehat{f}(j)  \label{10AA}
\end{equation}%
\begin{equation*}
=\left\{
\begin{array}{ll}
\frac{\lambda _{k}M_{\left\vert \alpha _{k}\right\vert }^{1/p-1}}{\lambda },
& j\in \left\{ M_{\left\vert \alpha _{k}\right\vert },...,\text{ ~}%
M_{\left\vert \alpha _{k}\right\vert +1}-1\right\} ,\text{ }k\in \mathbb{N}%
_{+}, \\
0, & \text{\thinspace }j\notin \bigcup\limits_{k=1}^{\infty }\left\{
M_{\left\vert \alpha _{k}\right\vert },...,\text{ ~}M_{\left\vert \alpha
_{k}\right\vert +1}-1\right\} .%
\end{array}%
\right.
\end{equation*}%
Let $M_{\left\vert \alpha _{l-1}\right\vert +1}\leq j\leq M_{\left\vert
\alpha _{l}\right\vert },$\ $l\in \mathbb{N}_{+}.$ Then%
\begin{equation}
S_{j}f=S_{M_{\left\vert \alpha _{l-1}\right\vert +1}}=\sum_{\eta =0}^{l-1}%
\frac{\lambda _{\eta }M_{_{\left\vert \alpha _{\eta }\right\vert }}^{1/p-1}}{%
\lambda }\left( D_{M_{_{\left\vert \alpha _{\eta }\right\vert
+1}}}-D_{M_{_{_{\left\vert \alpha _{\eta }\right\vert }}}}\right) .
\label{12AA}
\end{equation}%
Let $M_{\left\vert \alpha _{l}\right\vert }\leq j<M_{\left\vert \alpha
_{l}\right\vert +1},$ $l\in \mathbb{N}_{+}.$ Then%
\begin{equation}
S_{j}f=S_{M_{\left\vert \alpha _{l}\right\vert }}+\frac{\lambda
_{l}M_{\left\vert \alpha _{l}\right\vert }^{1/p-1}\psi _{M_{\left\vert
\alpha _{l}\right\vert }}D_{j-M_{_{\left\vert \alpha _{l}\right\vert }}}}{%
\lambda }  \label{11AA}
\end{equation}%
\begin{equation*}
=\sum_{\eta =0}^{l-1}\frac{\lambda _{\eta }M_{\left\vert \alpha _{\eta
}\right\vert }^{1/p-1}}{\lambda }\left( D_{M_{_{\left\vert \alpha _{\eta
}\right\vert +1}}}-D_{M_{_{\left\vert \alpha _{\eta }\right\vert }}}\right) +%
\frac{\lambda _{l}M_{\left\vert \alpha _{l}\right\vert }^{1/p-1}\psi
_{M_{\left\vert \alpha _{l}\right\vert }}D_{j-M_{\left\vert \alpha
_{l}\right\vert }}}{\lambda }.
\end{equation*}
\end{example}

{\bf Proof}:
The proof is quite analogously to the proof of Example \ref{example2.5}.
Hence, we leave out the details.
\QED

For our further investigations it is convenient to fix sequence $\left\{
\lambda _{k}:k\in \mathbb{N}\right\} .$ In many cases such appropriate choice
can be the following sequence
\begin{equation*}
\left\{ \lambda /\alpha _{k}^{1/2}:k\in \mathbb{N}\right\}.
\end{equation*}

\begin{example}
\label{example2.5.1}\bigskip Let $0<p\leq q\leq 1,$ $\lambda =\sup_{n\in
\mathbb{N}}m_{n}$ and $\left\{ \alpha _{k}:k\in \mathbb{N}\right\} $ be an
increasing sequence of positive integers such that\qquad
\begin{equation}
\sum_{k=0}^{\infty }\frac{1}{\alpha _{k}^{p/2}}<\infty ,  \label{1t2}
\end{equation}%
\begin{equation}
\sum_{\eta =0}^{k-1}\frac{M_{2\alpha _{\eta }}^{1/p}}{\alpha _{\eta }^{1/2}}<%
\frac{M_{2\alpha _{k}}^{1/p}}{\alpha _{k}^{1/2}},  \label{1t3}
\end{equation}%
and
\begin{equation}
\frac{32\lambda M_{2\alpha _{k-1}}^{1/p}}{\alpha _{k-1}^{1/2}}<\left\{
\begin{array}{ll}
\frac{M_{\alpha _{k}}^{1/p}}{\alpha _{k}^{3/2}} & p=q \\
\frac{M_{\alpha _{k}}^{1/p-1/q}}{\alpha _{k}^{3/2}} & \text{\thinspace }p<q%
\end{array}%
\right. \text{\ }  \label{1t4}
\end{equation}%

For the given sequence $\left\{ \alpha _{k}:k\in \mathbb{N}\right\} $, which
satisfy conditions (\ref{1t2})-(\ref{1t4}) we define a martingale $%
\,f=\left( f^{\left( n\right) }:\text{ }n\in \mathbb{N}\right) ,$
where\qquad
\begin{equation*}
f^{\left( n\right) }:=\sum_{\left\{ k:\text{ }2\alpha _{k}<n\right\}
}\lambda _{k}a_{k},
\end{equation*}%
and
\begin{equation*}
\lambda _{k}=\frac{\lambda }{\alpha _{k}^{1/2}}\ \ \ \text{and \ \ }a_{k}:=%
\frac{M_{2\alpha _{k}}^{1/p-1}}{\lambda }\left( D_{M_{2\alpha
_{k}+1}}-D_{M_{_{2\alpha _{k}}}}\right) .
\end{equation*}%

Then $\,f\in H_{p}$ and%
\begin{equation}
\widehat{f}(j)  \label{1t6}
\end{equation}%
\begin{equation*}
=\left\{
\begin{array}{ll}
\frac{M_{2\alpha _{k}}^{1/p-1}}{\alpha _{k}^{1/2}}, & j\in \left\{
M_{2\alpha _{k}},...,\text{ ~}M_{2\alpha _{k}+1}-1\right\} ,\text{ }k\in
\mathbb{N}_{+}, \\
0, & \text{\thinspace }j\notin \bigcup\limits_{k=1}^{\infty }\left\{
M_{2\alpha _{k}},...,\text{ ~}M_{2\alpha _{k}+1}-1\right\} .%
\end{array}%
\right.
\end{equation*}%

For every $l\in \mathbb{N}$ we have the following identities%
\begin{equation}
S_{j}f=\left\{
\begin{array}{ll}
S_{M_{2\alpha _{l-1}+1}}, & M_{2\alpha _{l-1}+1}\leq j\leq M_{2\alpha _{l}},%
\text{ } \\
S_{M_{2\alpha _{l}}}+\frac{M_{2\alpha _{l}}^{1/p-1}\psi _{M_{_{2\alpha
_{l}}}}D_{j-M_{_{2\alpha _{l}}}}}{\alpha _{l}^{1/2}}, & M_{2\alpha _{l}}\leq
j<M_{2\alpha _{l}+1}.%
\end{array}%
\right.  \label{1t6.1}
\end{equation}%
Moreover,
\begin{equation}
\left\vert S_{j}f\right\vert \leq \frac{2\lambda M_{2\alpha _{k-1}}^{p}}{%
\alpha _{k-1}^{1/2}},\text{ \ for all }0\leq j\leq M_{2\alpha _{k}}\text{
and }0<p\leq q\leq 1  \label{sn101}
\end{equation}%
and%
\begin{equation}
\left\vert S_{M_{2\alpha _{k}+1}}f\right\vert \geq \frac{M_{2\alpha
_{k}}^{1/p-1}}{4\alpha _{k}^{1/2}},\text{ }0<p<q<1.  \label{sn102}
\end{equation}
\end{example}

{\bf Proof}:
First, we note that such an increasing sequence $\left\{ \alpha
_{k}:k\in \mathbb{N}\right\} $ which satisfies conditions (\ref{1t2})-(\ref%
{1t4}), can be constructed.

By applying (\ref{1t2}) we conclude that condition (\ref{2AA}) is satisfied
and we can invoke Lemma \ref{1t2} to obtain that the martingale $\,f=\left(
f^{\left( n\right) }:n\in \mathbb{N}\right) \in H_{p}.$

Since
\begin{equation*}
\frac{\lambda _{k}M_{2\alpha _{k}}^{1/p-1}}{\lambda }=\frac{M_{2\alpha
_{k}}^{1/p-1}}{\alpha _{k}^{1/2}},
\end{equation*}
according to (\ref{fc1}) and (\ref{fc2}), we find that Fourier coefficients
are given by equality (\ref{1t6}).

Let $M_{2\alpha _{l-1}+1}+1\leq $ $j\leq M_{2\alpha _{l}},$ where $1\leq
l\leq k.$ In view of (\ref{1t6}) we can conclude that
\begin{equation}
S_{j}f=S_{M_{2\alpha _{l-1}+1}}f=\sum_{v=0}^{M_{2\alpha _{l-1}+1}-1}\widehat{%
f}(v)\psi _{v}  \label{1t6.2}
\end{equation}%
\begin{equation*}
=\sum_{\eta =0}^{l-1}\sum_{v=M_{2\alpha _{\eta }}}^{M_{2\alpha _{\eta }+1}-1}%
\frac{M_{2\alpha _{\eta }}^{1/p-1}}{\alpha _{\eta }^{1/2}}\psi _{v}
\end{equation*}%
\begin{equation*}
=\sum_{\eta =0}^{l-1}\frac{M_{2\alpha _{\eta }}^{1/p-1}}{\alpha _{\eta
}^{1/2}}\left( D_{M_{_{2\alpha _{\eta }+1}}}-D_{M_{_{2\alpha _{\eta
}}}}\right) .
\end{equation*}

Let $M_{2\alpha _{l}}\leq j<M_{2\alpha _{l}+1},$ where $1\leq l\leq k.$ We
apply (\ref{1t6}) again and invoke (\ref{dn21}) in Lemma \ref{1t6} to obtain
that%
\begin{equation}
S_{j}f=S_{M_{2\alpha _{l}}}f+\sum_{v=M_{2\alpha _{l}}}^{j-1}\widehat{f}%
(v)\psi _{v}  \label{1t6.3}
\end{equation}%
\begin{equation*}
=S_{M_{2\alpha _{l}}}f+\frac{M_{2\alpha _{l}}}{\alpha _{l}^{1/2}}%
\sum_{v=M_{2\alpha _{l}}}^{j-1}\psi _{v}
\end{equation*}%
\begin{equation*}
=S_{M_{2\alpha _{l}}}f+\frac{M_{2\alpha _{l}}}{\alpha _{l}^{1/2}}\left(
D_{j}-D_{M_{_{2\alpha _{l}}}}\right)
\end{equation*}%
\begin{equation*}
=S_{M_{2\alpha _{l}}}f+\frac{M_{2\alpha _{l}}\psi _{2\alpha
_{l}}D_{j-M_{_{2\alpha _{l}}}}}{\alpha _{l}^{1/2}}.
\end{equation*}

By combining (\ref{1t6.2}) and (\ref{1t6.3}) we get the identities given by (\ref{1t6.1}).

Moreover, if we take (\ref{1t6.2}) obtained by letting $j=M_{2\alpha _{l}}$ we can
rewrite (\ref{1t6.3}) as
\begin{equation}
S_{j}f=\sum_{\eta =0}^{l-1}\frac{M_{2\alpha _{\eta }}^{1/p-1}}{\alpha _{\eta
}^{1/2}}\left( D_{M_{_{2\alpha _{\eta }+1}}}-D_{M_{_{2\alpha _{\eta
}}}}\right)  \label{1t8}
\end{equation}%
\begin{equation}
=\sum_{\eta =0}^{l-1}\frac{M_{2\alpha _{\eta }}^{1/p-1}}{\alpha _{\eta
}^{1/2}}\left( D_{M_{_{2\alpha _{\eta }+1}}}-D_{M_{_{2\alpha _{\eta
}}}}\right) +\frac{M_{2\alpha _{l}}^{1/p-1}}{\alpha _{l}^{1/2}}\left(
D_{j}-D_{M_{_{2\alpha _{l}}}}\right) .  \notag
\end{equation}

Let $M_{2\alpha _{l-1}+1}-1\leq $ $j\leq M_{2\alpha _{l}},$ where $1\leq
l\leq k.$ Then, according to (\ref{1t3}) and (\ref{1t8}), we get that
\begin{equation}
\left\vert S_{j}f\right\vert \leq \sum_{\eta =0}^{l-1}\frac{\lambda
M_{2\alpha _{\eta -1}}^{1/p}}{\alpha _{\eta -1}^{1/2}}  \label{sn100}
\end{equation}%
\begin{equation*}
=\sum_{\eta =0}^{l-2}\frac{\lambda M_{2\alpha _{\eta -1}}^{1/p}}{\alpha
_{\eta -1}^{1/2}}+\frac{\lambda M_{2\alpha _{l-1}}^{1/p}}{\alpha _{l-1}^{1/2}%
}\leq \frac{2\lambda M_{2\alpha _{l-1}}^{1/p}}{\alpha _{l-1}^{1/2}}\leq
\frac{2\lambda M_{2\alpha _{k-1}}^{1/p}}{\alpha _{k-1}^{1/2}}.
\end{equation*}

Let $M_{2\alpha _{l}}\leq j<M_{2\alpha _{l}+1},$ where $1\leq l\leq k-1.$ In
view of (\ref{1t3}) and (\ref{1t8}) we can conclude that
\begin{equation}
\left\vert S_{j}f\right\vert \leq \sum_{\eta =0}^{l-1}\frac{\lambda
M_{2\alpha _{\eta -1}}^{1/p}}{\alpha _{\eta -1}^{1/2}}+\frac{M_{2\alpha
_{l}}^{1/p-1}j}{\alpha _{l}^{1/2}}  \label{sn99}
\end{equation}%
\begin{equation*}
\leq \frac{\lambda M_{2\alpha _{l}}^{1/p}}{\alpha _{l}^{1/2}}+\frac{\lambda
M_{2\alpha _{l}}^{1/p}}{\alpha _{l}^{1/2}}=\frac{\lambda M_{2\alpha
_{l}}^{1/p}}{\alpha _{l}^{1/2}}\leq \frac{2\lambda M_{2\alpha _{k-1}}^{1/p}}{%
\alpha _{k-1}^{1/2}}.
\end{equation*}

Estimates (\ref{sn99}) and (\ref{sn100}) insure us that (\ref{sn101}) is
satisfied.

Let $0<p<q<1$. By combining (\ref{1t4}) and (\ref{1t6}) with the expression (%
\ref{sn100}) obtained by taking $\ j=M_{2\alpha _{k}}$we see that
\begin{equation*}
\left\vert S_{M_{2\alpha _{k}+1}}f\right\vert =\left\vert \frac{M_{2\alpha
_{k}}^{1/p-1}}{\alpha _{k}^{1/2}}\psi _{M_{2\alpha _{k}}}+S_{M_{2\alpha
_{k}}}f\right\vert
\end{equation*}%
\begin{equation*}
=\left\vert \frac{M_{2\alpha _{k}}^{1/p-1}}{\alpha _{k}^{1/2}}\psi
_{M_{2\alpha _{k}}}+S_{M_{2\alpha _{k-1}+1}}f\right\vert
\end{equation*}%
\begin{equation}
\geq \left\vert \frac{M_{2\alpha _{k}}^{1/p-1}}{\alpha _{k}^{1/2}}\psi
_{M_{2\alpha _{k}}}\right\vert -\left\vert S_{M_{2\alpha
_{k-1}+1}}f\right\vert  \notag
\end{equation}%
\begin{equation*}
\geq \frac{M_{2\alpha _{k}}^{1/p-1}}{\alpha _{k}^{1/2}}-\frac{2\lambda
M_{2\alpha _{k-1}}^{1/p}}{\alpha _{k-1}^{1/2}}\geq \frac{M_{2\alpha
_{k}}^{1/p-1}}{\alpha _{k}^{1/2}}-\frac{M_{\alpha _{k}}^{1/p-1/q}}{16\alpha
_{k}^{3/2}}
\end{equation*}%
\begin{equation*}
\geq \frac{M_{2\alpha _{k}}^{1/p-1}}{\alpha _{k}^{1/2}}-\frac{M_{\alpha
_{k}}^{1/p-1}}{16\alpha _{k}^{3/2}}\geq \frac{M_{2\alpha _{k}}^{1/p-1}}{%
4\alpha _{k}^{1/2}}.
\end{equation*}

The proof is complete.
\QED

\begin{remark}
\label{remark2.5.1}It is obvious that if we cancel condition (\ref{1t4}) in
Example \ref{example2.5.1} we still get that $\,f=\left( f^{\left( n\right)
}:n\in \mathbb{N}\right) \in H_{p}$ and equality (\ref{1t6}) and estimate (%
\ref{sn101}) is valid. Moreover, \ (\ref{1t8})-(\ref{sn100}) are still
fulfilled.
\end{remark}

In the next chapter we many times will use a martingale, with the
following fixed sequence $\left\{ \lambda _{k}:k\in \mathbb{N}\right\}: $

\begin{example}
\label{example2.5.2}Let $0<p\leq 1$, \ $\lambda =\sup_{n\in \mathbb{N}%
}m_{n}\ $and $\left\{ \Phi _{k}:k\in \mathbb{N}\right\} $ be sequence of
non-negative and nondecreasing numbers, such that%
\begin{equation}
\sum_{k=1}^{\infty }\frac{1}{\Phi _{M_{_{2\alpha _{k}}}}^{p/4}}<\infty .
\label{psin}
\end{equation}%

For this given sequence $\left\{ \alpha _{k}:k\in \mathbb{N}\right\} $,
which satisfies condition (\ref{psin}) we define a martingale $%
\,f=\left( f^{\left( n\right) }:\text{ }n\in \mathbb{N}\right) ,$ where
\begin{equation*}
f^{\left( n\right) }:=\sum_{\left\{ k:\text{ }2\alpha _{k}<n\right\} }\lambda
_{k}a_{k}
\end{equation*}%
and
\begin{equation*}
\lambda _{k}=\frac{\lambda }{\Phi _{M_{2_{\alpha _{k}}}}^{1/4}}\text{,\ \ }%
a_{k}=\frac{M_{2\alpha _{k}}^{1/p-1}}{\lambda }\left( D_{M_{2\alpha
_{k}+1}}-D_{M_{_{2\alpha _{k}}}}\right).
\end{equation*}%

Then $\,f=\left( f^{\left( n\right) }:\text{ }n\in \mathbb{N}\right) \in
H_{p}.$ Moreover,%
\begin{equation}
\widehat{f}(j)=\left\{
\begin{array}{ll}
\frac{M_{2\alpha _{k}}^{1/p-1}}{\Phi _{M_{2_{\alpha _{k}}}}^{1/4}}, & j\in
\left\{ M_{2\alpha _{k}},...,\text{ ~}M_{2\alpha _{k}+1}-1\right\} ,\text{ }%
k\in \mathbb{N}_{+}, \\
0, & j\notin \bigcup\limits_{k=1}^{\infty }\left\{ M_{2\alpha _{k}},...,%
\text{ ~}M_{2\alpha _{k}+1}-1\right\} .\text{ }%
\end{array}%
\right.  \label{3.7}
\end{equation}
\end{example}

{\bf Proof}:
The proof is quite analogously to the proof of Example \ref{example2.5}. Hence, we leave out the details.
\QED

Sometimes an appropriate choice of sequence will be the following
\begin{equation*}
\left\{ \lambda _{k}=\lambda /M_{k}^{1/p}:k\in \mathbb{N}\right\} .
\end{equation*}

\begin{example}
\label{example2.6.1}Let $0<p\leq 1,$ $\lambda =\sup_{n\in \mathbb{N}}m_{n}\ $%
and \qquad
\begin{equation*}
f^{\left( n\right) }=\sum_{\left\{ k:2M_{k}<n\right\} }\lambda _{k}a_{k},
\end{equation*}%
where
\begin{equation*}
\lambda _{k}:=\frac{\lambda }{M_{k}^{1/p}}\text{ \ \ and \ \ }a_{k}:=\frac{%
M_{2M_{k}}^{1/p-1}}{\lambda }\left( D_{M_{2M_{k}+1}}-D_{M_{2M_{k}}}\right) .
\end{equation*}

Then $f=\left( f^{\left( n\right) }:n\in \mathbb{N}\right) \in H_{p},$ for
all $0<p\leq 1.$ Moreover,%
\begin{equation}
\widehat{f}(j)=\left\{
\begin{array}{ll}
\frac{M_{2M_{i}}^{1/p-1}}{M_{i}^{1/p}}, & \text{\thinspace }j\in \left\{
M_{2M_{i}},...,M_{2M_{i}+1}-1\right\} ,\text{ }i\in \mathbb{N}_{+}, \\
0, & \text{\thinspace }j\notin \bigcup\limits_{i=0}^{\infty }\left\{
M_{2M_{i}},...,M_{2M_{i}+1}-1\right\} ,\text{ }%
\end{array}%
\right.  \label{13}
\end{equation}%
and
\begin{equation}
\omega _{H_{p}}\left( \frac{1}{M_{n}},f\right) =O\left( \frac{1}{n^{1/p}}%
\right) ,\text{ \ \ when \ \ }n\rightarrow \infty .  \label{13m}
\end{equation}%
Let $M_{2M_{i}}\leq j<M_{2M_{i}+1},$ $l\in \mathbb{N}_{+}.$ Then%
\begin{equation}
S_{j}f=S_{M_{2M_{i}}}+\frac{M_{_{2M_{i}}}^{1/p-1}\psi
_{M_{2M_{i}}}D_{j-M_{2M_{i}}}}{M_{i}^{1/p}}.  \label{13s}
\end{equation}
\end{example}

{\bf Proof}:
This example is a simple consequence of Example \ref{example2.5} if we take $%
\alpha _{k}=M_{k}$ and $\lambda _{k}=\lambda /M_{k}^{1/p}.$ Since
\begin{equation*}
\sum_{k=0}^{\infty }\left\vert \lambda _{k}\right\vert ^{p}\leq
c\sum\limits_{k=0}^{\infty }\frac{1}{M_{k}}\leq c\sum\limits_{k=0}^{\infty }%
\frac{1}{2^{k}}<c<\infty ,
\end{equation*}%
we obtain that $f=\left( f^{\left( n\right) }:n\in \mathbb{N}\right) \in
H_{p},$ for all $0<p\leq 1.$

In view of (\ref{2aa0}) we can conclude that
\begin{equation*}
\omega _{H_{p}}\left( \frac{1}{M_{n}},f\right) \backsim \left( \sum_{\left\{
k;\text{ }2M_{k}\geq n\right\} }\frac{1}{M_{k}}\right) ^{1/p}\leq \left(
\sum_{\left\{ k;\text{ }2^{k+1}\geq n\right\} }\frac{1}{2^{k}}\right) ^{1/p}
\end{equation*}%
\begin{equation*}
\leq \left( \sum_{\left\{ k;\text{ }k\geq \log _{2}^{n}-1\right\} }\frac{1}{%
2^{k}}\right) ^{1/p}\leq \left( \frac{1}{2^{\log _{2}^{n}-1}}\right) ^{1/p}
\end{equation*}%
\begin{equation*}
=O\left( \frac{1}{n^{1/p}}\right) \rightarrow \infty ,\text{ when \ }%
n\rightarrow \infty .
\end{equation*}

If we apply (\ref{5CC}) obtained by taking $\alpha _{k}=M_{k}$ we immediately get
equality (\ref{13s}). So the proof is complete.
\QED

\newpage

\section{\textbf{Fourier coefficients and partial sums of Vilenkin-Fourier
series on Martingale Hardy spaces}\protect\bigskip}

\vspace{0.5cm}

\subsection{Some classical results on the Vilenkin-Fourier coefficients and
partial sums of Vilenkin-Fourier series}

According to the  Riemann-Lebesgue lemma (for details see e.g. the book \cite{sws}) we obtain that
\begin{equation*}
\widehat{f}\left( k\right) \rightarrow 0,\text{ \ \ when \ \ \ }k\rightarrow
\infty ,
\end{equation*}%
for each $f\in L_{1}.$

Moreover,
\begin{equation*}
\left\vert \widehat{f}\left( n\right) \right\vert \leq \frac{1}{2}\omega
_{1}\left( \frac{1}{M_{N}},f\right) ,\text{ }
\end{equation*}%
where $M_{N}\leq n\leq M_{N+1}$ \ and \ $f\in L_{1}.$

It is well-known (see e.g. the books \ \cite{AVD}\textit{ }and \cite{sws}
) that if $%
f\in L_{1}$ and the Vilenkin series
\begin{equation*}
T\left( x\right) =\sum_{j=0}^{\infty }c_{j}\psi _{j}\left( x\right)
\end{equation*}%
convergences to $f$ \ in $L_{1}$-norm, then
\begin{equation*}
c_{j}=\int_{G_{m}}f\overline{\psi }_{j}d\mu :=\widehat{f}\left( j\right) ,
\end{equation*}%
i.e. in this case the approximation series must be a Vilenkin-Fourier
series. An analogous result is true also if the Vilenkin series convergences
uniformly on $G_{m}$ to an integrable function $f$ .

Since $H_{1}\subset L_{1}$ it yields that $\widehat{f}\left( k\right)\rightarrow 0$ when $k\rightarrow \infty ,$ for every $f\in H_{1}.$ The classical
inequality of Hardy type is well known in the trigonometric as well as in
the Vilenkin-Fourier analysis. Namely,
\begin{equation*}
\overset{\infty }{\underset{k=1}{\sum }}\frac{\left\vert \widehat{f}\left(
k\right) \right\vert }{k}\leq c\left\Vert f\right\Vert _{H_{1}},
\end{equation*}%
where the function $f$ belongs to the Hardy space $H_{1}$ and $c$ is an
absolute constant. This was proved in the trigonometric case by Hardy and
Littlewood \cite{hl} (see also the book \cite{cw}) and for the
Walsh system it was proved in the book \cite{sws}.

Weisz \cite{We4,We1} generalized this result for the bounded Vilenkin systems and proved
that there is an absolute constant $c_{p},$ depending only on $p,$ such that
\begin{equation}
\overset{\infty }{\underset{k=1}{\sum }}\frac{\left\vert \widehat{f}\left(
k\right) \right\vert ^{p}}{k^{2-p}}\leq c_{p}\left\Vert f\right\Vert
_{H_{p}},\text{ \ }\left( f\in H_{p},\text{ }0<p\leq 2\right)  \label{3.1}
\end{equation}

Paley \cite{11} proved that the Walsh-Fourier coefficients of a function $%
f\in L_{p}$ $\left( 1<p<2\right) $ satisfy the condition
\begin{equation*}
\overset{\infty }{\underset{k=1}{\sum }}\left\vert \widehat{f}\left(
2^{k}\right) \right\vert ^{2}<\infty .
\end{equation*}

This results fails to hold for $p=1.$ However, it in fact holds for
functions $f\in L_{1},$ such that the maximal function $f^{\ast }$ belongs $L_{1},$ i.e $f\in
H_{1}$ (see e.g the book \cite{cw})$.$

\bigskip Let \ $0<p\leq 1$ and $f\in H_{p}.$ For Vilenkin systems the
following inequality was proved by Simon and Weisz \cite{we5}:
\begin{equation}
\left( \underset{k=1}{\overset{\infty }{\sum }}M_{k}^{2-2/p}\underset{j=1}{%
\overset{m_{k}-1}{\sum }}\left\vert \widehat{f}\left( jM_{k}\right)
\right\vert ^{2}\right) ^{1/2}\leq c_{p}\left\Vert f\right\Vert _{H_{p}}.
\label{3.2}
\end{equation}

\bigskip

\bigskip By using the Lebesgue constants we easily obtain that $S_{n_{k}}f$
\ convergence to $f$ \ in $L_{1}$-norm, for every integrable function $f$,
if and only if $\sup_{k}L_{n_{k}}\leq c<\infty .$ There are various results
when $p>1.$ It is also well-known that (see e.g. the book
\ \cite{sws})
\begin{equation*}
\left\Vert S_{n}f\right\Vert _{p}\leq c_{p}\left\Vert f\right\Vert _{p},%
\text{ \ when \ }p>1,
\end{equation*}%
but it can be proved also a more stronger result (see e.g. the book \ \cite{sws}):%
\begin{equation*}
\left\Vert S^{\ast }f\right\Vert _{p}\leq c_{p}\left\Vert f\right\Vert _{p},%
\text{ \ when \ }f\in L_{p},\text{ \ }p>1.
\end{equation*}

Uniform and point-wise convergence and some approximation properties of
partial sums in $L_{1}$ norms was investigated by Goginava \cite{gog1} (see
also \cite{gog2}), Shneider \cite{sn} and Avdispahi\'c and Memi\'c \cite{am}.
Fine \cite{fi} derived sufficient conditions for the uniform convergence,
which are in complete analogy with the Dini-Lipschitz conditions. Guli\'cev
\cite{9} estimated the rate of uniform convergence of a Walsh-Fourier series
using Lebesgue constants and modulus of continuity. Uniform convergence of
subsequences of partial sums was investigated also in Goginava and
Tkebuchava \cite{gt}. This corresponding problem was considered for Vilenkin
groups $G_{m}$ by Fridli \cite{4} and G\'at \cite{5}.

It is known (for details see e.g. the books \ \cite{sws}
and \cite{We4,We1}) that the subsequence $S_{M_{n}}$ of the partial
sums is bounded from the martingale Hardy space $H_{p}$ to the Lebesgue
space $L_{p},$ for all $p>0.$ However, (see Tephnadze \cite{tep7}) there
exists a martingale $f\in H_{p}$ $\left( 0<p<1\right) ,$ such that
\begin{equation*}
\underset{n\in \mathbb{N}}{\sup }\left\Vert S_{M_{n}+1}f\right\Vert
_{weak-L_{p}}=\infty .
\end{equation*}%
The reason of the divergence of $S_{M_{n}+1}f$ \ is that when $0<p<1$ the
Fourier coefficients of $f\in H_{p}$ are not uniformly bounded (see
Tephnadze \cite{tep6}). However, Simon \cite{si1} proved that there exists
an absolute constant $c_{p},$ depending only on $p,$ such that
\begin{equation*}
\overset{\infty }{\underset{k=1}{\sum }}\frac{\left\Vert S_{k}f\right\Vert
_{p}^{p}}{k^{2-p}}\leq c_{p}\left\Vert f\right\Vert _{H_{p}}^{p},\text{ \ \
\ }\left( 0<p<1\right)
\end{equation*}%
for all $f\in H_{p}.$

It is also well-known that Vilenkin systems do not form bases in the
space $L_{1}.$ Moreover, there is a function in the Hardy space $H_{1},$
such that the partial sums of $f$ \ are not bounded in the $L_{1}$-norm.
However, in G\'{a}t \cite{Ga1} the following strong convergence result was
obtained:%
\begin{equation*}
\underset{n\rightarrow \infty }{\lim }\frac{1}{\log n}\overset{n}{\underset{%
k=1}{\sum }}\frac{\left\Vert S_{k}f-f\right\Vert _{1}}{k}=0
\end{equation*}%
for all $f\in H_{1}.$

For the trigonometric analogue see Smith \cite{sm} and for the Walsh-Paley
system see Simon \cite{Si3}.

Moreover, there exists an absolute constant $c,$ such that%
\begin{equation*}
\frac{1}{\log n}\overset{n}{\underset{k=1}{\sum }}\frac{\left\Vert
S_{k}f\right\Vert _{1}}{k}\leq c\left\Vert f\right\Vert _{H_{1}} \text{ \ }%
\left( n=2,3...\right)
\end{equation*}%
and
\begin{equation*}
\underset{n\rightarrow \infty }{\lim }\frac{1}{\log n}\overset{n}{\underset{%
k=1}{\sum }}\frac{\left\Vert S_{k}f\right\Vert _{1}}{k}=\left\Vert
f\right\Vert _{H_{1}},
\end{equation*}%
for all $f\in H_{1}$ $.$

\subsection{Estimations of the Vilenkin-Fourier coefficients on
martingale Hardy spaces}

In this section we prove that the Fourier coefficients of $\ f\in
H_{p}, $ $0<p<1,$ are not uniformly bounded. Moreover, we
prove a new estimate of the Vilenkin-Fourier coefficients (see Theorem \ref%
{theorem3.1}). By applying this estimate we can obtain an independent proof
of (\ref{3.1}) and (\ref{3.2}).

The next theorem can be found in Tephnadze \cite{tep6}:

\begin{theorem}
\label{theorem3.1}a)Let $0<p<1$ and $f\in H_{p}.$ Then there exists an
absolute constant $c_{p},$ depending only on $p,$ such that
\begin{equation*}
\left\vert \widehat{f}\left( n\right) \right\vert \leq
c_{p}n^{1/p-1}\left\Vert f\right\Vert _{H_{p}}.\text{\ }
\end{equation*}%
b) Let $0<p<1$ and $\left\{ \Phi _{n}:n\in \mathbb{N}\right\} $ be any
non-decreasing sequence, satisfying the condition
\begin{equation}
\underset{n\rightarrow \infty }{\overline{\lim }}\frac{n^{1/p-1}}{\Phi _{n}}%
=\infty .  \label{3.0}
\end{equation}%
Then there exists a martingale $f\in H_{p},$ such that
\begin{equation*}
\underset{n\rightarrow \infty }{\overline{\lim }}\frac{\left\vert \widehat{f}%
\left( n\right) \right\vert }{\Phi _{n}}=\infty .
\end{equation*}
\end{theorem}

{\bf Proof}:
Since
\begin{equation*}
\left( \sum_{k=0}^{\infty }\left\vert \lambda _{k}\right\vert \right)
^{p}\leq \sum_{k=0}^{\infty }\left\vert \lambda _{k}\right\vert ^{p},\text{
\ }\left( 0<p\leq 1\right)
\end{equation*}%
by applying Lemma \ref{lemma2.1} we obtain that the proof of part a) will
be complete if we show that%
\begin{equation*}
\frac{\left\vert \widehat{a}\left( n\right) \right\vert }{\left( n+1\right)
^{1/p-1}}\leq c<\infty \text{ },\text{ when }0<p<1,
\end{equation*}%
for every $p$-atom $a,$ where $I$ denotes the support of the atom.

Let $a$ be an arbitrary $p$-atom with support$\ I$ and $\mu \left( I\right)
=M_{N}^{-1}.$ We may assume that $I=I_{N}.$ It is easy to see that $\widehat{%
a}\left( n\right) =0$ when $n\leq M_{N}.$ Therefore, we can suppose that $%
n>M_{N}.$

According to the estimate $\left\Vert a\right\Vert _{\infty }\leq M_{N}^{1/p}$
we can write that
\begin{equation}
\frac{\left\vert \widehat{a}\left( n\right) \right\vert }{\left( n+1\right)
^{1/p-1}}\leq \frac{1}{\left( n+1\right) ^{1/p-1}}\int\limits_{I_{N}}\left%
\Vert a\right\Vert _{\infty }d\mu  \label{3.3}
\end{equation}%
\begin{equation*}
\leq \frac{1}{\left( n+1\right) ^{1/p-1}}\int\limits_{I_{N}}M_{N}^{1/p}d\mu =%
\frac{M_{N}^{1/p}\left\vert I_{N}\right\vert }{\left( n+1\right) ^{1/p-1}}
\end{equation*}%
\begin{equation*}
\leq \frac{M_{N}^{1/p-1}}{\left( n+1\right) ^{1/p-1}}\ \leq c<\infty ,\text{
\ when }0<p<1.
\end{equation*}%
It follows that%
\begin{equation*}
\left\vert \widehat{f}\left( n\right) \right\vert \leq
c_{p}n^{1/p-1}\left\Vert f\right\Vert _{H_{p}}
\end{equation*}%
and the proof of the part a) is complete.

Let $0<p<1$ and $\left\{ \Phi _{n}:n\in \mathbb{N}\right\} $ be any
non-decreasing, non-negative sequence, satisfying condition (\ref{3.0}).
Then, for every $0<p<1,$ there exists an increasing sequence $\left\{ \alpha
_{k}:k\in \mathbb{N}\right\} $ of positive integers such that%
\begin{equation*}
\lim_{k\rightarrow \infty }\frac{M_{2\alpha _{k}}^{\left( 1/p-1\right) /2}}{%
\Phi _{M_{2\alpha _{k}}}^{1/2}}=\infty
\end{equation*}%
and
\begin{equation}
\sum_{k=0}^{\infty }\frac{\Phi _{M_{2\alpha _{k}}}^{p/2}}{M_{2\alpha
_{k}}^{\left( 1-p\right) /2}}<\infty .  \label{3.4}
\end{equation}

Let $f=\left( f^{\left( n\right) },\text{ }n\in \mathbb{N}\right) $ be from
Example \ref{example2.5}, where
\begin{equation}
\lambda _{k}=\frac{\lambda \Phi _{M_{2\alpha _{k}}}^{1/2}}{M_{2\alpha
_{k}}^{\left( 1/p-1\right) /2}}.  \label{3.5.0}
\end{equation}%
In view of (\ref{3.4}) we conclude that (\ref{2AA}) is satisfied and by using Example \ref{example2.5} we obtain that\ $f\in H_{p}.$

By now using (\ref{sn100}) with $\lambda _{k}$ defined by (\ref{3.5.0})
we readily see that

\begin{equation}
\widehat{f}(j)=  \label{3.5}
\end{equation}%
\begin{equation*}
=\left\{
\begin{array}{ll}
M_{2\alpha _{k}}^{\left( 1/p-1\right) /2}\Phi _{M_{2\alpha _{k}}}^{1/2},\,\,
& \text{\thinspace }j\in \left\{ M_{2\alpha _{k}},...,\text{ ~}M_{2\alpha
_{k}+1}-1\right\} ,\text{ }k\in \mathbb{N}, \\
0, & \text{\thinspace }j\notin \bigcup\limits_{k=1}^{\infty }\left\{
M_{2\alpha _{k}},...,\text{ ~}M_{2\alpha _{k}+1}-1\right\} \text{.}%
\end{array}%
\right.
\end{equation*}

By combining (\ref{3.5.0}) and (\ref{3.5}) we find that
\begin{equation*}
\underset{n\rightarrow \infty }{\overline{\lim }}\frac{\widehat{f}(n)}{\Phi
_{n}}\geq \lim_{k\rightarrow \infty }\frac{\widehat{f}(M_{2\alpha _{k}})}{%
\Phi _{M_{2\alpha _{k}}}}
\end{equation*}%
\begin{equation*}
\geq \underset{k\rightarrow \infty }{\lim }\frac{M_{2\alpha _{k}}^{\left(
1/p-1\right) /2}\Phi _{M_{2\alpha _{k}}}^{1/2}}{\Phi _{M_{2\alpha _{k}}}}
\end{equation*}%
\begin{equation*}
\geq \underset{k\rightarrow \infty }{\lim }\frac{M_{2\alpha _{k}}^{\left(
1/p-1\right) /2}}{\Phi _{M_{2\alpha _{k}}}^{1/2}}=\infty .
\end{equation*}

The proof is complete.
\QED

\subsection{Hardy and Paley type inequalities on martingale Hardy spaces}

By using Theorem \ref{theorem3.1} we prove Hardy type inequality (\ref{3.1})
for $0<p\leq 1$. Moreover, we also show (see Tephnadze \cite{tep4}) that the
sequence $\left\{ 1/k^{2-p}:k\geq 0\right\} $ in the inequality below can
not be improved:

\begin{theorem}
\label{theorem3.2}a)Let $0<p\leq 1$ and $f\in H_{p}.$ Then there exists an
absolute constant $c_{p},$ depending only on $p,$ such that
\begin{equation*}
\overset{\infty }{\underset{k=1}{\sum }}\frac{\left\vert \widehat{f}\left(
k\right) \right\vert ^{p}}{k^{2-p}}\leq c_{p}\left\Vert f\right\Vert
_{H_{p}}^{p}.\text{ \ }
\end{equation*}%
b) Let $\left\{ \Phi _{n}:n\in \mathbb{N}\right\} $ be any non-decreasing
sequence, satisfying the condition
\begin{equation*}
\underset{n\rightarrow \infty }{\overline{\lim }}\Phi _{n}=+\infty .
\end{equation*}%
Then there exists a martingale $f\in H_{p},$ such that
\begin{equation*}
\text{ }\underset{k=1}{\overset{\infty }{\sum }}\frac{\left\vert \widehat{f}%
\left( k\right) \right\vert ^{p}\Phi _{k}}{k^{2-p}}=\infty ,\text{ }\left(
0<p\leq 2\right) .
\end{equation*}
\end{theorem}

{\bf Proof}:
By applying Lemma \ref{lemma2.1} we obtain that the proof of the part a)
will be complete, if we show that%
\begin{equation*}
\overset{\infty }{\underset{k=1}{\sum }}\frac{\left\vert \widehat{a}\left(
k\right) \right\vert ^{p}}{k^{2-p}}<c<\infty .
\end{equation*}%
for every $p$-atom $a$ $\left( 0<p\leq 1\right) ,$ where $I$ denotes the
support of the atom$.$ Let $a$ be an arbitrary $p$-atom with support$\ I$
and $\mu \left( I\right) =M_{N}^{-1}.$ Analogously we may assume that $%
I=I_{N}$ \ and $n>M_{N}.$

First, we consider the case $0<p<1.$ According to (\ref{3.3}) in Theorem %
\ref{theorem3.1} we readily get that
\begin{equation*}
\overset{\infty }{\underset{k=M_{N}+1}{\sum }}\frac{\left\vert \widehat{a}%
\left( k\right) \right\vert ^{p}}{k^{2-p}}
\end{equation*}%
\begin{equation*}
\leq \overset{\infty }{\underset{k=M_{N}+1}{\sum }}\frac{1}{k}\left\vert
\frac{\widehat{a}\left( k\right) }{k^{1/p-1}}\right\vert ^{p}
\end{equation*}%
\begin{equation*}
\leq \overset{\infty }{\underset{k=M_{N}+1}{\sum }}\frac{1}{k}\frac{%
M_{N}^{1-p}}{k^{1-p}}
\end{equation*}%
\begin{equation*}
\leq M_{N}^{1-p}\overset{\infty }{\underset{k=M_{N}+1}{\sum }}\frac{1}{%
k^{2-p}}\leq c_{p}\,<\infty .
\end{equation*}

Let $p=1.$ Hence, Schwatz and Bessel inequalities imply that%
\begin{equation}
\overset{\infty }{\underset{k=M_{N}+1}{\sum }}\frac{\left\vert \widehat{a}%
\left( k\right) \right\vert }{k}  \label{3.05}
\end{equation}%
\begin{equation*}
\leq \left( \overset{\infty }{\underset{k=M_{N}+1}{\sum }}\left\vert
\widehat{a}\left( k\right) \right\vert ^{2}\right) ^{1/2}\left( \overset{%
\infty }{\underset{k=M_{N}+1}{\sum }}\frac{1}{k^{2}}\right) ^{1/2}
\end{equation*}%
\begin{equation*}
\leq \frac{1}{M_{N}^{1/2}}\left\Vert a\right\Vert _{2}
\end{equation*}%
\begin{equation*}
=\frac{1}{M_{N}^{1/2}}\left( \int\limits_{I_{N}}\left\vert a\right\vert
^{2}d\mu \right) ^{1/2}
\end{equation*}%
\begin{equation*}
=\frac{1}{M_{N}^{1/2}}\left( \int\limits_{I_{N}}M_{N}^{2}d\mu \right) ^{1/2}
\end{equation*}%
\begin{equation*}
=\frac{1}{M_{N}^{1/2}}M_{N}^{1/2}<c<\infty .
\end{equation*}

The proof of part a) is complete.

Next we note that for every $0<p<1$ there exists an increasing sequence $\left\{ \alpha
_{k}\geq 2:k\in \mathbb{N}\right\} $ of positive integers such that%
\begin{equation*}
\lim_{k\rightarrow \infty }\Phi _{M_{_{2\alpha _{k}}}}^{1/2}=\infty
\end{equation*}%
and
\begin{equation}
\sum_{k=1}^{\infty }\frac{1}{\Phi _{M_{_{2\alpha _{k}}}}^{p/4}}<\infty.
\label{3.6}
\end{equation}%

Let $f=\left( f^{\left( n\right) },\text{ }n\in \mathbb{N}\right) $ be the martingale defined in
Example \ref{example2.5.2}. According to (\ref{3.6}) we conclude that $f\in
H_{p}.$

Moreover, by using (\ref{3.7}) from Example \ref{example2.5.2} and obvious estimates we find
that
\begin{equation*}
\underset{l=1}{\overset{M_{2\alpha _{k}+1}-1}{\sum }}\frac{\left\vert
\widehat{f}\left( l\right) \right\vert ^{p}\Phi _{l}}{l^{2-p}}
\end{equation*}%
\begin{equation*}
=\underset{n=1}{\overset{k}{\sum }}\underset{l=M_{2\alpha _{n}}}{\overset{%
M_{2\alpha _{n}+1}-1}{\sum }}\frac{\left\vert \widehat{f}\left( l\right)
\right\vert ^{p}\Phi _{l}}{l^{2-p}}
\end{equation*}%
\begin{equation*}
\geq \underset{l=M_{2\alpha _{k}}}{\overset{M_{2\alpha _{k}+1}-1}{\sum }}%
\frac{\left\vert \widehat{f}\left( l\right) \right\vert ^{p}\Phi _{l}}{%
l^{2-p}}
\end{equation*}%
\begin{equation*}
\geq c\Phi _{M_{2\alpha _{k}}}\underset{l=M_{2\alpha _{k}}}{\overset{%
M_{2\alpha _{k}+1}-1}{\sum }}\frac{\left\vert \widehat{f}\left( l\right)
\right\vert ^{p}}{l^{2-p}}
\end{equation*}%
\begin{equation*}
\geq c\Phi _{M_{2\alpha _{k}}}\frac{M_{2\alpha _{k}}^{1-p}}{\Phi
_{M_{2\alpha _{k}}}^{p/4}}\underset{l=M_{2\alpha _{k}}}{\overset{M_{2\alpha
_{k}+1}-1}{\sum }}\frac{1}{l^{2-p}}
\end{equation*}%
\begin{equation*}
\geq c\Phi _{M_{2\alpha _{k}}}^{1/2}M_{2\alpha _{k}}^{1-p}\underset{%
l=M_{2\alpha _{k}}}{\overset{M_{2\alpha _{k}+1}-1}{\sum }}\frac{1}{%
M_{2\alpha _{k}+1}^{2-p}}
\end{equation*}%
\begin{equation*}
\geq c\Phi _{M_{2\alpha _{k}}}^{1/2}M_{2\alpha _{k}}^{1-p}\frac{1}{%
M_{2\alpha _{k}+1}^{1-p}}
\end{equation*}%
\begin{equation*}
\geq c\Phi _{M_{2\alpha _{k}}}^{1/2}\rightarrow \infty ,\text{ \ when \ \ }%
k\rightarrow \infty .
\end{equation*}

The proof is complete.
\QED

By using Theorem \ref{theorem3.1} we can also prove a Paley type inequality (%
\ref{3.2}) for $0<p\leq 1$. Moreover, we show (see Tephnadze \cite{tep4}%
) that the sequence $\left\{ 1/M_{k}^{2-2/p}:k\in \mathbb{N}\right\} $ in
the inequality below can not be improved:

\begin{theorem}
\label{theorem3.3}a)Let $0<p\leq 1$ and $f\in H_{p}.$ Then there exists an
absolute constant $c_{p},$ depending only on $p,$ such that
\begin{equation*}
\left( \underset{k=1}{\overset{\infty }{\sum }}M_{k}^{2-2/p}\underset{j=1}{%
\overset{m_{k}-1}{\sum }}\left\vert \widehat{f}\left( jM_{k}\right)
\right\vert ^{2}\right) ^{1/2}\leq c_{p}\left\Vert f\right\Vert _{H_{p}}.
\end{equation*}%
b) Let $0<p\leq 1$ and $\left\{ \Phi _{n}:n\in \mathbb{N}\right\} $ be any
non-decreasing sequence satisfying the condition
\begin{equation*}
\underset{n\rightarrow \infty }{\overline{\lim }}\Phi _{n}=+\infty .
\end{equation*}%
Then there exists a martingale $f\in H_{p},$ such that
\begin{equation*}
\underset{k=1}{\overset{\infty }{\sum }}\frac{\Phi _{M_{k}}}{M_{k}^{2/p-2}}%
\underset{j=1}{\overset{m_{k}-1}{\sum }}\left\vert \widehat{f}\left(
jM_{k}\right) \right\vert ^{2}=\infty .
\end{equation*}
\end{theorem}

{\bf Proof}:
First, we consider case $0<p<1.$ By applying Lemma \ref{lemma2.1} for the
proof of the part a) it suffices to show that%
\begin{equation*}
\underset{k=N+1}{\overset{\infty }{\sum }}M_{k}^{2-2/p}\underset{j=1}{%
\overset{m_{k}-1}{\sum }}\left\vert \widehat{a}\left( jM_{k}\right)
\right\vert ^{2}\leq c<\infty
\end{equation*}%
for every $p$-atom $a$ with support$\ I=I_{N}.$

According to (\ref{3.3}) in Theorem \ref{theorem3.1} we readily get that%
\begin{equation*}
\underset{k=N+1}{\overset{\infty }{\sum }}M_{k}^{2-2/p}\underset{j=1}{%
\overset{m_{k}-1}{\sum }}\left\vert \widehat{a}\left( jM_{k}\right)
\right\vert ^{2}
\end{equation*}%
\begin{equation*}
\leq \underset{k=N+1}{\overset{\infty }{\sum }}m_{k}^{2-2/p}\underset{j=1}{%
\overset{m_{k}-1}{\sum }}\left\vert \frac{\widehat{a}\left( jM_{k}\right) }{%
\left( jM_{k}\right) ^{1/p-1}}\right\vert ^{2}
\end{equation*}%
\begin{equation*}
\leq \underset{k=N+1}{\overset{\infty }{\sum }}m_{k}^{2-2/p}\underset{j=1}{%
\overset{m_{k}-1}{\sum }}\left( \frac{M_{N}^{1/p-1}}{\left( jM_{k}\right)
^{1/p-1}}\right) ^{2}
\end{equation*}%
\begin{equation*}
\leq M_{N}^{2/p-2}\underset{k=N+1}{\overset{\infty }{\sum }}m_{k}^{3}\frac{1%
}{M_{k}^{2/p-2}}\leq c_{p}<\infty .
\end{equation*}

Analogously to (\ref{3.05}) we can prove this result also in the case when $p=1$,
so the proof of part a) is complete for every $0<p\leq 1.$

On the other hand, to prove part b) we use the construction of martingale
from Example \ref{example2.5.2} (see also part b) of Theorem \ref{theorem3.2}%
). According to (\ref{3.7}) in Example \ref{example2.5.2} we get that
\begin{equation*}
\underset{l=1}{\overset{k}{\sum }}M_{2\alpha _{l}}^{2-2/p}\Phi _{M_{2\alpha
_{l}}}\underset{j=1}{\overset{m_{\alpha _{l}}-1}{\sum }}\left\vert \widehat{f%
}\left( jM_{2\alpha _{l}}\right) \right\vert ^{2}
\end{equation*}%
\begin{equation*}
\geq M_{2\alpha _{k}}^{2-2/p}\Phi _{M_{2\alpha _{k}}}\underset{j=1}{\overset{%
m_{\alpha _{k}}-1}{\sum }}\left\vert \widehat{f}\left( jM_{2\alpha
_{k}}\right) \right\vert ^{2}
\end{equation*}%
\begin{equation*}
\geq cM_{2\alpha _{k}}^{2-2/p}\Phi _{M_{2\alpha _{k}}}\underset{j=1}{\overset%
{m_{\alpha _{k}}-1}{\sum }}\frac{M_{2\alpha _{k}}^{2/p-2}}{\Phi _{M_{2\alpha
_{k}}}^{1/2}}
\end{equation*}%
\begin{equation*}
\geq c\Phi _{M_{2\alpha _{k}}}^{1/2}\rightarrow \infty ,\text{ \qquad when \
\ \ }k\rightarrow \infty .
\end{equation*}

The proof is complete.
\QED

\subsection{Maximal operators of partial sums of Vilenkin-Fourier series on martingale Hardy spaces}

In this section we will consider weighted Maximal operators of partial
sums of Vilenkin-Fourier series and prove $(H_{p},L_{p})$ and $%
(H_{p},weak-L_{p})$ type inequalities. We also show sharpness of these
theorems in a special sense.

The next theorem can be found in Tephnadze \cite{tep7}.

\begin{theorem}
\label{theorem4.1}a) Let $0<p<1.$ Then the \bigskip maximal operator
\begin{equation*}
\widetilde{S}_{p}^{\ast }f:=\sup_{n\in \mathbb{N}}\frac{\left\vert
S_{n}f\right\vert }{\left( n+1\right) ^{1/p-1}}
\end{equation*}%
is bounded from the Hardy space $H_{p}$\ to the Lebesgue space $L_{p}.$

b) Let $0<p<1$ and $\left\{ \Phi _{n}:n\in \mathbb{N}\right\} $ be any
non-decreasing sequence satisfying the condition%
\begin{equation}
\overline{\lim_{n\rightarrow \infty }}\frac{n^{1/p-1}}{\Phi _{n}}=+\infty .
\label{6l}
\end{equation}%
Then%
\begin{equation*}
\sup_{k\in \mathbb{N}}\frac{\left\Vert \frac{S_{M_{2n_{k}}+1}f_{k}}{\Phi
_{M_{2n_{k}}+2}}\right\Vert _{weak-L_{p}}}{\left\Vert f_{k}\right\Vert
_{H_{p}}}=\infty .
\end{equation*}
\end{theorem}

{\bf Proof}:
Let $0<p<1.$\textbf{\ }Since $\overset{\sim }{S}_{p}^{\ast }$ is bounded
from $L_{\infty }$ to $L_{\infty }$ by Lemma \ref{lemma2.2} we obtain that
the proof of part a) will be complete if we show that%
\begin{equation*}
\int\limits_{\overline{I}}\left\vert \overset{\sim }{S}_{p}^{\ast
}a\right\vert ^{p}d\mu \leq c<\infty
\end{equation*}%
for every $p$-atom $a,$ where $I$ denotes the support of the atom$.$

Let $a$ be a $p$-atom with support$\ I$ and $\mu \left( I\right) =M_{N}.$ We
may assume that $I=I_{N}.$ It is easy to see that $S_{n}\left( a\right) =0 $
when $n\leq M_{N}.$ Therefore, we can suppose that $n>M_{N}.$

Since $\left\Vert a\right\Vert _{\infty }\leq M_{N}^{1/p}$ it yields that
\begin{equation*}
\left\vert S_{n}\left( a\right) \right\vert \leq \int_{I_{N}}\left\vert
a\left( t\right) \right\vert \left\vert D_{n}\left( x-t\right) \right\vert
d\mu \left( t\right)
\end{equation*}%
\begin{equation*}
\leq \left\Vert a\right\Vert _{\infty }\int_{I_{N}}\left\vert D_{n}\left(
x-t\right) \right\vert d\mu \left( t\right)
\end{equation*}%
\begin{equation*}
\leq M_{N}^{1/p}\int_{I_{N}}\left\vert D_{n}\left( x-t\right) \right\vert
d\mu \left( t\right) .
\end{equation*}

Let $0<p<1$ and $x\in I_{s}\backslash I_{s+1}.$ From Lemma \ref{dn2.6} we
get that
\begin{equation}
\frac{\left\vert S_{n}a\left( x\right) \right\vert }{\left( n+1\right)
^{1/p-1}}\leq \frac{cM_{N}^{1/p-1}M_{s}}{\left( n+1\right) ^{1/p-1}}.
\label{13AA}
\end{equation}

Since $n>M_{N}$ we can conclude that
\begin{equation}
\frac{\left\vert S_{n}a\left( x\right) \right\vert }{\left( n+1\right)
^{1/p-1}}\leq cM_{s},\text{\ \ for \ }x\in I_{s}\backslash I_{s+1},\text{ \ }%
0\leq s\leq N-1.  \label{13AAM}
\end{equation}

The expression on the right-hand side of (\ref{13AAM}) does not depend on $n.$
Therefore,
\begin{equation}
\left\vert \widetilde{S}_{p}^{\ast }a\left( x\right) \right\vert \leq cM_{s}.
\label{13AAMA}
\end{equation}

By combining (\ref{1.1}) and (\ref{13AAMA}) we obtain that
\begin{equation*}
\int_{\overline{I_{N}}}\left\vert \widetilde{S}_{p}^{\ast }a\left( x\right)
\right\vert ^{p}d\mu \left( x\right)
\end{equation*}%
\begin{equation*}
=\overset{N-1}{\underset{s=0}{\sum }}\int_{I_{s}\backslash
I_{s+1}}\left\vert \widetilde{S}_{p}^{\ast }a\left( x\right) \right\vert
^{p}d\mu \left( x\right)
\end{equation*}%
\begin{equation*}
\leq c_{p}\overset{N-1}{\underset{s=0}{\sum }}\frac{M_{s}^{p}}{M_{s}}%
<c_{p}<\infty .
\end{equation*}%
The proof of part a) is complete.

Let $0<p<1.$ Under condition (\ref{6l}) there exists positive integers $%
\left\{ n_{k}:k\in \mathbb{N}\right\} $ such that
\begin{equation*}
\lim_{k\rightarrow \infty }\frac{\left( M_{2n_{k}}+2\right) ^{1/p-1}}{\Phi
_{M_{2n_{k}}+2}}=\infty ,\text{ \ \ }0<p<1.
\end{equation*}%
To prove part b) we will apply the $p$-atoms from Example \ref{example2.2}.
By combining (\ref{13l}) and (\ref{14l}) we obtain that
\begin{equation*}
\frac{\left\vert S_{M_{2n_{k}}+1}f_{k}\right\vert }{\Phi _{M_{2n_{k}}+2}}=%
\frac{\left\vert D_{M_{2n_{k}}+1}-D_{M_{2n_{k}}}\right\vert }{\Phi
_{M_{2n_{k}}+2}}
\end{equation*}%
\begin{equation*}
=\frac{\left\vert \psi _{M_{2n_{k}}}\right\vert }{\Phi _{M_{2n_{k}}+2}}=%
\frac{1}{\Phi _{M_{2n_{k}}+2}}.
\end{equation*}

Hence,
\begin{equation}
\mu \left\{ x\in G_{m}:\frac{\left\vert S_{M_{2n_{k}}+1}f_{k}\left( x\right)
\right\vert }{\Phi _{M_{2n_{k}}+2}}\geq \frac{1}{\Phi _{M_{2n_{k}}+2}}%
\right\} =1.  \label{16}
\end{equation}

By combining (\ref{15l}) and (\ref{16}) we have that%
\begin{equation*}
\frac{\frac{1}{\Phi _{M_{2n_{k}}+2}}\left( \mu \left\{ x\in G_{m}:\frac{%
\left\vert S_{M_{2n_{k}}+1}f_{k}\right\vert }{\Phi _{M_{2n_{k}}+2}}\geq
\frac{1}{\Phi _{M_{2n_{k}}+2}}\right\} \right) ^{1/p}}{\left\Vert
f_{k}\right\Vert _{H_{p}}}
\end{equation*}%
\begin{equation*}
\geq \frac{1}{\Phi _{M_{2n_{k}}+2}M_{_{2n_{k}}}^{1-1/p}}
\end{equation*}%
\begin{equation*}
=\frac{\left( M_{_{2n_{k}}}+2\right) ^{1/p-1}}{\Phi _{M_{2n_{k}}+2}}%
\rightarrow \infty ,\text{ when \ }k\rightarrow \infty .
\end{equation*}

The proof is complete.
\QED

The next corollary is very important for our further investigation:

\begin{corollary}
\label{corollary4.2} Let $0<p<1$ and $f\in H_{p}.$ Then there exists an
absolute constant $c_{p},$ depending only on $p,$ such that%
\begin{equation*}
\left\Vert S_{n}f\right\Vert _{p}\leq c_{p}\left( n+1\right)
^{1/p-1}\left\Vert f\right\Vert _{H_{p}},\text{ \ }n\in \mathbb{N}_{+}.
\end{equation*}
\end{corollary}

{\bf Proof}:
According to part a) of Theorem \ref{theorem4.1} we conclude that%
\begin{equation*}
\left\Vert \frac{S_{n}f}{\left( n+1\right) ^{1/p-1}}\right\Vert _{p}
\end{equation*}%
\begin{equation*}
\leq \left\Vert \sup_{n\in \mathbb{N}}\frac{\left\vert S_{n}f\right\vert }{%
\left( n+1\right) ^{1/p-1}}\right\Vert _{p}\leq c_{p}\left\Vert f\right\Vert
_{H_{p}},\text{ \ }n\in \mathbb{N}_{+}.
\end{equation*}%
The proof is complete.
\QED

We also mention the following consequences of Theorem \ref{theorem4.1}:
\begin{corollary}
\label{corollary4.2a} Let $0<p<1$ and $\left\{ \Phi _{n}:n\in \mathbb{N}%
\right\} $ be any non-decreasing sequence, satisfying the condition (\ref{6l}%
). Then there exists a martingale $f\in H_{p},$ such that%
\begin{equation*}
\sup_{n\in \mathbb{N}}\left\Vert \frac{S_{n}f}{\Phi _{n}}\right\Vert
_{weak-L_{p}}=\infty .
\end{equation*}
\end{corollary}

\begin{corollary}
\label{corollary4.2aa} Let $0<p<1$ and $\left\{ \Phi _{n}:n\in \mathbb{N}%
\right\} $ be any non-decreasing sequence, satisfying the condition (\ref{6l}%
). Then there exists a martingale $f\in H_{p},$ such that the following maximal operator%
\begin{equation*}
\sup_{n\in \mathbb{N}}\frac{\left\vert S_{n}f\right\vert }{\Phi _{n}}
\end{equation*}%
\textit{is not bounded from the Hardy space }$H_{p}$\textit{\ to the space }$%
weak-L_{p}.$
\end{corollary}

The next result can be found in Tephnadze \cite{tep7}.

\begin{theorem}
\label{theorem4.2} a) Let $f\in H_{1}.$ Then the \bigskip maximal operator
\begin{equation*}
\widetilde{S}^{\ast }f:=\sup_{n\in \mathbb{N}_{+}}\frac{\left\vert
S_{n}f\right\vert }{\log \left( n+1\right) }\text{\ }
\end{equation*}%
is bounded from the Hardy space $H_{1}$ to the space $L_{1}.$ \noindent

b) Let $q_{n}=M_{2n}+M_{2n-2}+...+M_{0}$ and $\left\{ \Phi _{n}:n\in \mathbb{%
N}\right\} $ be any non-decreasing sequence, satisfying the condition%
\begin{equation}
\overline{\lim_{n\rightarrow \infty }}\frac{\log n}{\Phi _{n}}=+\infty .
\label{bb}
\end{equation}%
Then%
\begin{equation*}
\sup_{k}\frac{\left\Vert \frac{S_{q_{n_{k}}}f_{k}}{\Phi _{q_{n_{k}}+1}}%
\right\Vert _{1}}{\left\Vert f_{k}\right\Vert _{H_{1}}}=\infty .
\end{equation*}
\end{theorem}

{\bf Proof}:
\textbf{\ }Since $\overset{\sim }{S}^{\ast }$ is bounded from $L_{\infty }$
to $L_{\infty },$ according Lemma \ref{lemma2.2} it is sufficies to show that%
\begin{equation*}
\int\limits_{\overline{I}}\left\vert \overset{\sim }{S}^{\ast }a\right\vert
d\mu \leq c<\infty
\end{equation*}%
for every $p$-atom $a,$ where $I$ denotes the support of the atom$.$ We may
assume that $I=I_{N}.$ Since $S_{n}\left( a\right) =0$ when $n\leq M_{N},$
we can suppose that $n>M_{N}.$

Since $\left\Vert a\right\Vert _{\infty }\leq M_{N}$ we have that
\begin{equation*}
\left\vert S_{n}\left( a\right) \right\vert \leq M_{N}\int_{I_{N}}\left\vert
D_{n}\left( x-t\right) \right\vert d\mu \left( t\right) .
\end{equation*}

Let $x\in I_{s}\backslash I_{s+1}.$ From Lemma \ref{dn2.6} we get that
\begin{equation*}
\frac{\left\vert S_{n}a\left( x\right) \right\vert }{\log \left( n+1\right) }%
\leq \frac{M_{s}}{\log \left( n+1\right) }.
\end{equation*}

Since $n>M_{N}$ we can conclude that
\begin{equation}
\frac{\left\vert S_{n}a\left( x\right) \right\vert }{\log \left( n+1\right) }%
\leq \frac{cM_{s}}{N}.  \label{13AAAM}
\end{equation}

The expression on the right-hand side of (\ref{13AAAM}) does not depend on $n.$ Thus,
\begin{equation}
\left\vert \widetilde{S}^{\ast }a\left( x\right) \right\vert \leq \frac{%
cM_{s}}{N},\text{ \ \ for \ }x\in I_{s}\backslash I_{s+1},\text{ \ }0\leq
s\leq N-1.  \label{13AAAMA}
\end{equation}

By combining (\ref{1.1}) and (\ref{13AAAMA}) we obtain that
\begin{equation*}
\int_{\overline{I_{N}}}\left\vert \widetilde{S}^{\ast }a\left( x\right)
\right\vert d\mu \left( x\right)
\end{equation*}%
\begin{equation*}
=\overset{N-1}{\underset{s=0}{\sum }}\int_{I_{s}\backslash
I_{s+1}}\left\vert \widetilde{S}^{\ast }a\left( x\right) \right\vert d\mu
\left( x\right)
\end{equation*}%
\begin{equation*}
\leq \frac{c}{N}\overset{N-1}{\underset{s=0}{\sum }}\frac{M_{s}}{M_{s}}\leq
\frac{cN}{N}<c<\infty .
\end{equation*}

The proof of part a) is complete.

Next, we note that under condition (\ref{bb}) there exists $\left\{
n_{k}:k\in \mathbb{N}\right\} ,$ such that
\begin{equation*}
\lim_{k\rightarrow \infty }\frac{\log q_{n_{k}+1}}{\Phi _{q_{n_{k}+1}}}%
=\infty .
\end{equation*}

Let $f_{k}$ be 1-atom from Example \ref{example2.2}. By applying
Corollary \ref{dn5} with (\ref{15l}) we have that%
\begin{equation*}
\frac{\left\Vert \frac{S_{q_{n_{k}}}f_{k}\left( x\right) }{\Phi
_{q_{n_{k}+1}}}\right\Vert _{1}}{\left\Vert f_{k}\left( x\right) \right\Vert
_{H_{1}}}
\end{equation*}%
\begin{equation*}
\geq \frac{1}{\left\Vert f_{k}\right\Vert _{H_{1}}}\left( \left\Vert \frac{%
D_{q_{n_{k}}}}{\Phi _{q_{n_{k}+1}}}\right\Vert _{1}-\left\Vert \frac{%
D_{M_{2n_{k}}}}{\Phi _{q_{n_{k}+1}}}\right\Vert _{1}\right)
\end{equation*}%
\begin{equation*}
\geq \frac{c}{\Phi _{q_{n_{k}+1}}}\left( \log q_{n_{k}}-1\right)
\end{equation*}%
\begin{equation*}
\geq \frac{c\log q_{n_{k}+1}}{\Phi _{q_{n_{k}+1}}}\rightarrow \infty ,\text{
when }k\rightarrow \infty .
\end{equation*}

Hence, also part b) is proved so the proof is complete.
\QED

The next corollary is very important for our further investigation:

\begin{corollary}
\label{corollary4.3} Let $f\in H_{1}.$ Then there exists an absolute
constant $c,$ such that%
\begin{equation*}
\left\Vert S_{n}f\right\Vert _{1}\leq c\log \left( n+1\right) \left\Vert
f\right\Vert _{H_{1}},\text{ \ }n\in \mathbb{N}_{+}.
\end{equation*}
\end{corollary}

{\bf Proof}:
According to part a) of Theorem \ref{theorem4.2} we readily conclude that%
\begin{equation*}
\left\Vert \frac{S_{n}f}{\log \left( n+1\right) }\right\Vert _{1}
\end{equation*}%
\begin{equation*}
\leq \left\Vert \sup_{n\in \mathbb{N}}\frac{\left\vert S_{n}f\right\vert }{%
\log \left( n+1\right) }\right\Vert _{1}\leq c\left\Vert f\right\Vert
_{H_{1}},\text{ \ }n\in \mathbb{N}_{+}.
\end{equation*}%
The proof is complete.
\QED

We also point out two more consequences of Theorem \ref{theorem4.2}:

\begin{corollary}
\label{corollary4.3a} Let $\left\{ \Phi _{n}:n\in \mathbb{N}\right\} $ be
any non-decreasing sequence, satisfying the condition (\ref{bb}). Then there
exists a martingale $f\in H_{1}$ such that%
\begin{equation*}
\sup_{n\in \mathbb{N}}\left\Vert \frac{S_{n}f}{\Phi _{n}}\right\Vert
_{1}=\infty .
\end{equation*}
\end{corollary}

\begin{corollary}
\label{corollary4.3aa} Let $0<p<1$ and $\left\{ \Phi _{n}:n\in \mathbb{N}%
\right\} $ be any non-decreasing sequence, satisfying the condition (\ref{bb}%
). Then there exists a martingale $f\in H_{1}$ such that the following maximal operator%
\begin{equation*}
\sup_{n\in \mathbb{N}}\frac{\left\vert S_{n}f\right\vert }{\Phi _{n}}
\end{equation*}%
\textit{is not bounded from the Hardy space }$H_{1}$\textit{\ to the Lebesgue space }$%
L_{1}.$
\end{corollary}

\qquad Now, we formulate a theorem about boundedness of a restricted  maximal operator of
partial sums. This result is presented for the first
time.

\begin{theorem}
\label{theorem4.4}a) Let $0<p\leq 1$ and $\alpha _{k}$ be a subsequence of
positive numbers, such that
\begin{equation}
\sup_{k\in \mathbb{N}}\rho \left( \alpha _{k}\right) =\varkappa <\infty .  \label{ak0}
\end{equation}%
Then the maximal operator
\begin{equation*}
\widetilde{S}^{\ast ,\vartriangle }f:=\sup_{k\in \mathbb{N}}\left\vert S_{\alpha
_{k}}f\right\vert
\end{equation*}%
is bounded from the Hardy space $H_{p}$ to the space $L_{p}.$

b) Let $0<p<1\ \ $and $\left\{ \alpha _{k}:k\in \mathbb{N}\right\} $ be a
subsequence of positive numbers satisfying the condition
\begin{equation}
\sup_{k\in \mathbb{N}}\rho \left( \alpha _{k}\right) =\infty .  \label{ak}
\end{equation}%
Then there exists a martingale $f\in H_{p}$ such that%
\begin{equation*}
\sup_{k\in \mathbb{N}}\left\Vert S_{\alpha _{k}}f\right\Vert
_{weak-L_{p}}=\infty ,\text{ \ \ }\left( 0<p<1\right) .
\end{equation*}
\end{theorem}

{\bf Proof}:
By using (\ref{2dn}) in Lemma \ref{dn2.1} and Corollary \ref{dn2.3} we
easily conclude that if condition (\ref{ak0}) holds, then
\begin{equation*}
\left\Vert D_{\alpha _{k}}\right\Vert _{1}\leq c\overset{\left\vert \alpha
_{k}\right\vert }{\underset{j=\left\langle \alpha _{k}\right\rangle }{\sum }}%
\left\Vert D_{M_{j}}\right\Vert _{1}
\end{equation*}%
\begin{equation*}
\leq \overset{\left\vert \alpha _{k}\right\vert }{\underset{j=\left\langle
\alpha _{k}\right\rangle }{\sum }}1=\rho \left( \alpha _{k}\right) +1\leq
c<\infty .
\end{equation*}

It follows that $\widetilde{S}^{\ast ,\vartriangle }$ is bounded from $%
L_{\infty }$ to $L_{\infty }.$ By Lemma \ref{lemma2.2} we obtain that the
proof of part a) will be complete if we show that%
\begin{equation*}
\int\limits_{\overline{I_{N}}}\left\vert \widetilde{S}^{\ast ,\vartriangle
}a\right\vert ^{p}d\mu \leq c<\infty
\end{equation*}%
for every $p$-atom $a$ with support $I=I_{N}$. Since $S_{\alpha _{k}}\left(
a\right) =0$ when $\alpha _{k}\leq M_{N},$ we can suppose that $\alpha
_{k}>M_{N}.$

Let \ $t\in I_{N}$ and $x\in I_{s}\backslash I_{s+1},$ $1\leq s\leq
\left\langle \alpha _{k}\right\rangle -1.$ By using (\ref{2dn}) in Lemma \ref%
{dn2.1} we get that%
\begin{equation*}
D_{\alpha _{k}}\left( x-t\right) =0
\end{equation*}%
and
\begin{equation}
\left\vert S_{\alpha _{k}}\left( a\right) \right\vert =0.  \label{dnsn0}
\end{equation}

Let $0<p\leq 1,$ \ $t\in I_{N}$ and $x\in I_{s}\backslash I_{s+1},$ $%
\left\langle \alpha _{k}\right\rangle \leq s\leq N-1.$ By applying the fact
that $\left\Vert a\right\Vert _{\infty }\leq M_{N}$ and Lemma \ref{dn2.6} we
find that
\begin{equation}
\left\vert S_{\alpha _{k}}\left( a\right) \right\vert  \label{dnsn}
\end{equation}%
\begin{equation*}
\leq M_{N}^{1/p}\int_{I_{N}}\left\vert D_{\alpha _{k}}\left( x-t\right)
\right\vert d\mu \left( t\right) \leq c_{p}M_{N}^{1/p-1}M_{s}.
\end{equation*}%
Set
\begin{equation*}
\varrho :=\min_{k\in \mathbb{N}}\left\langle \alpha _{k}\right\rangle .
\end{equation*}
Then, in view of (\ref{dnsn0}) and (\ref{dnsn}) we can conclude that
\begin{equation}
\left\vert \widetilde{S}^{\ast ,\vartriangle }a\left( x\right) \right\vert
=0,\text{ \ \ for }\ x\in I_{s}\backslash I_{s+1},\ 0\leq s\leq \varrho
\label{dnsnM0}
\end{equation}%
and
\begin{equation}
\left\vert \widetilde{S}^{\ast ,\vartriangle }a\left( x\right) \right\vert
\leq c_{p}M_{N}^{1/p-1}M_{s},\text{ for \ }x\in I_{s}\backslash I_{s+1},%
\text{ \ }\varrho \leq s\leq N-1.  \label{dnsnM}
\end{equation}

By the definition of $\varrho$ there exists at least one index $%
k_{i_{0}}\in \mathbb{N}_{+}$ such that $\varrho =\left\langle \alpha
_{k_{i_{0}}}\right\rangle .$ It follows that
\begin{equation}
N-\varrho =N-\left\langle \alpha _{k_{i_{0}}}\right\rangle \leq \left\vert
\alpha _{k_{i_{0}}}\right\vert -\left\langle \alpha _{k_{i_{0}}}\right\rangle
\label{dnsnM1}
\end{equation}%
\begin{equation*}
\leq \sup_{k\in \mathbb{N}}\rho \left( \alpha _{k}\right) =\varkappa <c<\infty
\end{equation*}%
and

\begin{equation}
\frac{M_{N}^{1-p}}{M_{\varrho }^{1-p}}\leq \lambda ^{\left( N-\varrho
\right) \left( 1-p\right) }\leq \lambda ^{\varkappa \left( 1-p\right)
}<c<\infty ,  \label{dnsnM2}
\end{equation}%
where \ \ $\lambda =\sup_{k}m_{k}.$

Let $p=1.$ We combine (\ref{dnsnM0})-(\ref{dnsnM1}) and invoke identity (\ref%
{1.1}) to obtain that
\begin{equation*}
\int_{\overline{I_{N}}}\left\vert \widetilde{S}^{\ast ,\vartriangle }a\left(
x\right) \right\vert d\mu \left( x\right)
\end{equation*}%
\begin{equation*}
=\overset{N-1}{\underset{s=0}{\sum }}\int_{I_{s}\backslash
I_{s+1}}\left\vert \widetilde{S}^{\ast ,\vartriangle }a\left( x\right)
\right\vert d\mu \left( x\right)
\end{equation*}%
\begin{equation*}
=\overset{N-1}{\underset{s=\varrho }{\sum }}\int_{I_{s}\backslash
I_{s+1}}\left\vert \widetilde{S}^{\ast ,\vartriangle }a\left( x\right)
\right\vert d\mu \left( x\right)
\end{equation*}%
\begin{equation*}
\leq c\overset{N-1}{\underset{s=\varrho +1}{\sum }}\frac{M_{s}}{M_{s}}=c%
\overset{N-1}{\underset{s=\varrho +1}{\sum }}1
\end{equation*}%
\begin{equation*}
\leq \varkappa <c<\infty .
\end{equation*}

Let $0<p<1.$ According to (\ref{1.1}) by using (\ref{dnsnM0}), (\ref%
{dnsnM}) and (\ref{dnsnM2}) we obtain that
\begin{equation*}
\int_{\overline{I_{N}}}\left\vert \widetilde{S}^{\ast ,\vartriangle }a\left(
x\right) \right\vert ^{p}d\mu \left( x\right)
\end{equation*}%
\begin{equation*}
=\overset{N-1}{\underset{s=\varrho }{\sum }}\int_{I_{s}\backslash
I_{s+1}}\left\vert \widetilde{S}^{\ast ,\vartriangle }a\left( x\right)
\right\vert ^{p}d\mu \left( x\right)
\end{equation*}%
\begin{equation*}
\leq c_{p}M_{N}^{1-p}\overset{N-1}{\underset{s=\varrho +1}{\sum }}\frac{%
M_{s}^{p}}{M_{s}}=c_{p}M_{N}^{1-p}\overset{N-1}{\underset{s=\varrho +1}{\sum
}}\frac{1}{M_{s}^{1-p}}
\end{equation*}%
\begin{equation*}
\leq \frac{c_{p}M_{N}^{1-p}}{M_{\varrho }^{1-p}}<c<\infty .
\end{equation*}

The proof of part a) is proved.

Under condition (\ref{ak}) we readily get that
\begin{equation*}
\sup_{\text{ }k\in \mathbb{N}}\frac{M_{\left\vert n_{k}\right\vert }}{%
M_{\left\langle n_{k}\right\rangle }}=\infty .
\end{equation*}

Moreover, there exists a sequence $\left\{ \alpha _{k}:\text{ }k\in \mathbb{N%
}\right\} \subset \left\{ n_{k}:\text{ }k\in \mathbb{N}\right\} $ such that $%
\alpha _{0}\geq 3$ and%
\begin{equation*}
\lim_{k\rightarrow \infty }\frac{M_{\left\vert \alpha _{k}\right\vert
}^{\left( 1-p\right) /2}}{M_{\left\langle \alpha _{k}\right\rangle }^{\left(
1-p\right) /2}}=\infty
\end{equation*}%
and
\begin{equation}
\sum_{k=0}^{\infty }\frac{M_{\left\langle \alpha _{k}\right\rangle }^{\left(
1-p\right) /2}}{M_{\left\vert \alpha _{k}\right\vert }^{\left( 1-p\right) /2}%
}<c<\infty .  \label{charpsn2}
\end{equation}

Let $f=\left( f^{\left( n\right) }:n\in \mathbb{N}\right) $ be the martingale defined in Example %
\ref{example2.6}, where
\begin{equation}
\lambda _{k}=\frac{\lambda M_{\left\langle \alpha _{k}\right\rangle
}^{\left( 1/p-1\right) /2}}{M_{\left\vert \alpha _{k}\right\vert }^{\left(
1/p-1\right) /2}}.  \label{charpsn3}
\end{equation}

Under condition (\ref{charpsn2}) we can conclude that $f\in H_{p}.$

By now using (\ref{10AA}) with $\lambda _{k}$ defined by (\ref{charpsn3})
we readily see that

\begin{equation*}
\widehat{f}(j)
\end{equation*}%
\begin{equation*}
=\left\{
\begin{array}{ll}
M_{\left\langle \alpha _{k}\right\rangle }^{\left( 1/p-1\right)
/2}M_{\left\vert \alpha _{k}\right\vert }^{\left( 1/p-1\right) /2}, & \text{%
\thinspace }j\in \left\{ M_{\left\vert \alpha _{k}\right\vert },...,\text{ ~}%
M_{\left\vert \alpha _{k}\right\vert +1}-1\right\} ,\text{ }k\in \mathbb{N},
\\
0, & j\notin \bigcup\limits_{k=0}^{\infty }\left\{ M_{\left\vert \alpha
_{k}\right\vert },...,\text{ ~}M_{\left\vert \alpha _{k}\right\vert
+1}-1\right\}.%
\end{array}%
\right.
\end{equation*}

Let $M_{\left\vert \alpha _{k}\right\vert }<j<\alpha _{k}.$ By using (\ref%
{11AA}) obtained by $l=k$ in the case when $\lambda _{k}$ are given by
expression (\ref{charpsn3}) we immediately get that
\begin{equation*}
S_{\alpha _{k}}f=S_{M_{\left\vert \alpha _{k}\right\vert }}f+M_{\left\langle
\alpha _{k}\right\rangle }^{\left( 1/p-1\right) /2}M_{\left\vert \alpha
_{k}\right\vert }^{\left( 1/p-1\right) /2}\psi _{M_{\left\vert \alpha
_{k}\right\vert }}D_{\alpha _{k}-M_{\left\vert \alpha _{k}\right\vert }}
\end{equation*}%
\begin{equation*}
:=I+II.
\end{equation*}

According to part a) of Theorem \ref{theorem4.4} for $I$ we have that%
\begin{equation*}
\left\Vert I\right\Vert _{weak-L_{p}}^{p}\leq \left\Vert S_{M_{\left\vert
\alpha _{k}\right\vert }}f\right\Vert _{weak-L_{p}}^{p}\leq c_{p}\left\Vert
f\right\Vert _{H_{p}}^{p}<\infty .
\end{equation*}

Under condition (\ref{ak}) we can conclude that
\begin{equation*}
\left\langle \alpha _{k}\right\rangle \neq \left\vert \alpha _{k}\right\vert
\text{ \ \ \ \ \ \ and \ \ \ \ \ }\left\langle \alpha _{k}-M_{\left\vert
\alpha _{k}\right\vert }\right\rangle =\left\langle \alpha _{k}\right\rangle
.
\end{equation*}

Let $%
x\in I_{\left\langle \alpha _{k}\right\rangle }\backslash I_{\left\langle
\alpha _{k}\right\rangle +1}.$ By using Lemma \ref{dn2.6.2} we obtain that
\begin{equation*}
\left\vert D_{\alpha _{k}-M_{\left\vert \alpha _{k}\right\vert }}\right\vert
\geq cM_{\left\langle \alpha _{k}\right\rangle }
\end{equation*}%
and%
\begin{equation*}
\left\vert II\right\vert =M_{\left\langle \alpha _{k}\right\rangle }^{\left(
1/p-1\right) /2}M_{\left\vert \alpha _{k}\right\vert }^{\left( 1/p-1\right)
/2}\left\vert D_{\alpha _{k}-M_{\left\vert \alpha _{k}\right\vert
}}\right\vert \geq M_{\left\langle \alpha _{k}\right\rangle }^{\left(
1/p+1\right) /2}M_{\left\vert \alpha _{k}\right\vert }^{\left( 1/p-1\right)
/2}.
\end{equation*}

It follows that%
\begin{equation*}
\left\Vert II\right\Vert _{weak-L_{p}}^{p}
\end{equation*}%
\begin{equation*}
\geq c_{p}\left( M_{\left\langle \alpha _{k}\right\rangle }^{\left(
1/p+1\right) /2}M_{\left\vert \alpha _{k}\right\vert }^{\left( 1/p-1\right)
/2}\right) ^{p}\mu \left\{ x\in G_{m}:\text{ }\left\vert II\right\vert \geq
c_{p}M_{\left\langle \alpha _{k}\right\rangle }^{\left( 1/p+1\right)
/2}M_{\left\vert \alpha _{k}\right\vert }^{\left( 1/p-1\right) /2}\right\}
\end{equation*}%
\begin{equation*}
\geq c_{p}M_{\left\vert \alpha _{k}\right\vert }^{\left( 1-p\right)
/2}M_{\left\langle \alpha _{k}\right\rangle }^{\left( 1+p\right) /2}\mu
\left\{ I_{_{\left\langle \alpha _{k}\right\rangle +1}}^{\left\langle \alpha
_{k}\right\rangle -1,\left\langle \alpha _{k}\right\rangle }\right\} \geq
\frac{c_{p}M_{\left\vert \alpha _{k}\right\vert }^{\left( 1-p\right) /2}}{%
M_{\left\langle \alpha _{k}\right\rangle }^{\left( 1-p\right) /2}}.
\end{equation*}

Hence, for large $k$,
\begin{equation*}
\left\Vert S_{\alpha _{k}}f\right\Vert _{weak-L_{p}}^{p}\geq \left\Vert
II\right\Vert _{weak-L_{p}}^{p}-\left\Vert I\right\Vert _{weak-L_{p}}^{p}
\end{equation*}%
\begin{equation*}
\geq \frac{1}{2}\left\Vert II\right\Vert _{weak-L_{p}}^{p}\geq \frac{%
c_{p}M_{\left\vert \alpha _{k}\right\vert }^{\left( 1-p\right) /2}}{%
M_{\left\langle \alpha _{k}\right\rangle }^{\left( 1-p\right) /2}}%
\rightarrow \infty ,\ \text{when}\ \ k\rightarrow \infty .
\end{equation*}

The proof is complete.
\QED

We also mention the following consequences:

\begin{corollary}
\label{corollary4.43} Let $p>0$ and $f\in H_{p}$. Then the maximal operator
\begin{equation*}
\widetilde{S}_{\#}^{\ast }f:=\sup_{n\in \mathbb{N}}\left\vert S_{M_{n}}f\right\vert
\end{equation*}%
is bounded from the Hardy space $H_{p}$ to the space $L_{p}.$
\end{corollary}

\begin{corollary}
\label{corollary4.44} Let $p>0$ and $f\in H_{p}$. Then the maximal operator
\begin{equation*}
\sup_{n\in \mathbb{N_{+}}}\left\vert S_{M_{n}+M_{n-1}}f\right\vert
\end{equation*}%
is bounded from the Hardy space $H_{p}$ to the space $L_{p}.$
\end{corollary}

\begin{corollary}
\label{corollary4.45} Let $p>0$ and $f\in H_{p}$. Then the maximal operator
\begin{equation*}
\sup_{n\in \mathbb{N_{+}}}\left\vert S_{M_{n}+1}f\right\vert
\end{equation*}%
is not bounded from the Hardy space $H_{p}$ to the space $L_{p}.$
\end{corollary}

\subsection{Norm convergence of partial sums of Vilenkin-Fourier series on martingale Hardy spaces}

By applying Corollaries \ref{corollary4.2} and \ref{corollary4.3} we find
necessary and sufficient conditions for the modulus of continuity of
martingale Hardy spaces $H_{p},$ for which the partial sums of
Vilenkin-Fourier series convergence in $L_{p}$-norm. We also study sharpness
of these results. All results in this section can be found in Tephnadze
\cite{tep7}.

\begin{theorem}
\label{theorem4.7bfejera}Let $0<p<1$ and $f\in H_{p}.$ Then there exists an
absolute constant $c_{p}$ depending only on $p$ such that
\begin{equation*}
\left\Vert S_{n}f\right\Vert _{H_{p}}\leq c_{p}n^{1/p-1}\left\Vert
f\right\Vert _{H_{p}}.
\end{equation*}
\end{theorem}

\begin{remark}

We note that the asymptotic behaviour of the sequence $\left\{ n^{1/p-1}:n\in \mathbb{N}\right\} $ in Theorem  \ref{theorem4.7bfejera} can not be improved (c.f. part b) of Theorem \ref{theorem4.1}).

\end{remark}

{\bf Proof}:
According to Corollaries \ref{corollary4.2} \ and  \ref%
{corollary4.43} if we invoke Example \ref{example002} we can conclude
that
\begin{equation}
\left\Vert S_{n}f\right\Vert _{H_{p}}\leq \left\Vert S
^{\#}f\right\Vert _{p}+\left\Vert S_{n}f\right\Vert _{p}  \label{fnmn1partial}
\end{equation}%
\begin{equation*}
\leq \left\Vert f\right\Vert _{H_{p}}+c_{p}n^{1/p-1}\left\Vert f\right\Vert
_{H_{p}}
\end{equation*}%
\begin{equation*}
\leq c_{p}n^{1/p-1}\left\Vert f\right\Vert _{H_{p}}.
\end{equation*}

The proof is complete.
\QED

\begin{theorem}
\label{theorem4.7afejera}Let $f\in H_{1}.$ Then there exists an absolute
constant $c$ such that
\begin{equation*}
\left\Vert S_{n}f\right\Vert _{H_{1}}\leq c\log n\left\Vert
f\right\Vert _{H_{1}}.
\end{equation*}
\end{theorem}

\begin{remark}

We note that the asymptotic behaviour of the sequence $\left\{ \log n:n\in \mathbb{N}\right\} $ in Theorem  \ref{theorem4.7afejera} can not be improved (c.f. part b) of Theorem \ref{theorem4.2}).

\end{remark}

{\bf Proof}:
According to Corollaries \ref{corollary4.3} and \ref%
{corollary4.43} if we invoke Example \ref{example002} we can write that%
\begin{equation*}
\left\Vert S_{n}f\right\Vert _{H_{1}}\leq \left\Vert S
^{\#}f\right\Vert _{1}+\left\Vert S_{n}f\right\Vert _{1}
\end{equation*}%
\begin{equation*}
\leq \left\Vert f\right\Vert _{H_{1}}+c\log n\left\Vert f\right\Vert
_{H_{1}}
\end{equation*}%
\begin{equation*}
\leq c\log n\left\Vert f\right\Vert _{H_{1}}.
\end{equation*}

The proof is complete.
\QED

\begin{theorem}
\label{theorem4.7fejer2}Let $p>0$ and $f\in H_{p}$. Then there exists an
absolute constant $c_{p}$ depending only on $p$ such that
\begin{equation*}
\left\Vert S_{M_{n}}f\right\Vert _{H_{p}}\leq c_{p}\left\Vert
f\right\Vert _{H_{p}}.
\end{equation*}
\end{theorem}

{\bf Proof}:
In the view of Corollary \ref{corollary4.43} and Example \ref%
{example002} we can conclude that
\begin{equation*}
\left\Vert S_{M_{n}}f\right\Vert _{H_{p}}\leq \left\Vert \sup_{0\leq
l\leq n}\left\vert S_{M_{l}}f\right\vert \right\Vert _{p}\leq
\left\Vert f\right\Vert _{H_{p}}.
\end{equation*}

The proof is complete.
\QED

\begin{theorem}
\label{theorem4.5}Let $0<p<1,$ $f\in H_{p}$ and $M_{k}<n\leq M_{k+1}$. Then
there is an absolute constant $c_{p}$ depending only on $p$ such that
\begin{equation*}
\left\Vert S_{n}f -f\right\Vert _{H_{p}}\leq c_{p}n^{1/p-1}\omega _{H_{p}}\left(
\frac{1}{M_{k}},f\right) .
\end{equation*}
\end{theorem}

{\bf Proof}:
Let $0<p<1$ and $M_{k}<n\leq M_{k+1}.$ By using Corollary \ref{corollary4.2}
we immediately get that%
\begin{equation*}
\left\Vert S_{n}f-f\right\Vert _{{H_{p}}}^{p}
\end{equation*}%
\begin{equation*}
\leq \left\Vert S_{n}f-S_{M_{k}}f\right\Vert _{H_{p}}^{p}+\left\Vert
S_{M_{k}}f-f\right\Vert _{{H_{p}}}^{p}
\end{equation*}%
\begin{equation*}
=\left\Vert S_{n}\left( S_{M_{k}}f-f\right) \right\Vert _{H_{p}}^{p}+\left\Vert
S_{M_{k}}f-f\right\Vert _{{H_{p}}}^{p}
\end{equation*}%
\begin{equation*}
\leq c_{p}\left( n^{1-p}+1\right) \left\Vert S_{M_{k}}f-f\right\Vert _{{H_{p}}}^{p}
\end{equation*}%
\begin{equation*}
c_{p}\left( n^{1-p}+1\right)\left\Vert S_{M_{k}}f-f\right\Vert _{{H_{p}}}^{p}
\end{equation*}%
\begin{equation*}
\leq c_{p}n^{1-p}\omega _{H_{p}}^{p}\left( \frac{1}{M_{k}},f\right) .
\end{equation*}

The proof is complete.
\QED

\begin{theorem}
\label{theorem4.6}Let $f\in H_{1}$ and $M_{k}<n\leq M_{k+1}$. Then there is
an absolute constant $c$ such that
\begin{equation*}
\left\Vert S_{n}f-f\right\Vert _{H_{1}}\leq c\lg n\omega _{H_{1}}\left( \frac{1}{%
M_{k}},f\right) .
\end{equation*}
\end{theorem}

{\bf Proof}:
Let $M_{k}<n\leq M_{k+1}.$ Then, by using Corollary \ref{corollary4.3},\ we
immediately get that
\begin{equation*}
\left\Vert S_{n}f-f\right\Vert _{H_{1}}
\end{equation*}%
\begin{equation*}
\leq \left\Vert S_{n}\left( S_{M_{k}}f-f\right) \right\Vert _{H_{1}}+\left\Vert
S_{M_{k}}f-f\right\Vert _{H_{1}}
\end{equation*}%
\begin{equation*}
\leq c\left( \log n+1\right) \left\Vert S_{M_{k}}f-f\right\Vert _{H_{1}}
\end{equation*}%
\begin{equation*}
\leq c\log n\omega _{H_{1}}\left( \frac{1}{M_{k}},f\right) .
\end{equation*}%
The proof is complete.
\QED

\begin{theorem}
\label{theorem4.7}a) Let $0<p<1,$ $f\in H_{p}$ and%
\begin{equation*}
\omega _{H_{p}}\left( \frac{1}{M_{n}},f\right) =o\left( \frac{1}{%
M_{n}^{1/p-1}}\right) ,\text{ when }n\rightarrow \infty .
\end{equation*}%
Then%
\begin{equation*}
\left\Vert S_{k}f-f\right\Vert _{H_{p}}\rightarrow 0,\,\,\,\text{when \ }%
\,\,\,k\rightarrow \infty .
\end{equation*}

b) For every $0<p<1$ there exists a martingale $f\in H_{p}$, for which
\begin{equation*}
\omega _{H_{p}}\left( \frac{1}{M_{n}},f\right) =O\left( \frac{1}{%
M_{n}^{1/p-1}}\right) ,\text{ \ when \ \ \ }n\rightarrow \infty
\end{equation*}%
\textit{and}
\begin{equation*}
\left\Vert S_{k}f-f\right\Vert _{weak-L_{p}}\nrightarrow 0,\,\,\,\text{%
when\thinspace \thinspace \thinspace }k\rightarrow \infty .
\end{equation*}
\end{theorem}

{\bf Proof}:
Let $0<p<1,$ $\ f\in H_{p}$ and%
\begin{equation*}
\omega _{H_{p}}\left( \frac{1}{M_{n}},f\right) =o\left( \frac{1}{%
M_{n}^{1/p-1}}\right) ,\text{ when }n\rightarrow \infty .
\end{equation*}%
By using Theorem \ref{theorem4.5} we immediately get that%
\begin{equation*}
\left\Vert S_{n}f-f\right\Vert _{H_{p}}\rightarrow 0,\text{ when }n\rightarrow
\infty .
\end{equation*}

The proof of part a) is complete.

Let $f=\left( f^{\left( n\right) }:\text{ }n\in \mathbb{N}\right) $ be the martingale defined in Example \ref{example2.5}, where
\begin{equation}
\lambda _{k}=\frac{\lambda }{M_{2\alpha _{k}}^{1/p-1}},\text{ where }\lambda
=\sup_{n\in \mathbb{N}}m_{n}.  \label{charpsn3.1}
\end{equation}

Since
\begin{equation*}
\overset{\infty }{\sum_{k=0}}\frac{\lambda ^{p}}{M_{2\alpha _{k}}^{1-p}}%
<c<\infty
\end{equation*}

we conclude that $f\in H_{p}.$ Moreover, according to (\ref{4AA}) we find
that
\begin{equation*}
\omega _{H_{p}}(\frac{1}{M_{n}},f)=O\left( \overset{\infty }{\sum_{k=n}}%
\frac{1}{M_{k}^{1/p-1}}\right)
\end{equation*}%
\begin{equation*}
=O\left( \frac{1}{M_{n}^{1/p-1}}\right) ,\ when\ \ n\rightarrow \infty .
\end{equation*}%

By now using (\ref{5AA}) with $\lambda _{k}$ defined by (\ref{charpsn3.1})
we readily see that

\begin{equation}
\widehat{f}(j)=\left\{
\begin{array}{ll}
1, & j\in \left\{ M_{2\alpha _{k}},...,\text{ ~}M_{2\alpha _{k}+1}-1\right\}
,\text{ }k\in \mathbb{N}_{+}, \\
0, & j\notin \bigcup\limits_{k=1}^{\infty }\left\{ M_{2\alpha _{k}},...,%
\text{ ~}M_{2\alpha _{k}+1}-1\right\} .\text{ }%
\end{array}%
\right.  \label{29f}
\end{equation}

According to (\ref{29f}) and Corollary \ref{dn2.3} we can establish
that
\begin{equation*}
\limsup\limits_{k\rightarrow \infty }\Vert S_{M_{2\alpha _{k}}+1}f-f\Vert
_{weak-L_{p}}
\end{equation*}%
\begin{equation*}
\geq\limsup\limits_{k\rightarrow \infty }\Vert S_{M_{2\alpha _{k}}}f+\widehat{f}%
(M_{2\alpha _{k}})\psi _{M_{2\alpha _{k}}}-f\Vert _{weak-L_{p}}
\end{equation*}%
\begin{equation*}
\geq \limsup\limits_{k\rightarrow \infty }\left( \Vert \psi _{M_{2\alpha
_{k}}}\Vert _{weak-L_{p}}-\Vert f-S_{M_{2\alpha _{k}}}f\Vert
_{weak-L_{p}}\right)
\end{equation*}%
\begin{equation*}
\geq \limsup\limits_{k\rightarrow \infty }\left( 1-o\left( 1\right) \right)
=1.
\end{equation*}

This completes the proof of the part b) and the proof is complete.
\QED

\begin{theorem}
\label{theorem4.8} a) Let $f\in H_{1}\ $and
\begin{equation*}
\omega _{H_{1}}\left( \frac{1}{M_{n}},f\right) =o\left( \frac{1}{n}\right),%
\text{ when }n\rightarrow \infty .
\end{equation*}%
Then%
\begin{equation*}
\left\Vert S_{k}f-f\right\Vert _{H_{1}}\rightarrow 0,\,\,\,\text{when\thinspace
\thinspace \thinspace }k\rightarrow \infty .
\end{equation*}

b) There exists a martingale $f\in H_{1},$ \ for which%
\begin{equation*}
\omega _{H_{1}}\left( \frac{1}{M_{2M_{n}}},f\right) =O\left( \frac{1}{M_{n}}%
\right) ,\ \ \ \text{when}\ \ n\rightarrow \infty
\end{equation*}%
and%
\begin{equation*}
\left\Vert S_{k}f-f\right\Vert _{1}\nrightarrow 0,\,\,\,\text{when\thinspace
\thinspace \thinspace }k\rightarrow \infty .
\end{equation*}
\end{theorem}

{\bf Proof}:
Let $f\in H_{1}$ and%
\begin{equation*}
\omega _{H_{1}}\left( \frac{1}{M_{n}},f\right) =o\left( \frac{1}{n}\right)
,\ \ when\ \ n\rightarrow \infty .
\end{equation*}

By using Corollary \ref{corollary4.3} we immediately get that%
\begin{equation*}
\left\Vert S_{n}f-f\right\Vert _{H_{1}}\rightarrow \infty ,\text{ when }%
n\rightarrow \infty .
\end{equation*}

For the proof of part b) we use the martingale defined in Example \ref%
{example2.6.1}, for $p=1.$ Then,$\ $ $f\in H_{1}\ $ and by applying (\ref{13m}%
) we can conclude that
\begin{equation*}
\omega _{H_{1}}(\frac{1}{M_{n}},f)=O\left( \frac{1}{n}\right) ,\ \ when\ \
n\rightarrow \infty .
\end{equation*}

By combining (\ref{13}) and (\ref{13s}) with Corollaries \ref{dn5} and \ref%
{dn2.3} we find that
\begin{equation*}
\limsup\limits_{k\rightarrow \infty }\Vert S_{q_{M_{k}}}f-f\Vert
_{1}
\end{equation*}%
\begin{equation*}
=\limsup\limits_{k\rightarrow \infty }\Vert \frac{\psi
_{M_{M_{k}}}D_{q_{M_{k}}}}{M_{2k}}+S_{M_{M_{k}}}f-f\Vert _{1}
\end{equation*}%
\begin{equation*}
-\limsup\limits_{k\rightarrow \infty }\left( \frac{1}{M_{2k}}\Vert
D_{q_{M_{k}}}\Vert _{1}+\Vert S_{M_{M_{k}}}f-f\Vert _{1}\right)
\end{equation*}%
\begin{equation*}
\geq c-\limsup\limits_{k\rightarrow \infty }\left( \overset{\infty }{%
\sum_{i=k+1}}\frac{1}{M_{2i}}-\frac{1}{M_{2k}}\right) =c>0.
\end{equation*}

The proof is complete.
\QED

\subsection{Strong convergence of partial sums of Vilenkin-Fourier series on martingale Hardy spaces}

In this section we prove a Hardy type inequality for partial sums with respect to
Vilenkin systems. This result was first proved by Simon \cite{si1}. Here we
use some new estimations and present a simpler proof, which is due to in
Tephnadze \cite{tep7}. We also show sharpness of this result in the special
sense (see Tephnadze \cite{tep4}).

\begin{theorem}
\label{theorem4.10}(Simon \cite{si1}) a) Let $0<p<1\ $and $f\in H_{p}.$ Then
there is an absolute constant $c_{p},$ depending only on $p,$ such that
\begin{equation*}
\overset{\infty }{\underset{k=1}{\sum }}\frac{\left\Vert S_{k}f\right\Vert
_{p}^{p}}{k^{2-p}}\leq c_{p}\left\Vert f\right\Vert _{H_{p}}^{p}.
\end{equation*}%
b) Let $0<p<1$ and $\left\{ \Phi _{n}:n\in \mathbb{N}\right\} $ be any
non-decreasing sequence satisfying the condition
\begin{equation}
\underset{n\rightarrow \infty }{\overline{\lim }}\Phi _{n}=+\infty .
\label{2ca}
\end{equation}%
Then there exists a martingale $f\in H_{p}$ such that
\begin{equation*}
\text{ }\underset{k=1}{\overset{\infty }{\sum }}\frac{\left\Vert
S_{k}f\right\Vert _{weak-L_{p}}^{p}\Phi _{k}}{k^{2-p}}=\infty.
\end{equation*}
\end{theorem}

{\bf Proof}:
Let $0<p<1.$ By applying (\ref{1.1}) with (\ref{13AA}) in Theorem\ \ref%
{theorem4.1} we have that%
\begin{equation*}
\overset{\infty }{\underset{k=M_{N}}{\sum }}\frac{\left\Vert
S_{k}a\right\Vert _{p}^{p}}{k^{2-p}}
\end{equation*}%
\begin{equation*}
\leq \overset{\infty }{\underset{k=M_{N}}%
{\sum }}\frac{1}{k}\int_{\overline{I_{N}}}\left\vert \frac{S_{k}a\left(
x\right) }{k^{1/p-1}}\right\vert ^{p}d\mu \left( x\right)
\end{equation*}%
\begin{equation*}
=\overset{\infty }{\underset{k=M_{N}}{\sum }}\frac{1}{k}\overset{N-1}{%
\underset{s=0}{\sum }}\int_{I_{s}\backslash I_{s+1}}\left\vert \frac{%
S_{k}a\left( x\right) }{k^{1/p-1}}\right\vert ^{p}d\mu \left( x\right)
\end{equation*}%
\begin{equation*}
\leq c\overset{\infty }{\underset{k=M_{N}}{\sum }}\frac{1}{k}\overset{N-1}{%
\underset{s=0}{\sum }}\int_{I_{s}\backslash I_{s+1}}\left\vert \frac{%
M_{N}^{1/p-1}M_{s}}{k^{1/p-1}}\right\vert ^{p}d\mu
\end{equation*}%
\begin{equation*}
\leq c_{p}M_{N}^{1-p}\overset{\infty }{\underset{k=M_{N}}{\sum }}\frac{1}{%
k^{2-p}}\overset{N-1}{\underset{s=0}{\sum }}\int_{I_{s}\backslash
I_{s+1}}M_{s}^{p}
\end{equation*}%
\begin{equation*}
\leq c_{p}M_{N}^{1-p}\overset{\infty }{\underset{k=M_{N}}{\sum }}\frac{1}{%
k^{2-p}}\overset{N-1}{\underset{s=0}{\sum }}M_{s}^{p-1}d\mu
\end{equation*}%
\begin{equation*}
+c_{p}M_{N}^{1-p}\overset{\infty }{\underset{k=M_{N}}{\sum }}\frac{1}{k^{2-p}%
}\leq c_{p}<\infty .
\end{equation*}

This completes the proof of the part a).

Next, we note that under condition (\ref{2ca}) there exists an increasing
sequence $\left\{ \alpha _{k}\geq 2:k\in \mathbb{N}_{+}\right\} $ of
positive integers such that%
\begin{equation}
\sum_{k=1}^{\infty }\frac{1}{\Phi _{M_{_{2\alpha _{k}}}}^{p/4}}<\infty .
\label{36.1}
\end{equation}

Let $f=\left( f^{\left( k\right) }:k\in \mathbb{N}\right) $ be the martingale defined in Example %
\ref{example2.5.2}. By using (\ref{36.1}) we conclude that the martingale $%
\,f\in H_{p}.$

Let $M_{\alpha _{k}}\leq j<M_{\alpha _{k}+1}.$ By using (\ref{3.7}) if we
repeat the steps which leads to (\ref{11AA}) obtained by $l=k$ in the case
when $\lambda _{k}=1/\Phi _{M_{2\alpha _{\eta }}}^{1/4}$ we obtain that%
\begin{equation*}
S_{j}f=\sum_{\eta =0}^{k-1}\frac{M_{2\alpha _{\eta }}^{1/p-1}}{\Phi
_{M_{2\alpha _{\eta }}}^{1/4}}\left( D_{M_{_{2\alpha _{\eta
}+1}}}-D_{M_{2_{\alpha _{\eta }}}}\right)
\end{equation*}%
\begin{equation*}
+\frac{M_{2\alpha _{k}}^{1/p-1}\psi _{M_{2\alpha _{k}}}D_{j-M_{_{2\alpha
_{k}}}}}{\Phi _{M_{2\alpha _{k}}}^{1/4}}
\end{equation*}%
\begin{equation*}
:=I+II.
\end{equation*}

We calculate each term separately. By using Corollary \ref%
{dn2.3} with condition $\alpha _{n}\geq 2$ $\left( n\in \mathbb{N}\right) $
for $I$ we find that%
\begin{equation*}
D_{M_{\alpha _{n}}}=0,\text{ for }x\in G_{m}\backslash I_{1}
\end{equation*}%
and%
\begin{equation*}
I=0,\text{ \qquad for }x\in G_{m}\backslash I_{1}.
\end{equation*}

Denote by $\mathbb{N}_{0}$ the subset of positive integers $\mathbb{N}_{+},$
for which $\left\langle n\right\rangle =0.$ Then every $n\in \mathbb{N}%
_{0},$ $M_{k}<n<$ $M_{k+1}$ $\left( k>1\right) $ can be written as
\begin{equation*}
n=n_{0}M_{0}+\sum_{j=1}^{k}n_{j}M_{j},
\end{equation*}%
where $n_{0}\in \left\{ 1,...,m_{0}-1\right\} $ and $n_{j}\in \left\{
0,...,m_{j}-1\right\} ,$ $~(j\in \mathbb{N}_{+}).$

Let $j\in \mathbb{N}_{0}$ and $x\in G_{m}\backslash I_{1}=I_{0}\backslash
I_{1}.$ According to Lemma \ref{dn2.6.2} we find that
\begin{equation*}
\left\vert D_{\alpha _{k}-M_{\alpha _{k}}}\right\vert \geq cM_{0}\geq c>0
\end{equation*}%
and
\begin{equation*}
\left\vert II\right\vert =\frac{M_{2\alpha _{k}}^{1/p-1}}{\Phi _{M_{2\alpha
_{k}}}^{1/4}}\left\vert D_{_{j-M_{_{2\alpha _{k}}}}}\left( x\right)
\right\vert =\frac{M_{2\alpha _{k}}^{1/p-1}}{\Phi _{M_{2\alpha _{k}}}^{1/4}}.
\end{equation*}%

It follows that
\begin{equation*}
\left\vert S_{j}f\left( x\right) \right\vert =\left\vert II\right\vert =%
\frac{M_{2\alpha _{k}}^{1/p-1}}{\Phi _{M_{2\alpha _{k}}}^{1/4}},\text{
\qquad for }x\in G_{m}\backslash I_{1},
\end{equation*}%
and%
\begin{equation}
\left\Vert S_{j}f\right\Vert _{weak-L_{p}}  \label{13f}
\end{equation}%
\begin{equation}
\geq \frac{M_{2\alpha _{k}}^{1/p-1}}{2\Phi _{M_{2\alpha _{k}}}^{1/4}}\mu
\left( x\in G_{m}\backslash I_{1}:\left\vert S_{j}f\left( x\right)
\right\vert >\frac{M_{2\alpha _{k}}^{1/p-1}}{2\Phi _{M_{2\alpha _{k}}}^{1/4}}%
\right) ^{1/p}  \notag
\end{equation}%
\begin{equation*}
=\frac{M_{2\alpha _{k}}^{1/p-1}}{2\Phi _{M_{2\alpha _{k}}}^{1/4}}\left\vert
G_{m}\backslash I_{1}\right\vert \geq \frac{cM_{2\alpha _{k}}^{1/p-1}}{\Phi
_{M_{2\alpha _{k}}}^{1/4}}.
\end{equation*}

Since
\begin{equation*}
\underset{\left\{ n\in \mathbb{N}_{0}:M_{k}\leq n\leq M_{k+1},\right\} }{%
\sum }1\geq cM_{k},
\end{equation*}%
where $c$ is an absolute constant, by applying (\ref{13f}) we obtain that%
\begin{equation*}
\underset{j=1}{\overset{M_{2\alpha _{k}+1}-1}{\sum }}\frac{\left\Vert
S_{j}f\right\Vert _{weak-L_{p}}^{p}\Phi _{j}}{j^{2-p}}
\end{equation*}
\begin{equation*}
\geq \underset{%
j=M_{2\alpha _{k}}}{\overset{M_{2\alpha _{k}+1}-1}{\sum }}\frac{\left\Vert
S_{j}f\right\Vert _{weak-L_{p}}^{p}\Phi _{j}}{j^{2-p}}
\end{equation*}%
\begin{equation*}
\geq \Phi _{M_{2\alpha _{k}}}\underset{\left\{ j\in \mathbb{N}%
_{0}:M_{2\alpha _{k}}\leq j\leq M_{2\alpha _{k}+1}\right\} }{\sum }\frac{%
\left\Vert S_{j}f\right\Vert _{weak-L_{p}}^{p}}{j^{2-p}}
\end{equation*}%
\begin{equation*}
\geq c\Phi _{M_{2\alpha _{k}}}\frac{M_{2\alpha _{k}}^{1-p}}{\Phi
_{M_{2\alpha _{k}}}^{p/4}}\underset{\left\{ j\in \mathbb{N}_{0}:M_{2\alpha
_{k}}\leq j\leq M_{2\alpha _{k}+1}\right\} }{\sum }\frac{1}{j^{2-p}}
\end{equation*}%
\begin{equation*}
\geq c\Phi _{M_{\alpha _{k}}}^{3/4}M_{\alpha _{k}}^{1-p}\underset{\left\{
j\in \mathbb{N}_{0}:M_{2\alpha _{k}}\leq j\leq M_{2\alpha _{k}+1}\right\} }{%
\sum }\frac{1}{M_{2\alpha _{k}+1}^{2-p}}
\end{equation*}%
\begin{equation*}
\geq c\frac{\Phi _{M_{2\alpha _{k}}}^{3/4}}{M_{2\alpha _{k}+1}}\underset{%
\left\{ j\in \mathbb{N}_{0}:M_{2\alpha _{k}}\leq j\leq M_{2\alpha
_{k}+1}\right\} }{\sum }1
\end{equation*}%
\begin{equation*}
\geq c\Phi _{M_{2\alpha _{k}}}^{3/4}\rightarrow \infty ,\ \text{when \ \ \ \
}k\rightarrow \infty .
\end{equation*}

Hence, also the part b) is proved so the proof is complete.
\QED

\subsection{{An application concerning estimations of Vilenkin-Fourier
coefficients}\protect\bigskip\protect\bigskip\protect\bigskip}

The following inequalities follow from our results:

\begin{corollary}
\label{corollary3.1} Let $0<p<1$ and $f\in H_{p}.$ Then there exists an
absolute constant $c_{p},$ depending only on $p,$ such that
\begin{equation*}
\left\vert \widehat{f}\left( n\right) \right\vert \leq
c_{p}n^{1/p-1}\left\Vert f\right\Vert _{H_{p}},\text{\ }
\end{equation*}%
\begin{equation*}
\overset{\infty }{\underset{k=1}{\sum }}\frac{\left\vert \widehat{f}\left(
k\right) \right\vert ^{p}}{k^{2-p}}\leq c_{p}\left\Vert f\right\Vert
_{H_{p}}^{p},\text{ \ }\left( f\in H_{p},\text{ }0<p<1\right)
\end{equation*}%
and%
\begin{equation*}
\left( \underset{k=1}{\overset{\infty }{\sum }}M_{k}^{2-2/p}\underset{j=1}{%
\overset{m_{k}-1}{\sum }}\left\vert \widehat{f}\left( jM_{k}\right)
\right\vert ^{2}\right) ^{1/2}\leq c_{p}\left\Vert f\right\Vert _{H_{p}}.
\end{equation*}
\end{corollary}

\begin{remark}
The first inequality was first proved by Tephnadze \cite{tep9}, the second
by Weisz \cite{We4,We1} and the third also by Weisz \cite{we5}. The proofs
are different and simpler then the original ones.
\end{remark}

{\bf Proof}:
Let $0<p<1.$ Since
\begin{equation}
\left\vert \widehat{f}\left( n\right) \right\vert =\left\vert
S_{n+1}f-S_{n}f\right\vert \leq \left\vert S_{n+1}f\right\vert +\left\vert
S_{n}f\right\vert  \label{appl}
\end{equation}%
and
\begin{equation*}
\frac{\left\vert \widehat{f}\left( n\right) \right\vert }{\left( n+1\right)
^{1/p-1}}
\end{equation*}%
\begin{equation*}
\leq \frac{\left\vert S_{n+1}f\right\vert }{\left( n+1\right)
^{1/p-1}}+\frac{\left\vert S_{n}f\right\vert }{n^{1/p-1}}.
\end{equation*}%
By using part a) of Theorem \ref{theorem4.1} we can write that%
\begin{equation*}
\frac{\left\vert \widehat{f}\left( n\right) \right\vert }{\left( n+1\right)
^{1/p-1}}
\end{equation*}

\begin{equation*}
\leq \left\Vert \frac{S_{n+1}f}{\left( n+1\right) ^{1/p-1}}%
\right\Vert _{p}+\left\Vert \frac{S_{n}f}{n^{1-p}}\right\Vert _{p}
\end{equation*}%
\begin{equation*}
\leq c_{p}\left\Vert f\right\Vert _{H_{p}}.
\end{equation*}%
It follows that%
\begin{equation*}
\left\vert \widehat{f}\left( n\right) \right\vert \leq
c_{p}n^{1/p-1}\left\Vert f\right\Vert _{H_{p}}
\end{equation*}%
and the first inequality is proved.

To prove the second inequality we use (\ref{appl}) again. We find that
\begin{equation}
\frac{\left\vert \widehat{f}\left( n\right) \right\vert }{\left( n+1\right)
^{2/p-1}} \label{appl3}
\end{equation}
\begin{equation*}
\leq \frac{\left\vert S_{n+1}f\right\vert }{\left( n+1\right)
^{2/p-1}}+\frac{\left\vert S_{n}f\right\vert }{n^{2/p-1}}.
\end{equation*}%

By combining (\ref{appl3}) and Theorem \ref{theorem4.1} we get that%
\begin{equation*}
\overset{\infty }{\underset{n=1}{\sum }}\frac{\left\vert \widehat{f}\left(
n\right) \right\vert ^{p}}{n^{2-p}}
\end{equation*}
\begin{equation*}
\leq \overset{\infty }{\underset{n=1}{%
\sum }}\left( \left\Vert \frac{S_{n+1}f}{\left( n+1\right) ^{2/p-1}}%
\right\Vert _{p}+\left\Vert \frac{S_{n}f}{n^{2/p-1}}\right\Vert _{p}\right)
^{p}
\end{equation*}%
\begin{equation*}
\leq \overset{\infty }{\underset{n=1}{\sum }}\left( \frac{\left\Vert
S_{n+1}f\right\Vert _{p}^{p}}{\left( n+1\right) ^{2-p}}+\frac{\left\Vert
S_{n}f\right\Vert _{p}^{p}}{n^{2-p}}\right)
\end{equation*}%
\begin{equation*}
\leq 2\overset{\infty }{\underset{n=1}{\sum }}\frac{\left\Vert
S_{n}f\right\Vert _{p}^{p}}{n^{2-p}}\leq c_{p}\left\Vert f\right\Vert
_{H_{p}}^{p}.
\end{equation*}%
and also the second inequality is proved. The third inequality can be proved
analogously so we leave the details.

The proof is complete.
\QED

\newpage

\section{\textbf{Vilenkin-Fej\'er means on martingale Hardy spaces}\protect\bigskip}

\vspace{0.5cm}

\subsection{Some classical results on Vilenkin-Fej\'er means}

In the one-dimensional case Yano \cite{Yano} proved that
\begin{equation*}
\left\Vert K_{n}\right\Vert \leq 2\ \ \ (n\in \mathbb{N}).
\end{equation*}%
Consequently,
\begin{equation*}
\left\Vert \sigma _{n}f-f\right\Vert _{p}\rightarrow 0,\text{ \ \ \ when \ \ \
\ }n\rightarrow \infty ,\text{ \ }(f\in L_{p},\text{ \ }1\leq p\leq \infty ).
\end{equation*}%
However (see \cite{JOO, sws}) the rate of convergence can not be better then
$O\left( n^{-1}\right) $ $\left( n\rightarrow \infty \right) $ for
non-constant functions. a.e, if $f\in L_{p},$ $1\leq p\leq \infty $ and
\begin{equation*}
\left\Vert \sigma _{M_{n}}f-f\right\Vert _{p}=o\left( \frac{1}{M_{n}}\right)
,\text{ \ when \ \ }n\rightarrow \infty ,
\end{equation*}%
then \textit{\ }$f$ \ is a constant function.

Fridli \cite{FR} used dyadic modulus of continuity to characterize the set
of functions in the space $L_{p}$, whose Vilenkin-Fej\'{e}r means converge
at a given rate. It is also known that (see e.g books \cite{AVD} and
\cite{sws})%
\begin{equation*}
\left\Vert \sigma _{n}f-f\right\Vert _{p}
\end{equation*}%
\begin{equation*}
\leq c_{p}\omega _{p}\left( \frac{1}{M_{N}},f\right) +c_{p}\sum_{s=0}^{N-1}%
\frac{M_{s}}{M_{N}}\omega _{p}\left( \frac{1}{M_{s}},f\right) ,\text{ \ }%
\left( 1\leq p\leq \infty ,\text{ \ }n\in \mathbb{N}\right) .
\end{equation*}

By applying this estimate, we immediately obtain that if $f\in lip\left(
\alpha ,p\right) ,$ i.e.,%
\begin{equation*}
\omega _{p}\left( \frac{1}{M_{n}},f\right) =O\left( \frac{1}{M_{n}^{\alpha }}%
\right) ,\text{ \ \ }n\rightarrow \infty ,
\end{equation*}%
then%
\begin{equation*}
\left\Vert \sigma _{n}f-f\right\Vert _{p}=\left\{
\begin{array}{ll}
O\left( \frac{1}{M_{N}}\right) , & \text{if }\alpha >1, \\
O\left( \frac{N}{M_{N}}\right) , & \text{if }\alpha =1, \\
O\left( \frac{1}{M_{N}^{\alpha }}\right) , & \text{if }\alpha <1\text{. }%
\end{array}%
\right.
\end{equation*}

Weisz \cite{We2} considered the norm convergence of Fej\'er means of
Vilenkin-Fourier series and proved that
\begin{equation}
\left\Vert \sigma _{k}f\right\Vert _{p}\leq c_{p}\left\Vert f\right\Vert
_{H_{p}},\text{ \ \ }p>1/2\text{ \ \ and \ \ \ }f\in H_{p}.  \label{f100}
\end{equation}

This result implies that%
\begin{equation*}
\frac{1}{n^{2p-1}}\overset{n}{\underset{k=1}{\sum }}\frac{\left\Vert \sigma
_{k}f\right\Vert _{p}^{p}}{k^{2-2p}}\leq c_{p}\left\Vert f\right\Vert
_{H_{p}}^{p},\text{ \ \ \ }\left( 1/2<p<\infty \right) .
\end{equation*}

If (\ref{f100}) hold for $0<p\leq 1/2,$ then we would have that
\begin{equation}
\frac{1}{\log ^{\left[ 1/2+p\right] }n}\overset{n}{\underset{k=1}{\sum }}%
\frac{\left\Vert \sigma _{k}f\right\Vert _{p}^{p}}{k^{2-2p}}\leq
c_{p}\left\Vert f\right\Vert _{H_{p}}^{p},\text{ \ \ \ }\left( 0<p\leq
1/2\right) .  \label{2cc}
\end{equation}

However, in Tephnadze \cite{tep1} it was shown that the assumption $p>1/2$ in (\ref%
{f100}) is essential. In particular, is was proved that there exists a
martingale $f\in H_{1/2}$ such that
\begin{equation*}
\sup_{n\in \mathbb{N}}\left\Vert \sigma _{n}f\right\Vert _{1/2}=+\infty .
\end{equation*}

For Vilenkin systems in \cite{tep5} it was proved that (\ref{2cc}) holds,
though inequality (\ref{f100}) is not true for $0<p\leq 1/2.$

In the one-dimensional case the weak type inequality
\begin{equation*}
\mu \left( \sigma ^{\ast }f>\lambda \right) \leq \frac{c}{\lambda }%
\left\Vert f\right\Vert _{1},\text{ \qquad }\left( f\in L_{1},\text{ \ \ \ }%
\lambda >0\right)
\end{equation*}%
can be found in Zygmund \cite{13} for the trigonometric series, in Schipp
\cite{Sc} for Walsh series and in P\'al, Simon \cite{PS} for bounded Vilenkin
series. Fujji \cite{Fu} and Simon \cite{Si2} verified that $\sigma ^{\ast }$
is bounded from $H_{1}$ to $L_{1}$. Weisz \cite{We2} generalized this result
and proved the boundedness of $\sigma ^{\ast }$ from the martingale space $%
H_{p}$ to the Lebesgue space $L_{p}$ for $p>1/2$. Simon \cite{Si1} gave a
counterexample, which shows that boundedness does not hold for $0<p<1/2.$
The counterexample for $p=1/2$ due to Goginava \cite{GoAMH}, (see also \cite%
{BGG} and \cite{BGG2}). Weisz \cite{we4} proved that $\sigma ^{\ast }$ is
bounded from the Hardy space $H_{1/2}$ to the space $weak-L_{1/2}$. In \cite%
{tep2} and \cite{tep3} (for Walsh system see \cite{GoSzeged}) it was proved
that the maximal operator $\widetilde{\sigma }_{p}^{\ast }$ with respect to
Vilenkin systems defined by
\begin{equation*}
\widetilde{\sigma }_{p}^{\ast }:=\sup_{n\in \mathbb{N}}\frac{\left\vert
\sigma _{n}\right\vert }{\left( n+1\right) ^{1/p-2}\log ^{2\left[ 1/2+p%
\right] }\left( n+1\right) },
\end{equation*}%
where $0<p\leq 1/2$ and $\left[ 1/2+p\right] $ denotes integer part of $%
1/2+p,$ is bounded from the Hardy space $H_{p}$ to the Lebesgue space $L_{p}.$
Moreover, the order of deviant behavior of the $n$-th Fej\'er mean was given
exactly. As a corollary we get that
\begin{equation*}
\left\Vert \sigma _{n}f\right\Vert _{p}\leq c_{p}\left( n+1\right)
^{1/p-2}\log ^{2\left[ 1/2+p\right] }\left( n+1\right) \left\Vert
f\right\Vert _{H_{p}}.
\end{equation*}

For Walsh-Kaczmarz system analogical theorems were proved in \cite{GNCz} and
\cite{tep4}.

For the one-dimensional Vilenkin-Fourier series Weisz \cite{We2} proved that
the maximal operator
\begin{equation*}
\sigma ^{\#}f=\sup_{n\in \mathbb{N}}\left\vert \sigma _{M_{n}}f\right\vert
\end{equation*}%
is bounded from the martingale Hardy space $H_{p}$ to the Lebesgue space $L_{p}$ for $%
p>0.$ For the Walsh-Fourier series Goginava \cite{gog1} proved that the
operator $\left\vert \sigma _{2^{n}}f\right\vert $ is not bounded from the
space $H_{p}$ to the space $H_{p},$ for $0<p\leq 1.$

\subsection{Divergence of Vilenkin-Fej\'er means on martingale Hardy spaces}

\begin{theorem}
\label{theorem4fejer}a)There exist a martingale $f\in H_{1/2}$ such that
\begin{equation*}
\sup\limits_{n\in \mathbb{N}}\left\Vert \sigma _{n}f\right\Vert
_{1/2}=+\infty .
\end{equation*}%
b) Let $0<p<1/2.$ There exist a martingale $f\in H_{p},$ such that%
\begin{equation*}
\sup\limits_{n\in \mathbb{N}}\left\Vert \sigma _{n}f\right\Vert
_{weak-L_{p}}=+\infty .
\end{equation*}
\end{theorem}

\begin{remark}
This result for $p=1/2$ can be found in Tephnadze \cite{tep1}. By
interpolation automatically follows the second part of this Theorem. Here,
we make a more concrete proof of this fact by even constructing a martingale
$f\in H_{p},$ for which the Vilenkin-Fej\'er means are not bounded in the
space $weak-L_{p}.$
\end{remark}

{\bf Proof}:
Let $f=\left( f^{\left( n\right) }:n\in \mathbb{N}\right) $ be the martingale defined in Example %
\ref{example2.5.1} in the case when $p=q=1/2.$ We can write that%
\begin{equation*}
\sigma _{q_{\alpha _{k}}}f=\frac{1}{q_{\alpha _{k}}}\sum_{j=1}^{M_{2\alpha
_{k}}}S_{j}f+\frac{1}{q_{\alpha _{k}}}\sum_{j=M_{2\alpha _{k}}+1}^{q_{\alpha
_{k}}}S_{j}f
\end{equation*}%
\begin{equation*}
:=I+II.
\end{equation*}

According to (\ref{sn101}) in Example \ref{example2.5.1} we can
conclude that
\begin{equation}
\left\vert I\right\vert \leq \frac{1}{q_{\alpha _{k}}}\sum_{j=1}^{M_{2\alpha
_{k}}}\left\vert S_{i}f\right\vert  \label{1t10}
\end{equation}%
\begin{equation*}
\leq \frac{2\lambda M_{2\alpha _{k-1}}^{2}}{\alpha _{k-1}^{1/2}}\frac{%
M_{2\alpha _{k}}}{q_{\alpha _{k}}}\leq \frac{2\lambda M_{2\alpha _{k-1}}^{2}%
}{\alpha _{k-1}^{1/2}}.
\end{equation*}

By applying (\ref{1t6.1}) obtained by letting $l=k$ we can rewrite $II$ as%
\begin{equation*}
II=\frac{\left( q_{\alpha _{k}}-M_{2\alpha _{k}}\right) S_{M_{2\alpha _{k}}}%
}{q_{\alpha _{k}}}
\end{equation*}%
\begin{equation*}
+\frac{M_{2\alpha _{k}}\psi _{M_{2\alpha _{k}}}}{\alpha _{k}^{1/2}q_{\alpha
_{k}}}\sum_{j=M_{2\alpha _{k}}+1}^{q_{\alpha _{k}}}D_{j-M_{_{2\alpha _{k}}}}
\end{equation*}%
\begin{equation*}
:=II_{1}+II_{2}.
\end{equation*}%
According to (\ref{sn101}) for $j=M_{2\alpha _{k}}$ and $p=1/2$  we find that%
\begin{equation}
\left\vert II_{1}\right\vert \leq \left\vert S_{M_{2\alpha _{k}}}\right\vert
\leq \frac{2\lambda M_{2\alpha _{k-1}}^{2}}{\alpha _{k-1}^{1/2}}.
\label{1t9}
\end{equation}

In view of (\ref{1t10}) and (\ref{1t9}) we invoke estimate (\ref{1t4})
for $\ $ $p=q=1/2$ to conclude that%
\begin{equation*}
\left\vert II_{1}\right\vert \leq \frac{2\lambda
M_{2\alpha _{k-1}}^{2}}{\alpha _{k-1}^{1/2}}\leq \frac{M_{\alpha _{k}}^{2}}{%
16\alpha _{k}^{3/2}}
\end{equation*}
and
\begin{equation*}
\left\vert I\right\vert \leq \frac{2\lambda
M_{2\alpha _{k-1}}^{2}}{\alpha _{k-1}^{1/2}}\leq \frac{M_{\alpha _{k}}^{2}}{%
16\alpha _{k}^{3/2}}.
\end{equation*}

Let%
\begin{equation*}
x\in I_{2\alpha _{k}}^{2\eta ,2\left( \eta +1\right) },\,\,\,\,\,\eta =\left[
\frac{\alpha _{k}}{2}\right] +1,...,\alpha _{k}-3.
\end{equation*}

By applying Lemma \ref{lemma3} we have that%
\begin{equation*}
4q_{\alpha _{k}-1}\left\vert K_{q_{\alpha _{k}-1}}\left( x\right)
\right\vert \geq M_{2\eta }M_{2\left( \eta +1\right) }\geq M_{2\eta }^{2}.
\end{equation*}%

Hence, for $II_{2}$ we readily get that
\begin{equation*}
\left\vert II_{2}\right\vert =\frac{M_{2\alpha _{k}}}{\alpha
_{k}^{1/2}q_{\alpha _{k}}}\left\vert \psi _{M_{2\alpha
_{k}}}\sum_{j=1}^{q_{\alpha _{k}-1}}D_{j}\right\vert
\end{equation*}%
\begin{equation*}
=\frac{M_{2\alpha _{k}}}{q_{\alpha _{k}}}\frac{q_{\alpha _{k}-1}}{\alpha
_{k}^{1/2}}\left\vert K_{q_{\alpha _{k}}-1}\right\vert
\end{equation*}%
\begin{equation*}
\geq \frac{q_{\alpha _{k}-1}\left\vert K_{q_{\alpha _{k}}-1}\right\vert }{%
2\alpha _{k}^{1/2}}\geq \frac{M_{2\eta }^{2}}{8\alpha _{k}^{1/2}}.
\end{equation*}

Since $M_{2\eta }\geq M_{\alpha _{k}}$ we obtain that%
\begin{equation}
\left\vert \sigma _{q_{\alpha _{k}}}f\left( x\right) \right\vert \geq
\left\vert II_{2}\right\vert -\left( \left\vert I\right\vert +\left\vert
II_{1}\right\vert \right)  \label{1t11}
\end{equation}%
\begin{equation*}
\geq \frac{1}{8\alpha _{k}^{1/2}}\left( M_{2\eta }^{2}-\frac{M_{\alpha
_{k}}^{2}}{2\alpha _{k}}\right) .
\end{equation*}

From (\ref{1t11}) it follows that
\begin{equation*}
\left\vert \sigma _{q_{\alpha _{k}}}f\left( x\right) \right\vert \geq \frac{%
cM_{2\eta }^{2}}{\alpha _{k}^{1/2}},\,\,x\in I_{2\alpha _{k}}^{2\eta
,2\left( \eta +1\right) },
\end{equation*}%
where%
\begin{equation*}
\eta =\left[ \frac{\alpha _{k}}{2}\right] +1,...,\alpha _{k}-3,\,
\end{equation*}

Thus,
\begin{equation*}
\int_{G_{m}}\left\vert \sigma _{q_{\alpha _{k}}}f\left( x\right) \right\vert
^{1/2}d\mu \left( x\right)
\end{equation*}%
\begin{equation*}
\geq \sum_{\eta =\left[ \alpha _{k}/2\right] +1}^{\alpha
_{k}-3}\sum_{x_{2\eta +3}=0}^{m_{2\eta +3}-1}...\sum_{x_{_{2\alpha
_{k}-1}}=0}^{m_{2\alpha _{k}-1}-1}\int_{I_{2\alpha _{k}}^{2\eta ,2\left(
\eta +1\right) }}\left\vert \sigma _{q_{\alpha _{k}}}f\left( x\right)
\right\vert ^{1/2}d\mu \left( x\right)
\end{equation*}%
\begin{equation*}
\geq \frac{c}{\alpha _{k}^{1/4}}\sum_{\eta =\left[ \alpha _{k}/2\right]
+1}^{\alpha _{k}-3}\sum_{x_{2\eta +3}=0}^{m_{2\eta
+3}-1}...\sum_{x_{_{2\alpha _{k}-1}}=0}^{m_{2\alpha _{k}-1}-1}\left\vert
I_{2\alpha _{k}}^{2\eta ,2\left( \eta +1\right) }\right\vert M_{2\eta }
\end{equation*}%
\begin{equation*}
\geq \frac{c}{\alpha _{k}^{1/4}}\sum_{\eta =\left[ \alpha _{k}/2\right]
+1}^{\alpha _{k}-3}\frac{m_{2\eta +3...}m_{2\alpha _{k}-1}}{M_{2\alpha _{k}}}%
M_{2\eta }
\end{equation*}%
\begin{equation*}
\geq \frac{c}{\alpha _{k}^{1/4}}\sum_{\eta =\left[ \alpha _{k}/2\right]
+1}^{\alpha _{k}-3}\frac{M_{2\eta }}{M_{2\eta +2}}
\end{equation*}%
\begin{equation*}
\geq \frac{c}{\alpha _{k}^{1/4}}\sum_{\eta =\left[ \alpha _{k}/2\right]
+1}^{\alpha _{k}-3}1
\end{equation*}%
\begin{equation*}
\geq \alpha _{k}^{3/4}\rightarrow \infty ,\ \ when\ \ k\longrightarrow \infty .
\end{equation*}

The proof of part a) is complete so we turn to the proof of part b). Let $%
0<p<1/2.$ Let $f=\left( f^{\left( n\right) }:n\in \mathbb{N}\right) $ be the martingale defined in Example \ref{example2.5.1} in the case when $0<p<q=1/2.$ We can write
that
\begin{equation*}
\sigma _{M_{2\alpha _{k}}+1}f
\end{equation*}%
\begin{equation*}
=\frac{1}{M_{2\alpha _{k}}+1}\sum_{j=0}^{M_{2\alpha _{k}}}S_{j}f+\frac{%
S_{M_{2\alpha _{k}}+1}f}{M_{2\alpha _{k}}+1}
\end{equation*}%
\begin{equation*}
:=III+IV.
\end{equation*}

We combine (\ref{sn101}) and (\ref{sn102}) and invoke (\ref{1t4}) in the
case when $p<q=1/2$ to obtain the following estimates:
\begin{equation*}
\left\vert III\right\vert \leq \frac{M_{2\alpha _{k}}}{M_{2\alpha _{k}}+1}%
\frac{2\lambda M_{2\alpha _{k-1}}^{1/p}}{\alpha _{k-1}^{1/2}}
\end{equation*}%
\begin{equation*}
\leq \frac{2\lambda M_{2\alpha _{k-1}}^{1/p}}{\alpha _{k-1}^{1/2}}\leq \frac{%
M_{\alpha _{k}}^{1/p-2}}{16\alpha _{k}^{3/2}}
\end{equation*}%
and
\begin{equation*}
\left\vert IV\right\vert \geq \frac{\left\vert S_{M_{2\alpha
_{k}}+1}f\right\vert }{M_{2\alpha _{k}}+1}\geq \frac{M_{2\alpha _{k}}^{1/p-2}%
}{2\alpha _{k}}
\end{equation*}

Let $x\in G_{m}.$ We conclude that
\begin{equation*}
\left\vert \sigma _{M_{2\alpha _{k}}+1}f\left( x\right) \right\vert \geq
\left\vert IV\right\vert -\left\vert III\right\vert
\end{equation*}%
\begin{equation*}
\geq \frac{M_{2\alpha _{k}}^{1/p-2}}{2\alpha _{k}^{1/2}}-\frac{M_{\alpha
_{k}}^{1/p-2}}{16\alpha _{k}^{3/2}}\geq \frac{M_{2\alpha _{k}}^{1/p-2}}{%
4\alpha _{k}^{1/2}}.
\end{equation*}

It follows that%
\begin{equation*}
\frac{M_{2\alpha _{k}}^{1/p-2}}{4\alpha _{k}^{1/2}}\left( \mu \left\{ x\in
G_{m}:\left\vert \sigma _{M_{_{2\alpha _{k}}}+1}f\left( x\right) \right\vert
\geq \frac{M_{2\alpha _{k}}^{1/p-2}}{4\alpha _{k}^{1/2}}\right\} \right)
^{1/p}
\end{equation*}

\begin{equation*}
\geq \frac{M_{2\alpha _{k}}^{1/p-2}}{4\alpha _{k}^{1/2}}\rightarrow \infty ,%
\text{ when }k\rightarrow \infty .
\end{equation*}

Hence, also part b) is proved so the proof is complete.
\QED

\subsection{Maximal operators of Vilenkin-Fej\'er means on martingale
Hardy spaces}

In this subsection we consider weighted Maximal operators of Fej\'er means of
Vilenkin-Fourier series and prove $(H_{p},L_{p})$ and $(H_{p},weak-L_{p})$
type inequalities. In all cases we also show sharpness in a special sense.

First we note the following consequence of Theorem \ref{theorem4fejer}:

\begin{corollary}
a) There exist a martingale $f\in H_{1/2}$ such that
\begin{equation*}
\left\Vert \sigma ^{\ast }f\right\Vert _{1/2}=+\infty .
\end{equation*}%
b) Let $0<p<1/2.$ There exist a martingale $f\in H_{p}$ such that
\begin{equation*}
\left\Vert \sigma ^{\ast }f\right\Vert _{weak-L_{p}}=+\infty .
\end{equation*}
\end{corollary}

The next theorem can be found in Tephnadze \cite{tep3}, but here we will
give a simpler proof of part b).

\begin{theorem}
\label{theorem3fejermax2}a) Let $0<p<1/2.$ Then the \bigskip maximal operator%
\begin{equation*}
\widetilde{\sigma }_{p}^{\ast }f:=\sup_{n\in \mathbb{N}}\frac{\left\vert
\sigma _{n}f\right\vert }{\left( n+1\right) ^{1/p-2}}
\end{equation*}%
is bounded from the Hardy martingale space $H_{p}$ to the Lebesgue space $L_{p}.$

b) Let $\left\{ \Phi _{n}:n\in \mathbb{N}\right\} $ be any non-decreasing
sequence satisfying the condition%
\begin{equation}
\overline{\lim_{n\rightarrow \infty }}\frac{\left( n+1\right) ^{1/p-2}}{\Phi
_{n}}=+\infty .  \label{cond2}
\end{equation}%
Then
\begin{equation*}
\sup_{k\in \mathbb{N}}\frac{\left\Vert \frac{\sigma _{M_{_{2n_{k}}}+1}f_{k}}{%
\Phi _{M_{_{2n_{k}}}+1}}\right\Vert _{weal-L_{p}}}{\left\Vert
f_{k}\right\Vert _{H_{p}}}=\infty .
\end{equation*}
\end{theorem}

{\bf Proof}:
First place we note that $\sigma _{n}$ is bounded from $L_{\infty }$
to $L_{\infty }$ (see (\ref{fn4}) in Corollary \ref{lemma7kn}). Hence, by
Lemma \ref{lemma2.2} the proof of Theorem \ref{theorem3fejermax2} will be
complete if we show that%
\begin{equation*}
\int\limits_{\overline{I}_{N}}\left( \underset{n\in \mathbb{N}}{\sup }\frac{%
\left\vert \sigma _{n}a\right\vert }{\left( n+1\right) ^{1/p-2}}\right)
^{p}d\mu \leq c<\infty
\end{equation*}%
for every p-atom $a,$ where $I$ denotes the support of the atom$.$

Let $a$ be an arbitrary p-atom with support$\ I$ and $\mu \left( I\right)
=M_{N}^{-1}.$ We may assume that $I=I_{N}.$ It is easy to see that $\sigma
_{n}\left( a\right) =0$ when $n\leq M_{N}$. Therefore we can suppose that $%
n>M_{N}$.

Since $\left\Vert a\right\Vert _{\infty }\leq M_{N}^{1/p}$ it follows that
\begin{equation*}
\frac{\left\vert \sigma _{n}\left( a\right) \right\vert }{\left( n+1\right)
^{1/p-2}}
\end{equation*}%
\begin{equation*}
\leq \frac{1}{\left( n+1\right) ^{1/p-2}}\int_{I_{N}}\left\vert a\left(
t\right) \right\vert \left\vert K_{n}\left( x-t\right) \right\vert d\mu
\left( t\right)
\end{equation*}%
\begin{equation*}
\leq \frac{\left\Vert a\right\Vert _{\infty }}{\left( n+1\right) ^{1/p-2}}%
\int_{I_{N}}\left\vert K_{n}\left( x-t\right) \right\vert d\mu \left(
t\right)
\end{equation*}%
\begin{equation*}
\leq \frac{cM_{N}^{1/p}}{\left( n+1\right) ^{1/p-2}}\int_{I_{N}}\left\vert
K_{n}\left( x-t\right) \right\vert d\mu \left( t\right) .
\end{equation*}

Let $x\in I_{N}^{k,l},\,0\leq k<l\leq N.$ From Corollary \ref{lemma5aa} we
can deduce that
\begin{equation}
\frac{\left\vert \sigma _{n}\left( a\right) \right\vert }{\left( n+1\right)
^{1/p-2}}  \label{1112}
\end{equation}%
\begin{equation*}
\leq \frac{c_{p}M_{N}^{1/p}}{M_{N}^{1/p-2}}\frac{M_{l}M_{k}}{M_{N}^{2}}%
=c_{p}M_{l}M_{k}.
\end{equation*}

The expression on the right-hand side of (\ref{1112}) does not depend on $n.$ Thus,
\begin{equation}
\left\vert \widetilde{\sigma }_{p}^{\ast }a\left( x\right) \right\vert \leq
c_{p}M_{l}M_{k},\text{ \ for \ }x\in I_{N}^{k,l},\,0\leq k<l\leq N.
\label{1112M}
\end{equation}%

By using (\ref{1112M}) with identity (\ref{1.1}) we obtain that
\begin{equation*}
\int_{\overline{I_{N}}}\left\vert \widetilde{\sigma }_{p}^{\ast }a\left(
x\right) \right\vert ^{p}d\mu \left( x\right)
\end{equation*}%
\begin{equation*}
=\overset{N-2}{\underset{k=0}{\sum }}\overset{N-1}{\underset{l=k+1}{\sum }}%
\sum\limits_{x_{j}=0,j\in
\{l+1,...,N-1\}}^{m_{j-1}}\int_{I_{N}^{k,l}}\left\vert \widetilde{\sigma }%
_{p}^{\ast }a\left( x\right) \right\vert ^{p}d\mu \left( x\right)
\end{equation*}%
\begin{equation*}
+\overset{N-1}{\underset{k=0}{\sum }}\int_{I_{N}^{k,N}}\left\vert \widetilde{%
\sigma }_{p}^{\ast }a\left( x\right) \right\vert ^{p}d\mu \left( x\right)
\end{equation*}%
\begin{equation*}
\leq c_{p}\overset{N-2}{\underset{k=0}{\sum }}\overset{N-1}{\underset{l=k+1}{%
\sum }}\frac{m_{l+1}...m_{N-1}}{M_{N}}M_{l}^{p}M_{k}^{p}
\end{equation*}%
\begin{equation*}
+c_{p}\overset{N-1}{\underset{k=0}{\sum }}\frac{1}{M_{N}}M_{N}^{p}M_{k}^{p}
\end{equation*}%
\begin{equation*}
\leq c_{p}\overset{N-2}{\underset{k=0}{\sum }}\overset{N-1}{\underset{l=k+1}{%
\sum }}\frac{M_{l}^{p}M_{k}^{p}}{M_{l}}+c_{p}\overset{N-1}{\underset{k=0}{%
\sum }}\frac{M_{k}^{p}}{M_{N}^{1-p}}
\end{equation*}%
\begin{equation*}
:=I+II.
\end{equation*}

We estimate each term separately.
\begin{equation}
I=c_{p}\overset{N-2}{\underset{k=0}{\sum }}\overset{N-1}{\underset{l=k+1}{%
\sum }}\frac{1}{M_{l}^{1-2p}}\frac{M_{l}^{p}M_{k}^{p}}{M_{l}^{2p}} \label{nec3}
\end{equation}%
\begin{equation*}
\leq c_{p}\overset{N-2}{\underset{k=0}{\sum }}\overset{N-1}{\underset{l=k+1}{%
\sum }}\frac{1}{M_{l}^{1-2p}}
\end{equation*}%
\begin{equation*}
\leq c_{p}\overset{N-2}{\underset{k=0}{\sum }}\overset{N-1}{\underset{l=k+1}{%
\sum }}\frac{1}{2^{l\left( 1-2p\right) }}
\end{equation*}%
\begin{equation*}
\leq c_{p}\overset{N-2}{\underset{k=0}{\sum }}\frac{1}{2^{k\left(
1-2p\right) }}<c_{p}<\infty .
\end{equation*}

It is obvious that
\begin{equation}
II\leq \frac{c_{p}}{M_{N}^{1-2p}}\overset{N-1}{\underset{k=0}{\sum }}\frac{%
M_{k}^{p}}{M_{N}^{p}}  \label{nec4}
\end{equation}%
\begin{equation*}
\leq \frac{c_{p}}{M_{N}^{1-2p}}<c<\infty .
\end{equation*}

The proof of the part a) is complete by combining the estimates above.

Let $0<p<1/2.$ Under condition (\ref{cond2}) there exists an increasing
sequence of positive integers $\left\{ \lambda _{k}:\text{ }k\in \mathbb{N}%
\right\} $ such that
\begin{equation*}
\lim_{k\rightarrow \infty }\frac{\lambda _{k}^{1/p-2}}{\Phi _{\lambda _{k}}}%
=\infty .
\end{equation*}%

It is evident that for every $\lambda _{k}$ there exists a positive integers
$m_{k}^{,}$ such that $q_{_{m_{k}^{^{\prime }}}}<\lambda
_{k}<2q_{_{m_{k}^{,}}}.$ Since $\left\{ \Phi _{n}:n\in \mathbb{N}\right\} $
is a non-decreasing function we have that
\begin{equation*}
\overline{\underset{k\rightarrow \infty }{\lim }}\frac{M_{2m_{k}^{,}}^{1/p-2}%
}{\Phi _{M_{2m_{k}^{,}}+1}}
\end{equation*}%
\begin{equation*}
\geq \frac{1}{2}\overline{\underset{k\rightarrow \infty }{\lim }}\frac{%
\left( M_{2m_{k}^{,}}+1\right) ^{1/p-2}}{\Phi _{M_{2m_{k}^{,}}+1}}
\end{equation*}%
\begin{equation*}
\geq c\lim_{k\rightarrow \infty }\frac{\lambda _{k}^{1/p-2}}{\Phi _{\lambda
_{k}}}=\infty .
\end{equation*}

Let$\ \left\{ n_{k}:\text{ }k\in \mathbb{N}\right\} \subset \left\{
m_{k}^{,}:\text{ }k\in \mathbb{N}\right\} $ be a sequence of positive numbers
such that
\begin{equation*}
\lim_{k\rightarrow \infty }\frac{M_{2n_{k}}^{1/p-2}}{\Phi _{M_{_{2n_{k}}}+1}}%
=\infty
\end{equation*}%
and $f_{k}$ be the atom defined in Example \ref{example2.2}.

By combining (\ref{13l}) and (\ref{14l}) in Example \ref{example2.2} we find
that
\begin{equation*}
\frac{\left\vert \sigma _{M_{_{2n_{k}}}+1}f_{k}\right\vert }{\Phi
_{M_{_{2n_{k}}}+1}}
\end{equation*}%
\begin{equation*}
=\frac{1}{\Phi _{M_{_{2n_{k}}}+1}\left( M_{_{2n_{k}}}+1\right) }\left\vert
\overset{M_{_{2n_{k}}}+1}{\underset{j=0}{\sum }}S_{j}f_{k}\right\vert
\end{equation*}%
\begin{equation*}
=\frac{1}{\Phi _{M_{_{2n_{k}}}+1}\left( M_{_{2n_{k}}}+1\right) }\left\vert
S_{M_{_{2n_{k}}}+1}f_{k}\right\vert
\end{equation*}%
\begin{equation*}
=\frac{1}{\Phi _{M_{_{2n_{k}}}+1}\left( M_{_{2n_{k}}}+1\right) }\left\vert
D_{M_{_{2n_{k}}}+1}-D_{M_{_{2n_{k}}}}\right\vert
\end{equation*}%
\begin{equation*}
=\frac{1}{\Phi _{M_{_{2n_{k}}}+1}\left( M_{_{2n_{k}}}+1\right) }\left\vert
\psi _{M_{_{2n_{k}}}}\right\vert
\end{equation*}%
\begin{equation*}
=\frac{1}{\Phi _{M_{_{2n_{k}}}+1}\left( M_{_{2n_{k}}}+1\right) }
\end{equation*}%
\begin{equation*}
\geq \frac{c}{M_{_{_{_{2n_{k}}}}}\Phi _{M_{_{2n_{k}}}+1}}.
\end{equation*}

Hence,%
\begin{equation*}
\mu \left\{ x\in G_{m}:\frac{\left\vert \sigma _{M_{_{2n_{k}}}+1}f_{k}\left(
x\right) \right\vert }{\Phi _{M_{_{2n_{k}}}+1}}\geq \frac{c}{%
M_{_{_{_{2n_{k}}}}}\Phi _{M_{_{2n_{k}}}+1}}\right\}
\end{equation*}%
\begin{equation*}
\geq \mu \left( G_{m}\right) =1.\qquad
\end{equation*}

Therefore, by using (\ref{15l}) in Example \ref{example2.2} we get that

\begin{equation*}
\frac{\frac{c}{M_{_{_{_{2n_{k}}}}}\Phi _{M_{_{2n_{k}}}+1}}\left( \mu \left\{
x\in G_{m}:\frac{\left\vert \sigma _{M_{_{2n_{k}}}+1}f_{k}\left( x\right)
\right\vert }{\Phi _{M_{_{2n_{k}}}+1}}\geq \frac{c}{M_{_{_{_{2n_{k}}}}}\Phi
_{M_{_{2n_{k}}}+1}}\right\} \right) ^{1/p}}{\left\Vert f_{k}\right\Vert
_{H_{p}}}
\end{equation*}%
\begin{equation*}
\geq \frac{c}{M_{_{_{_{2n_{k}}}}}\Phi _{M_{_{2n_{k}}}+1}}%
M_{_{2n_{k}}}^{1/p-1}=\frac{cM_{_{2n_{k}}}^{1/p-2}}{\Phi _{M_{_{2n_{k}}}+1}}%
\rightarrow \infty \text{,\qquad when }k\rightarrow \infty .
\end{equation*}

Thus, also part b) is proved so the proof is complete.
\QED

We also point out the following consequence of Theorem \ref{theorem3fejermax2}, which we need later on.

\begin{corollary}
\label{corollary3fejermax1} a) Let $0<p<1/2$ and $f\in H_{p}.$ Then there
exists an absolute constant $c_{p},$ depending only on $p$, such that
\begin{equation*}
\left\Vert \sigma _{n}f\right\Vert _{p}\leq c_{p}\left( n+1\right)
^{1/p-2}\left\Vert f\right\Vert _{H_{p}},\text{ \ }n\in \mathbb{N}_{+}.
\end{equation*}
\end{corollary}

{\bf Proof}:
According to part a) of Theorem \ref{theorem3fejermax2} we conclude that%
\begin{equation*}
\left\Vert \frac{\sigma _{n}f}{\left( n+1\right) ^{1/p-2}}\right\Vert _{p}
\end{equation*}%
\begin{equation*}
\leq \left\Vert \sup_{n\in \mathbb{N}}\frac{\left\vert \sigma
_{n}f\right\vert }{\left( n+1\right) ^{1/p-2}}\right\Vert _{p}\leq
c\left\Vert f\right\Vert _{H_{p}},\text{ \ }n\in \mathbb{N}_{+}.
\end{equation*}%
The proof is complete.
\QED

We also mention the following corollaries:

\begin{corollary}
\label{corollary3fejermax1.2}Let $\left\{ \Phi _{n}:n\in \mathbb{N}\right\} $
be any non-decreasing sequence satisfying the condition (\ref{cond2}). Then
there exists a martingale $f\in H_{p}$ such that%
\begin{equation*}
\sup_{n\in \mathbb{N}}\left\Vert \frac{\sigma _{n}f}{\Phi _{n}}\right\Vert
_{weak-L_{p}}=\infty .
\end{equation*}
\end{corollary}

\begin{corollary}
\label{corollary3fejermax1.3}Let $\left\{ \Phi _{n}:n\in \mathbb{N}\right\} $
be any non-decreasing sequence satisfying the condition (\ref{cond2}). Then the following maximal operator%
\begin{equation*}
\sup_{n\in \mathbb{N}}\frac{\left\vert \sigma _{n}f\right\vert }{\Phi _{n}}
\end{equation*}%
\textit{is not bounded from the Hardy space }$H_{p}$\textit{\ to the space }$%
weak-L_{p}.$
\end{corollary}

The next theorem can be found in Tephnadze \cite{tep2}.

\begin{theorem}
\label{theorem2fejermax}a) The maximal operator%
\begin{equation*}
\overset{\sim }{\sigma }^{\ast }f:=\sup_{n\in \mathbb{N}}\frac{\left\vert
\sigma _{n}f\right\vert }{\log ^{2}\left( n+1\right) }
\end{equation*}%
\textit{is bounded from the Hardy space }$H_{1/2}$\textit{\ to the Lebesgue space }$%
L_{1/2}.$

\textit{b) }Let $\left\{ \Phi _{n}:n\in \mathbb{N}\right\} $ be any
non-decreasing sequence satisfying the condition%
\begin{equation}
\overline{\lim_{n\rightarrow \infty }}\frac{\log ^{2}\left( n+1\right) }{%
\Phi _{n}}=+\infty .  \label{cond1}
\end{equation}%
Then
\begin{equation*}
\sup_{k\in \mathbb{N}}\frac{\left\Vert \frac{\sigma _{q_{n_{k}}}f_{k}}{\Phi
_{q_{n_{k}}}}\right\Vert _{1/2}}{\left\Vert f_{k}\right\Vert _{H_{1/\text{ }%
2}}}=\infty .
\end{equation*}%
\end{theorem}

{\bf Proof}:
First we note that $\sigma _{n}$ is bounded from $L_{\infty }$ to $L_{\infty
}$ (see (\ref{fn4}) in Corollary \ref{lemma7kn}). Hence, according to Lemma %
\ref{lemma2.2}, to prove part a) it suffices to show that
\begin{equation*}
\int\limits_{\overline{I}_{N}}\left( \underset{n\in \mathbb{N}}{\sup }\frac{%
\left\vert \sigma _{n}a\right\vert }{\log ^{2}\left( n+1\right) }\right)
^{1/2}d\mu \leq c<\infty
\end{equation*}%
for every 1/2-atom $a,$ where $I$ denotes the support of the atom$.$

Let $a$ be an arbitrary 1/2-atom with support$\ I$ and $\mu \left( I\right)
=M_{N}^{-1}.$ We may assume that $I=I_{N}.$ It is easy to see that $\sigma
_{n}\left( a\right) =0$ when $n\leq M_{N}$. Therefore we can suppose that $%
n>M_{N}$.

Since $\left\Vert a\right\Vert _{\infty }\leq M_{N}^{2}$ we obtain that
\begin{equation*}
\frac{\left\vert \sigma _{n}\left( a\right) \right\vert }{\log ^{2}\left(
n+1\right) }
\end{equation*}%
\begin{equation*}
\leq \frac{1}{\log ^{2}\left( n+1\right) }\int_{I_{N}}\left\vert a\left(
t\right) \right\vert \left\vert K_{n}\left( x-t\right) \right\vert d\mu
\left( t\right)
\end{equation*}%
\begin{equation*}
\leq \frac{\left\Vert a\right\Vert _{\infty }}{\log ^{2}\left( n+1\right) }%
\int_{I_{N}}\left\vert K_{n}\left( x-t\right) \right\vert d\mu \left(
t\right)
\end{equation*}%
\begin{equation*}
\leq \frac{M_{N}^{2}}{\log ^{2}\left( n+1\right) }\int_{I_{N}}\left\vert
K_{n}\left( x-t\right) \right\vert d\mu \left( t\right) .
\end{equation*}

Let $x\in I_{N}^{k,l},\,0\leq k<l\leq N.$ Then, from Corollary \ref{lemma5aa}
it follows that
\begin{equation}
\frac{\left\vert \sigma _{n}\left( a\right) \right\vert }{\log ^{2}\left(
n+1\right) }  \label{2t12}
\end{equation}%
\begin{equation*}
\leq \frac{cM_{N}^{2}}{N^{2}}\frac{M_{l}M_{k}}{M_{N}^{2}}=\frac{cM_{l}M_{k}}{%
N^{2}}.
\end{equation*}%
The expression on the right-hand side of (\ref{2t12}) does not depend on $n.$ Therefore,
\begin{equation}
\left\vert \widetilde{\sigma }^{\ast }a\left( x\right) \right\vert \leq
\frac{cM_{l}M_{k}}{N^{2}},\text{ \ for \ }x\in I_{N}^{k,l},\,0\leq k<l\leq N.
\label{2t12M}
\end{equation}%
By applying (\ref{2t12M}) with identity (\ref{1.1}) we obtain that%
\begin{equation*}
\int_{\overline{I_{N}}}\left\vert \widetilde{\sigma }^{\ast }a\left(
x\right) \right\vert ^{1/2}d\mu \left( x\right)
\end{equation*}%
\begin{equation*}
=\overset{N-2}{\underset{k=0}{\sum }}\overset{N-1}{\underset{l=k+1}{\sum }}%
\sum\limits_{x_{j}=0,j\in
\{l+1,...,N-1\}}^{m_{j}-1}\int_{I_{N}^{k,l}}\left\vert \overset{\sim }{%
\sigma }^{\ast }a\left( x\right) \right\vert ^{1/2}d\mu \left( x\right)
\end{equation*}%
\begin{equation*}
+\overset{N-1}{\underset{k=0}{\sum }}\int_{I_{N}^{k,N}}\left\vert \widetilde{%
\sigma }^{\ast }a\left( x\right) \right\vert ^{1/2}d\mu \left( x\right)
\end{equation*}%
\begin{equation*}
\leq c\overset{N-2}{\underset{k=0}{\sum }}\overset{N-1}{\underset{l=k+1}{%
\sum }}\frac{m_{l+1}\ldots m_{N-1}}{M_{N}}\frac{M_{l}^{1/2}M_{k}^{1/2}}{N}
\end{equation*}%
\begin{equation*}
+c\overset{N-1}{\underset{k=0}{\sum }}\frac{1}{M_{N}}\frac{%
M_{N}^{1/2}M_{k}^{1/2}}{N}:=I+II.
\end{equation*}%
We estimate each term separately.%
\begin{equation*}
I\leq \frac{c}{N}\overset{N-2}{\underset{k=0}{\sum }}\overset{N-1}{\underset{%
l=k+1}{\sum }}\frac{M_{k}^{1/2}}{M_{l}^{1/2}}\leq \frac{c}{N}\overset{N-2}{%
\underset{k=0}{\sum }}1\leq c<\infty .
\end{equation*}%

Moreover,
\begin{equation*}
II\leq \frac{1}{M_{N}^{1/2}N}\overset{N-1}{\underset{k=0}{\sum }}%
M_{k}^{1/2}\leq \frac{c}{N}<c<\infty .
\end{equation*}%
The proof of part a) is complete.

Let$\ \left\{ \lambda _{k}:\text{ }k\in \mathbb{N}\right\} $ be an
increasing sequence of positive integers such that
\begin{equation*}
\lim_{k\rightarrow \infty }\frac{\log ^{2}\left( \lambda _{k}+1\right) }{%
\Phi _{\lambda _{k}}}=\infty .
\end{equation*}

It is evident that for every $\lambda _{k}$ there exists a positive integers
$m_{k}^{,}$ such that $q_{m_{k}^{^{\prime }}}\leq \lambda
_{k}<q_{m_{k}^{\prime }+1}<c$ $q_{m_{k}^{^{\prime }}}.$ Since $\Phi _{n}$ is
a non-decreasing function we have that
\begin{equation*}
\overline{\underset{k\rightarrow \infty }{\lim }}\frac{\left(
m_{k}^{^{\prime }}\right) ^{2}}{\Phi _{q_{m_{k}^{,}}}}\geq
c\lim_{k\rightarrow \infty }\frac{\log ^{2}\left( \lambda _{k}+1\right) }{%
\Phi _{\lambda _{k}}}=\infty .
\end{equation*}

Let $\ \left\{ n_{k}:k\in \mathbb{N}\right\} \subset \left\{ m_{k}^{\prime
}:k\in \mathbb{N}\right\} $ be a subsequence of positive numbers $\mathbb{N}%
_{+}$ such that
\begin{equation*}
\lim_{k\rightarrow \infty }\frac{n_{k}^{2}}{\Phi _{q_{n_{k}}}}=\infty
\end{equation*}%
and let $f_{k}$ be the atom defined in Example \ref{example2.2}. We combine (\ref{13l}%
) and (\ref{14l}) in Example \ref{example2.2} and invoke (\ref{dn21}) in
Lemma \ref{dn1} to obtain that
\begin{equation*}
\frac{\left\vert \sigma _{q_{n_{k}}}f_{k}\right\vert }{\Phi _{q_{n_{k}}}}
\end{equation*}%
\begin{equation*}
=\frac{1}{\Phi _{q_{n_{k}}}q_{n_{k}}}\left\vert \overset{q_{n_{k}}}{\underset%
{j=1}{\sum }}S_{j}f_{k}\right\vert
\end{equation*}%
\begin{equation*}
=\frac{1}{\Phi _{q_{n_{k}}}q_{n_{k}}}\left\vert \overset{q_{n_{k}}}{\underset%
{j=M_{_{2n_{k}}}+1}{\sum }}S_{j}f_{k}\right\vert
\end{equation*}%
\begin{equation*}
=\frac{1}{\Phi _{q_{n_{k}}}q_{n_{k}}}\left\vert \overset{q_{n_{k}}}{\underset%
{j=M_{_{2n_{k}}}+1}{\sum }}\left( D_{j}-D_{M_{_{2n_{k}}}}\right) \right\vert
\end{equation*}%
\begin{equation*}
=\frac{1}{\Phi _{q_{n_{k}}}q_{n_{k}}}\left\vert \overset{q_{n_{k}-1}}{%
\underset{j=1}{\sum }}\left( D_{j+M_{_{2n_{k}}}}-D_{M_{_{2n_{k}}}}\right)
\right\vert
\end{equation*}%
\begin{equation*}
=\frac{1}{\Phi _{q_{n_{k}}}q_{n_{k}}}\left\vert \overset{q_{n_{k}-1}}{%
\underset{j=1}{\sum }}D_{j}\right\vert
\end{equation*}%
\begin{equation*}
=\frac{1}{\Phi _{q_{n_{k}}}}\frac{q_{n_{k}-1}}{q_{n_{k}}}\left\vert
K_{q_{n_{k}-1}}\right\vert .
\end{equation*}

Let $x\in $ $I_{_{2n_{k}}}^{2s,2l}$. Then, in view of Lemma \ref{lemma3}, we
can conclude that
\begin{equation*}
\frac{\left\vert \sigma _{q_{n_{k}}}f_{k}(x)\right\vert }{\Phi _{q_{n_{k}}}}%
\geq \frac{cM_{2s}M_{2l}}{M_{_{2n_{k}}}\Phi _{q_{n_{k}}}}.
\end{equation*}%
Hence,
\begin{equation*}
\int_{G_{m}}\left\vert \frac{\sigma _{q_{n_{k}}}f_{k}}{\Phi _{q_{n_{k}}}}%
\right\vert ^{1/2}d\mu
\end{equation*}%
\begin{equation*}
\geq \text{ }\overset{n_{k}-3}{\underset{s=0}{\sum }}\overset{n_{k}-1}{%
\underset{l=s+1}{\sum }}\overset{m_{2l+1}}{\underset{x_{2l+1=0}}{\sum }}...%
\overset{m_{2n_{k}-1}}{\underset{x_{2n_{k}-1}=0}{\sum }}%
\int_{I_{_{2n_{k}}}^{2s,2l}}\left\vert \frac{\sigma _{q_{n_{k}}}f_{k}}{\Phi
_{q_{n_{k}}}}\right\vert ^{1/2}d\mu
\end{equation*}%
\begin{equation*}
\geq c\overset{n_{k}-3}{\underset{s=0}{\sum }}\overset{n_{k}-1}{\underset{%
l=s+1}{\sum }}\frac{m_{_{2l+1}}...m_{2n_{k}-1}}{M_{2n_{k}}}\frac{%
M_{2s}^{1/2}M_{2l}^{1/2}}{\Phi _{q_{n_{k}}}^{1/2}M_{_{2n_{k}}}^{1/2}}
\end{equation*}%
\begin{equation*}
\geq c\overset{n_{k}-3}{\underset{s=0}{\sum }}\overset{n_{k}-1}{\underset{%
l=s+1}{\sum }}\frac{M_{2s}^{1/2}}{M_{2l}^{1/2}M_{_{2n_{k}}}^{1/2}\Phi
_{q_{n_{k}}}^{1/2}}
\end{equation*}%
\begin{equation*}
\geq \frac{cn_{k}}{M_{_{2n_{k}}}^{1/2}\Phi _{q_{n_{k}}}^{1/2}}.
\end{equation*}%
From (\ref{15l}) in Example \ref{example2.2} we have that
\begin{equation*}
\frac{\left( \int_{G_{m}}\left\vert \frac{\sigma _{q_{n_{k}}}f_{k}}{\Phi
_{q_{n_{k}}}}\right\vert ^{1/2}d\mu \right) ^{2}}{\left\Vert
f_{k}\right\Vert _{H_{1/\text{ }2}}}
\end{equation*}%
\begin{equation*}
\geq \frac{cn_{k}^{2}}{M_{_{2n_{k}}}\Phi _{q_{n_{k}}}}M_{2n_{k}}
\end{equation*}%
\begin{equation*}
\geq \frac{cn_{k}^{2}}{\Phi _{q_{n_{k}}}}\rightarrow \infty ,\qquad \text{%
when \qquad }k\rightarrow \infty .
\end{equation*}

Thus, also part b) is proved so the proof is complete.

\QED

We also point out the following consequence of Theorem \ref{theorem2fejermax}, which we need later on.

\begin{corollary}
\label{corollary2fejermax} Let $f\in H_{1/2}.$ Then there exists an absolute
constant $c$ such that
\begin{equation*}
\left\Vert \sigma _{n}f\right\Vert _{1/2}\leq c\log ^{2}\left( n+1\right)
\left\Vert f\right\Vert _{H_{1/2}},\text{ \ }n\in \mathbb{N}_{+}.
\end{equation*}
\end{corollary}

{\bf Proof}:
According to part a) of Theorem \ref{theorem2fejermax} we readily conclude
that%
\begin{equation*}
\left\Vert \frac{\sigma _{n}f}{\log ^{2}\left( n+1\right) }\right\Vert _{1/2}
\end{equation*}%
\begin{equation*}
\leq \left\Vert \sup_{n\in \mathbb{N}}\frac{\left\vert \sigma
_{n}f\right\vert }{\log ^{2}\left( n+1\right) }\right\Vert _{1/2}\leq
c\left\Vert f\right\Vert _{H_{1/2}},\text{ \ }n\in \mathbb{N}_{+}.
\end{equation*}%
The proof is complete.
\QED

We also mention the following corollaries:

\begin{corollary}
\label{corollary3fejermax1.2.2}Let $\left\{ \Phi _{n}:n\in \mathbb{N}%
\right\} $ be any non-decreasing sequence, satisfying the condition (\ref%
{cond1}). Then there exists a martingale $f\in H_{1/2},$ such that%
\begin{equation*}
\sup_{n\in \mathbb{N}}\left\Vert \frac{\sigma _{n}f}{\Phi _{n}}\right\Vert
_{1/2}=\infty .
\end{equation*}
\end{corollary}

\begin{corollary}
\label{corollary3fejermax1.3.2}Let $\left\{ \Phi _{n}:n\in \mathbb{N}%
\right\} $ be any non-decreasing sequence, satisfying the condition (\ref%
{cond1}). Then the following maximal operator%
\begin{equation*}
\sup_{n\in \mathbb{N}}\frac{\left\vert \sigma _{n}f\right\vert }{\Phi _{n}}
\end{equation*}%
\textit{is not bounded from the Hardy space }$H_{1/2}$\textit{\ to the Lebesgue space
}$L_{1/2}.$
\end{corollary}

The next results was proved by Tephnadze \cite{pt2}.

\begin{theorem}
\label{theorem0fejermax}a) Let $0<p\leq 1/2$ and $\left\{ n_{k}:k\in \mathbb{%
N}\right\} $ be a subsequence of positive numbers such that
\begin{equation}
\sup_{k\in \mathbb{N}}\rho \left( n_{k}\right) =\varkappa <c<\infty .  \label{fenk}
\end{equation}%
Then the maximal operator
\begin{equation*}
\widetilde{\sigma }^{\ast ,\vartriangle }f:=\sup_{k\in \mathbb{N}}\left\vert \sigma
_{n_{k}}f\right\vert
\end{equation*}%
is bounded from the Hardy space $H_{p}$ to the Lebesgue space $L_{p}.$

b) Let $0<p<1/2\ \ $and $\left\{ n_{k}:k\in \mathbb{N}\right\} $ be a
subsequence of positive numbers satisfying the condition
\begin{equation}
\sup_{k\in \mathbb{N}}\rho \left( n_{k}\right) =\infty .  \label{fenk1}
\end{equation}%
Then there exists an martingale $f\in H_{p}$ such that%
\begin{equation*}
\sup_{k\in \mathbb{N}}\left\Vert \sigma _{n_{k}}f\right\Vert
_{weak-L_{p}}=\infty ,\text{ \ \ }\left( 0<p<1/2\right) .
\end{equation*}
\end{theorem}

{\bf Proof}:
Since $\widetilde{\sigma }^{\ast ,\vartriangle }$ is bounded from $L_{\infty
}$ to $L_{\infty },$ by using Lemma \ref{lemma2.2}, we obtain that the proof
of part a) is complete if we show
that%
\begin{equation*}
\int_{\overline{I_{N}}}\left\vert \widetilde{\sigma }^{\ast ,\vartriangle
}a\left( x\right) \right\vert <c<\infty ,
\end{equation*}%
for every $p$-atom $a$ with support$\ I_{N}$ and $\mu \left( I_{N}\right)
=M_{N}^{-1}.$ Analogously to Theorem \ref{theorem2fejermax} we may assume that
$n_{k}>M_{N}$.

Since $\left\Vert a\right\Vert _{\infty }\leq M_{N}^{1/p}$ we find that%
\begin{equation}
\left\vert \sigma _{n_{k}}\left( a\right) \right\vert \leq
\int_{I_{N}}\left\vert a\left( t\right) \right\vert \left\vert
K_{n_{k}}\left( x-t\right) \right\vert d\mu \left( t\right)  \label{400z}
\end{equation}%
\begin{equation*}
\leq \left\Vert a\right\Vert _{\infty }\int_{I_{N}}\left\vert
K_{n_{k}}\left( x-t\right) \right\vert d\mu \left( t\right)
\end{equation*}%
\begin{equation*}
\leq M_{N}^{1/p}\int_{I_{N}}\left\vert K_{n_{k}}\left( x-t\right)
\right\vert d\mu \left( t\right).
\end{equation*}

Let $x\in I_{N}^{i,j}\ $and $i<j<\left\langle n_{k}\right\rangle .$ Then $%
x-t\in I_{N}^{i,j}$ for $t\in I_{N}\ \ $and, according to Lemma \ref{lemma2}
we obtain that
\begin{equation*}
\left\vert K_{M_{l}}\left( x-t\right) \right\vert =0,\text{ \ for all }%
\left\langle n_{k}\right\rangle \leq \text{ }l\leq \left\vert
n_{k}\right\vert .
\end{equation*}

By applying (\ref{400z}) and (\ref{fn5}) in Lemma \ref{lemma7kn} we get
that
\begin{eqnarray*}
\left\vert \sigma _{n_{k}}a\left( x\right) \right\vert &\leq &M_{N}^{1/p}%
\overset{\left\vert n_{k}\right\vert }{\underset{l=\left\langle
n_{k}\right\rangle }{\sum }}\int_{I_{N}}\left\vert K_{M_{l}}\left(
x-t\right) \right\vert d\mu \left( t\right) =0,\text{ } \\
\text{for }x &\in &I_{N}^{i,j},\text{ \ }0\leq i<j<\left\langle
n_{k}\right\rangle \leq l\leq \left\vert n_{k}\right\vert .  \notag
\end{eqnarray*}

Let $x\in I_{N}^{i,j},\,$where $\left\langle n_{k}\right\rangle \leq j\leq
N. $ Then, in view of Corollary \ref{lemma5aa}, we have that
\begin{equation*}
\int_{I_{N}}\left\vert K_{n_{k}}\left( x-t\right) \right\vert d\mu \left(
t\right) \leq \frac{cM_{i}M_{j}}{M_{N}^{2}}.
\end{equation*}

By using again (\ref{400z}) we obtain that
\begin{equation*}
\left\vert \sigma _{n_{k}}a\left( x\right) \right\vert \leq
c_{p}M_{N}^{1/p-2}M_{i}M_{j}.
\end{equation*}%
Set
\begin{equation*}
\varrho :=\min_{k\in \mathbb{N}}\left\langle \alpha _{k}\right\rangle .
\end{equation*}
Then%
\begin{equation}
\left\vert \widetilde{\sigma }^{\ast ,\vartriangle }a\left( x\right)
\right\vert =0,\text{ \ \ for }\ x\in I_{N}^{i,j},\ 0\leq i<j\leq \varrho
\label{403zM0}
\end{equation}%
and
\begin{equation}
\left\vert \widetilde{\sigma }^{\ast ,\vartriangle }a\left( x\right)
\right\vert \leq c_{p}M_{N}^{1/p-2}M_{i}M_{j},\text{ for \ }x\in I_{N}^{i,j},%
\text{ \ }i<\varrho \leq j\leq N-1.  \label{403zM1}
\end{equation}

Analogously to (\ref{dnsnM1}) we can conclude that
\begin{equation}
N-\varrho \leq \varkappa  \label{403zM20}
\end{equation}%
and
\begin{equation}
\frac{M_{N}^{1-2p}}{M_{\varrho }^{1-2p}}\leq \lambda ^{\left( N-\varrho
\right) \left( 1-2p\right) }\leq \lambda ^{\varkappa \left( 1-2p\right)
}<c<\infty ,  \label{403zM2}
\end{equation}%
where \ \ $\lambda =\sup_{k}m_{k}.$

Let $p=1/2.$ By combining (\ref{403zM0})-(\ref{403zM20}) with (\ref{1.1}) we
get that
\begin{equation*}
\int_{\overline{I_{N}}}\left\vert \widetilde{\sigma }^{\ast ,\vartriangle
}a\right\vert ^{1/2}d\mu
\end{equation*}%
\begin{equation*}
=\overset{N-2}{\underset{i=0}{\sum }}\overset{N-1}{\underset{j=i+1}{\sum }}%
\int_{I_{N}^{i,j}}\left\vert \widetilde{\sigma }^{\ast ,\vartriangle
}a\right\vert ^{1/2}d\mu
\end{equation*}%
\begin{equation*}
+\overset{N-1}{\underset{i=0}{\sum }}\int_{I_{N}^{i,N}}\left\vert \widetilde{%
\sigma }^{\ast ,\vartriangle }a\right\vert ^{1/2}d\mu
\end{equation*}%
\begin{equation*}
\leq \overset{\varrho -1}{\underset{i=0}{\sum }}\overset{N-1}{\underset{%
j=\varrho }{\sum }}\int_{I_{N}^{i,j}}\left\vert \widetilde{\sigma }^{\ast
,\vartriangle }a\right\vert ^{1/2}d\mu
\end{equation*}%
\begin{equation*}
+\overset{N-2}{\underset{i=\varrho }{\sum }}\overset{N-1}{\underset{j=i+1}{%
\sum }}\int_{I_{N}^{i,j}}\left\vert \widetilde{\sigma }^{\ast ,\vartriangle
}a\right\vert ^{1/2}d\mu +\overset{N-1}{\underset{i=0}{\sum }}%
\int_{I_{N}^{i,N}}\left\vert \widetilde{\sigma }^{\ast ,\vartriangle
}a\right\vert ^{1/2}d\mu
\end{equation*}%
\begin{equation*}
\leq c\overset{\varrho }{\underset{i=0}{\sum }}M_{i}^{1/2}\overset{N-1}{%
\underset{j=\varrho +1}{\sum }}\frac{1}{M_{j}^{1/2}}
\end{equation*}%
\begin{equation*}
+\overset{N-2}{\underset{i=\varrho }{\sum }}M_{i}^{1/2}\overset{N-1}{%
\underset{j=i+1}{\sum }}\frac{1}{M_{j}^{1/2}}+c_{p}\overset{N-1}{\underset{%
i=0}{\sum }}\frac{M_{i}^{1/2}}{M_{N}^{1/2}}
\end{equation*}%
\begin{equation*}
\leq cM_{\varrho }^{1/2}\frac{1}{M_{\varrho }^{1/2}}+c\overset{N-2}{\underset%
{i=\varrho }{\sum }}M_{i}^{1/2}\frac{1}{M_{i}^{1/2}}+c
\end{equation*}%
\begin{equation*}
\leq N-\varrho +c\leq c<\infty .
\end{equation*}

Let $0<p<1/2.$ By combining (\ref{403zM0}), (\ref{403zM1}) and (\ref{403zM2}%
) with (\ref{1.1}) we have that
\begin{equation*}
\int_{\overline{I_{N}}}\left\vert \widetilde{\sigma }^{\ast ,\vartriangle
}a\right\vert ^{p}d\mu
\end{equation*}%
\begin{equation*}
=\overset{N-2}{\underset{i=0}{\sum }}\overset{N-1}{\underset{j=i+1}{\sum }}%
\int_{I_{N}^{i,j}}\left\vert \widetilde{\sigma }^{\ast ,\vartriangle
}a\right\vert ^{p}d\mu
\end{equation*}%
\begin{equation*}
+\overset{N-1}{\underset{i=0}{\sum }}\int_{I_{N}^{k,N}}\left\vert \widetilde{%
\sigma }^{\ast ,\vartriangle }a\right\vert ^{p}d\mu
\end{equation*}%
\begin{equation*}
\leq \overset{\varrho -1}{\underset{i=0}{\sum }}\overset{N-1}{\underset{%
j=\varrho }{\sum }}\int_{I_{N}^{i,j}}\left\vert \widetilde{\sigma }^{\ast
,\vartriangle }a\right\vert ^{p}d\mu
\end{equation*}%
\begin{equation*}
+\overset{N-2}{\underset{i=\varrho }{\sum }}\overset{N-1}{\underset{j=i+1}{%
\sum }}\int_{I_{N}^{i,j}}\left\vert \widetilde{\sigma }^{\ast ,\vartriangle
}a\right\vert ^{p}d\mu +\overset{N-1}{\underset{i=0}{\sum }}%
\int_{I_{N}^{i,N}}\left\vert \widetilde{\sigma }^{\ast ,\vartriangle
}a\right\vert ^{p}d\mu
\end{equation*}%
\begin{equation*}
\leq c_{p}M_{N}^{1-2p}\overset{\varrho }{\underset{i=0}{\sum }}M_{i}^{p}%
\overset{N-1}{\underset{j=\varrho +1}{\sum }}\frac{1}{M_{j}^{1-p}}
\end{equation*}%
\begin{equation*}
+M_{N}^{1-2p}\overset{N-2}{\underset{i=\varrho }{\sum }}M_{i}^{p}\overset{N-1%
}{\underset{j=i+1}{\sum }}\frac{1}{M_{j}^{1-p}}+c_{p}\overset{N-1}{\underset{%
i=0}{\sum }}\frac{M_{i}^{p}}{M_{N}^{p}}
\end{equation*}%
\begin{equation*}
\leq \frac{c_{p}M_{N}^{1-2p}}{M_{\varrho }^{1-2p}}+c_{p}\leq c_{p}\lambda
^{\left( \left\vert n_{k}\right\vert -\left\langle n_{k}\right\rangle
\right) \left( 1-2p\right) }+c_{p}<\infty .
\end{equation*}

The proof of part a) is complete.

Let $\left\{ n_{k}:k\in \mathbb{N}\right\} $ be a sequence of positive
numbers satisfying condition (\ref{fenk1}). Then
\begin{equation}
\sup_{k\in \mathbb{N}}\frac{M_{\left\vert n_{k}\right\vert }}{M_{\left\langle
n_{k}\right\rangle }}=\infty .  \label{12h}
\end{equation}

Under condition (\ref{12h}) there exists a sequence $\left\{ \alpha _{k}:%
\text{ }k\in \mathbb{N}\right\} \subset \left\{ n_{k}:\text{ }k\in \mathbb{N}%
\right\} $ such that
\begin{equation*}
\sum_{k=0}^{\infty }\frac{M_{\left\langle \alpha _{k}\right\rangle }^{\left(
1-2p\right) /2}}{M_{\left\vert \alpha _{k}\right\vert }^{\left( 1-2p\right)
/2}}<c<\infty .
\end{equation*}

Let $f=\left( f^{\left( n\right) }:\text{ }n\in \mathbb{N}\right) $ be
the martingae defined in Example \ref{example2.6}, where
\begin{equation}
\lambda _{k}=\frac{\lambda M_{\left\langle \alpha _{k}\right\rangle
}^{\left( 1/p-2\right) /2}}{M_{\left\vert \alpha _{k}\right\vert }^{\left(
1/p-2\right) /2}}\text{, \ \ \ \ }\lambda =\sup_{n\in \mathbb{N}}m_{n}.  \label{12hhaa}
\end{equation}

By applying (\ref{lemma2.1}) we can conclude that $f\in H_{p}.$

By now using (\ref{10AA}) with $\lambda _{k}$ defined by (\ref{12hhaa}) we obtain that

\begin{equation}
\widehat{f}(j)  \label{6aacharp}
\end{equation}%
\begin{equation*}
=\left\{
\begin{array}{ll}
M_{\left\vert \alpha _{k}\right\vert }^{1/2p}M_{\left\langle \alpha
_{k}\right\rangle }^{\left( 1/p-2\right) /2}, & j\in \left\{ M_{\left\vert
\alpha _{k}\right\vert },...,\text{ ~}M_{\left\vert \alpha _{k}\right\vert
+1}-1\right\} ,\text{ }k\in \mathbb{N}_{+}, \\
0, & \text{\thinspace }j\notin \bigcup\limits_{k=0}^{\infty }\left\{
M_{\left\vert \alpha _{k}\right\vert },...,\text{ ~}M_{\left\vert \alpha
_{k}\right\vert +1}-1\right\} .%
\end{array}%
\right.
\end{equation*}

Therefore,

\begin{equation*}
\sigma _{_{\alpha _{k}}}f=\frac{1}{\alpha _{k}}\sum_{j=1}^{M_{\left\vert
\alpha _{k}\right\vert }}S_{j}f+\frac{1}{\alpha _{k}}\sum_{j=M_{\left\vert
\alpha _{k}\right\vert }+1}^{\alpha _{k}}S_{j}f
\end{equation*}%
\begin{equation*}
:=I+II.
\end{equation*}%

Let $M_{\left\vert \alpha _{k}\right\vert }<j\leq \alpha _{k}.$ Then,
according to (\ref{6aacharp}) we have that
\begin{equation}
S_{j}f=S_{M_{\left\vert \alpha _{k}\right\vert }}f+M_{\left\vert \alpha
_{k}\right\vert }^{1/2p}M_{\left\langle \alpha _{k}\right\rangle }^{\left(
1/p-2\right) /2}D_{j-M_{\left\vert \alpha _{k}\right\vert }}.  \label{8aafn}
\end{equation}

By applying (\ref{8aafn}) we can rewrite $II$ as%
\begin{equation*}
II=\frac{\alpha _{k}-M_{\left\vert \alpha _{k}\right\vert }}{\alpha _{k}}%
S_{M_{\left\vert \alpha _{k}\right\vert }}f
\end{equation*}%
\begin{equation*}
+\frac{M_{\left\vert \alpha _{k}\right\vert }^{1/2p}M_{\left\langle \alpha
_{k}\right\rangle }^{\left( 1/p-2\right) /2}\psi _{M_{\left\vert \alpha
_{k}\right\vert }}}{\alpha _{k}}\sum_{j=M_{\left\vert \alpha _{k}\right\vert
}+1}^{\alpha _{k}}D_{j-M_{\left\vert \alpha _{k}\right\vert }}
\end{equation*}%
\begin{equation*}
:=II_{1}+II_{2}.
\end{equation*}

In view of Corollary \ref{corollary4.43} for $II_{1}$ we find that
\begin{equation*}
\left\Vert II_{1}\right\Vert _{weak-L_{p}}^{p}
\end{equation*}%
\begin{equation*}
\leq \left( \frac{\alpha _{k}-M_{\left\vert \alpha _{k}\right\vert }}{\alpha
_{k}}\right) ^{p}\left\Vert S_{M_{\left\vert \alpha _{k}\right\vert
}}f\right\Vert _{weak-L_{p}}^{p}
\end{equation*}%
\begin{equation*}
\leq c_{p}\left\Vert f\right\Vert _{H_{p}}^{p}<\infty .
\end{equation*}

By using part a) of Theorem \ref{theorem0fejermax} for $I$ we obtain that
\begin{equation*}
\left\Vert I\right\Vert _{weak-L_{p}}^{p}=\left( \frac{M_{\left\vert \alpha
_{k}\right\vert }}{\alpha _{k}}\right) ^{p}\left\Vert \sigma _{M_{\left\vert
\alpha _{k}\right\vert }}f\right\Vert _{weak-L_{p}}^{p}
\end{equation*}%
\begin{equation*}
\leq c\left\Vert f\right\Vert _{H_{p}}^{p}\leq c_{p}<\infty .
\end{equation*}

Under condition (\ref{fenk1}) we can conclude that
\begin{equation*}
\left\langle \alpha _{k}\right\rangle \neq \left\vert \alpha _{k}\right\vert
\text{ \ \ \ \ \ \ \ and \ \ \ \ \ \ }\left\langle \alpha _{k}-M_{\left\vert
\alpha _{k}\right\vert }\right\rangle =\left\langle \alpha _{k}\right\rangle.
\end{equation*}%
Let $x\in $ $I_{_{\left\langle \alpha _{k}\right\rangle +1}}^{\left\langle
\alpha _{k}\right\rangle -1,\left\langle \alpha _{k}\right\rangle }.$
According to Lemma \ref{lemma8ccc} we get that
\begin{equation*}
\left\vert II_{2}\right\vert =\frac{M_{\left\vert \alpha _{k}\right\vert
}^{1/2p}M_{\left\langle \alpha _{k}\right\rangle }^{\left( 1/p-2\right) /2}}{%
\alpha _{k}}\left\vert \sum_{j=1}^{\alpha _{k}-M_{\left\vert \alpha
_{k}\right\vert }}D_{j}\right\vert
\end{equation*}%
\begin{equation*}
=\frac{M_{\left\vert \alpha _{k}\right\vert }^{1/2p}M_{\left\langle \alpha
_{k}\right\rangle }^{\left( 1/p-2\right) /2}}{\alpha _{k}}\left( \alpha
_{k}-M_{\left\vert \alpha _{k}\right\vert }\right) \left\vert K_{\alpha
_{k}-M_{\left\vert \alpha _{k}\right\vert }}\right\vert
\end{equation*}%
\begin{equation*}
\geq cM_{\left\vert \alpha _{k}\right\vert }^{1/2p-1}M_{\left\langle \alpha
_{k}\right\rangle }^{\left( 1/p-2\right) /2}\left( \alpha _{k}-M_{\left\vert
\alpha _{k}\right\vert }\right) \left\vert K_{\alpha _{k}-M_{\left\vert
\alpha _{k}\right\vert }}\right\vert
\end{equation*}%
\begin{equation*}
\geq cM_{\left\vert \alpha _{k}\right\vert }^{1/2p-1}M_{\left\langle \alpha
_{k}\right\rangle }^{\left( 1/p+2\right) /2}.
\end{equation*}

It follows that%
\begin{equation*}
\left\Vert II_{2}\right\Vert _{weak-L_{p}}^{p}
\end{equation*}%
\begin{equation*}
\geq c_{p}\left( M_{\left\vert \alpha _{k}\right\vert }^{\left( 1/p-2\right)
/2}M_{\left\langle \alpha _{k}\right\rangle }^{\left( 1/p+2\right)
/2}\right) ^{p}\mu \left\{ x\in G_{m}:\text{ }\left\vert II_{2}\right\vert
\geq c_{p}M_{\left\vert \alpha _{k}\right\vert }^{\left( 1/p-2\right)
/2}M_{\left\langle \alpha _{k}\right\rangle }^{\left( 1/p+2\right)
/2}\right\} -c_{p}
\end{equation*}%
\begin{equation*}
\geq c_{p}M_{\left\vert \alpha _{k}\right\vert }^{1/2-p}M_{\left\langle
\alpha _{k}\right\rangle }^{1/2+p}\mu \left\{ I_{_{\left\langle \alpha
_{k}\right\rangle +1}}^{\left\langle \alpha _{k}\right\rangle
-1,\left\langle \alpha _{k}\right\rangle }\right\}
\end{equation*}%
\begin{equation*}
\geq \frac{c_{p}M_{\left\vert \alpha _{k}\right\vert }^{1/2-p}}{%
M_{\left\langle \alpha _{k}\right\rangle }^{1/2-p}}\rightarrow \infty ,\text{
\ when \ \ }k\rightarrow \infty .
\end{equation*}%

By now combining the estimates above we can conclude that
\begin{equation*}
\left\Vert \sigma _{\alpha _{k}}f\right\Vert _{weak-L_{p}}^{p}
\end{equation*}%
\begin{equation*}
\geq \left\Vert II_{2}\right\Vert _{weak-L_{p}}^{p}-\left\Vert
II_{1}\right\Vert _{weak-L_{p}}^{p}-\left\Vert I\right\Vert _{weak-L_{p}}^{p}
\end{equation*}%
\begin{equation*}
\geq \frac{1}{2}\left\Vert II_{2}\right\Vert _{weak-L_{p}}^{p}
\end{equation*}%
\begin{equation*}
\geq \frac{c_{p}M_{\left\vert \alpha _{k}\right\vert }^{1/2-p}}{%
M_{\left\langle \alpha _{k}\right\rangle }^{1/2-p}}\rightarrow \infty ,\text{
when }k\rightarrow \infty .
\end{equation*}

The proof is complete.
\QED

From the part a) of Theorem \ref{theorem0fejermax}  follows immediately the
following well-known results of Weisz \protect\cite{We2}:

\begin{corollary}
\label{Corollary0fejermax}Let $p>0.$ Then the maximal operator
\begin{equation*}
\sigma ^{\#}f:=\underset{n\in
\mathbb{N}
}{\sup }\left\vert \sigma _{M_{n}}f\right\vert
\end{equation*}%
is bounded from the Hardy space $H_{p}$ to the Lebesgue space $L_{p}.$
\end{corollary}

{\bf Proof}:
It is obvious that
\begin{equation*}
\left\vert M_{n}\right\vert =\left[ M_{n}\right] =n
\end{equation*}%
and
\begin{equation*}
\sup_{n\in \mathbb{N}}\rho \left( M_{n}\right) =0<c<\infty .
\end{equation*}

It follows that condition (\ref{fenk}) is satisfied and the proof is
complete by just applying part a) of Theorem \ref{theorem0fejermax}.
\QED

Our next results reads:

\begin{theorem}
\label{Theorem1fejermax}Let $0<p\leq 1.$ Then the operator $\left\vert
\sigma _{M_{n}}f\right\vert $ is not bounded from the martingale Hardy space
$H_{p}$ to the martingale Hardy space $H_{p}.$
\end{theorem}

\begin{remark}
This result for the Walsh system can be found in Goginava \cite{gog1} and for bounded Vilenkin systems in the paper of Persson and Tephnadze \cite{pt1}.
\end{remark}

{\bf Proof}:
Let $f_{k}$ be the martingale from Example \ref{example2.2}. By combining (%
\ref{dn21}) in Lemma \ref{dn1} and (\ref{13l}), (\ref{14l}) in Example \ref%
{example2.2} we find that
\begin{equation*}
\sigma _{M_{2n_{k}+1}}f_{k}
\end{equation*}%
\begin{equation*}
=\frac{1}{M_{2n_{k}+1}}\overset{M_{2n_{k}+1}}{\underset{j=1}{\sum }}%
S_{j}f_{k}
\end{equation*}%
\begin{equation*}
=\frac{1}{M_{2n_{k}+1}}\overset{M_{2n_{k}+1}}{\underset{j=M_{2n_{k}}+1}{\sum
}}S_{j}f_{k}
\end{equation*}%
\begin{equation*}
=\frac{1}{M_{2n_{k}+1}}\overset{M_{2n_{k}+1}}{\underset{j=M_{2n_{k}}+1}{\sum
}}\left( D_{j}-D_{M_{2n_{k}}}\right)
\end{equation*}%
\begin{equation*}
=\frac{1}{M_{2n_{k}+1}}\overset{M_{2n_{k}}}{\underset{j=1}{\sum }}\left(
D_{j+M_{2n_{k}}}-D_{M_{2n_{k}}}\right)
\end{equation*}%
\begin{equation*}
=\frac{\psi _{M_{2n_{k}}}}{M_{2n_{k}+1}}\overset{M_{2n_{k}}}{\underset{j=1}{%
\sum }}D_{j}=\frac{\psi _{M_{2n_{k}}}}{m_{2n_{k}}}K_{M_{2n_{k}}}.
\end{equation*}

It is evident that%
\begin{equation*}
S_{M_{N}}\left( \left\vert \sigma _{M_{2n_{k}+1}}f_{k}\left( x\right)
\right\vert \right)
\end{equation*}%
\begin{equation*}
=\int_{G_{m}}\left\vert \sigma _{M_{2n_{k}}}f_{k}\left( t\right) \right\vert
D_{M_{N}}\left( x-t\right) d\mu \left( t\right)
\end{equation*}%
\begin{equation*}
\geq \int_{I_{2n_{k}}}\left\vert \sigma _{M_{2n_{k}}}f_{k}\left( t\right)
\right\vert D_{M_{N}}\left( x-t\right) d\mu \left( t\right)
\end{equation*}%
\begin{equation*}
\geq c\int_{I_{2n_{k}}}\left\vert K_{M_{2n_{k}+1}}\left( t\right)
\right\vert D_{M_{N}}\left( x-t\right) d\mu \left( t\right)
\end{equation*}%
\begin{equation*}
\geq cM_{2n_{k}}\int_{I_{2n_{k}}}D_{M_{N}}\left( x-t\right) d\mu \left(
t\right) .
\end{equation*}

Since
\begin{equation*}
\psi _{j}\left( t\right) =1\text{, \ \ \ for \ \ }t\in I_{2n_{k}},\text{ \ \
}j=0,1,...,M_{2n_{k}}-1
\end{equation*}%
we obtain that%
\begin{equation*}
S_{M_{N}}\left( \left\vert \sigma _{M_{2n_{k}+1}}f_{k}\left( x\right)
\right\vert \right)
\end{equation*}%
\begin{equation*}
\geq cM_{2n_{k}}\frac{1}{M_{2n_{k}}}D_{M_{N}}\left( x\right) ,\text{ \ }%
N=0,1,...,M_{N}-1
\end{equation*}%
and%
\begin{equation*}
\sup_{1\leq N<M_{2n_{k}}}S_{M_{N}}\left( \left\vert \sigma
_{M_{2n_{k}+1}}f_{k}\left( x\right) \right\vert \right)
\end{equation*}%
\begin{equation*}
\geq \sup_{1\leq N<M_{2n_{k}}}S_{M_{N}}\left( \left\vert \sigma
_{M_{2n_{k}+1}}f_{k}\left( x\right) \right\vert \right)
\end{equation*}%
\begin{equation*}
\geq \sup_{1\leq N<M_{2n_{k}}}D_{M_{N}}\left( x\right) .
\end{equation*}%
Let $x\in I_{s}/I_{s+1},$ for some $s=0,1,\cdots ,2n_{k}-1.$ Then, from
Corollary \ref{dn2.3} it follows that%
\begin{equation*}
\sup_{1\leq N<M_{2n_{k}}}S_{M_{N}}\left( \left\vert \sigma
_{M_{2n_{k}+1}}f_{k}\left( x\right) \right\vert \right) \geq cM_{s}.
\end{equation*}

According to Lemma \ref{lemma2.3.4} for every $0<p\leq 1$ it yields that
\begin{equation}
\left\Vert \left\vert \sigma _{M_{2n_{k}+1}}f_{k}\right\vert \right\Vert
_{H_{p}}^{p}  \label{667}
\end{equation}%
\begin{equation*}
=\left\Vert \sup_{1\leq N<2n_{k}}S_{M_{N}}\left( \left\vert \sigma
_{M_{2n_{k}+1}}f_{k}\left( x\right) \right\vert \right) \right\Vert _{p}^{p}
\end{equation*}%
\begin{equation*}
\geq \int_{G_{m}}\left( \sup_{1\leq N<M_{2n_{k}}}S_{M_{N}}\left( \left\vert
\sigma _{M_{2n_{k}+1}}f_{k}\left( x\right) \right\vert \right) \right)
^{p}d\mu \left( x\right)
\end{equation*}%
\begin{equation}
\geq \overset{2n_{k}-1}{\underset{s=1}{\sum }}\int_{I_{s}\backslash
I_{s+1}}\left( \sup_{1\leq N<M_{2n_{k}}}S_{M_{N}}\left( \left\vert \sigma
_{M_{2n_{k}+1}}f_{k}\left( x\right) \right\vert \right) \right) ^{p}d\mu
\left( x\right)  \notag
\end{equation}%
Let $0<p<1.$ Then%
\begin{equation*}
\left\Vert \left\vert \sigma _{M_{2n_{k}+1}}f_{k}\right\vert \right\Vert
_{H_{p}}^{p}\geq c\overset{2n_{k}-1}{\underset{s=1}{\sum }}\frac{M_{s}^{p}}{%
M_{s}}=c_{p}>0.
\end{equation*}

Let $p=1.$ We find that%
\begin{equation}
\left\Vert \left\vert \sigma _{M_{2n_{k}+1}}f_{k}\right\vert \right\Vert
_{H_{1}}\geq c\overset{2n_{k}-1}{\underset{s=1}{\sum }}\frac{M_{s}}{M_{s}}%
\geq cn_{k}.  \label{777}
\end{equation}

By combining (\ref{15l}) and (\ref{667}) with (\ref{777}) in Example \ref%
{example2.2} we can conclude that%
\begin{equation*}
\frac{\left\Vert \left\vert \sigma _{M_{2n_{k}+1}}f_{k}\right\vert
\right\Vert _{H_{p}}}{\left\Vert f_{k}\right\Vert _{H_{p}}}\geq \frac{c_{p}}{%
M_{2n_{k}}^{1-1/p}}\rightarrow \infty ,\text{ \ \ when \ \ \ }k\rightarrow
\infty ,\text{ \ \ }0<p<1
\end{equation*}%
and
\begin{equation*}
\frac{\left\Vert \left\vert \sigma _{M_{2n_{k}+1}}f_{k}\right\vert
\right\Vert _{H_{1}}}{\left\Vert f_{k}\right\Vert _{H_{1}}}\geq
cn_{k}\rightarrow \infty ,\text{ \ \ when \ \ \ }k\rightarrow \infty .
\end{equation*}

The proof is complete.
\QED

\begin{corollary}
\label{Corollary0fejermax2}Let $p>0.$ Then the maximal operator
\begin{equation*}
\sigma ^{\#}:=\underset{n\in
\mathbb{N}
}{\sup }\left\vert \sigma _{M_{n}}f\right\vert
\end{equation*}%
is bounded from the Hardy space $H_{p}$ to the Lebesgue space $L_{p}$, but it is not
bounded from the Hardy space $H_{p}$ to the Hardy space $H_{p}.$
\end{corollary}

{\bf Proof}:
Since%
\begin{equation}
\left\vert \sigma _{M_{n}}f\right\vert \leq \sigma ^{\#}f  \label{modsup}
\end{equation}%
this result follows immediately from Corollary \ref{Corollary0fejermax} and
Theorem \ref{Theorem1fejermax} so the proof is complete.
\QED

\begin{corollary}
\label{Corollary0fejermax1}Let $p>0.$ Then the operator $\left\vert
\sigma _{M_{n}}f\right\vert $ is bounded from the Hardy space $H_{p}$ to the
space Lebesgue  $L_{p},$ but is not bounded from the Hardy space $H_{p}$ to the Hardy space $%
H_{p}.$
\end{corollary}

{\bf Proof}:
By using (\ref{modsup}) this result immediately follows Corollary \ref%
{Corollary0fejermax} and Theorem \ref{Theorem1fejermax}, so the proof is
complete.
\QED

\subsection{Norm convergence of Vilenkin-Fej\'er means on martingale Hardy
spaces}

In this section we find necessary and sufficient conditions for the
modulus of continuity of martingale Hardy spaces $H_{p}$ for which the
partial sums of Vilenkin-Fourier series convergence in $L_{p}$-norm. We also
study sharpness of these results.

Theorems \ref{theorem4.7bfejer}-\ref{theorem4.7fejer2} can be found in
Persson and Tephnadze \cite{pt1} and Theorem \ref{theorem9fejer}-\ref%
{theorem9fejer1} are proved in Tephnadze \cite{tep6}.

\begin{theorem}
\label{theorem4.7bfejer}Let $0<p<1/2$ and $f\in H_{p}.$ Then there exists an
absolute constant $c_{p}$, depending only on $p,$ such that
\begin{equation*}
\left\Vert \sigma _{n}f\right\Vert _{H_{p}}\leq c_{p}n^{1/p-2}\left\Vert
f\right\Vert _{H_{p}},\ \ \ \left( 0<p<1\right).
\end{equation*}
\end{theorem}

\begin{remark}

We note that the asymptotic behaviour of the sequence $\left\{ n^{1/p-2}:n\in \mathbb{N}\right\} $ in Theorem  \ref{theorem4.7bfejer} can not be improved (c.f. part b) of Theorem \ref{theorem3fejermax2}).

\end{remark}

{\bf Proof}:
According to Theorem \ref{theorem4.7bfejer} \ and Corollary \ref%
{corollary2fejermax} if we invoke Example \ref{example003} we can conclude
that
\begin{equation}
\left\Vert \sigma _{n}f\right\Vert _{H_{p}}\leq \left\Vert \sigma
^{\#}f\right\Vert _{p}+\left\Vert \sigma _{n}f\right\Vert _{p}  \label{fnmn1}
\end{equation}%
\begin{equation*}
\leq \left\Vert f\right\Vert _{H_{p}}+c_{p}n^{1/p-2}\left\Vert f\right\Vert
_{H_{p}}
\end{equation*}%
\begin{equation*}
\leq c_{p}n^{1/p-2}\left\Vert f\right\Vert _{H_{p}}.
\end{equation*}

The proof is complete.
\QED

\begin{theorem}
\label{theorem4.7afejer}Let $f\in H_{1/2}.$ Then there exists an absolute
constant $c$ such that
\begin{equation*}
\left\Vert \sigma _{n}f\right\Vert _{H_{1/2}}\leq c\log ^{2}n\left\Vert
f\right\Vert _{H_{1/2}}.
\end{equation*}
\end{theorem}

\begin{remark}

We note that the asymptotic behaviour of the sequence of sequence $\left\{ \log ^{2}n:n\in \mathbb{N}\right\} $ in Theorem  \ref{theorem4.7bfejer} can not be improved (c.f. part b) of Theorem \ref{theorem2fejermax}).

\end{remark}

{\bf Proof}:
According to Theorem \ref{theorem4.7bfejer} \ and Corollary \ref%
{corollary2fejermax} if we invoke Example \ref{example003} we can write that%
\begin{equation*}
\left\Vert \sigma _{n}f\right\Vert _{H_{1/2}}\leq \left\Vert \sigma
^{\#}f\right\Vert _{1/2}+\left\Vert \sigma _{n}f\right\Vert _{1/2}
\end{equation*}%
\begin{equation*}
\leq \left\Vert f\right\Vert _{H_{1/2}}+c\log ^{2}n\left\Vert f\right\Vert
_{H_{1/2}}
\end{equation*}%
\begin{equation*}
\leq c\log ^{2}n\left\Vert f\right\Vert _{H_{1/2}}.
\end{equation*}

The proof is complete.
\QED

\begin{theorem}
\label{theorem4.7fejer2}Let $0<p<1$ and $f\in H_{p}$. Then there exists an
absolute constant $c_{p}$, depending only on $p,$ such that
\begin{equation*}
\left\Vert \sigma _{M_{n}}f\right\Vert _{H_{p}}\leq c_{p}\left\Vert
f\right\Vert _{H_{p}},\ \ \ \left( 0<p<1\right).
\end{equation*}
\end{theorem}

{\bf Proof}:
In the view of Corollary \ref{Corollary0fejermax2} and Example \ref%
{example003} we can conclude that
\begin{equation*}
\left\Vert \sigma _{M_{n}}f\right\Vert _{H_{p}}\leq \left\Vert \sup_{0\leq
l\leq n}\left\vert \sigma _{M_{l}}f\right\vert \right\Vert _{p}\leq
\left\Vert f\right\Vert _{H_{p}}.
\end{equation*}

The proof is complete.
\QED

\begin{theorem}
\label{theorem9fejer}a) Let $0<p<1/2,$ $f\in H_{p}$ and
\begin{equation}
\omega _{p}\left( \frac{1}{M_{n}},f\right) =o\left( \frac{1}{M_{n}^{1/p-2}}%
\right) \text{when \ }n\rightarrow \infty .  \label{modcond}
\end{equation}%
Then
\begin{equation*}
\left\Vert \sigma _{n}f-f\right\Vert_{H_{p}}\rightarrow 0,\text{ when }%
n\rightarrow \infty .
\end{equation*}
\end{theorem}

\textit{b) Let }$0<p<1/2$. \textit{\ Then there exists a martingale }$f\in
H_{p}$\textit{\ \ for which}\textbf{\ }%
\begin{equation*}
\omega \left( \frac{1}{M_{n}},f\right) _{H_{p}}=O\left( \frac{1}{%
M_{n}^{1/p-2}}\right) \text{ \ when \ }n\rightarrow \infty
\end{equation*}%
\textit{and}
\begin{equation*}
\left\Vert \sigma _{n}f-f\right\Vert _{weak-L_{p}}\nrightarrow 0,\,\,\,\text{%
when\thinspace \thinspace \thinspace }n\rightarrow \infty .
\end{equation*}

{\bf Proof}:
Let $f\in H_{p},$ $0<p<1/2,$ and $M_{N}<n\leq M_{N+1}.$ By a routine
calculation we immediately get that
\begin{equation}
\sigma _{n}S_{M_{N}}f-S_{M_{N}}f  \label{knsn1}
\end{equation}%
\begin{equation*}
=\frac{1}{n}\sum_{k=0}^{M_{N}}S_{k}S_{M_{N}}f+\frac{1}{n}%
\sum_{k=M_{N}+1}^{n}S_{k}S_{M_{N}}f-S_{M_{N}}f
\end{equation*}%
\begin{equation*}
=\frac{1}{n}\sum_{k=0}^{M_{N}}S_{k}f+\frac{1}{n}%
\sum_{k=M_{N}+1}^{n}S_{M_{N}}f-S_{M_{N}}f
\end{equation*}%
\begin{equation*}
=\frac{1}{n}\sum_{k=0}^{M_{N}}S_{k}f+\frac{n-M_{N}}{n}S_{M_{N}}f-S_{M_{N}}f
\end{equation*}%
\begin{equation*}
=\frac{M_{N}}{n}\sigma _{M_{N}}f-\frac{M_{N}}{n}S_{M_{N}}f
\end{equation*}%
\begin{equation*}
=\frac{M_{N}}{n}\left( S_{M_{N}}\sigma _{M_{N}}f-S_{M_{N}}f\right)
\end{equation*}%
\begin{equation*}
=\frac{M_{N}}{n}S_{M_{N}}\left( \sigma _{M_{N}}f-f\right) .
\end{equation*}

Hence,
\begin{equation*}
\left\Vert \sigma _{n}f-f\right\Vert _{H_{p}}^{p}\leq \left\Vert \sigma
_{n}f-\sigma _{n}S_{M_{N}}f\right\Vert _{H_{p}}^{p}
\end{equation*}%
\begin{equation*}
+\left\Vert \sigma _{n}S_{M_{N}}f-S_{M_{N}}f\right\Vert
_{H_{p}}^{p}+\left\Vert S_{M_{N}}f-f\right\Vert _{H_{p}}^{p}
\end{equation*}%
\begin{equation*}
=\left\Vert \sigma _{n}\left( S_{M_{N}}f-f\right) \right\Vert
_{H_{p}}^{p}+\left\Vert S_{M_{N}}f-f\right\Vert _{H_{p}}^{p}+\left\Vert
\sigma _{n}S_{M_{N}}f-S_{M_{N}}f\right\Vert _{H_{p}}^{p}
\end{equation*}%
\begin{equation*}
\leq c_{p}\left( n^{1-2p}+1\right) \omega _{H_{p}}^{p}\left( \frac{1}{M_{N}}%
,f\right) +\left\Vert S_{M_{N}}\left( \sigma _{M_{N}}f-f\right) \right\Vert
_{p}^{p}.
\end{equation*}

Let $p>0$ and $\ f\in H_{p}.$ According to Corollary \ref{corollary4.43}
and Theorem \ref{theorem4.7fejer2} we get that
\begin{equation}
\Vert S_{M_{k}}f-f\Vert _{H_{p}}\rightarrow 0,\text{ \ when }k\rightarrow
\infty ,\text{ \ \ }\left( p>0\right)  \label{sncon}
\end{equation}%
and%
\begin{equation}
\Vert \sigma _{M_{k}}f-f\Vert _{H_{p}}\rightarrow 0,\text{ \ when }%
k\rightarrow \infty ,\text{ \ \ }\left( p>0\right).  \label{fecon}
\end{equation}%
Hence,%
\begin{equation*}
\left\Vert \sigma _{n}S_{M_{N}}f-S_{M_{N}}f\right\Vert _{H_{p}}^{p}
\end{equation*}%
\begin{equation*}
=\frac{M_{N}^{p}}{n^{p}}\left\Vert S_{M_{N}}\left( \sigma _{M_{N}}f-f\right)
\right\Vert _{H_{p}}^{p}
\end{equation*}%
\begin{equation*}
\leq \left\Vert S_{M_{N}}\left( \sigma _{M_{N}}f-f\right) \right\Vert
_{H_{p}}^{p}
\end{equation*}%
\begin{equation*}
\leq \left\Vert \sigma _{M_{N}}f-f\right\Vert _{H_{p}}^{p}\rightarrow 0,%
\text{ when }k\rightarrow \infty \text{.}
\end{equation*}

It follows that under condition (\ref{modcond}) we have that%
\begin{equation*}
\left\Vert \sigma _{n}f-f\right\Vert _{H_{p}}\rightarrow 0,\text{ when }%
n\rightarrow \infty .
\end{equation*}%
This completes the proof of part a) so we turn to the part b).

First, we consider case $0<p<1/2.$ Let $f=\left( f^{n}:n\in \mathbb{N}%
\right) $ be the martingale defined in Lemma \ref{example2.5}, where
\begin{equation*}
\lambda _{k}=\frac{\lambda }{M_{2\alpha _{k}}^{1/p-2}}.\text{ \ }
\end{equation*}

By combining (\ref{2AA}) and (\ref{4AA}) in Lemma \ref{example2.5} we
conclude that the martingale $f\in H_{p}\ \ $and%
\begin{equation*}
\omega _{H_{p}}(\frac{1}{M_{n}},f)\leq \left( \sum_{\left\{ k,2\alpha
_{k}\geq n\right\} }\frac{1}{M_{2\alpha _{k}}^{1-2p}}\right) ^{1/p}
\end{equation*}%
\begin{equation*}
=O\left( \frac{1}{M_{n}^{1/p-2}}\right) ,\text{ when \ }n\rightarrow \infty .
\end{equation*}

A simple calculation gives that%
\begin{equation}
\widehat{f}(j)=\left\{
\begin{array}{ll}
M_{2\alpha _{i}}, & \text{\thinspace }j\in \left\{ M_{2\alpha
_{i}},...,M_{2\alpha _{i}+1}-1\right\} ,\text{ }i\in \mathbb{N}_{+}, \\
0, & j\notin \bigcup\limits_{i=0}^{\infty }\left\{ M_{2\alpha
_{i}},...,M_{2\alpha _{i}+1}-1\right\} .%
\end{array}%
\right.  \label{29}
\end{equation}

By using (\ref{29}) we obtain that%
\begin{equation*}
\Vert \sigma _{M_{2\alpha _{k}}+1}f-f\Vert _{weak-L_{p}}
\end{equation*}%
\begin{equation*}
=\Vert \frac{M_{2\alpha _{k}}\sigma _{M_{2\alpha _{k}}}f}{M_{2\alpha _{k}}+1}%
+\frac{S_{M_{2\alpha _{k}}}f}{M_{2\alpha _{k}}+1}+\frac{M_{2\alpha _{k}}\psi
_{M_{2\alpha _{k}}+1}}{M_{2\alpha _{k}}+1}-\frac{M_{2\alpha _{k}}f}{%
M_{2\alpha _{k}}+1}-\frac{f}{M_{2\alpha _{k}}+1}\Vert _{weak-L_{p}}
\end{equation*}%
\begin{equation*}
\geq \frac{M_{2\alpha _{k}}}{M_{2\alpha _{k}}+1}\Vert \psi _{M_{2\alpha
_{k}}}\Vert _{weak-L_{p}}
\end{equation*}%
\begin{equation*}
-\frac{M_{2\alpha _{k}}}{M_{2\alpha _{k}}+1}\Vert \sigma _{M_{2\alpha
_{k}}}f-f\Vert _{weak-L_{p}}
\end{equation*}%
\begin{equation*}
-\frac{1}{M_{2\alpha _{k}}+1}\Vert S_{M_{2\alpha _{k}}}f-f\Vert _{weak-L_{p}}
\end{equation*}%
\begin{equation*}
\geq \frac{M_{2\alpha _{k}}}{M_{2\alpha _{k}}+1}-o\left( 1\right) \text{, \
\ when \ }k\rightarrow \infty .
\end{equation*}%

Therefore,
\begin{equation*}
\limsup\limits_{k\rightarrow \infty }\Vert \sigma _{M_{2\alpha
_{k}}+1}f-f\Vert _{weak-L_{p}}
\end{equation*}%
\begin{equation*}
\geq \limsup\limits_{k\rightarrow \infty }\frac{M_{2\alpha _{k}}}{M_{2\alpha
_{k}}+1}\geq c>0.
\end{equation*}

Hence also part b) is proved so the proof is complete.
\QED

\begin{theorem}
\label{theorem9fejer1}a) Let $f\in H_{1/2}$ and
\begin{equation}
\omega _{H_{1/2}}\left( \frac{1}{M_{n}},f\right) =o\left( \frac{1}{n^{2}}%
\right) ,\text{ when \ }n\rightarrow \infty.  \label{fecon2}
\end{equation}%
Then
\begin{equation*}
\left\Vert \sigma _{n}f-f\right\Vert_{H_{1/2}}\rightarrow 0,\text{ when }%
n\rightarrow \infty .
\end{equation*}%
b) There exists a martingale $f\in H_{1/2}$\ \ for which
\begin{equation*}
\omega _{H_{1/2}}\left( \frac{1}{M_{n}},f\right) =O\left( \frac{1}{n^{2}}%
\right) ,\text{ \ when \ }n\rightarrow \infty
\end{equation*}%
and
\begin{equation*}
\left\Vert \sigma _{n}f-f\right\Vert _{1/2}\nrightarrow 0,\,\,\,\text{%
when\thinspace \thinspace \thinspace }n\rightarrow \infty .
\end{equation*}
\end{theorem}

{\bf Proof}:
Let $f\in H_{1/2}$ and $M_{N}<n\leq M_{N+1}.$ By using identity (\ref{knsn1}%
) we find that%
\begin{equation*}
\left\Vert \sigma _{n}f-f\right\Vert _{H_{1/2}}^{1/2}\leq \left\Vert \sigma
_{n}f-\sigma _{n}S_{M_{N}}f\right\Vert _{H_{1/2}}^{1/2}
\end{equation*}%
\begin{equation*}
+\left\Vert \sigma _{n}S_{M_{N}}f-S_{M_{N}}f\right\Vert
_{H_{1/2}}^{1/2}+\left\Vert S_{M_{N}}f-f\right\Vert _{H_{1/2}}^{1/2}
\end{equation*}%
\begin{equation*}
=\left\Vert \sigma _{n}\left( S_{M_{N}}f-f\right) \right\Vert
_{H_{1/2}}^{1/2}+\left\Vert S_{M_{N}}f-f\right\Vert
_{H_{1/2}}^{1/2}+\left\Vert \sigma _{n}S_{M_{N}}f-S_{M_{N}}f\right\Vert
_{H_{1/2}}^{1/2}
\end{equation*}%
\begin{equation*}
\leq c\left( \log ^{2}n+1\right) \omega _{H_{1/2}}^{1/2}\left( \frac{1}{M_{N}%
},f\right) +\left\Vert \sigma _{n}S_{M_{N}}f-S_{M_{N}}f\right\Vert
_{H_{1/2}}^{1/2}.
\end{equation*}%
Let $p>0$ and $f\in H_{p}.$ By combining (\ref{sncon}), (\ref{fecon}) and
equality (\ref{knsn1}) we have that
\begin{equation}
\left\Vert \sigma _{n}S_{M_{N}}f-S_{M_{N}}f\right\Vert _{H_{1/2}}^{1/2} \label{10q}
\end{equation}
\begin{equation*}
=\frac{M_{N}}{n}\left\Vert S_{M_{N}}\left( \sigma _{M_{N}}f-f\right)
\right\Vert _{H_{1/2}}^{1/2}\rightarrow 0 \text{ when }k\rightarrow \infty
\text{.}
\end{equation*}

It follows that under condition (\ref{fecon2}) we obtain that%
\begin{equation*}
\left\Vert \sigma _{n}f-f\right\Vert _{H_{1/2}}\rightarrow 0,\text{ when }%
n\rightarrow \infty .
\end{equation*}%
and the proof of part a) is complete.

Let $\ f=\left( f^{\left( n\right) }:n\in
\mathbb{N}
\right) $ be a martingale from the Lemma \ref{example2.6.1} where $p=1/2$.
Then $\ f\in H_{1/2}$ and%
\begin{equation*}
\omega _{H_{1/2}}(\frac{1}{M_{n}},f)=O\left( \frac{1}{n^{2}}\right) ,\text{
\ when }n\rightarrow \infty .
\end{equation*}%
Hence,%
\begin{equation}
\sigma _{q_{M_{k}}}f-f  \label{nn}
\end{equation}%
\begin{equation*}
=\frac{M_{2M_{k}}\sigma _{M_{2M_{k}}}f}{q_{M_{k}}}+\frac{1}{q_{M_{k}}}%
\sum_{j=M_{2M_{k}}+1}^{q_{M_{k}}}S_{j}f
\end{equation*}%
\begin{equation*}
-\frac{M_{2M_{k}}f}{q_{M_{k}}}-\frac{q_{M_{k}-1}f}{q_{M_{k}}}.
\end{equation*}%
Let $M_{2M_{k}}<j\leq q_{M_{k}}.$ By combining (\ref{13}) and (\ref{13s}) we
have that
\begin{equation*}
\frac{1}{q_{M_{k}}}\sum_{j=M_{_{2M_{k}}}+1}^{q_{M_{k}}}S_{j}f
\end{equation*}%
\begin{equation*}
=\frac{q_{M_{k}-1}S_{M_{2M_{k}}}f}{q_{M_{k}}}+\frac{M_{2M_{k}}\psi
_{M_{_{2M_{k}}}}}{q_{M_{k}}M_{k}^{2}}\sum_{j=1}^{q_{M_{k}-1}}D_{j}
\end{equation*}%
\begin{equation*}
=\frac{q_{M_{k}-1}S_{M_{2M_{k}}}f}{q_{M_{k}}}+\frac{M_{2M_{k}}\psi
_{M_{_{2M_{k}}}}q_{M_{k}-1}K_{q_{M_{k}-1}}}{q_{M_{k}}M_{k}^{2}}.
\end{equation*}

By using (\ref{nn}) we get that%
\begin{equation}
\Vert \sigma _{q_{M_{k}}}f-f\Vert _{1/2}^{1/2}  \label{a11}
\end{equation}%
\begin{equation*}
\geq \frac{c}{M_{k}}\Vert q_{M_{k}-1}K_{q_{M_{k}-1}}\Vert _{1/2}^{1/2}
\end{equation*}%
\begin{equation*}
-\left( \frac{M_{2M_{k}}}{q_{M_{k}}}\right) ^{1/2}\Vert \sigma
_{M_{_{2M_{k}}}}f-f\Vert _{1/2}^{1/2}
\end{equation*}%
\begin{equation*}
-\left( \frac{q_{M_{k}-1}}{q_{M_{k}}}\right) ^{1/2}\Vert
S_{M_{2M_{k}}}f-f\Vert _{1/2}^{1/2}.
\end{equation*}

Let \
\begin{equation*}
x\in I_{2M_{k}}^{2s,2\eta },\,\,\,\,\,\,s=\eta +2,\eta +3,M_{k}-2.
\end{equation*}

By applying Lemma \ref{lemma3} we find that%
\begin{equation*}
q_{M_{k}-1}\left\vert K_{q_{M_{k}-1}}\left( x\right) \right\vert \geq \frac{%
M_{2\eta }M_{2s}}{4}.
\end{equation*}

Hence,%
\begin{equation}
\int_{G_{m}}\left\vert q_{M_{k}-1}K_{q_{M_{k}-1}}\right\vert ^{1/2}d\mu
\label{33}
\end{equation}%
\begin{equation}
\geq c\sum_{\eta =1}^{M_{k}-4}\sum_{s=\eta
+2}^{M_{k}-2}\sum_{x_{2s+1}=0}^{m_{2s+1}-1}...\sum_{x_{_{2\alpha
_{k}-1}}=0}^{m_{2M_{k}-1}-1}\int_{I_{2M_{k}}^{2s,2\eta }}\left\vert
q_{M_{k}-1}K_{q_{M_{k}-1}}\right\vert ^{1/2}d\mu  \notag
\end{equation}%
\begin{equation*}
\geq c\sum_{\eta =1}^{M_{k}-4}\sum_{s=\eta +2}^{M_{k}-2}\frac{1}{2M_{2s}}%
M_{2s}^{1/2}M_{2\eta }^{1/2}\geq cM_{k}.
\end{equation*}

By combining (\ref{a11}) and (\ref{33}) we obtain that%
\begin{equation*}
\underset{k\rightarrow \infty }{\lim \sup }\Vert \sigma _{q_{M_{k}}}f-f\Vert
_{1/2}\geq c>0.
\end{equation*}

Thus also the part b) is proved so the proof is complete.
\QED

\subsection{{Strong convergence of Vilenkin-Fej\'er means on  martingale
Hardy spaces}\protect\bigskip\protect\bigskip\protect\bigskip}

The first result in this section is due to Blahota and Tephnadze \cite{bt1}:

\begin{theorem}
\label{theorem2fejerstrong}a) Let $0<p<1/2$ and $f\in H_{p}.$ Then there
exists an absolute constant $c_{p},$ depending only on $p$, such that
\end{theorem}

\begin{equation*}
\overset{\infty }{\underset{k=1}{\sum }}\frac{\left\Vert \sigma
_{k}f\right\Vert _{p}^{p}}{k^{2-2p}}\leq c_{p}\left\Vert f\right\Vert
_{H_{p}}^{p}.
\end{equation*}%
\textit{b) Let }$0<p<1/2$\textit{\ and }$\left\{ \Phi _{k}\ :k\in \mathbb{N}%
\right\} \ $be\textit{\ any non-decreasing sequence satisfying the
conditions }$\Phi _{n}\uparrow \infty $\textit{\ and}
\begin{equation}
\overline{\underset{k\rightarrow \infty }{\lim }}\frac{k^{2-2p}}{\Phi _{k}}%
=\infty .  \label{cond000}
\end{equation}%
\textit{Then there exists a martingale }$f\in H_{p}$\textit{\ such that }%
\begin{equation*}
\underset{k=1}{\overset{\infty }{\sum }}\frac{\left\Vert \sigma
_{k}f\right\Vert _{weak-L_{p}}^{p}}{\Phi _{k}}=\infty .
\end{equation*}

{\bf Proof}:
According to Lemma \ref{lemma2.1} it suffices to show that%
\begin{equation*}
\overset{\infty }{\underset{m=1}{\sum }}\frac{\left\Vert \sigma
_{m}a\right\Vert _{p}^{p}}{m^{2-2p}}\leq c<\infty
\end{equation*}%
for every $p$-atom $a$ with support$\ I$, $\mu \left( I\right) =M_{N}^{-1}.$
We may assume that $I=I_{N}$ and $n>M_{N}$.

Let $x\in I_{N}.$ Since $\sigma _{n}$ is bounded from $L_{\infty }$ to $%
L_{\infty }$ (see (\ref{fn4}) in Corollary \ref{lemma7kn}) and $\left\Vert
a\right\Vert _{\infty }\leq M_{N}^{1/p}$ we obtain that
\begin{equation*}
\int_{I_{N}}\left\vert \sigma _{m}a\right\vert ^{p}d\mu \leq c\frac{%
\left\Vert a\right\Vert _{\infty }^{p}}{M_{N}}\leq c<\infty ,\text{ \ }%
0<p<1/2.
\end{equation*}%
Hence,
\begin{equation*}
\overset{\infty }{\underset{m=1}{\sum }}\frac{\int_{I_{N}}\left\vert \sigma
_{m}a\right\vert ^{p}d\mu }{m^{2-2p}}\leq \overset{\infty }{\underset{m=1}{%
\sum }}\frac{1}{m^{2-2p}}\leq c<\infty .
\end{equation*}

It is easy to see that
\begin{equation*}
\left\vert \sigma _{m}a\left( x\right) \right\vert
\end{equation*}%
\begin{equation*}
\leq \int_{I_{N}}\left\vert a\left( t\right) \right\vert \left\vert
K_{m}\left( x-t\right) \right\vert d\mu \left( t\right)
\end{equation*}%
\begin{equation*}
\leq \left\Vert a\left( x\right) \right\Vert _{\infty
}\int_{I_{N}}\left\vert K_{m}\left( x-t\right) \right\vert d\mu \left(
t\right)
\end{equation*}%
\begin{equation*}
\leq M_{N}^{1/p}\int_{I_{N}}\left\vert K_{m}\left( x-t\right) \right\vert
d\mu \left( t\right) .
\end{equation*}

Let $x\in I_{N}^{k,l},\,0\leq k<l\leq N.$ Then, from Corollary \ref{lemma5aa}
it follows that
\begin{equation}
\left\vert \sigma _{m}a\left( x\right) \right\vert \leq
c_{p}M_{l}M_{k}M_{N}^{1/p-2}.  \label{12q}
\end{equation}

By combining (\ref{1.1}) and (\ref{12q}) if we invoke also estimates (\ref%
{nec3}) and (\ref{nec4}) we obtain that%
\begin{equation}
\int_{\overline{I_{N}}}\left\vert \sigma _{m}a\left( x\right) \right\vert
^{p}d\mu \left( x\right)  \label{7aaa}
\end{equation}%
\begin{equation*}
=\overset{N-2}{\underset{k=0}{\sum }}\overset{N-1}{\underset{l=k+1}{\sum }}%
\sum\limits_{x_{j}=0,\text{ }j\in \{l+1,\dots
,N-1\}}^{m_{j-1}}\int_{I_{N}^{k,l}}\left\vert \sigma _{m}a\left( x\right)
\right\vert ^{p}d\mu \left( x\right)
\end{equation*}%
\begin{equation*}
+c_{p}\overset{N-1}{\underset{k=0}{\sum }}\int_{I_{N}^{k,N}}\left\vert
\sigma _{m}a\left( x\right) \right\vert ^{p}d\mu \left( x\right)
\end{equation*}%
\begin{equation}
\leq c_{p}\overset{N-2}{\underset{k=0}{\sum }}\overset{N-1}{\underset{l=k+1}{%
\sum }}\frac{m_{l+1}\dotsm m_{N-1}}{M_{N}}M_{l}^{p}M_{k}^{p}M_{N}^{1-2p}
\notag
\end{equation}%
\begin{equation*}
+c_{p}\overset{N-1}{\underset{k=0}{\sum }}\frac{1}{M_{N}}M_{k}^{p}M_{N}^{1-p}
\end{equation*}%
\begin{equation*}
\leq c_{p}M_{N}^{1-2p}\overset{N-2}{\underset{k=0}{\sum }}\overset{N-1}{%
\underset{l=k+1}{\sum }}\frac{M_{l}^{p}M_{k}^{p}}{M_{l}}+c_{p}\overset{N-1}{%
\underset{k=0}{\sum }}\frac{M_{k}^{p}}{M_{N}^{p}}
\end{equation*}%
\begin{equation*}
\leq c_{p}M_{N}^{1-2p}.
\end{equation*}

Let $0<p<1/2.$ By using (\ref{7aaa}) we get that%
\begin{equation*}
\overset{\infty }{\underset{m=M_{N}+1}{\sum }}\frac{\int_{\overline{I_{N}}%
}\left\vert \sigma _{m}a\left( x\right) \right\vert ^{p}d\mu \left( x\right)
}{m^{2-2p}}
\end{equation*}%
\begin{equation*}
\leq \overset{\infty }{\underset{m=M_{N}+1}{\sum }}\frac{cM_{N}^{1-2p}}{%
m^{2-p}}<c<\infty .
\end{equation*}
and the proof of part a) is complete.

Under condition (\ref{cond000}) there exists an increasing numbers $\left\{
\alpha _{k}:k\in \mathbb{N}\right\} $, such that
\begin{equation*}
\underset{k\rightarrow \infty }{\lim }\frac{cM_{\left\vert \alpha
_{k}\right\vert +1}^{2-2p}}{\Phi _{M_{\left\vert \alpha _{k}\right\vert +1}}}%
=\infty .
\end{equation*}

There exists a sequence $\left\{ \alpha _{k}:\text{ }k\in \mathbb{N}\right\}
\subset \left\{ n_{k}:\text{ }k\in \mathbb{N}\right\} $ such that%
\begin{equation*}
\left\vert \alpha _{k}\right\vert >2,\text{ for \ }k\in \mathbb{N}
\end{equation*}%
and%
\begin{equation}
\sum_{\eta =0}^{\infty }\frac{\Phi _{M_{\left\vert \alpha _{\eta
}\right\vert +1}}^{1/2}}{M_{\left\vert \alpha _{\eta }\right\vert }^{1-p}}
\label{121}
\end{equation}%
\begin{equation*}
=m_{\left\vert \alpha _{\eta }\right\vert }^{1-p}\sum_{\eta =0}^{\infty }%
\frac{\Phi _{M_{\left\vert \alpha _{\eta }\right\vert +1}}^{1/2}}{%
M_{\left\vert \alpha _{\eta }\right\vert +1}^{1-p}}<c<\infty .
\end{equation*}

Let $f$ be the martingale defined in Example \ref{example2.6} in the case when
\begin{equation}
\lambda _{k}=\frac{\lambda \Phi _{M_{\left\vert \alpha _{k}\right\vert
+1}}^{1/2p}}{M_{\left\vert \alpha _{k}\right\vert }^{1/p-1}}.  \label{121g}
\end{equation}

By applying (\ref{121}) we can conclude that $f\in H_{p}.$

By now using (\ref{121g}) with $\lambda _{k}$ defined by (\ref{10AA})
we readily get that
\begin{equation}
\widehat{f}(j)  \label{6aa}
\end{equation}%
\begin{equation*}
=\left\{
\begin{array}{ll}
\Phi _{M_{\left\vert \alpha _{k}\right\vert +1}}^{1/2p}, & \text{\thinspace }%
j\in \left\{ M_{\left\vert \alpha _{k}\right\vert },\dots ,\text{ ~}%
M_{\left\vert \alpha _{k}\right\vert +1}-1\right\} ,\text{ }k\in \mathbb{N}%
\text{,} \\
0, & \text{\thinspace }j\notin \bigcup\limits_{k=0}^{\infty }\left\{
M_{\left\vert \alpha _{k}\right\vert },\dots ,\text{ ~}M_{\left\vert \alpha
_{k}\right\vert +1}-1\right\} .\text{ }%
\end{array}%
\right.
\end{equation*}

Hence,%
\begin{equation}
\sigma _{_{\alpha _{k}}}f=\frac{1}{\alpha _{k}}\sum_{j=1}^{M_{\left\vert
\alpha _{k}\right\vert }}S_{j}f+\frac{1}{\alpha _{k}}\sum_{j=M_{\left\vert
\alpha _{k}\right\vert }+1}^{\alpha _{k}}S_{j}f  \label{7aa}
\end{equation}%
\begin{equation*}
:=III+IV.
\end{equation*}%

Let $M_{\left\vert n_{k}\right\vert }<j<\alpha _{k}.$ According to (\ref{6aa}%
) and (\ref{12AA}), we can write that%
\begin{equation}
S_{j}f  \label{8aa}
\end{equation}%
\begin{equation}
=\sum_{\eta =0}^{k-1}\Phi _{M_{\left\vert \alpha _{\eta }\right\vert
+1}}^{1/2p}\left( D_{M_{\left\vert \alpha _{\eta }\right\vert
+1}}-D_{M_{\left\vert \alpha _{\eta }\right\vert }}\right)  \notag
\end{equation}%
\begin{equation*}
+\Phi _{M_{\left\vert \alpha _{k}\right\vert +1}}^{1/2p}\psi _{M_{\left\vert
\alpha _{k}\right\vert }}D_{j-M_{\left\vert \alpha _{k}\right\vert }}.
\end{equation*}

By using (\ref{8aa}) in $IV$ we obtain that
\begin{equation}
IV=\frac{\alpha _{k}-M_{\left\vert n_{k}\right\vert }}{\alpha _{k}}%
\sum_{\eta =0}^{k-1}\Phi _{M_{\left\vert \alpha _{\eta }\right\vert
+1}}^{1/2p}\left( D_{M_{\left\vert \alpha _{\eta }\right\vert
+1}}-D_{M_{\left\vert \alpha _{\eta }\right\vert }}\right)  \label{9aa}
\end{equation}%
\begin{equation*}
+\frac{\Phi _{M_{\left\vert \alpha _{k}\right\vert +1}}^{1/2p}\psi
_{M_{\left\vert \alpha _{k}\right\vert }}}{\alpha _{k}}\sum_{j=M_{\left\vert
\alpha _{k}\right\vert }+1}^{\alpha _{k}}D_{j-M_{\left\vert \alpha
_{k}\right\vert }}
\end{equation*}%
\begin{equation*}
:=IV_{1}+IV_{2}.
\end{equation*}

We calculate each term separately. First
we define a set of positive numbers $n\in \mathbb{N}_{1,}$, for which $%
\left\langle n\right\rangle =1,$ that is,
\begin{equation*}
\mathbb{N}_{1}:=\left\{ n\in \mathbb{N}:\text{ }n=n_{1}M_{1}+\sum_{i=2}^{%
\left\vert n\right\vert }n_{i}M_{i}\right\} ,
\end{equation*}
where $n_{1}\in \left\{ 1,\ldots ,m_{1}-1\right\} $ and $n_{i}\in \left\{
0,\ldots ,m_{i}-1\right\} ,$ for $i\geq 2.$

Let $\alpha _{k}\in \mathbb{N}_{1}$ and $x\in I_{2}^{0,1}$. Since $\alpha
_{k}-M_{\left\vert \alpha _{k}\right\vert }\in \mathbb{N}_{2}\ \ $and $%
\left\vert \alpha _{k}\right\vert \neq \left\langle \alpha _{k}\right\rangle
$ for every $\left\vert \alpha _{k}\right\vert >2,$ by applying Lemma \ref%
{lemma8ccc}, we find that
\begin{equation*}
\left\vert IV_{2}\right\vert =\frac{c_{p}\Phi _{M_{\left\vert \alpha
_{k}\right\vert +1}}^{1/2p}}{\alpha _{k}}\left\vert \sum_{j=1}^{\alpha
_{k}-M_{\left\vert \alpha _{k}\right\vert }}D_{_{j}}\left( x\right)
\right\vert
\end{equation*}%
\begin{equation*}
=\frac{c_{p}\Phi _{M_{\left\vert \alpha _{k}\right\vert +1}}^{1/2p}}{\alpha
_{k}}\left\vert \left( \alpha _{k}-M_{\left\vert \alpha _{k}\right\vert
}\right) K_{\alpha _{k}-M_{\left\vert \alpha _{k}\right\vert }}\left(
x\right) \right\vert
\end{equation*}%
\begin{equation*}
\geq \frac{c_{p}\Phi _{M_{\left\vert \alpha _{k}\right\vert +1}}^{1/2p}}{%
\alpha _{k}}.
\end{equation*}

Let $x\in I_{2}^{0,1},$ $n\geq 2$ and $1\leq s_{n}\leq m_{n}-1.$ By
combining Corollary \ref{dn2.3}, Lemma \ref{lemma2} and (\ref{mag}) in Lemma %
\ref{lemma6kn} we have that%
\begin{equation*}
K_{s_{n}M_{n}}\left( x\right) =D_{M_{n}}\left( x\right) =0,\text{ \ \ \ for
\ \ }n\geq 2.
\end{equation*}%
Since $\left\vert \alpha _{k}\right\vert >2,$ $k\in \mathbb{%
\mathbb{N}
}$, we obtain that%
\begin{equation}
IV_{1}=0,\text{ \ for }x\in I_{2}^{0,1}.  \label{10aaa}
\end{equation}%
Moreover, if we invoke (\ref{12AA}) and (\ref{11AA}) with (\ref{121g}) we
get that%
\begin{equation*}
S_{j}f=\left\{
\begin{array}{c}
\Phi ^{1/2p}\left( M_{\left\vert \alpha _{s}\right\vert +1}\right) \psi
_{M_{\left\vert \alpha _{s}\right\vert }}D_{_{j-M_{\left\vert \alpha
_{s}\right\vert }}}, \\
M_{\left\vert \alpha _{s}\right\vert }<j\leq M_{\left\vert \alpha
_{s}\right\vert +1},\text{ \ \ \ }s\in \mathbb{N}_{+} \\
0, \\
M_{\left\vert \alpha _{s}\right\vert +1}<j\leq M_{\left\vert \alpha
_{s+1}\right\vert },\text{ }s\in \mathbb{N}_{+}%
\end{array}%
\right.
\end{equation*}
and
\begin{equation}
III=\frac{1}{n}\sum_{\eta =0}^{k-1}\Phi ^{1/2p}\left( M_{\left\vert \alpha
_{\eta }\right\vert +1}\right) \psi _{M_{\left\vert \alpha _{\eta
}\right\vert }}\sum_{v=M_{\left\vert \alpha _{\eta }\right\vert
}+1}^{M_{\left\vert \alpha _{\eta }\right\vert +1}}D_{v-M_{\left\vert \alpha
_{\eta }\right\vert }}  \label{30ao}
\end{equation}%
\begin{equation*}
=\frac{1}{n}\sum_{\eta =0}^{k-1}\Phi ^{1/2p}\left( M_{\left\vert \alpha
_{\eta }\right\vert +1}\right) \psi _{M_{\left\vert \alpha _{\eta
}\right\vert }}\sum_{v=1}^{\left( m_{\left\vert \alpha _{\eta }\right\vert
}-1\right) M_{\left\vert \alpha _{\eta }\right\vert }}D_{v-M_{\left\vert
\alpha _{\eta }\right\vert }}
\end{equation*}%
\begin{equation*}
=\frac{1}{n}\sum_{\eta =0}^{k-1}\Phi ^{1/2p}\left( M_{\left\vert \alpha
_{\eta }\right\vert +1}\right) \left( m_{\left\vert \alpha _{\eta
}\right\vert }-1\right) M_{\left\vert \alpha _{\eta }\right\vert }K_{\left(
m_{\left\vert \alpha _{\eta }\right\vert }-1\right) M_{\left\vert \alpha
_{\eta }\right\vert }}\left( x\right) =0.
\end{equation*}

Let $0<p<1/2,$ $n\in \mathbb{N}_{2}$ and $M_{\left\vert \alpha
_{k}\right\vert }<n<M_{\left\vert \alpha _{k}\right\vert +1}.$ By combining (%
\ref{7aa})-(\ref{30ao}) we find that%
\begin{equation*}
\left\Vert \sigma _{n}f\right\Vert _{weak-L_{p}}^{p}
\end{equation*}%
\begin{equation*}
\geq \frac{c_{p}\Phi _{M_{\left\vert \alpha _{k}\right\vert +1}}^{1/2}}{%
\alpha _{k}^{p}}\mu \left\{ x\in I_{2}^{0,1}:\text{ }\left\vert
IV_{2}\right\vert \geq \frac{c_{p}\Phi _{M_{\left\vert \alpha
_{k}\right\vert +1}}^{1/2p}}{\alpha _{k}}\right\}
\end{equation*}%
\begin{equation*}
\geq \frac{c_{p}\Phi _{M_{\left\vert \alpha _{k}\right\vert +1}}^{1/2}}{%
\alpha _{k}^{p}}\mu \left\{ I_{2}^{0,1}\right\} \geq \frac{c_{p}\Phi
_{M_{\left\vert \alpha _{k}\right\vert +1}}^{1/2}}{M_{\left\vert \alpha
_{k}\right\vert +1}^{p}}.
\end{equation*}

Since
\begin{equation*}
\underset{\left\{ n\in \mathbb{N}_{1}:M_{k}\leq n\leq M_{k+1}\right\} }{\sum
}1\geq cM_{k},
\end{equation*}%
we obtain that
\begin{equation*}
\underset{n=1}{\overset{\infty }{\sum }}\frac{\left\Vert \sigma
_{n}f\right\Vert _{weak-L_{p}}^{p}}{\Phi _{n}}
\end{equation*}%
\begin{equation*}
\geq \underset{\left\{ n\in \mathbb{N}_{1}:\text{ }M_{\left\vert \alpha
_{k}\right\vert }<n<M_{\left\vert \alpha _{k}\right\vert +1}\right\} }{\sum }%
\frac{\left\Vert \sigma _{n}f\right\Vert _{weak-L_{p}}^{p}}{\Phi _{n}}
\end{equation*}%
\begin{equation*}
\geq \frac{c_{p}}{\Phi _{M_{\left\vert \alpha _{k}\right\vert +1}}^{1/2}}%
\underset{\left\{ \mathbb{N}_{1}:\text{ }M_{\left\vert \alpha
_{k}\right\vert }<n<M_{\left\vert \alpha _{k}\right\vert +1}\right\} }{\sum }%
\frac{1}{M_{\left\vert \alpha _{k}\right\vert +1}^{p}}
\end{equation*}%
\begin{equation*}
\geq \frac{c_{p}M_{\left\vert \alpha _{k}\right\vert +1}^{1-p}}{\Phi
_{M_{\left\vert \alpha _{k}\right\vert +1}}^{1/2}}\rightarrow \infty ,\text{
when \ \ }k\rightarrow \infty .
\end{equation*}%
Therefore also part b) is proved so the proof is complete.
\QED

\begin{corollary}
\label{corollary1fejerstrong12}Let $0<p<1/2$ and $f\in H_{p}.$ Then there
exists an absolute constant $c_{p},$ depending only on $p$, such that
\end{corollary}

\begin{equation*}
\overset{\infty }{\underset{k=1}{\sum }}\frac{\left\Vert \sigma
_{k}f\right\Vert _{H_{p}}^{p}}{k^{2-2p}}\leq c_{p}\left\Vert f\right\Vert
_{H_{p}}^{p},
\end{equation*}

{\bf Proof}:
By combining (\ref{fnmn1}) and Theorem \ref{theorem2fejerstrong} we have that%
\begin{equation*}
\overset{\infty }{\underset{k=1}{\sum }}\frac{\left\Vert \sigma
_{k}f\right\Vert _{H_{p}}^{p}}{k^{2-2p}}
\end{equation*}%
\begin{equation*}
\leq \overset{\infty }{\underset{k=1}{\sum }}\frac{\left\Vert \sigma
_{k}f\right\Vert _{p}^{p}+\left\Vert \sigma ^{\#}f\right\Vert _{p}^{p}}{%
k^{2-2p}}
\end{equation*}%
\begin{equation*}
\leq c\overset{\infty }{\underset{k=1}{\sum }}\frac{\left\Vert \sigma
_{k}f\right\Vert _{p}^{p}+\left\Vert f\right\Vert _{p}^{p}}{k^{2-2p}}
\end{equation*}%
\begin{equation*}
\leq c_{p}\overset{\infty }{\underset{k=1}{\sum }}\frac{\left\Vert \sigma
_{k}f\right\Vert _{p}^{p}}{k^{2-2p}}+c_{p}\overset{\infty }{\underset{k=1}{%
\sum }}\frac{\left\Vert f\right\Vert _{p}^{p}}{k^{2-2p}}
\end{equation*}%
\begin{equation*}
\leq c\left\Vert f\right\Vert _{H_{p}}^{p}
\end{equation*}%
and the proof is complete.
\QED

\begin{corollary}
\label{corollary1fejerstrong1}Let $0<p<1/2$ and $f\in H_{p}.$ Then there
exists an absolute constant $c_{p},$ depending only on $p$, such that

\begin{equation*}
\frac{1}{n}\overset{n}{\underset{k=1}{\sum }}\frac{\left\Vert \sigma
_{k}f\right\Vert _{H_{p}}^{p}}{k^{1-2p}}\leq c_{p}\left\Vert f\right\Vert
_{H_{p}}^{p},
\end{equation*}%
\begin{equation*}
\frac{1}{n}\overset{n}{\underset{k=1}{\sum }}\frac{\left\Vert \sigma
_{k}f-f\right\Vert _{H_{p}}^{p}}{k^{1-2p}}=0,
\end{equation*}%
and

\begin{equation*}
\frac{1}{n}\overset{n}{\underset{k=1}{\sum }}\frac{\left\Vert \sigma
_{k}f\right\Vert _{H_{p}}^{p}}{k^{1-2p}}=\left\Vert f\right\Vert
_{H_{p}}^{p}.
\end{equation*}
\end{corollary}

In Blahota and Tephnadze \cite{bt1} also the endpoint case $p=1/2$ was considered and
the following result was proved:

\begin{theorem}
\label{theorem1fejerstrong}Let $f\in H_{1/2}.$ Then there exists an absolute
constant $c$ such that
\end{theorem}

\begin{equation*}
\frac{1}{\log n}\overset{n}{\underset{k=1}{\sum }}\frac{\left\Vert \sigma
_{k}f\right\Vert _{1/2}^{1/2}}{k}\leq c\left\Vert f\right\Vert
_{H_{1/2}}^{1/2}.
\end{equation*}

{\bf Proof}:
In view of Lemma \ref{lemma2.1} it suffices to show that%
\begin{equation*}
\frac{1}{\log n}\overset{n}{\underset{m=1}{\sum }}\frac{\left\Vert \sigma
_{m}a\right\Vert _{1/2}^{1/2}}{m}\leq c<\infty
\end{equation*}%
for every $p$-atom $a$ with support$\ I$, $\mu \left( I\right) =M_{N}^{-1}.$
Analogously to the previous Theorems we may assume that $I=I_{N}$ and $%
n>M_{N}$.

Let $x\in I_{N}.$ Since $\sigma _{n}$ is bounded from $L_{\infty }$ to $%
L_{\infty }$ (see (\ref{fn4}) in Corollary \ref{lemma7kn}) and $\left\Vert
a\right\Vert _{\infty }\leq M_{N}^{2}$ we obtain that
\begin{equation*}
\int_{I_{N}}\left\vert \sigma _{m}a\right\vert ^{1/2}d\mu \leq \frac{%
\left\Vert a\right\Vert _{\infty }^{1/2}}{M_{N}}\leq c<\infty .
\end{equation*}%
Hence,
\begin{equation*}
\frac{1}{\log n}\overset{n}{\underset{m=1}{\sum }}\frac{\int_{I_{N}}\left%
\vert \sigma _{m}a\right\vert ^{1/2}d\mu }{m}
\end{equation*}%
\begin{equation*}
\leq \frac{c}{\log n}\overset{n}{\underset{m=1}{\sum }}\frac{1}{m}\leq
c<\infty .
\end{equation*}

It is easy to see that
\begin{equation*}
\left\vert \sigma _{m}a\left( x\right) \right\vert
\end{equation*}%
\begin{equation*}
\leq \int_{I_{N}}\left\vert a\left( t\right) \right\vert \left\vert
K_{m}\left( x-t\right) \right\vert d\mu \left( t\right)
\end{equation*}%
\begin{equation*}
\leq \left\Vert a\left( x\right) \right\Vert _{\infty
}\int_{I_{N}}\left\vert K_{m}\left( x-t\right) \right\vert d\mu \left(
t\right)
\end{equation*}%
\begin{equation*}
\leq M_{N}^{2}\int_{I_{N}}\left\vert K_{m}\left( x-t\right) \right\vert d\mu
\left( t\right) .
\end{equation*}

Let $x\in I_{N}^{k,l},\,0\leq k<l<N.$ Then, from Lemma \ref{lemma5} it
follows that
\begin{equation}
\left\vert \sigma _{m}a\left( x\right) \right\vert \leq \frac{%
cM_{l}M_{k}M_{N}}{m}.  \label{12q1}
\end{equation}

Let $x\in I_{N}^{k,N},\,0\leq k<N.$ Then, according to Lemma \ref{lemma5},
we have that%
\begin{equation}
\left\vert \sigma _{m}a\left( x\right) \right\vert \leq cM_{N}^{2}\frac{M_{k}%
}{M_{N}}\leq cM_{k}M_{N}.  \label{12aq1}
\end{equation}

By combining (\ref{1.1}) with (\ref{12q1}) and (\ref{12aq1}) we find that%
\begin{equation*}
\int_{\overline{I_{N}}}\left\vert \sigma _{m}a\left( x\right) \right\vert
^{1/2}d\mu \left( x\right)
\end{equation*}%
\begin{equation*}
=\overset{N-2}{\underset{k=0}{\sum }}\overset{N-1}{\underset{l=k+1}{\sum }}%
\sum\limits_{x_{j}=0,\text{ }j\in \{l+1,\dots
,N-1\}}^{m_{j-1}}\int_{I_{N}^{k,l}}\left\vert \sigma _{m}a\left( x\right)
\right\vert ^{1/2}d\mu \left( x\right)
\end{equation*}%
\begin{equation*}
+\overset{N-1}{\underset{k=0}{\sum }}\int_{I_{N}^{k,N}}\left\vert \sigma
_{m}a\left( x\right) \right\vert ^{1/2}d\mu \left( x\right)
\end{equation*}%
\begin{equation}
\leq c\overset{N-2}{\underset{k=0}{\sum }}\overset{N-1}{\underset{l=k+1}{%
\sum }}\frac{m_{l+1}\dotsm m_{N-1}}{M_{N}}\frac{%
M_{l}^{1/2}M_{k}^{1/2}M_{N}^{1/2}}{m^{1/2}}  \notag
\end{equation}%
\begin{equation*}
+c\overset{N-1}{\underset{k=0}{\sum }}\frac{1}{M_{N}}M_{k}^{1/2}M_{N}^{1/2}
\end{equation*}%
\begin{equation*}
\leq cM_{N}^{1/2}\overset{N-2}{\underset{k=0}{\sum }}\overset{N-1}{\underset{%
l=k+1}{\sum }}\frac{M_{l}^{1/2}M_{k}^{1/2}}{m^{1/2}M_{l}}
\end{equation*}%
\begin{equation*}
+c\overset{N-1}{\underset{k=0}{\sum }}\frac{M_{k}^{1/2}}{M_{N}^{1/2}}\leq
\frac{cM_{N}^{1/2}N}{m^{1/2}}+c.
\end{equation*}

It follows that%
\begin{equation*}
\frac{1}{\log n}\overset{n}{\underset{m=M_{N}+1}{\sum }}\frac{\int_{%
\overline{I_{N}}}\left\vert \sigma _{m}a\left( x\right) \right\vert
^{1/2}d\mu \left( x\right) }{m}
\end{equation*}%
\begin{equation*}
\leq \frac{1}{\log n}\overset{n}{\underset{m=M_{N}+1}{\sum }}\left( \frac{%
cM_{N}^{1/2}N}{m^{3/2}}+\frac{c}{m}\right) <c<\infty .
\end{equation*}

The proof is complete.
\QED

\begin{corollary}
\label{corollary1fejerstrong2}Let $f\in H_{1/2}.$ Then
\begin{equation*}
\frac{1}{\log n}\overset{n}{\underset{k=1}{\sum }}\frac{\left\Vert \sigma
_{k}f\right\Vert _{H_{1/2}}^{1/2}}{k}\leq c\left\Vert f\right\Vert
_{H_{1/2}}^{1/2},\text{ }
\end{equation*}%
\begin{equation*}
\lim_{n\rightarrow \infty }\frac{1}{\log n}\overset{n}{\underset{k=1}{\sum }}%
\frac{\left\Vert \sigma _{k}f-f\right\Vert _{H_{1/2}}^{1/2}}{k}=0
\end{equation*}%
and

\begin{equation*}
\lim_{n\rightarrow \infty }\frac{1}{\log n}\overset{n}{\underset{k=1}{\sum }}%
\frac{\left\Vert \sigma _{k}f\right\Vert _{H_{1/2}}^{1/2}}{k}=\left\Vert
f\right\Vert _{H_{1/2}}^{1/2}.
\end{equation*}
\end{corollary}

{\bf Proof}:
According to Theorem \ref{theorem1fejerstrong} we have that%
\begin{equation*}
\frac{1}{\log n}\overset{n}{\underset{k=1}{\sum }}\frac{\left\Vert \sigma
_{k}f\right\Vert _{H_{1/2}}^{1/2}}{k}
\end{equation*}%
\begin{equation*}
\leq \frac{1}{\log n}\overset{n}{\underset{k=1}{\sum }}\frac{\left\Vert
\sigma _{k}f\right\Vert _{1/2}^{1/2}+\left\Vert \sigma ^{\#}f\right\Vert
_{1/2}^{1/2}}{k}
\end{equation*}%
\begin{equation*}
\leq \frac{c}{\log n}\overset{n}{\underset{k=1}{\sum }}\frac{\left\Vert
\sigma _{k}f\right\Vert _{1/2}^{1/2}+\left\Vert f\right\Vert _{1/2}^{1/2}}{k}
\end{equation*}%
\begin{equation*}
\leq \frac{c}{\log n}\overset{n}{\underset{k=1}{\sum }}\frac{\left\Vert
\sigma _{k}f\right\Vert _{1/2}^{1/2}}{k}+\frac{c}{\log n}\overset{n}{%
\underset{k=1}{\sum }}\frac{\left\Vert f\right\Vert _{1/2}^{1/2}}{k}
\end{equation*}%
\begin{equation*}
\leq c\left\Vert f\right\Vert _{H_{1/2}}^{1/2}.
\end{equation*}

The first inequality is proved. Analogously we can prove the second and third
statement. We leave out the details.
\QED

\newpage

\section{\protect\bigskip\textbf{Vilenkin-N\"orlund means on  martingale
Hardy spaces}}

\vspace{0.5cm}

\subsection{Some classical results on Vilenkin-N\"orlund means}

It is well-known in the literature that the so-called N\"orlund means are
generalizations of the Fej\'er, Ces\`aro and logarithmic means. The N\"orlund
summation is a general summability method. Therefore it is of prior
interest to study the behavior of operators related to N\"orlund means of
Fourier series with respect to orthonormal systems.

In \cite{tep4} it was proved that there exists a martingale $f\in H_{p},$ $%
(0<p\leq 1),$ such that the maximal operator of N\"orlund logarithmic means $%
L^{\ast }$ is not bounded in the Lebesgue space $L_{p}.$ Riesz logarithmic means with
respect to the trigonometric system was studied by a lot of authors.We
mention, for instance, the paper by Szasz \cite{Sz} and Yabuta \cite{Ya}.
These means with respect to the Walsh and Vilenkin systems were investigated
by Simon \cite{Si1} and G\`{a}t \cite{Ga1}. Blahota and G\'at \cite{bg}
considered norm summability of N\"orlund logarithmic means and showed that
Riesz logarithmic means $R_{n}$ have better approximation properties on some
unbounded Vilenkin groups than the Fej\'er means. Moreover, in \cite{tep10}
it was proved that the maximal operator of Riesz`s means is bounded from the
Hardy space $H_{p}$ to the Lebesgue space $L_{p}$ for $p>1/2$ but not
when $0<p\leq 1/2.$

M\'oricz and Siddiqi \cite{Mor} investigate the approximation properties of
some special N\"orlund means of Walsh-Fourier series of $L_{p}$ functions in
norm. The case when $\left\{ q_{k}=1/k:k\in \mathbb{N}\right\} $ was
excluded, since the methods of M\'oricz and Siddiqi are not applicable to N\"o%
rlund logarithmic means. In \cite{Ga2} G\'{a}t and Goginava proved some
convergence and divergence properties of the N\"orlund logarithmic means of
functions in the class of continuous functions and in the Lebesgue space $%
L_{1}.$ In particular, they gave a negative answer to the question of M\'oricz
and Siddiqi \cite{Mor}. G\'{a}t and Goginava \cite{Ga3} proved that for each
measurable function satisfying
\begin{equation*}
\phi \left( u\right) =o\left( u\log ^{1/2}u\right) ,\text{ \ \ \ \ \ when \ \
\ }u\rightarrow \infty ,
\end{equation*}%
there exists an integrable function $f$ such that
\begin{equation*}
\int_{G_{m}}\phi \left( \left\vert f\left( x\right) \right\vert \right) d\mu
\,\left( x\right) <\infty
\end{equation*}%
and that there exists a set with positive measure such that the
Walsh-logarithmic means of the function diverges on this set.
Fridli, Manchanda and Siddiqi \cite{FMS} improved and extended results of M\'oricz and Siddiqi \cite{Mor} to dyadic homogeneous Banach spaces and Martingale Hardy spaces.

In \cite{gog8} Goginava investigated the behavior of Ces\`aro means of
Walsh-Fourier series in detail. In the two-dimensional case approximation
properties of N\"orlund and Ces\`aro means were considered by Nagy (see \cite{na}%
, \cite{n} and \cite{nagy}). The maximal operator $\sigma ^{\alpha ,\ast }$ $%
\left( 0<\alpha <1\right) $ of the $\left( C,\alpha \right) $ means of
Vilenkin systems was investigated by Weisz \cite{we6}. In this paper Weisz
proved that $\sigma ^{\alpha ,\ast }$ is bounded from the martingale space $%
H_{p}$ to the Lebesgue space $L_{p}$ for $p>1/\left( 1+\alpha \right) .$ Goginava
\cite{gog4} gave a counterexample which shows that boundedness does not
hold for $0<p\leq 1/\left( 1+\alpha \right) .$ Weisz and Simon \cite{sw}
showed that the maximal operator $\sigma ^{\alpha ,\ast }$ is bounded from
the Hardy space $H_{1/\left( 1+\alpha \right) }$ to the space $%
weak-L_{1/\left( 1+\alpha \right) }$.

\subsection{Maximal operators of N\"orlund means on martingale Hardy spaces%
}

In this section we first state our main result concerning the maximal
operator of the N\"orlund summation method (see (\ref{1.2})). We also show
that this result is in a sense sharp. The proof can be found in Persson,
Tephnadze and Wall \cite{ptw}.

\begin{theorem}
\label{theorem1norlund}a) The maximal operator $t^{\ast }$ of \ the
summability method (\ref{1.2}) with non-decreasing sequence $\{q_{k}:k\in
\mathbb{N}\}$ is bounded from the Hardy space $H_{1/2}$ to the space $%
weak-L_{1/2}.$

The statement in a) is sharp in the following sense:

b) Let $0<p<1/2$ and $\{q_{k}:k\in \mathbb{N}\}$ be a non-decreasing sequence
satisfying the condition
\begin{equation}
\frac{q_{0}}{Q_{n}}\geq \frac{c}{n},\text{ \ \ }\left( c>0\right) .
\label{cond1nor}
\end{equation}%
Then there exists a martingale $f\in H_{p},$ such that \textit{\ }
\begin{equation*}
\underset{n\in
\mathbb{N}
}{\sup }\left\Vert t_{n}f\right\Vert _{weak-L_{p}}=\infty .
\end{equation*}
\end{theorem}

{\bf Proof}:
Let the sequence $\{q_{k}:k\in \mathbb{N}\}$ be non-decreasing. By combining
(\ref{2b}) and (\ref{2bb}) and using Abel transformation we get that%
\begin{equation*}
\left\vert t_{n}f\right\vert \leq \left\vert \frac{1}{Q_{n}}\overset{n}{%
\underset{j=1}{\sum }}q_{n-j}S_{j}f\right\vert
\end{equation*}%
\begin{equation*}
\leq \frac{1}{Q_{n}}\left( \overset{n-1}{\underset{j=1}{\sum }}\left\vert
q_{n-j}-q_{n-j-1}\right\vert j\left\vert \sigma _{j}f\right\vert
+q_{0}n\left\vert \sigma _{n}f\right\vert \right)
\end{equation*}%
\begin{equation*}
\leq \frac{c}{Q_{n}}\left( \overset{n-1}{\underset{j=1}{\sum }}\left(
q_{n-j}-q_{n-j-1}\right) j+q_{0}n\right) \sigma ^{\ast }f\leq \sigma ^{\ast
}f
\end{equation*}%
so that
\begin{equation}
t^{\ast }f\leq \sigma ^{\ast }f.  \label{12aaaa}
\end{equation}

In view of (\ref{12aaaa}) we can conclude that the maximal operators $%
t^{\ast }$ is bounded from the Hardy space $H_{1/2}$ to the space $%
weak-L_{1/2}.$ The proof of part a) is complete.

Let $f:=\left( f^{\left( n\right) },n\in
\mathbb{N}
\right) $ be the martingale defined in Example \ref{example2.5.1} in the case when $0<p<q=1/2$%
. We can write that%
\begin{equation*}
t_{M_{2\alpha _{k}}+1}f=\frac{1}{Q_{M_{2\alpha _{k}}+1}}\sum_{j=0}^{M_{2%
\alpha _{k}}}q_{M_{2\alpha _{k}}+1-j}S_{j}f+\frac{q_{0}}{Q_{M_{2\alpha
_{k}}+1}}S_{M_{2\alpha _{k}}+1}f
\end{equation*}%
\begin{equation*}
:=I+II.
\end{equation*}

According to (\ref{sn101}) in Example \ref{example2.5.1} we can conclude
that
\begin{equation}
\left\vert I\right\vert \leq \frac{1}{Q_{M_{2\alpha _{k}}+1}}%
\sum_{j=0}^{M_{2\alpha _{k}}}q_{j}\left\vert S_{M_{2\alpha
_{k}}+1-j}f\right\vert  \label{10nor}
\end{equation}%
\begin{equation*}
\leq \frac{2\lambda M_{2\alpha _{k-1}}^{1/p}}{\alpha _{k-1}^{1/2}}\frac{1}{%
Q_{M_{2\alpha _{k}}+1}}\sum_{j=0}^{M_{2\alpha _{k}}}q_{M_{2\alpha _{k}}+1-j}
\end{equation*}%
\begin{equation*}
\leq \frac{2\lambda M_{2\alpha _{k-1}}^{1/p}}{\alpha _{k-1}^{1/2}}\leq \frac{%
M_{\alpha _{k}}^{1/p-2}}{16\alpha _{k}^{3/2}}.
\end{equation*}

If we now apply (\ref{sn102}) for $II$ we find that
\begin{equation}
\left\vert II\right\vert =\frac{q_{0}}{Q_{M_{2\alpha _{k}}+1}}\left\vert
S_{M_{2\alpha _{k}}+1}f\right\vert \geq \frac{q_{0}}{Q_{M_{2\alpha _{k}}+1}}%
\frac{M_{2\alpha _{k}}^{1/p-1}}{4\alpha _{k}^{1/2}}.  \label{100nor}
\end{equation}

Without lost the generality we may assume that $c=1$ in (\ref{cond1nor}). By
combining (\ref{10nor}) and (\ref{100nor}) we obtain that
\begin{equation*}
\left\vert t_{M_{2\alpha _{k}}+1}f\right\vert \geq \left\vert II\right\vert
-\left\vert I\right\vert
\end{equation*}%
\begin{equation*}
\geq \frac{q_{0}}{Q_{M_{2\alpha _{k}}+1}}\frac{M_{2\alpha _{k}}^{1/p-1}}{%
4\alpha _{k}^{1/2}}-\frac{4\lambda M_{\alpha _{k}}^{1/p-2}}{\alpha _{k}^{3/2}%
}
\end{equation*}%
\begin{equation*}
\geq \frac{M_{2\alpha _{k}}^{1/p-2}}{4\alpha _{k}^{1/2}}-\frac{4\lambda
M_{\alpha _{k}}^{1/p-2}}{\alpha _{k}^{3/2}}\geq \frac{M_{2\alpha
_{k}}^{1/p-2}}{8\alpha _{k}^{1/2}}.
\end{equation*}

On the other hand,%
\begin{equation}
\mu \left\{ x\in G_{m}:\left\vert t_{M_{2\alpha _{k}}+1}f\left( x\right)
\right\vert \geq \frac{M_{2\alpha _{k}}^{1/p-2}}{8\alpha _{k}^{1/2}}\right\}
=\mu \left( G_{m}\right) =1.  \label{88nor}
\end{equation}

Let $0<p<1/2.$ Then
\begin{equation}
\frac{M_{2\alpha _{k}}^{1/p-2}}{8\alpha _{k}^{1/2}}\cdot \mu \left\{ x\in
G_{m}:\left\vert t_{M_{2\alpha _{k}}+1}f\left( x\right) \right\vert \geq
\frac{M_{2\alpha _{k}}^{1/p-2}}{8\alpha _{k}^{1/2}}\right\}  \label{99nor}
\end{equation}%
\begin{equation*}
=\frac{M_{2\alpha _{k}}^{1/p-2}}{8\alpha _{k}^{1/2}}\rightarrow \infty ,%
\text{ \ when }k\rightarrow \infty .
\end{equation*}

\textbf{\ }The proof is complete.
\QED

The next Theorems (Theorems \ref{theorem3fejermax2a}-\ref{theorem2fejermaxaa}%
) are due to Persson, Tephnadze and Wall \cite{ptw2}:

\begin{theorem}
\label{theorem3fejermax2a}a) Let $0<p<1/2$ and the sequence $\left\{ q_{k}:k\in
\mathbb{N}\right\} $ be non-decreasing. Then the \bigskip maximal operator%
\begin{equation*}
\widetilde{t}_{p,1}^{\ast }f:=\sup_{n\in \mathbb{N}}\frac{\left\vert
t_{n}f\right\vert }{\left( n+1\right) ^{1/p-2}}
\end{equation*}%
is bounded from the Hardy martingale space $H_{p}$ to the Lebesgue space $L_{p}.$
\end{theorem}

\begin{remark}
\label{Remark2nor}Since the Fej\'er means are examples of N\"orlund means with
non-decreasing sequence $\left\{ q_{k}:k\in \mathbb{N}\right\} $ we
immediately obtain from part b) of Theorem \ref{theorem3fejermax2} that the
asymptotic behaviour of the sequence of weights
\begin{equation*}
\left\{ 1/\left( k+1\right) ^{1/p-2}:k\in \mathbb{N}\right\}
\end{equation*}
in N\"orlund means in Theorem \ref{theorem3fejermax2a} can not be improved.
\end{remark}

{\bf Proof}:
The idea of proof is similar to that of part a) of Theorem \ref{theorem3fejermax2}, but in more general situation so we give the details.

First we note that $t_{n}$ is bounded from $L_{\infty }$ to $L_{\infty }$
(see Corollary \ref{Corollary3nn}). Let $a$ be an arbitrary p-atom, with
support$\ I$ and $\mu \left( I\right) =M_{N}^{-1}.$ We may assume that $%
I=I_{N}.$ It is easy to see that $S_{n}\left( a\right) =t_{n}\left( a\right)
=0,$ when $n\leq M_{N}$. Therefore, we can suppose that $n>M_{N}$. Hence,

\begin{equation}
t_{n}\left( a\right) =\frac{1}{Q_{n}}\overset{n}{\underset{k=1}{\sum }}%
q_{n-k}S_{k}\left( a\right)  \label{nom}
\end{equation}%
\begin{equation*}
=\frac{1}{Q_{n}}\overset{n}{\underset{k=M_{N}}{\sum }}q_{n-k}S_{k}\left(
a\right)
\end{equation*}%
\begin{equation*}
\frac{1}{Q_{n}}\overset{n}{\underset{k=M_{N}}{\sum }}q_{n-k}\int_{I_{N}}a%
\left( x\right) D_{k}\left( x-t\right) d\mu \left( t\right) .
\end{equation*}

Since $\left\Vert a\right\Vert _{\infty }\leq M_{N}^{1/p}$ it follows that
\begin{equation*}
\frac{\left\vert t_{n}\left( a\right) \right\vert }{\left( n+1\right)
^{1/p-2}}\leq \frac{\left\Vert a\right\Vert _{\infty }}{\left( n+1\right)
^{1/p-2}}\int_{I_{N}}\left\vert \frac{1}{Q_{n}}\overset{n}{\underset{k=M_{N}}%
{\sum }}q_{n-k}D_{k}\left( x-t\right) \right\vert d\mu \left( t\right)
\end{equation*}%
\begin{equation*}
\leq \frac{M_{N}^{1/p}}{\left( n+1\right) ^{1/p-2}}\int_{I_{N}}\left\vert
\frac{1}{Q_{n}}\overset{n}{\underset{k=M_{N}}{\sum }}q_{n-k}D_{k}\left(
x-t\right) \right\vert d\mu \left( t\right) .
\end{equation*}

Let $x\in I_{N}^{k,l},\,0\leq k<l\leq N.$ From Lemma \ref{lemma5aaa} we can
deduce that
\begin{equation}
\frac{\left\vert t_{n}\left( a\right) \right\vert }{\left( n+1\right)
^{1/p-2}}\leq \frac{c_{p}M_{N}^{1/p}}{M_{N}^{1/p-2}}\frac{M_{l}M_{k}}{%
M_{N}^{2}}=c_{p}M_{l}M_{k}.  \label{1112a}
\end{equation}

The expression on the right-hand side of (\ref{1112a}) does not depend on $n.$ Therefore,
\begin{equation}
\left\vert \widetilde{t}_{p,1}^{\ast }a\left( x\right) \right\vert \leq
c_{p}M_{l}M_{k},\text{ \ \ \ for \ \ \ }x\in I_{N}^{k,l},\,0\leq k<l\leq N.
\label{1112aM}
\end{equation}

By combining (\ref{1.1}) with (\ref{1112aM}) we obtain that
\begin{equation*}
\int_{\overline{I_{N}}}\left\vert \widetilde{t}_{p,1}^{\ast }a\left(
x\right) \right\vert ^{p}d\mu \left( x\right)
\end{equation*}%
\begin{equation*}
=\overset{N-2}{\underset{k=0}{\sum }}\overset{N-1}{\underset{l=k+1}{\sum }}%
\sum\limits_{x_{j}=0,j\in
\{l+1,...,N-1\}}^{m_{j-1}}\int_{I_{N}^{k,l}}\left\vert \widetilde{t}%
_{p,1}^{\ast }a\left( x\right) \right\vert ^{p}d\mu \left( x\right)
\end{equation*}%
\begin{equation*}
+\overset{N-1}{\underset{k=0}{\sum }}\int_{I_{N}^{k,N}}\left\vert \widetilde{%
t}_{p,1}^{\ast }a\left( x\right) \right\vert ^{p}d\mu \left( x\right)
\end{equation*}%
\begin{equation*}
\leq c_{p}\overset{N-2}{\underset{k=0}{\sum }}\overset{N-1}{\underset{l=k+1}{%
\sum }}\frac{m_{l+1}...m_{N-1}}{M_{N}}M_{l}^{p}M_{k}^{p}+c_{p}\overset{N-1}{%
\underset{k=0}{\sum }}\frac{1}{M_{N}}M_{N}^{p}M_{k}^{p}
\end{equation*}%
\begin{equation*}
\leq c_{p}\overset{N-2}{\underset{k=0}{\sum }}\overset{N-1}{\underset{l=k+1}{%
\sum }}\frac{M_{l}^{p}M_{k}^{p}}{M_{l}}+c_{p}\overset{N-1}{\underset{k=0}{%
\sum }}\frac{M_{k}^{p}}{M_{N}^{1-p}}<c<\infty .
\end{equation*}

 The proof is complete by using this estimate and Lemma \ref{lemma2.2}.
\QED

\begin{theorem}
\label{theorem2fejermaxaa} Let the sequence $\left\{ q_{k}:k\in \mathbb{N}%
\right\} $ be non-decreasing. Then the maximal operator%
\begin{equation*}
\overset{\sim }{t}_{1}^{\ast }f:=\sup_{n\in \mathbb{N}}\frac{\left\vert
t_{n}f\right\vert }{\log ^{2}\left( n+1\right) }
\end{equation*}%
\textit{is bounded from the Hardy space }$H_{1/2}$\textit{\ to the Lebesgue space }$%
L_{1/2}.$
\end{theorem}

\begin{remark}
\label{Remark3nor}Since the Fej\'er means are examples of N\"orlund means with
non-decreasing sequence $\left\{ q_{k}:k\in \mathbb{N}\right\} $ we
immediately obtain from part b) of Theorem \ref{theorem2fejermax} that the
asymptotic behaviour of the sequence of \ weights
\begin{equation*}
\left\{ 1/\log ^{2}\left( n+1\right) :n\in \mathbb{N}\right\}
\end{equation*}
in N\"orlund means in Theorem \ref{theorem2fejermaxaa} can not be improved.
\end{remark}

{\bf Proof}: The idea of proof is similar as that of part a) of Theorem \ref{theorem2fejermax}, but since this case is more general we give the details.

Analogously to Theorem \ref{theorem3fejermax2a} we may assume that $n>M_{N}$
and $a$ be a  $p$-atom  with support $I=I_{N}.$ Since $\left\Vert a\right\Vert
_{\infty }\leq M_{N}^{2}$ if we apply (\ref{nom}) we obtain that
\begin{equation*}
\frac{\left\vert t_{n}\left( a\right) \right\vert }{\log ^{2}\left(
n+1\right) }
\end{equation*}%
\begin{equation*}
\leq \frac{\left\Vert a\right\Vert _{\infty }}{\log ^{2}\left( n+1\right) }%
\int_{I_{N}}\left\vert \frac{1}{Q_{n}}\overset{n}{\underset{k=M_{N}}{\sum }}%
q_{n-k}D_{k}\left( x-t\right) \right\vert d\mu \left( t\right)
\end{equation*}%
\begin{equation*}
\leq \frac{M_{N}^{2}}{\log ^{2}\left( n+1\right) }\int_{I_{N}}\left\vert
\frac{1}{Q_{n}}\overset{n}{\underset{k=M_{N}}{\sum }}q_{n-k}D_{k}\left(
x-t\right) \right\vert d\mu \left( t\right) .
\end{equation*}

Let $x\in I_{N}^{k,l},\,0\leq k<l\leq N.$ Then, from Lemma \ref{lemma5aaa}
it follows that
\begin{equation*}
\frac{\left\vert t_{n}\left( a\right) \right\vert }{\log ^{2}\left(
n+1\right) }\leq \frac{cM_{N}^{2}}{\log ^{2}\left(
n+1\right) }\frac{M_{l}M_{k}}{M_{N}^{2}}=\frac{%
cM_{l}M_{k}}{\log ^{2}\left(
n+1\right) }
\end{equation*}%
and%
\begin{equation}
\left\vert \widetilde{t_{1}}^{\ast }a\left( x\right) \right\vert \leq \frac{%
cM_{l}M_{k}}{N^{2}}.  \label{2t122M}
\end{equation}%
\

By combining (\ref{1.1}) with (\ref{2t122M}) we find that%
\begin{equation*}
\int_{\overline{I_{N}}}\left\vert \widetilde{t_{1}}^{\ast }a\left( x\right)
\right\vert ^{1/2}d\mu \left( x\right)
\end{equation*}%
\begin{equation*}
=\overset{N-2}{\underset{k=0}{\sum }}\overset{N-1}{\underset{l=k+1}{\sum }}%
\sum\limits_{x_{j}=0,j\in
\{l+1,...,N-1\}}^{m_{j}-1}\int_{I_{N}^{k,l}}\left\vert \overset{\sim }{t_{1}}%
^{\ast }a\left( x\right) \right\vert ^{1/2}d\mu \left( x\right)
\end{equation*}%
\begin{equation*}
+\overset{N-1}{\underset{k=0}{\sum }}\int_{I_{N}^{k,N}}\left\vert \widetilde{%
t_{1}}^{\ast }a\left( x\right) \right\vert ^{1/2}d\mu \left( x\right)
\end{equation*}%
\begin{equation*}
\leq c\overset{N-2}{\underset{k=0}{\sum }}\overset{N-1}{\underset{l=k+1}{%
\sum }}\frac{m_{l+1}...m_{N-1}}{M_{N}}\frac{M_{l}^{1/2}M_{k}^{1/2}}{N}
\end{equation*}%
\begin{equation*}
+c\overset{N-1}{\underset{k=0}{\sum }}\frac{1}{M_{N}}\frac{%
M_{N}^{1/2}M_{k}^{1/2}}{N}
\end{equation*}%
\begin{equation*}
\leq c\overset{N-2}{\underset{k=0}{\sum }}\overset{N-1}{\underset{l=k+1}{%
\sum }}\frac{M_{k}^{1/2}}{NM_{l}^{1/2}}
\end{equation*}%
\begin{equation*}
+c\overset{N-1}{\underset{k=0}{\sum }}\frac{1}{M_{N}^{1/2}}\frac{M_{k}^{1/2}%
}{N}\leq c<\infty .
\end{equation*}%
The proof is complete.
\QED

The next results deal with a N\"orlund means with non-increasing sequence $%
\left\{ q_{k}:k\in \mathbb{N}\right\} .$ We begin by stating a divergence
result for all such summability methods when $0<p<1/2$, which can be found
in Persson, Tephnade and Wall \cite{ptw}.

\begin{theorem}
\label{theorem2norlunda}Let $0<p<1/2.$ Then, for all N\"orlund means with
non-increasing sequence $\left\{ q_{k}:k\in \mathbb{N}\right\} $ there
exists a martingale $f\in H_{p}$ such that%
\begin{equation*}
\underset{n\in
\mathbb{N}
}{\sup }\left\Vert t_{n}f\right\Vert _{weak-L_{p}}=\infty .
\end{equation*}
\end{theorem}

{\bf Proof}:
For the proof we use the martingale defined in Example \ref{example2.5.1} (see
also the part b) of Theorem \ref{theorem1norlund}).

It is obvious that for every non-increasing sequence $\{q_{k}:k\in
\mathbb{N}
\}$ it automatically holds that
\begin{equation*}
\frac{q_{0}}{Q_{M_{2\alpha _{k}}+1}}\geq \frac{1}{M_{2\alpha _{k}}+1}.
\end{equation*}

Since
\begin{equation*}
t_{M_{2\alpha _{k}}+1}f=\frac{1}{Q_{M_{2\alpha _{k}}+1}}\sum_{j=0}^{M_{2%
\alpha _{k}}}q_{M_{2\alpha _{k}}+1-j}S_{j}f
\end{equation*}%
\begin{equation*}
+\frac{q_{0}}{Q_{M_{2\alpha _{k}}+1}}S_{M_{2\alpha _{k}}+1}f:=I+II.
\end{equation*}

by combining (\ref{10nor}) and (\ref{100nor}) we see that
\begin{equation*}
\left\vert t_{M_{2\alpha _{k}}+1}f\right\vert \geq \left\vert II\right\vert
-\left\vert I\right\vert \geq \frac{M_{2\alpha _{k}}^{1/p-2}}{8\alpha
_{k}^{1/2}}.
\end{equation*}

Analogously to\ (\ref{88nor}) and (\ref{99nor}) we get that%
\begin{equation*}
\sup_{k\in \mathbb{N}}\left\Vert t_{M_{2\alpha _{k}}+1}f\right\Vert _{weak-L_{p}}=\infty .
\end{equation*}

The proof is complete.
\QED

\begin{corollary}
\label{corollary2norlunda}Let $0<p<1/2\ \ $and $t_{n}$ be N\"orlund means with
non-increasing sequence $\{q_{k}:k\in \mathbb{N}\}.$ Then the maximal
operator $t^{\ast }$ is not bounded from the martingale Hardy space $H_{p}$
to the space $weak-L_{p},$ that is there exists a martingale $f\in H_{p},$
such that
\begin{equation*}
\underset{n\in
\mathbb{N}
}{\sup }\left\Vert t^{\ast }f\right\Vert _{weak-L_{p}}=\infty .
\end{equation*}
\end{corollary}

Next, we present necessary condition for the N\"orlund means with
non-increasing sequence $\left\{ q_{k}:k\in
\mathbb{N}
\right\} $, when $1/2\leq p<1,$ which can be found in Persson, Tephnade and
Wall \cite{ptw}.

\begin{theorem}
\label{theorem2norlund} a) Let $0<p<1/\left(1+\alpha \right),$ $0<\alpha
\leq 1,$ and $\{q_{k}:k\in \mathbb{N}\}$ be a non-increasing sequence
satisfying the condition
\begin{equation}
\overline{\lim_{n\rightarrow \infty }}\frac{n^{\alpha }}{Q_{n}}=c>0,\text{ }%
0<\alpha \leq 1.  \label{cond4}
\end{equation}%
Then there exists a martingale $f\in H_{p}$ such that%
\begin{equation}
\underset{n\in
\mathbb{N}
}{\sup }\left\Vert t_{n}f\right\Vert _{weak-L_{p}}=\infty . \label{cond444}
\end{equation}

b) Let $\{q_{k}:k\in
\mathbb{N}
\}$ be a non-increasing sequence satisfying the condition
\begin{equation}
\overline{\lim_{n\rightarrow \infty }}\frac{n^{\alpha }}{Q_{n}}=\infty ,%
\text{ \ \ \ }\left( 0<\alpha \leq 1\right) .  \label{cond22}
\end{equation}

Then there exists an martingale $f\in H_{1/\left( 1+\alpha \right) },$ such
that
\begin{equation}
\underset{n\in
\mathbb{N}
}{\sup }\left\Vert t_{n}f\right\Vert _{weak-L_{1/\left( 1+\alpha \right)
}}=\infty .\label{cond445}
\end{equation}
\end{theorem}

{\bf Proof}:
Under condition (\ref{cond4}) there exists an\ increasing sequence $\left\{
\alpha _{k}:k\in \mathbb{N}\right\} $ of positive integers such that
\begin{equation*}
\frac{M_{2\alpha _{k}}^{\alpha }}{Q_{M_{2\alpha _{k}}+1}}\geq c,\text{ \ }%
k\in \mathbb{N}
\end{equation*}%
and the estimates (\ref{1t2})-(\ref{1t4}) are satisfied. To prove part a) we use the martingale defined in Example \ref%
{example2.5.1} in the case when $0<p<q=1/\left( 1+\alpha \right) .$

Since
\begin{equation*}
t_{M_{2\alpha _{k}}+1}f=\frac{1}{Q_{M_{2\alpha _{k}}+1}}\sum_{j=0}^{M_{2%
\alpha _{k}}}q_{M_{2\alpha _{k}}+1-j}S_{j}f+\frac{1}{Q_{M_{2\alpha _{k}}+1}}%
S_{M_{2\alpha _{k}}+1}f
\end{equation*}%
\begin{equation*}
:=I+II,
\end{equation*}

by combining (\ref{10nor}) and (\ref{100nor}), we get that%
\begin{equation*}
\left\vert t_{M_{2\alpha _{k}}+1}f\right\vert \geq \left\vert II\right\vert
-\left\vert I\right\vert
\end{equation*}%
\begin{equation*}
=\frac{M_{2\alpha _{k}}^{1/p-1}}{4\alpha _{k}^{1/2}}\frac{1}{Q_{M_{2\alpha
_{k}}+1}}-\frac{2\lambda M_{\alpha _{k-1}}^{1/p}}{\alpha _{k-1}^{1/2}}.
\end{equation*}

Without lost the generality we may assume that $c=1$ in (\ref{cond4}). By
using (\ref{1t4}) we find that%
\begin{equation*}
\left\vert t_{M_{2\alpha _{k}}+1}f\right\vert
\end{equation*}%
\begin{equation*}
\geq \frac{M_{2\alpha _{k}}^{1/p-1-\alpha }}{4\alpha _{k}^{1/2}}-\frac{%
2\lambda M_{\alpha _{k-1}}^{1/p}}{\alpha _{k-1}^{1/2}}\geq \frac{M_{2\alpha
_{k}}^{1/p-1-\alpha }}{8\alpha _{k}^{1/2}}
\end{equation*}%
\begin{equation*}
\geq \frac{M_{2\alpha _{k}}^{1/p-1-\alpha }}{4\alpha _{k}^{1/2}}-\frac{%
\lambda M_{\alpha _{k}}^{1/p-1-\alpha }}{16\alpha _{k}^{3/2}}\geq \frac{%
M_{2\alpha _{k}}^{1/p-1-\alpha }}{8\alpha _{k}^{1/2}}
\end{equation*}%
and
\begin{equation*}
\frac{M_{2\alpha _{k}}^{1/p-1-\alpha }}{8\alpha _{k}}\cdot \mu \left\{ x\in
G_{m}:\left\vert t_{M_{2\alpha _{k}}+1}f\right\vert \geq \frac{M_{2\alpha
_{k}}^{1/p-1-\alpha }}{8\alpha _{k}}\right\}
\end{equation*}%
\begin{equation*}
=\frac{M_{2\alpha _{k}}^{1/p-1-\alpha }}{8\alpha _{k}}\rightarrow \infty ,%
\text{ \ when }k\rightarrow \infty .
\end{equation*}

Which means that \ref{cond444} holds and part (a) is proved.

Under condition (%
\ref{cond22}) there exists an increasing sequence $\left\{ \alpha _{k}:k\in
\mathbb{N}
\right\}$ which satisfies the conditions
\begin{equation}
\sum_{k=0}^{\infty }\frac{Q_{M_{2\alpha _{k}}+1}^{1/2\left( 1+\alpha \right)
}}{M_{2\alpha _{k}}^{\alpha /2\left( 1+\alpha \right) }}\leq c<\infty ,
\label{cond23}
\end{equation}%
\begin{equation}
\sum_{\eta =0}^{k-1}Q_{M_{2\alpha _{\eta }}+1}M_{2\alpha _{\eta }}^{\alpha
/2+1}\leq Q_{M_{2\alpha _{k}}+1}M_{2\alpha _{k}}^{\alpha /2+1}
\label{cond24}
\end{equation}%
and%
\begin{equation}
32\lambda Q_{M_{2\alpha _{k-1}}+1}M_{2\alpha _{k-1}}^{\alpha /2+1}<\frac{%
M_{2\alpha _{k}}^{\alpha /2}}{Q_{M_{2\alpha _{k}}+1}^{1/2}},  \label{cond25}
\end{equation}%
where $\lambda =\sup_{n}m_{n}.$

Let the martingale defined in the Example \ref{example2.5}, where $\lambda _{k}$
is defined by%
\begin{equation*}
\lambda _{k}=\frac{\lambda Q_{M_{2\alpha _{k}}+1}^{1/2}}{M_{2\alpha
_{k}}^{\alpha /2}}
\end{equation*}%
for which the sequence $\left\{ \alpha _{k}:k\in
\mathbb{N}
\right\} $ satisfies conditions (\ref{cond23})-(\ref{cond25}) and $a_{k}$
are given by (\ref{3AA}) for $p=1/\left( 1+\alpha \right) .$ If we apply (%
\ref{cond23}) analogously we can conclude that $\,f\in H_{1/\left( 1+\alpha
\right) }.$ Hence,%
\begin{equation}
\widehat{f}(j)=  \label{26nor}
\end{equation}%
\begin{equation*}
\left\{
\begin{array}{ll}
Q_{M_{2\alpha _{k}}+1}^{1/2}M_{2\alpha _{k}}^{\alpha /2}, & \text{\thinspace
}j\in \left\{ M_{2\alpha _{k}},...,M_{2\alpha _{k}+1}-1\right\} \text{ \ }%
k\in
\mathbb{N}
, \\
0, & \text{\thinspace }j\notin \bigcup\limits_{k=1}^{\infty }\left\{
M_{2\alpha _{k}},...,M_{2\alpha _{k}+1}-1\right\} ,%
\end{array}%
\right.
\end{equation*}%
and%
\begin{equation*}
t_{M_{2\alpha _{k}}+1}f
\end{equation*}%
\begin{equation*}
=\frac{1}{Q_{M_{2\alpha _{k}}+1}}\sum_{j=0}^{M_{2\alpha _{k}}}q_{M_{2\alpha
_{k}}+1-j}S_{j}f+\frac{1}{Q_{M_{2\alpha _{k}}+1}}S_{M_{2\alpha
_{k}}+1}f:=III+IV.
\end{equation*}

Let $M_{2\alpha _{s}}<$ $j\leq M_{2\alpha _{s}+1},$ where $s=0,...,k-1.$
From (\ref{26nor}) we have that%
\begin{equation}
\left\vert S_{j}f\right\vert \leq \left\vert \sum_{\eta
=0}^{s-1}Q_{M_{2\alpha _{\eta }}+1}^{1/2}M_{2\alpha _{\eta }}^{\alpha
/2}\left( D_{M_{_{2\alpha _{\eta }+1}}}-D_{M_{_{2\alpha _{\eta }}}}\right)
\right\vert  \label{89nor}
\end{equation}%
\begin{equation}
+Q_{M_{2\alpha _{s}}+1}^{1/2}M_{2\alpha _{s}}^{\alpha /2}\left\vert \left(
D_{j}-D_{M_{_{2\alpha _{s}}}}\right) \right\vert \leq 4\lambda Q_{M_{2\alpha
_{s}}+1}^{1/2}M_{2\alpha _{k-1}}^{\alpha /2+1}.  \notag
\end{equation}

Let $M_{\alpha _{s-1}+1}+1\leq $ $j\leq M_{\alpha _{s}},$ where $\
s=1,...,k. $ Then%
\begin{equation*}
\left\vert S_{j}f\right\vert =\left\vert \sum_{\eta =0}^{s-1}Q_{M_{2\alpha
_{\eta }}+1}^{1/2}M_{2\alpha _{\eta }}^{\alpha /2}\left( D_{M_{_{2\alpha
_{\eta }+1}}}-D_{M_{_{2\alpha _{\eta }}}}\right) \right\vert
\end{equation*}%
\begin{equation*}
\leq \lambda \sum_{\eta =0}^{s-1}Q_{M_{2\alpha _{\eta }}+1}M_{2\alpha _{\eta
}}^{\alpha /2+1}\leq 2\lambda Q_{M_{2\alpha _{s-1}}+1}M_{2\alpha
_{s-1}}^{\alpha /2+1}
\end{equation*}%
\begin{equation*}
\leq 2\lambda Q_{M_{2\alpha _{k-1}}+1}M_{2\alpha _{k-1}}^{\alpha /2+1}
\end{equation*}%
and%
\begin{equation}
\left\vert III\right\vert \leq \frac{1}{Q_{M_{2\alpha _{k}}+1}}%
\sum_{j=0}^{M_{2\alpha _{k}}}q_{M_{2\alpha _{k}}+1-j}\left\vert
S_{j}f\right\vert  \label{101nor}
\end{equation}%
\begin{equation*}
\leq 2\lambda Q_{M_{2\alpha _{k-1}}+1}M_{2\alpha _{k-1}}^{\alpha /2+1}\frac{1%
}{Q_{M_{2\alpha _{k}}+1}}\sum_{j=0}^{M_{2\alpha _{k}}}q_{M_{2\alpha
_{k}}+1-j}
\end{equation*}%
\begin{equation*}
\leq 2\lambda Q_{M_{2\alpha _{k-1}}+1}M_{2\alpha _{k-1}}^{\alpha /2+1}.
\end{equation*}

If we apply (\ref{26nor}) and (\ref{89nor}) we get that%
\begin{equation}
\left\vert IV\right\vert \geq Q_{M_{2\alpha _{k}}+1}^{1/2}M_{\alpha
_{k}}^{\alpha /2}\frac{q_{0}\left\vert D_{M_{2\alpha _{k}}+1}-D_{M_{2\alpha
_{k}}}\right\vert }{Q_{M_{2\alpha _{k}}+1}}  \label{102nor}
\end{equation}%
\begin{equation*}
-\frac{1}{Q_{M_{2\alpha _{k}}+1}}\left\vert S_{M_{2\alpha _{k}}}f\right\vert
\end{equation*}%
\begin{equation*}
\geq \frac{q_{0}M_{2\alpha _{k}}^{\alpha /2}}{Q_{M_{2\alpha _{k}}+1}^{1/2}}%
-2\lambda Q_{M_{2\alpha _{k-1}}+1}M_{2\alpha _{k-1}}^{\alpha /2+1}\geq \frac{%
M_{2\alpha _{k}}^{\alpha /2}}{4Q_{M_{2\alpha _{k}}+1}^{1/2}}.
\end{equation*}

By combining (\ref{101nor}) and (\ref{102nor}) we find that
\begin{equation*}
\left\vert t_{M_{2\alpha _{k}}+1}f\right\vert \geq \left\vert IV\right\vert
-\left\vert III\right\vert
\end{equation*}%
\begin{equation*}
\geq \frac{M_{2\alpha _{k}}^{\alpha /2}}{4Q_{M_{2\alpha _{k}}+1}^{1/2}}%
-2\lambda Q_{M_{2\alpha _{k-1}}+1}M_{2\alpha _{k-1}}^{\alpha /2+1}\geq \frac{%
M_{2\alpha _{k}}^{\alpha /2}}{8Q_{M_{2\alpha _{k}}+1}^{1/2}}.
\end{equation*}

Hence, it yields that%
\begin{equation*}
\frac{M_{2\alpha _{k}}^{\alpha /2\left( 1+\alpha \right) }}{8Q_{M_{2\alpha
_{k}}+1}^{1/2\left( 1+\alpha \right) }}\mu \left\{ x\in G_{m}:\left\vert
t_{M_{2\alpha _{k}}+1}f\left( x\right) \right\vert \geq \frac{M_{2\alpha
_{k}}^{\alpha /2}}{8Q_{M_{2\alpha _{k}}+1}^{1/2}}\right\}
\end{equation*}%
\begin{equation*}
=\frac{M_{2\alpha _{k}}^{\alpha /2}}{8Q_{M_{2\alpha _{k}}+1}^{1/2}}\mu
\left( G_{m}\right) =\frac{M_{2\alpha _{k}}^{\alpha /2}}{8Q_{M_{2\alpha
_{k}}+1}^{1/2}}\rightarrow \infty ,\text{ \ when \ }k\rightarrow \infty .
\end{equation*}

Which means that \ref{cond445} holds and the proof is complete.
\QED

\begin{corollary}
\label{theorem2norlund1000} Let $0<p<1/\left( 1+\alpha \right) ,$ $%
0<\alpha \leq 1$ and $\{q_{k}:k\in \mathbb{N}\}$ be a non-increasing
sequence satisfying the condition (\ref{cond4}). Then there exists a
martingale $f\in H_{p}$ such that%
\begin{equation*}
\left\Vert t^{\ast }f\right\Vert _{weak-L_{p}}=\infty .
\end{equation*}
\end{corollary}

\begin{corollary}
Let $\{q_{k}:k\in
\mathbb{N}
\}$ be a non-increasing sequence satisfying the condition (\ref{cond22}).
Then there exists an martingale $f\in H_{1/\left( 1+\alpha \right) }$ such
that
\begin{equation*}
\left\Vert t^{\ast }f\right\Vert _{weak-L_{1/\left( 1+\alpha \right)
}}=\infty .
\end{equation*}
\end{corollary}

Our next result reads:

\begin{theorem}
\label{Theorem3normax}a) The maximal operator $t^{\ast }$ of the N\"orlund
summability method with non-increasing sequence $\{q_{k}:k\in
\mathbb{N}
\},$ satisfying the condition (\ref{6a}) and (\ref{7a}) is bounded from the
Hardy space $H_{1/\left( 1+\alpha \right) }$ to the space $weak-L_{1/\left(
1+\alpha \right) },$ for $0<\alpha \leq 1.$

b) Let $0<\alpha \leq 1$ and $\{q_{k}:k\in \mathbb{N}\}$ be a non-increasing
sequence satisfying the conditions%
\begin{equation}
\overline{\lim_{n\rightarrow \infty }}\frac{n^{\alpha }}{Q_{n}}\geq
c_{\alpha }>0  \label{cond29}
\end{equation}%
and
\begin{equation}
\left\vert q_{n}-q_{n+1}\right\vert \geq c_{\alpha }n^{\alpha -2},\text{ \ }%
n\in
\mathbb{N}
\text{.}  \label{cond30}
\end{equation}

Then there exists a martingale $f\in H_{1/\left( 1+\alpha \right) }$ such
that
\begin{equation*}
\underset{n\in
\mathbb{N}
}{\sup }\left\Vert t_{n}f\right\Vert _{{1/\left( 1+\alpha \right)
}}=\infty .
\end{equation*}
\end{theorem}

\begin{remark}
Part a) of this result can be found in \cite{ptw}, but below we give a
different proof. Part b) of the theorem is new.
\end{remark}

{\bf Proof}:
By Lemma \ref{lemma2.3} the proof of part a) is complete if we show
that%
\begin{equation*}
t\mu \left\{ x\in \overline{I_{N}}:t^{\ast
}f\geq t^{1+\alpha }\right\} \leq c<\infty, \text{ \ \ \ \ }t\geq 0,
\end{equation*}%
for every $1/\left( 1+\alpha \right) $-atom $a.$ We may assume that $a$ is
an arbitrary $1/\left( 1+\alpha \right) $-atom with support$\ I,$ $\mu
\left( I\right) =M_{N}^{-1}$ and $I=I_{N}.$ It is easy to see that $%
t_{m}\left( a\right) =0,$ when $m\leq M_{N}.$ Therefore we can suppose that $%
m>M_{N}.$

Let $x\in I_{N}.$ Since $t_{m}$ is bounded from $L_{\infty }$ to $L_{\infty
} $ (the boundedness follows from Corollary \ref{corollary3n}) and $%
\left\Vert a\right\Vert _{\infty }\leq cM_{N}^{1/\left( 1+\alpha \right) }$
we obtain that
\begin{equation*}
\left\vert t_{m}a\left( x\right) \right\vert \leq \int_{I_{N}}\left\vert
a\left( t\right) \right\vert \left\vert F_{m}\left( x-t\right) \right\vert
d\mu \left( t\right)
\end{equation*}%
\begin{equation*}
\leq \left\Vert a\right\Vert _{\infty }\int_{I_{N}}\left\vert F_{m}\left(
x-t\right) \right\vert d\mu \left( t\right)
\end{equation*}%
\begin{equation*}
\leq M_{N}^{1+\alpha }\int_{I_{N}}\left\vert F_{m}\left( x-t\right)
\right\vert d\mu \left( t\right) .
\end{equation*}

Let $x\in I_{N}^{k,l},\,0\leq k<l\leq N$ and $m>M_{N}.$ From Corollary \ref%
{corollary3na} we get that
\begin{equation}
\left\vert t_{m}a\left( x\right) \right\vert \leq c_{\alpha
}M_{k}M_{l}^{\alpha }.  \label{saa11}
\end{equation}

The expression on the right-hand side of (\ref{saa11}) does not depend on $m.$ Hence,
\begin{equation}
\left\vert t^{\ast }a\left( x\right) \right\vert \leq c_{\alpha
}M_{k}M_{l}^{\alpha }  \label{saa110}
\end{equation}

Let $n\geq N.$ According to (\ref{saa110}) we conclude that
\begin{equation*}
\mu \left\{ x\in \overline{I_{N}}:t^{\ast }f\geq c_{\alpha }M_{n}^{1+\alpha
}\right\} =0.
\end{equation*}

Thus, we can suppose that $0<n<N.$ Let $\lambda =\sup_{n}m_{n}$ and $\left[ x%
\right] $ denotes the  integer part of $x.$ It is obvious that for fixed $\lambda
$ there exists a positive number $\theta $ so that $\lambda ^{1/\theta }\leq
2.$ Then, for every $k<n$, it yields that%
\begin{equation*}
M_{k}M_{n+\left[ \left( n-k\right) /\theta \right] }^{\alpha }\leq
M_{k}M_{n}^{\alpha }\lambda ^{\left[ \left( n-k-1\right) /\theta \right] }
\end{equation*}%
\begin{equation*}
\leq M_{k}M_{n}^{\alpha }\left( \lambda ^{1/\theta }\right) ^{n-k-1}\leq
M_{k}M_{n}^{\alpha }2^{n-k-1}
\end{equation*}%
\begin{equation*}
\leq M_{k}M_{n}^{\alpha }m_{k}m_{k+1}...m_{n-1}\leq M_{n}^{1+\alpha }.
\end{equation*}

It is obvious that if $n+\left[ \left( n-k-1\right) /\lambda \right] >N$,
for some $k<l<N$ we readily get that%
\begin{equation*}
c_{\alpha }M_{k}M_{l}^{\alpha }\leq c_{\alpha }M_{k}M_{N}^{\alpha }\leq
c_{\alpha }M_{k}M_{n+\left[ \left( n-k-1\right) /\lambda \right] }^{\alpha
}\leq c_{\alpha }M_{n}^{1+\alpha }.
\end{equation*}

It follows that for such $k<l<N$ we have the following estimate
\begin{equation*}
\left\vert t_{m}a\left( x\right) \right\vert \leq c_{\alpha }M_{n}^{1+\alpha
},\text{ for }x\in I_{N}^{k,l}
\end{equation*}%
and
\begin{equation*}
\mu \left\{ x\in I_{N}^{k,l}:t^{\ast }f\geq c_{\alpha }M_{n}^{1+\alpha
}\right\} =0.
\end{equation*}

Therefore, we may assume that $n+\left[ \left( n-k-1\right) /\lambda \right] \leq N.$

By combining (\ref{1.1}) and (\ref{saa110}) we obtain that%
\begin{equation*}
\mu \left\{ x\in \overline{I_{N}}:t^{\ast }f\geq c_{\alpha }M_{n}^{1+\alpha
}\right\}
\end{equation*}%
\begin{equation*}
\leq \overset{N-1}{\underset{k=n}{\sum }}\overset{N}{\underset{l=k+1}{\sum }}%
\sum\limits_{x_{j}=0,\text{ }j\in \{l+1,\dots ,N-1\}}^{m_{j-1}}\left\vert
I_{N}^{k,l}\right\vert
\end{equation*}%
\begin{equation*}
+\overset{n}{\underset{k=0}{\sum }}\overset{N}{\underset{l=n+\left[ \left(
n-k-1\right) /\lambda \right] }{\sum }}\sum\limits_{x_{j}=0,\text{ }j\in
\{l+1,\dots ,N-1\}}^{m_{j-1}}\left\vert I_{N}^{k,l}\right\vert
\end{equation*}%
\begin{equation*}
\leq \overset{N-1}{\underset{k=n}{\sum }}\overset{N}{\underset{l=k+1}{\sum }}%
\frac{1}{M_{l}}+\overset{n}{\underset{k=0}{\sum }}\overset{N}{\underset{l=n+%
\left[ \left( n-k\right) /\theta \right] }{\sum }}\frac{1}{M_{l}}\leq \frac{c%
}{M_{n}}.
\end{equation*}%
Hence,%
\begin{equation*}
\sup_{n\in \mathbb{N}}M_{n}\mu \left\{ x\in \overline{I_{N}}:t^{\ast }f\geq
M_{n}^{1+\alpha }\right\} \leq c<\infty .
\end{equation*}

Part a) is proved.

Under condition (\ref{cond29}) there exists an\ increasing sequence $%
\left\{ \alpha _{k}:k\in \mathbb{N}\right\}$ of positive integers such that
\begin{equation}
\frac{M_{2\alpha _{k}+1}^{\alpha }}{Q_{M_{2\alpha _{k}}+1}}>c_{\alpha }>0,%
\text{ \ }k\in \mathbb{N},  \label{cond32.1}
\end{equation}%
and estimates (\ref{1t2})-(\ref{1t4}) are satisfied. To prove part b) of
Theorem \ref{theorem2norlund} we use the martingale defined in Example \ref%
{example2.5.1} in the case when $p=q=1/\left( 1+\alpha \right).$

 In particular, as we proved there, the martingale belongs to the space $H_{1/\left( 1+\alpha
\right) }.$

Moreover,
\begin{equation*}
\widehat{f}(j)=\left\{
\begin{array}{ll}
\frac{M_{2\alpha _{k}}^{\alpha }}{\alpha _{k}^{1/2}}, & j\in \left\{
M_{2\alpha _{k}},...,M_{2\alpha _{k}+1}-1\right\} ,\text{ }k\in
\mathbb{N}
\text{,} \\
0, & \text{\thinspace \thinspace }j\notin \bigcup\limits_{k=1}^{\infty
}\left\{ M_{2\alpha _{k}},...,M_{2\alpha _{k}+1}-1\right\} .%
\end{array}%
\right.
\end{equation*}

We can write that%
\begin{equation*}
t_{M_{2\alpha _{k}}+M_{2s}}f=\frac{1}{Q_{M_{2\alpha _{k}}+M_{2s}}}%
\sum_{j=0}^{M_{2\alpha _{k}}}q_{M_{2\alpha _{k}}+M_{2s}-j}S_{j}f
\end{equation*}%
\begin{equation*}
+\frac{1}{Q_{M_{2\alpha _{k}}+M_{2s}}}\sum_{j=M_{2\alpha
_{k}}+1}^{M_{2\alpha _{k}}+M_{2s}}q_{M_{2\alpha _{k}}+M_{2s}-j}S_{j}f:=I+II.
\end{equation*}

According to (\ref{sn101}) we can conclude that%
\begin{equation}
\left\vert I\right\vert \leq \frac{1}{Q_{M_{2\alpha _{k}}+M_{2s}}}%
\sum_{j=0}^{M_{2\alpha _{k}}}q_{M_{2\alpha _{k}}+M_{2s}-j}\left\vert
S_{j}f\right\vert  \label{10anor}
\end{equation}%
\begin{equation*}
\leq \frac{2\lambda M_{2\alpha _{k-1}}^{\alpha +1}}{\alpha _{k-1}^{1/2}}%
\frac{1}{Q_{M_{2\alpha _{k}}+M_{2s}}}\sum_{j=0}^{M_{2\alpha
_{k}}}q_{M_{2\alpha _{k}}+M_{2s}-j}
\end{equation*}%
\begin{equation*}
\leq \frac{2\lambda M_{2\alpha _{k-1}}^{\alpha +1}}{\alpha _{k-1}^{1/2}}\leq
\frac{M_{\alpha _{k}}^{\alpha }}{16\alpha _{k}^{3/2}}.
\end{equation*}

Let $x\in I_{s}/I_{s+1}$ and $M_{2\alpha _{k}}+1\leq $ $j\leq M_{2\alpha
_{k}}+M_{2s}.$ In view of the second inequality of (\ref{1t6.1}) in the case
when $l=k$ and $p=1/\left( 1+\alpha \right)$ we can write that%
\begin{equation*}
S_{j}f=S_{M_{_{2\alpha _{k}}}}f+\frac{M_{2\alpha _{k}}^{\alpha }\psi
_{M_{_{2\alpha _{k}}}}D_{j-M_{_{2\alpha _{k}}}}}{\alpha _{k}^{1/2}}.
\end{equation*}

Hence, it yields that
\begin{equation*}
II=\frac{1}{Q_{M_{2\alpha _{k}}+M_{2s}}}\sum_{j=M_{2\alpha
_{k}}+1}^{M_{2\alpha _{k}}+M_{2s}}q_{M_{2\alpha _{k}}+M_{2s}-j}\frac{%
M_{2\alpha _{k}}^{\alpha }\psi _{M_{_{2\alpha _{k}}}}D_{j-M_{_{2\alpha
_{k}}}}}{\alpha _{k}^{1/2}}
\end{equation*}%
\begin{equation*}
+\frac{1}{Q_{M_{2\alpha _{k}}+M_{2s}}}\sum_{j=M_{2\alpha
_{k}}+1}^{M_{2\alpha _{k}}+M_{2s}}q_{M_{2\alpha
_{k}}+M_{2s}-j}S_{M_{_{2\alpha _{k}}}}f
\end{equation*}%
\begin{equation*}
:=II_{1}+II_{2}.
\end{equation*}

By using again (\ref{sn101}) we get that%
\begin{equation*}
\left\vert II_{2}\right\vert \leq \frac{2\lambda M_{2\alpha _{k-1}}^{\alpha
+1}}{\alpha _{k-1}^{1/2}}\frac{1}{Q_{M_{2\alpha _{k}}+M_{2s}}}%
\sum_{j=M_{2\alpha _{k}}+1}^{M_{2\alpha _{k}}+M_{2s}}q_{M_{2\alpha
_{k}}+M_{2s}-j}
\end{equation*}%
\begin{equation*}
\leq \frac{2\lambda M_{2\alpha _{k-1}}^{\alpha +1}}{\alpha _{k-1}^{1/2}}\leq
\frac{M_{\alpha _{k}}^{\alpha +1}}{16\alpha _{k}^{1/2}}.
\end{equation*}

Let $x\in I_{s}/I_{s+1},$ for $\left[ \alpha _{k}/2\right] <s\leq \alpha
_{k}.$ Then, according to (\ref{cond32.1}) we find that%
\begin{equation*}
\left\vert II_{1}\right\vert =\frac{1}{Q_{M_{2\alpha _{k}}+M_{2s}}}%
\left\vert \frac{\psi _{M_{2\alpha _{k}}}M_{2\alpha _{k}}^{\alpha }}{\alpha
_{k}^{1/2}}\sum_{j=1}^{M_{2s}}q_{M_{2s}-j}D_{j}\right\vert
\end{equation*}%
\begin{equation*}
=\frac{M_{2\alpha _{k}}^{\alpha }}{\alpha _{k}^{1/2}Q_{M_{2\alpha
_{k}}+M_{2s}}}\left\vert \sum_{j=1}^{M_{2s}}q_{M_{2s}-j}j\right\vert
\end{equation*}%
\begin{equation*}
\geq \frac{M_{2\alpha _{k}}^{\alpha }}{\alpha _{k}^{1/2}Q_{M_{2\alpha
_{k}+1}}}\left\vert \sum_{j=1}^{M_{2s}}q_{M_{2s}-j}j\right\vert
\end{equation*}%
\begin{equation*}
\geq \frac{c}{\lambda ^{\alpha }\alpha _{k}^{1/2}}\left\vert
\sum_{j=1}^{M_{2s}}q_{M_{2s}-j}j\right\vert .
\end{equation*}%
We invoke Abel transformation and apply (\ref{cond30}) to get that%
\begin{equation*}
\left\vert \sum_{j=1}^{M_{2s}}q_{M_{2s}-j}j\right\vert
\end{equation*}%
\begin{equation*}
=\left\vert \sum_{j=1}^{M_{2s}}\left( q_{M_{2s}-j}-q_{M_{2s}-j-1}\right)
\frac{j\left( j+1\right) }{2}\right\vert
\end{equation*}%
\begin{equation*}
\geq \frac{c_{\alpha }M_{2s}^{2}}{\alpha _{k}^{1/2}}\sum_{j=\left[ M_{2s}/2%
\right] }^{M_{2s}}\left\vert q_{M_{2s}-j}-q_{M_{2s}-j-1}\right\vert j^{2}
\end{equation*}%
\begin{equation*}
\geq c_{\alpha }M_{2s}^{2}\sum_{j=\left[ M_{2s}/2\right] }^{M_{2s}}\left%
\vert q_{M_{2s}-j}-q_{M_{2s}-j-1}\right\vert
\end{equation*}%
\begin{equation*}
\geq c_{\alpha }M_{2s}^{2}\sum_{j=1}^{\left[ M_{2s}/2\right] }\left\vert
q_{j}-q_{j+1}\right\vert
\end{equation*}%
\begin{equation*}
\geq c_{\alpha }M_{2s}^{2}\sum_{j=1}^{\left[ M_{2s}/2\right] }\frac{1}{%
j^{\alpha -2}}\geq c_{\alpha }M_{2s}^{\alpha -1}M_{2s}^{2}\geq c_{\alpha
}M_{2s}^{\alpha +1}.
\end{equation*}

Hence,%
\begin{equation*}
\left\vert II_{1}\right\vert \geq \frac{c_{\alpha }}{\alpha _{k}^{1/2}}%
\left\vert \sum_{j=1}^{M_{2s}}q_{M_{2s}-j}j\right\vert \geq \frac{c_{\alpha
}M_{2s}^{\alpha +1}}{\alpha _{k}^{1/2}}.
\end{equation*}

By now using the estimates above we obtain that
\begin{equation}
\int_{G_{m}}\left\vert t_{M_{2\alpha _{k}}+M_{2s}}f\right\vert ^{1/\left(
1+\alpha \right) }d\mu \geq \left\vert II_{1}\right\vert -\left\vert
II_{2}\right\vert -\left\vert I\right\vert  \label{100anor}
\end{equation}%
\begin{equation*}
\geq \frac{c_{\alpha }M_{2s}^{1+\alpha }}{\alpha _{k}^{1/2}}-\frac{4\lambda
M_{\alpha _{k}}^{\alpha +1}}{\alpha _{k}^{3/2}}\geq \frac{c_{\alpha
}M_{2s}^{1+\alpha }}{\alpha _{k}^{1/2}}.
\end{equation*}

By combining (\ref{10anor}) and (\ref{100anor}) we find that
\begin{equation*}
\int_{G_{m}}\left\vert t^{\ast }f\right\vert ^{1/\left( 1+\alpha \right)
}d\mu
\end{equation*}%
\begin{equation*}
\geq \sum_{s=\left[ \alpha _{k}/2\right] }^{\alpha
_{k}-1}\int_{I_{s}/I_{s+1}}\left\vert t_{M_{2\alpha
_{k}}+M_{2s}}f\right\vert ^{1/\left( 1+\alpha \right) }d\mu
\end{equation*}%
\begin{equation*}
\geq c_{\alpha }\sum_{s=\left[ \alpha _{k}/2\right] }^{\alpha _{k}-1}\frac{%
M_{2s}}{M_{2s}\alpha _{k}^{1/2\left( 1+\alpha \right) }}\geq c_{\alpha
}\sum_{s=\left[ \alpha _{k}/2\right] }^{\alpha _{k}-1}\frac{1}{\alpha
_{k}^{1/2\left( 1+\alpha \right) }}
\end{equation*}%
\begin{equation*}
\geq \frac{c_{\alpha }}{\alpha _{k}^{1/2\left( 1+\alpha \right) }}\sum_{s=%
\left[ \alpha _{k}/2\right] }^{\alpha _{k}-1}1
\end{equation*}%
\begin{equation*}
\geq \frac{c_{\alpha }\alpha _{k}}{\alpha _{k}^{1/2\left( 1+\alpha \right) }}%
\geq c_{\alpha }\alpha _{k}^{1/2}\rightarrow \infty ,\text{ when \ }%
k\rightarrow \infty .
\end{equation*}

\textbf{\ }The proof is complete.
\QED

Our next result reads:

\begin{theorem}
\label{Theorem1nor}Let $f\in H_{p},$ where $0$\thinspace $<p\,<1/\left( )
1+\alpha \right)$ for some $0<\alpha \leq 1,$ and $\{q_{k}:k\in \mathbb{N}%
\} $ be a sequence of non-increasing numbers satisfying conditions (\ref{6a}%
) and (\ref{7a}). Then the maximal operator%
\begin{equation*}
\overset{\sim }{t}_{p,\alpha }^{\ast }:=\frac{\left\vert t_{n}f\right\vert }{%
\left( n+1\right) ^{1/p-1-\alpha }}
\end{equation*}%
is bounded from the martingale Hardy space $H_{p}$ to the Lebesgue space $L_{p}.$

b) Let $\left\{ \Phi _{n}:n\in \mathbb{N}_{+}\right\} $ be any
non-decreasing sequence, satisfying the condition

\begin{equation}
\overline{\lim_{n\rightarrow \infty }}\frac{\left( n+1\right) ^{1/p-1-\alpha
}}{\Phi _{n}}=+\infty .  \label{6lbbb}
\end{equation}

Then there exists N\"orlund means with non-increasing sequence $\{q_{k}:k\in
\mathbb{N}\}$ satisfying the conditions (\ref{cond29}) and (\ref{cond30})  such
that
\begin{equation*}
\sup_{k\in \mathbb{N}}\frac{\left\Vert \frac{t_{M_{_{2n_{k}}}+1}f_{k}}{\Phi
_{M_{_{2n_{k}}}+1}}\right\Vert _{weal-L_{p}}}{\left\Vert f_{k}\right\Vert
_{H_{p}}}=\infty .
\end{equation*}

\end{theorem}

\begin{remark}
\label{remark1nor111}Part a) can be found in the paper Blahota and Tephnadze \cite%
{bt2}, while part b) has not been stated before for such a general case.
\end{remark}

{\bf Proof}:
Since the N\"orlund means $t_{n}$ are bounded from $L_{\infty }$ to $L_{\infty }$
(the boundedness follows from Corollary \ref{corollary3n}), according to
Lemma \ref{lemma2.2} it actually suffices to show that
\begin{equation*}
\int_{\overline{I_{N}}}\left\vert \overset{\sim }{t}_{p,\alpha }^{\ast
}a\right\vert ^{p}d\mu <c
\end{equation*}%
for some constant $c\ $and every $p$-atom $a$. We may assume that $a$ is an
arbitrary $p$-atom with support$\ I,$ $\mu \left( I\right) =M_{N}^{-1}$ and
$I=I_{N}.$ It is easy to see that $S_{n}\left( a\right) =t_{n}\left(
a\right) =0$ when $n\leq M_{N}.$ Therefore we can suppose that $n>M_{N}.$

Let $x\in I_{N}.$ Since $\left\Vert a\right\Vert _{\infty }\leq M_{N}^{1/p}$
we obtain that
\begin{equation*}
\left\vert t_{n}a\left( x\right) \right\vert \leq \int_{I_{N}}\left\vert
a\left( t\right) \right\vert \left\vert F_{n}\left( x-t\right) \right\vert
d\mu \left( t\right)
\end{equation*}%
\begin{equation*}
\leq \left\Vert a\right\Vert _{\infty }\int_{I_{N}}\left\vert F_{n}\left(
x-t\right) \right\vert d\mu \left( t\right) \leq
M_{N}^{1/p}\int_{I_{N}}\left\vert F_{n}\left( x-t\right) \right\vert d\mu
\left( t\right) .
\end{equation*}

Let $x\in I_{N}^{k,l},\,0\leq k<l<N.$ Then, from Lemma \ref{lemma3n} we get
that%
\begin{equation}
\left\vert t_{n}a\left( x\right) \right\vert \leq \frac{c_{\alpha
,p}M_{N}^{1/p-1}M_{l}^{\alpha }M_{k}}{n^{\alpha }}.  \label{12nor}
\end{equation}

Let $x\in I_{N}^{k,N},\,0\leq k<N.$ Then, according to Lemma \ref{lemma3n},
we have that
\begin{equation}
\left\vert t_{n}a\left( x\right) \right\vert \leq c_{\alpha
,p}M_{N}^{1/p-1}M_{k}.  \label{12anor}
\end{equation}

Let $x\in I_{N}^{k,l},\,0\leq k<l\leq N.$ Since $n>M_{N}$ we can conclude
that
\begin{equation}
\frac{\left\vert t_{n}a\left( x\right) \right\vert }{n^{1/p-1-\alpha }}\leq
c_{\alpha ,p}M_{l}^{\alpha }M_{k}.  \label{12anorM}
\end{equation}

The expression on the right-hand side of (\ref{12anorM}) does not depend on $n.$ Hence, we can conclude that%
\begin{equation}
\left\vert \overset{\sim }{t}_{p,\alpha }^{\ast }a\left( x\right)
\right\vert \leq c_{\alpha ,p}M_{l}^{\alpha }M_{k}  \label{12anorMA}
\end{equation}%
for $x\in I_{N}^{k,l},\,0\leq k<l\leq N.$

By combining (\ref{1.1}) and (\ref{12anorMA}) we obtain that%
\begin{equation*}
\int_{\overline{I_{N}}}\left\vert \overset{\sim }{t}_{p,\alpha }^{\ast
}a\right\vert ^{p}d\mu
\end{equation*}%
\begin{equation*}
=\overset{N-2}{\underset{k=0}{\sum }}\overset{N-1}{\underset{l=k+1}{\sum }}%
\sum\limits_{x_{j}=0,j\in \{l+1,\dots
,N-1\}}^{m_{j}-1}\int_{I_{N}^{k,l}}\left\vert \overset{\sim }{t}_{p,\alpha
}^{\ast }a\right\vert ^{p}d\mu
\end{equation*}%
\begin{equation*}
+\overset{N-1}{\underset{k=0}{\sum }}\int_{I_{N}^{k,N}}\left\vert \overset{%
\sim }{t}_{p,\alpha }^{\ast }a\right\vert ^{p}d\mu
\end{equation*}%
\begin{equation*}
\leq c_{\alpha ,p}\overset{N-2}{\underset{k=0}{\sum }}\overset{N-1}{\underset%
{l=k+1}{\sum }}\frac{1}{M_{l}}M_{l}^{\alpha p}M_{k}^{p}+c_{\alpha ,p}\overset%
{N-1}{\underset{k=0}{\sum }}\frac{1}{M_{N}}M_{N}^{1-p}M_{k}^{p}
\end{equation*}%
\begin{equation*}
\leq c_{\alpha ,p}\overset{N-2}{\underset{k=0}{\sum }}M_{k}^{p}\overset{N-1}{%
\underset{l=k+1}{\sum }}\frac{1}{M_{l}^{1-\alpha p}}+c_{\alpha ,p}\overset{%
N-1}{\underset{k=0}{\sum }}\frac{M_{k}^{p}}{M_{N}^{p}}\leq c_{\alpha
,p}<\infty .
\end{equation*}

The proof of part a) is complete.

Let $0<p<1/\left( 1+\alpha \right) .$ Under condition (\ref{6lbbb}) there
exists positive integers $n_{k}$ such that
\begin{equation*}
\lim_{k\rightarrow \infty }\frac{\left( M_{2n_{k}}+1\right) ^{1/p-1-\alpha }%
}{\Phi _{M_{2n_{k}}+1}}=\infty ,\text{ \ \ }0<p<1/\left( 1+\alpha \right) .
\end{equation*}%
To prove part b) we apply the $p$-atoms defined in Example \ref{example2.2}.
Under conditions (\ref{6a}) and (\ref{7a}), if we invoke (\ref{13l}) and (%
\ref{14l}) we find that
\begin{equation*}
\frac{\left\vert t_{M_{2n_{k}}+1}f_{k}\right\vert }{\Phi _{M_{2n_{k}}+1}}=%
\frac{\left\vert S_{M_{2n_{k}}+1}\right\vert }{Q_{M_{2n_{k}}+1}\Phi
_{M_{2n_{k}}+1}}
\end{equation*}%
\begin{equation*}
=\frac{q_{0}\left\vert D_{M_{2n_{k}}+1}-D_{M_{2n_{k}}}\right\vert }{%
Q_{M_{2n_{k}}+1}\Phi _{M_{2n_{k}}+1}}=\frac{q_{0}\left\vert \psi
_{M_{2n_{k}}}\right\vert }{Q_{M_{2n_{k}}+1}\Phi _{M_{2n_{k}}+1}}=\frac{%
c_{\alpha }}{M_{2n_{k}}^{\alpha }\Phi _{M_{2n_{k}}+1}}.
\end{equation*}

Hence,
\begin{equation}
\mu \left\{ x\in G_{m}:\frac{\left\vert t_{M_{2n_{k}}+1}f_{k}\left( x\right)
\right\vert }{\Phi _{M_{2n_{k}}+1}}\geq \frac{c_{\alpha }}{%
M_{2n_{k}}^{\alpha }\Phi _{M_{2n_{k}}+1}}\right\} =1.  \label{16bbb}
\end{equation}

By combining (\ref{15l}) and (\ref{16bbb}) we have that%
\begin{equation*}
\frac{\frac{c_{\alpha }}{M_{2n_{k}}^{\alpha }\Phi _{M_{2n_{k}}+1}}\left( \mu
\left\{ x\in G_{m}:\frac{\left\vert S_{M_{2n_{k}}+1}f_{k}\left( x\right)
\right\vert }{\Phi _{M_{2n_{k}}+1}}\geq \frac{q_{0}}{M_{2n_{k}}^{\alpha
}\Phi _{M_{2n_{k}}+1}}\right\} \right) ^{1/p}}{\left\Vert f_{k}\right\Vert
_{H_{p}}}
\end{equation*}%
\begin{equation*}
\geq \frac{c_{\alpha }M_{2n_{k}}^{1/p-1-\alpha }}{\Phi _{M_{2n_{k}}+1}}\geq
\frac{c_{\alpha }\left( M_{_{2n_{k}}}+1\right) ^{1/p-1-\alpha }}{\Phi
_{M_{2n_{k}}+1}}\rightarrow \infty ,\text{ when \ }k\rightarrow \infty .
\end{equation*}

The proof is complete.
\QED

\begin{corollary}
\label{corollary3fejermax1.5}  Let $0<p<1/\left( 1+\alpha \right) $ and $%
f\in H_{p}.$ Then there exists an absolute constant $c_{p,\alpha },$
depending only on $p$ and $\alpha $, such that
\begin{equation*}
\left\Vert t_{n}f\right\Vert _{p}\leq c_{p,\alpha }\left( n+1\right)
^{1/p-1-\alpha }\left\Vert f\right\Vert _{H_{p}},\text{ \ }n\in \mathbb{N}%
_{+}.
\end{equation*}
\end{corollary}

{\bf Proof}:
According to part a) of Theorem \ref{Theorem1nor} we conclude that%
\begin{equation*}
\left\Vert \frac{t_{n}f}{\left( n+1\right) ^{1/p-1-\alpha }}\right\Vert
_{p}\leq \left\Vert \sup_{n\in \mathbb{N}}\frac{\left\vert t_{n}f\right\vert
}{\left( n+1\right) ^{1/p-1-\alpha }}\right\Vert _{p}
\end{equation*}%
\begin{equation*}
\leq c_{p,\alpha }\left\Vert f\right\Vert _{H_{p}},\text{ \ }n\in \mathbb{N}%
_{+}.
\end{equation*}%
The proof is complete.
\QED

\begin{corollary}
\label{corollary3fejermax1.2.5}Let $\left\{ \Phi _{n}:n\in \mathbb{N}%
\right\} $ be any non-decreasing sequence satisfying the condition (\ref%
{6lbbb}). Then there exists a martingale $f\in H_{p}$ such that%
\begin{equation*}
\sup_{n\in \mathbb{N}}\left\Vert \frac{t_{n}f}{\Phi _{n}}\right\Vert
_{weak-L_{p}}=\infty.
\end{equation*}
\end{corollary}

\begin{corollary}
\label{corollary3fejermax1.3.5}Let $\left\{ \Phi _{n}:n\in \mathbb{N}%
\right\} $ be any non-decreasing sequence satisfying the condition (\ref%
{6lbbb}). Then the maximal operator%
\begin{equation*}
\sup_{n\in \mathbb{N}}\frac{\left\vert t_{n}f\right\vert }{\Phi _{n}}
\end{equation*}%
\textit{is not bounded from the Hardy space }$H_{p}$\textit{\ to the space }$%
weak-L_{p}.$
\end{corollary}

We now formulate our final result in this section.

\begin{theorem}
\label{Theorem1norin}Let $f\in H_{1/(1+\alpha )},$ where $0<\alpha \leq 1$
and $\{q_{k}:k\in \mathbb{N}\}$ be a sequence of non-increasing numbers
satisfying the conditions (\ref{6a}) and (\ref{7a}). Then there exists an absolute
constant $c_{\alpha }$ depending only on $\alpha$ such that the maximal
operator%
\begin{equation*}
\overset{\sim }{t}_{\alpha }^{\ast }:=\frac{\left\vert t_{n}f\right\vert }{%
\log ^{1+\alpha }\left( n+1\right) }
\end{equation*}%
is bounded from the martingale Hardy space $H_{1/(1+\alpha )}$ to the Lebesgue space $%
L_{1/(1+\alpha )}.$

b) Let $\left\{ \Phi _{n}:n\in \mathbb{N}_{+}\right\} $ be any
non-decreasing sequence satisfying the condition%
\begin{equation}
\overline{\lim_{n\rightarrow \infty }}\frac{\log ^{1+\alpha }\left(
n+1\right) }{\Phi _{n}}=+\infty .  \label{nom1}
\end{equation}%
Then there exists N\"orlund means with non-increasing sequence $\{q_{k}:k\in
\mathbb{N}\}$ satisfying the conditions (\ref{cond29}) and (\ref{cond30}) such
that
\begin{equation*}
\sup_{k\in \mathbb{N}}\frac{\left\Vert \sup_{n}\left\vert \frac{t_{n}f_{k}}{%
\Phi _{n}}\right\vert \right\Vert _{1/(1+\alpha )}}{\left\Vert f\right\Vert
_{H_{1/(1+\alpha )}}}=\infty .
\end{equation*}
\end{theorem}

\begin{remark}
Part a) of this result can be found in the paper Bhahota and  Tephnadze \cite%
{bpt1}. Part b) has not been stated before for such a general case.
\end{remark}

{\bf Proof}:
According to Lemma \ref{lemma2.2} the proof of part a) will be complete if
we show that%
\begin{equation*}
\int_{\overline{I_{N}}}\left\vert \overset{\sim }{t}_{\alpha }^{\ast
}a\right\vert ^{1/\left( 1+\alpha \right) }d\mu <\infty
\end{equation*}%
for every $1/\left( 1+\alpha \right) $-atom $a.$ We may assume that $a$ is
an arbitrary $1/\left( 1+\alpha \right) $-atom with support$\ I,$ $\mu
\left( I\right) =M_{N}^{-1}$ and $I=I_{N}.$ It is easy to see that $%
t_{m}\left( a\right) =0,$ when $m\leq M_{N}.$ Therefore we can suppose that $%
m>M_{N}.$

Let $x\in I_{N}.$ Since $t_{m}$ is bounded from $L_{\infty }$ to $L_{\infty
} $ (the boundedness follows from Corollary \ref{corollary3n}) and $%
\left\Vert a\right\Vert _{\infty }\leq M_{N}^{1/\left( 1+\alpha \right) }$
we obtain that
\begin{equation*}
\left\vert t_{m}a\left( x\right) \right\vert \leq \int_{I_{N}}\left\vert
a\left( t\right) \right\vert \left\vert F_{m}\left( x-t\right) \right\vert
d\mu \left( t\right)
\end{equation*}%
\begin{equation*}
\leq \left\Vert a\left( x\right) \right\Vert _{\infty
}\int_{I_{N}}\left\vert F_{m}\left( x-t\right) \right\vert d\mu \left(
t\right)
\end{equation*}%
\begin{equation*}
\leq M_{N}^{1+\alpha }\int_{I_{N}}\left\vert F_{m}\left( x-t\right)
\right\vert d\mu \left( t\right) .
\end{equation*}

Let $x\in I_{N}^{k,l},\,0\leq k<l<N.$ From Lemma \ref{lemma3n} we get that
\begin{equation}
\left\vert t_{m}a\left( x\right) \right\vert \leq \frac{c_{\alpha
}M_{k}M_{l}^{\alpha }M_{N}^{\alpha }}{m^{\alpha }}.  \label{saa}
\end{equation}

Let $x\in I_{N}^{k,N},\,0\leq k<N.$ In the view of Lemma \ref{lemma3n} we
have that%
\begin{equation}
\left\vert t_{m}a\left( x\right) \right\vert \leq c_{\alpha
}M_{k}M_{N}^{\alpha }.  \label{saaa}
\end{equation}

Let $x\in I_{N}^{k,l},\,0\leq k<l\leq N.$ Since $n>M_{N}$ we can conclude
that
\begin{equation}
\frac{\left\vert t_{n}a\left( x\right) \right\vert }{\log ^{1+\alpha }n}\leq
\frac{c_{\alpha }M_{k}M_{l}^{\alpha }}{N^{1+\alpha }}.  \label{saaaM}
\end{equation}

The expression on the right-hand side of (\ref{saaaM}) does not depend on $n.$ Hence,
\begin{equation}
\left\vert \overset{\sim }{t_{\alpha }}^{\ast }a\left( x\right) \right\vert
\leq \frac{c_{\alpha }M_{k}M_{N}^{\alpha }}{N^{1+\alpha }},  \label{saaaMA}
\end{equation}%
for $x\in I_{N}^{k,l},\,0\leq k<l\leq N.$

According to (\ref{1.1}) and (\ref{saaaMA}) we obtain that
\begin{equation*}
\int_{\overline{I_{N}}}\left\vert \overset{\sim }{t_{\alpha }}^{\ast
}a\right\vert ^{1/\left( 1+\alpha \right) }d\mu
\end{equation*}%
\begin{equation*}
=\overset{N-2}{\underset{k=0}{\sum }}\overset{N-1}{\underset{l=k+1}{\sum }}%
\sum\limits_{x_{j}=0,j\in \{l+1,\dots
,N-1\}}^{m_{j}-1}\int_{I_{N}^{k,l}}\left\vert \overset{\sim }{t_{\alpha }}%
^{\ast }a\right\vert ^{1/\left( 1+\alpha \right) }d\mu
\end{equation*}%
\begin{equation*}
+\overset{N-1}{\underset{k=0}{\sum }}\int_{I_{N}^{k,N}}\left\vert \overset{%
\sim }{t_{\alpha }}^{\ast }a\right\vert ^{1/\left( 1+\alpha \right) }d\mu
\end{equation*}%
\begin{equation*}
\leq \frac{c_{\alpha }}{N}\overset{N-2}{\underset{k=0}{\sum }}\overset{N-1}{%
\underset{l=k+1}{\sum }}\frac{1}{M_{l}}M_{l}^{\alpha /\left( 1+\alpha
\right) }M_{k}^{1/\left( 1+\alpha \right) }
\end{equation*}%
\begin{equation*}
+\frac{c_{\alpha }}{N}\overset{N-1}{\underset{k=0}{\sum }}\frac{1}{M_{N}}%
M_{N}^{\alpha /\left( 1+\alpha \right) }M_{k}^{1/\left( 1+\alpha \right) }
\end{equation*}%
\begin{equation*}
\leq \frac{c_{\alpha }}{N}\overset{N-2}{\underset{k=0}{\sum }}\overset{N-1}{%
\underset{l=k+1}{\sum }}\frac{M_{l}^{\alpha /\left( 1+\alpha \right)
}M_{k}^{1/\left( 1+\alpha \right) }}{M_{l}}
\end{equation*}%
\begin{equation*}
+\frac{c_{\alpha }}{N}\overset{N-1}{\underset{k=0}{\sum }}\frac{%
M_{k}^{1/\left( 1+\alpha \right) }}{M_{N}^{1/\left( 1+\alpha \right) }}\leq
c_{\alpha }<\infty .
\end{equation*}

The proof of part a) is complete.

Under condition (\ref{nom1}) there exists a positive integers $m_{k}^{,}$
such that $M_{2m_{k}^{,}+1}\leq \lambda _{k}<2M_{2m_{k}^{,}+1}.$ Since $\Phi
_{n}$ is a non-decreasing function we have that
\begin{equation*}
\overline{\underset{k\rightarrow \infty }{\lim }}\frac{\left(
m_{k}^{^{\prime }}\right) ^{1+\alpha }}{\Phi _{M_{2m_{k}^{,}+1}}}\geq
c\lim_{k\rightarrow \infty }\frac{\log ^{1+\alpha }\left( \lambda
_{k}\right) }{\Phi _{M_{2\lambda _{k}+1}}}=\infty .
\end{equation*}

Let $\ \left\{ n_{k}:k\in \mathbb{N}_{+}\right\} \subset \left\{
m_{k}^{\prime }:k\in \mathbb{N}_{+}\right\} $ be a sequence of positive
numbers such that
\begin{equation*}
\lim_{k\rightarrow \infty }\frac{n_{k}^{1+\alpha }}{\Phi _{M_{2n_{k}+1}}}%
=\infty.
\end{equation*}

To prove part b) we use the $1/\left( 1+\alpha \right) $-atoms defined in
Example \ref{example2.2}. If we apply (\ref{dn21}) in Lemma \ref{dn1} with (%
\ref{13l}) and (\ref{14l}) and invoke Abel transformation we get that%
\begin{equation*}
\frac{\left\vert t_{M_{2n_{k}}+M_{2s}}f_{k}\right\vert }{\Phi
_{M_{2n_{k}}+M_{2s}}}=\frac{1}{\Phi _{M_{2n_{k}}+M_{2s}}Q_{M_{2n_{k}}+M_{2s}}%
}\left\vert
\sum_{j=M_{2n_{k}}+1}^{M_{2n_{k}}+M_{2s}}q_{M_{2n_{k}}+M_{2s}-j}\left(
D_{j}-D_{M_{2n_{k}}}\right) \right\vert
\end{equation*}%
\begin{equation*}
=\frac{1}{\Phi _{M_{2n_{k}}+M_{2s}}Q_{M_{2n_{k}}+M_{2s}}}\left\vert
\sum_{j=1}^{M_{2s}}q_{M_{2s}-j}\left( D_{j+M_{2n_{k}}}-D_{M_{2n_{k}}}\right)
\right\vert
\end{equation*}%
\begin{equation*}
=\frac{1}{\Phi _{M_{2n_{k}}+M_{2s}}Q_{M_{2n_{k}}+M_{2s}}}\left\vert \psi
_{M_{2n_{k}}}\sum_{j=1}^{M_{2s}}q_{M_{2s}-j}D_{j}\right\vert
\end{equation*}%
\begin{equation*}
=\frac{1}{\Phi _{M_{2n_{k}}+M_{2s}}Q_{M_{2n_{k}}+M_{2s}}}\left\vert
\sum_{j=1}^{M_{2s}}q_{M_{2s}-j}j\right\vert .
\end{equation*}

Let $x\in I_{2s}/I_{2s+1},$ $s=\left[ n_{k}/2\right] ,...,n_{k}.$ If we
again use abel transformation, then under the conditions (\ref{cond29}) and (%
\ref{cond30}) we find that%
\begin{equation*}
\frac{\left\vert t_{M_{2n_{k}}+M_{2s}}f_{k}\right\vert }{\Phi
_{M_{2n_{k}}+M_{2s}}}\geq \frac{c}{\Phi
_{M_{2n_{k}}+M_{2s}}M_{2n_{k}}^{\alpha }}\left\vert
\sum_{j=1}^{M_{2s}}\left( q_{M_{2s}-j}-q_{M_{2s}-j-1}\right) j^{2}\right\vert
\end{equation*}%
\begin{equation*}
\geq \frac{c}{\Phi _{M_{2n_{k}}+M_{2s}}M_{2n_{k}}^{\alpha }}\sum_{j=\left[
M_{2s}\right] /2}^{M_{2s}}\left\vert q_{M_{2s}-j}-q_{M_{2s}-j-1}\right\vert
j^{2}
\end{equation*}%
\begin{equation*}
\geq \frac{cM_{2s}^{2}}{\Phi _{M_{2n_{k}}+M_{2s}}M_{2n_{k}}^{\alpha }}%
\sum_{j=\left[ M_{2s}\right] /2}^{M_{2s}}\left\vert
q_{M_{2s}-j}-q_{M_{2s}-j-1}\right\vert
\end{equation*}%
\begin{equation*}
\geq \frac{cM_{2s}^{2}}{\Phi _{M_{2n_{k}}+M_{2s}}M_{2n_{k}}^{\alpha }}%
\sum_{j=1}^{\left[ M_{2s}\right] /2}\left\vert q_{j}-q_{j+1}\right\vert
\end{equation*}%
\begin{equation*}
\geq \frac{cM_{2s}^{2}}{\Phi _{M_{2n_{k}}+M_{2s}}M_{2n_{k}}^{\alpha }}%
\sum_{j=1}^{\left[ M_{2s}\right] /2}j^{\alpha -2}\geq \frac{%
cM_{2s}^{2}M_{2s}^{\alpha -1}}{\Phi _{M_{2n_{k}}+M_{2s}}M_{2n_{k}}^{\alpha }}
\end{equation*}%
\begin{equation*}
\geq \frac{cM_{2s}^{\alpha +1}}{\Phi _{M_{2n_{k}}+M_{2s}}M_{2n_{k}}^{\alpha }%
}\geq \frac{cM_{2s}^{\alpha +1}}{\Phi _{M_{2n_{k}+1}}M_{2n_{k}}^{\alpha }}.
\end{equation*}

Hence,
\begin{equation*}
\int_{G_{m}}\left( \sup_{n}\left\vert \frac{t_{n}f_{k}}{\Phi _{n}}%
\right\vert \right) ^{1/\left( 1+\alpha \right) }d\mu
\end{equation*}%
\begin{equation*}
\geq \text{ }\overset{n_{k}-1}{\underset{s=\left[ n_{k}/2\right] }{\sum }}%
\int_{I_{2s}/I_{2s+1}}\left\vert \frac{t_{M_{2n_{k}}+M_{2s}}f_{k}}{\Phi
_{M_{2n_{k}}+M_{2s}}}\right\vert ^{1/\left( 1+\alpha \right) }d\mu
\end{equation*}%
\begin{equation*}
\geq \frac{c}{\Phi _{M_{2n_{k}+1}}^{1/\left( 1+\alpha \right) }}\text{ }%
\overset{n_{k}-1}{\underset{s=\left[ n_{k}/2\right] }{\sum }}%
\int_{I_{2s}/I_{2s+1}}\frac{M_{2s}}{M_{2n_{k}}^{\alpha /\left( 1+\alpha
\right) }}d\mu
\end{equation*}%
\begin{equation*}
\geq \frac{c}{\Phi _{M_{2n_{k}+1}}^{1/\left( 1+\alpha \right) }}\text{ }%
\overset{n_{k}-1}{\underset{s=\left[ n_{k}/2\right] }{\sum }}\frac{%
M_{2s}^{\alpha +1}}{M_{2n_{k}}^{\alpha /\left( 1+\alpha \right) }}\frac{1}{%
M_{2s}^{\alpha +1}}
\end{equation*}%
\begin{equation*}
\geq \frac{c}{\Phi _{M_{2n_{k}+1}}^{1/\left( 1+\alpha \right) }}\text{ }%
\overset{n_{k}-1}{\underset{s=\left[ n_{k}/2\right] }{\sum }}\frac{1}{%
M_{2n_{k}}^{\alpha /\left( 1+\alpha \right) }}\geq \frac{cn_{k}}{%
M_{2n_{k}}^{\alpha /\left( 1+\alpha \right) }\Phi _{M_{2n_{k}+1}}^{1/\left(
1+\alpha \right) }}.
\end{equation*}%
Therefore, by also using (\ref{15l}) for $p=1/\left( 1+\alpha \right) $ we
have that
\begin{equation*}
\frac{\left( \int_{G_{m}}\left( \sup_{n}\left\vert \frac{t_{n}f_{k}}{\Phi
_{n}}\right\vert \right) ^{1/\left( 1+\alpha \right) }d\mu \right)
^{1+\alpha }}{\left\Vert f_{k}\right\Vert _{H_{1/\left( 1+\alpha \right) }}}
\end{equation*}%
\begin{equation*}
\geq \frac{cn_{k}^{1+\alpha }}{M_{_{2n_{k}}}^{\alpha }\Phi _{M_{2n_{k}+1}}}%
M_{_{2n_{k}}}^{\alpha }\geq \frac{cn_{k}^{1+\alpha }}{\Phi _{M_{2n_{k}+1}}}%
\rightarrow \infty ,\text{when \qquad }k\rightarrow \infty .
\end{equation*}

The proof is complete.
\QED

\begin{corollary}
Let $\left\{ \Phi _{n}:n\in \mathbb{N}_{+}\right\} $ be any non-decreasing
sequence satisfying the condition (\ref{nom1}). Then the following maximal operator%
\begin{equation*}
\sup_{n}\left\vert \frac{t_{n}f}{\Phi _{n}}\right\vert
\end{equation*}%
is not bounded from the Hardy space $H_{1/\left( 1+\alpha \right) }$ to the
space $L_{1/\left( 1+\alpha \right) }.$
\end{corollary}

\subsection{Strong convergence of N\"orlund means on martingale Hardy
spaces}

The first result in this section is due to Persson, Tephnadze and Wall
\cite{ptw2}.

\begin{theorem}
\label{theorem2fejerstronga} Let $0<p<1/2$, $f\in H_{p}$ and the sequence $%
\left\{ q_{k}:k\in \mathbb{N}\right\} $ be non-decreasing$.$ Then there
exists an absolute constant $c_{p}$ depending only on $p$ such that
\end{theorem}

\begin{equation*}
\overset{\infty }{\underset{k=1}{\sum }}\frac{\left\Vert t_{k}f\right\Vert
_{p}^{p}}{k^{2-2p}}\leq c_{p}\left\Vert f\right\Vert _{H_{p}}^{p}.
\end{equation*}

\begin{remark}
\label{Remark4nor}Since the Fej\'er means are examples of N\"orlund means with
non-decreasing sequence $\left\{ q_{k}:k\in \mathbb{N}\right\} $ we
immediately obtain from part b) of Theorem \ref{theorem2fejerstrong} that
the asymptotic behaviour of the sequence of weights
\begin{equation*}
\left\{ 1/k^{2-2p}:k\in \mathbb{N}\right\}
\end{equation*}
in N\"orlund means in Theorem \ref{theorem2fejerstronga} can not be improved.
\end{remark}

{\bf Proof}:
The proof is similar to that of part a) of Theorem \ref{theorem2fejerstrong}, but since this case is more general we give the details.

By Lemma \ref{lemma2.1} the proof is complete if we show that%
\begin{equation}
\overset{\infty }{\underset{m=1}{\sum }}\frac{\left\Vert t_{m}a\right\Vert
_{p}^{p}}{m^{2-2p}}\leq c_{p}  \label{14ctn}
\end{equation}%
for every $p$-atom $a$ with support$\ I$, $\mu \left( I\right) =M_{N}^{-1}.$
We may assume that $I=I_{N}.$ It is easy to see that $S_{n}\left( a\right)
=t_{n}\left( a\right) =0,$ when $n\leq M_{N}$. Therefore, we can suppose
that $n>M_{N}$.

Let $x\in I_{N}.$ Since N\"orlund means $t_{n}$ with non-decreasing sequence $%
\left\{ q_{k}:k\in \mathbb{N}\right\} $ are bounded from $L_{\infty }$ to $%
L_{\infty }$ (the boundedness follows from Corollary \ref{Corollary3nn}) and
$\left\Vert a\right\Vert _{\infty }\leq M_{N}^{1/p}$ we obtain that
\begin{equation*}
\int_{I_{N}}\left\vert t_{m}a\right\vert ^{p}d\mu
\end{equation*}%
\begin{equation*}
\leq \frac{\left\Vert a\right\Vert _{\infty }^{p}}{M_{N}}\leq 1,\text{ \ }%
0<p\leq 1/2.
\end{equation*}%
Hence,
\begin{equation}
\overset{\infty }{\underset{m=1}{\sum }}\frac{\int_{I_{N}}\left\vert
t_{m}a\right\vert ^{p}d\mu }{m^{2-2p}}  \label{14btn}
\end{equation}%
\begin{equation*}
\leq \overset{\infty }{\underset{k=1}{\sum }}\frac{1}{m^{2-2p}}\leq c<\infty
.
\end{equation*}

It is easy to see that

\begin{equation*}
\left\vert t_{n}a\left( x\right) \right\vert
\end{equation*}%
\begin{equation*}
=\left\vert \int_{I_{N}}a\left( t\right) \frac{1}{Q_{n}}\overset{n}{\underset%
{k=M_{N}}{\sum }}q_{n-k}D_{k}\left( x-t\right) d\mu \left( t\right)
\right\vert
\end{equation*}%
\begin{equation*}
\leq \int_{I_{N}}\left\vert a\left( t\right) \right\vert \left\vert \frac{1}{%
Q_{n}}\overset{n}{\underset{k=M_{N}}{\sum }}q_{n-k}D_{k}\left( x-t\right)
\right\vert d\mu \left( t\right)
\end{equation*}%
\begin{equation*}
\leq \left\Vert a\right\Vert _{\infty }\int_{I_{N}}\left\vert \frac{1}{Q_{n}}%
\overset{n}{\underset{k=M_{N}}{\sum }}q_{n-k}D_{k}\left( x-t\right)
\right\vert d\mu \left( t\right)
\end{equation*}%
\begin{equation*}
\leq M_{N}^{1/p}\int_{I_{N}}\left\vert \frac{1}{Q_{n}}\overset{n}{\underset{%
k=M_{N}}{\sum }}q_{n-k}D_{k}\left( x-t\right) \right\vert d\mu \left(
t\right) .
\end{equation*}

Let $x\in I_{N}^{k,l},\,0\leq k<l\leq N.$ Then, in the view of Lemma \ref%
{lemma5aaa}, we get that
\begin{equation}
\left\vert t_{m}a\left( x\right) \right\vert \leq
c_{p}M_{l}M_{k}M_{N}^{1/p-2},\text{ for }0<p<1/2.  \label{12q1tn}
\end{equation}

According to (\ref{1.1}) with (\ref{12q1tn}) we find that%
\begin{equation}
\int_{\overline{I_{N}}}\left\vert t_{m}a\right\vert ^{p}d\mu  \label{7aaatn}
\end{equation}%
\begin{equation*}
=\overset{N-2}{\underset{k=0}{\sum }}\overset{N-1}{\underset{l=k+1}{\sum }}%
\sum\limits_{x_{j}=0,\text{ }j\in \{l+1,\dots
,N-1\}}^{m_{j-1}}\int_{I_{N}^{k,l}}\left\vert t_{m}a\right\vert ^{p}d\mu
\end{equation*}%
\begin{equation*}
+\overset{N-1}{\underset{k=0}{\sum }}\int_{I_{N}^{k,N}}\left\vert
t_{m}a\right\vert ^{p}d\mu
\end{equation*}%
\begin{equation}
\leq c_{p}\overset{N-2}{\underset{k=0}{\sum }}\overset{N-1}{\underset{l=k+1}{%
\sum }}\frac{m_{l+1}\dotsm m_{N-1}}{M_{N}}M_{l}^{p}M_{k}^{p}M_{N}^{1-2p}
\notag
\end{equation}%
\begin{equation*}
+c_{p}\overset{N-1}{\underset{k=0}{\sum }}\frac{1}{M_{N}}M_{k}^{p}M_{N}^{1-p}
\end{equation*}%
\begin{equation*}
\leq c_{p}M_{N}^{1-2p}\overset{N-2}{\underset{k=0}{\sum }}\overset{N-1}{%
\underset{l=k+1}{\sum }}\frac{M_{l}^{p}M_{k}^{p}}{M_{l}}
\end{equation*}%
\begin{equation*}
+c_{p}\overset{N-1}{\underset{k=0}{\sum }}\frac{M_{k}^{p}}{M_{N}^{p}}\leq
cM_{N}^{1-2p}.
\end{equation*}

Moreover, according to (\ref{7aaatn}), we get that%
\begin{equation*}
\overset{\infty }{\underset{m=M_{N}+1}{\sum }}\frac{\int_{\overline{I_{N}}%
}\left\vert t_{m}a\right\vert ^{p}d\mu }{m^{2-2p}}
\end{equation*}%
\begin{equation*}
\leq c_{p}\overset{\infty }{\underset{m=M_{N}+1}{\sum }}\frac{M_{N}^{1-2p}}{%
m^{2-2p}}<c<\infty ,\text{ \ }\left( 0<p<1/2\right) .
\end{equation*}%
Now, by combining this estimate with (\ref{14btn}) we obtain (\ref{14ctn})
so the proof is complete.
\QED

Also the next theorem is proved in Persson, Tephnadze and Wall \cite{ptw2}.

\begin{theorem}
\label{theorem1fejerstronga}Let $f\in H_{1/2}$ and the sequence $\left\{
q_{k}:k\in \mathbb{N}\right\} $ be non-decreasing satisfying condition (\ref%
{fn01}). Then there exists an absolute constant $c,$ such that
\begin{equation*}
\frac{1}{\log n}\overset{n}{\underset{k=1}{\sum }}\frac{\left\Vert
t_{k}f\right\Vert _{1/2}^{1/2}}{k}\leq c\left\Vert f\right\Vert
_{H_{1/2}}^{1/2}.
\end{equation*}
\end{theorem}

{\bf Proof}:
The proof is similar to that of part a) of Theorem \ref{theorem1fejerstrong}, but since this case is more general we present the details.

According to Lemma \ref{lemma2.1} it suffices to show that%
\begin{equation*}
\frac{1}{\log n}\overset{n}{\underset{m=1}{\sum }}\frac{\left\Vert
t_{m}a\right\Vert _{1/2}^{1/2}}{m}\leq c
\end{equation*}%
for every $p$-atom $a$ with support$\ I$, $\mu \left( I\right) =M_{N}^{-1}.$
We may assume that $I=I_{N}.$ It is easy to see that $S_{n}\left( a\right)
=t_{n}\left( a\right) =0$ when $n\leq M_{N}$. Therefore we can suppose
that $n>M_{N}$.

Let $x\in I_{N}.$ Since $t_{n}$ is bounded from $L_{\infty }$ to $L_{\infty
} $ (the boundedness follows from Corollary \ref{Corollary3nn}) and $%
\left\Vert a\right\Vert _{\infty }\leq M_{N}^{2}$ we obtain that
\begin{equation*}
\int_{I_{N}}\left\vert t_{m}a\right\vert ^{1/2}d\mu \leq \frac{\left\Vert
a\right\Vert _{\infty }^{1/2}}{M_{N}}\leq c<\infty.
\end{equation*}%
Hence,
\begin{equation}
\frac{1}{\log n}\overset{n}{\underset{m=1}{\sum }}\frac{\int_{I_{N}}\left%
\vert t_{m}a\right\vert ^{1/2}d\mu }{m}  \label{14b}
\end{equation}%
\begin{equation*}
\leq \frac{1}{\log n}\overset{n}{\underset{m=1}{\sum }}\frac{1}{m}\leq
c<\infty ,\text{ \ }n\in \mathbb{N}\text{.}
\end{equation*}

It is easy to see that

\begin{equation*}
\left\vert t_{m}a\left( x\right) \right\vert \leq \int_{I_{N}}\left\vert
a\left( t\right) \right\vert \left\vert F_{m}\left( x-t\right) \right\vert
d\mu \left( t\right)
\end{equation*}%
\begin{equation*}
\leq \left\Vert a\right\Vert _{\infty }\int_{I_{N}}\left\vert F_{m}\left(
x-t\right) \right\vert d\mu \left( t\right)
\end{equation*}%
\begin{equation*}
\leq M_{N}^{2}\int_{I_{N}}\left\vert F_{m}\left( x-t\right) \right\vert d\mu
\left( t\right) .
\end{equation*}

Let $x\in I_{N}^{k,l},\,0\leq k<l<N.$ Then, in the view of Lemma \ref%
{lemma5a}, we get that
\begin{equation}
\left\vert t_{m}a\left( x\right) \right\vert \leq \frac{cM_{l}M_{k}M_{N}}{m}.
\label{12qtn}
\end{equation}

Let $x\in I_{N}^{k,N}.$ Then, according to Lemma \ref{lemma5a}, we find that
\begin{equation}
\left\vert t_{m}a\left( x\right) \right\vert \leq cM_{k}M_{N}.  \label{12q2}
\end{equation}

By combining (\ref{12qtn}) and (\ref{12q2}) with (\ref{1.1})\ we can conclude
that%
\begin{equation*}
\int_{\overline{I_{N}}}\left\vert t_{m}a\left( x\right) \right\vert
^{1/2}d\mu \left( x\right)
\end{equation*}%
\begin{equation}
\leq c\overset{N-2}{\underset{k=0}{\sum }}\overset{N-1}{\underset{l=k+1}{%
\sum }}\frac{m_{l+1}\dotsm m_{N-1}}{M_{N}}\frac{%
M_{l}^{1/2}M_{k}^{1/2}M_{N}^{1/2}}{m^{1/2}}  \notag
\end{equation}%
\begin{equation*}
+c\overset{N-1}{\underset{k=0}{\sum }}\frac{1}{M_{N}}M_{k}^{1/2}M_{N}^{1/2}
\end{equation*}%
\begin{equation*}
\leq cM_{N}^{1/2}\overset{N-2}{\underset{k=0}{\sum }}\overset{N-1}{\underset{%
l=k+1}{\sum }}\frac{M_{l}^{1/2}M_{k}^{1/2}}{m^{1/2}M_{l}}+c\overset{N-1}{%
\underset{k=0}{\sum }}\frac{M_{k}^{1/2}}{M_{N}^{1/2}}
\end{equation*}%
\begin{equation*}
\leq \frac{cM_{N}^{1/2}N}{m^{1/2}}+c.
\end{equation*}

It follows that%
\begin{equation}
\frac{1}{\log n}\overset{n}{\underset{m=M_{N}+1}{\sum }}\frac{\int_{%
\overline{I_{N}}}\left\vert t_{m}a\left( x\right) \right\vert ^{1/2}d\mu
\left( x\right) }{m}  \label{15b}
\end{equation}%
\begin{equation*}
\leq \frac{1}{\log n}\overset{n}{\underset{m=M_{N}+1}{\sum }}\left( \frac{%
cM_{N}^{1/2}N}{m^{3/2}}+\frac{c}{m}\right) <c<\infty.
\end{equation*}

The proof is complete by just combining (\ref{14b}) and (\ref{15b}).
\QED

Next, we investigate N\"orlund means with non-increasing sequence $%
\{q_{k}:k\in \mathbb{N}\}.$ At first we consider the case $0$\thinspace $%
<p\,<1/\left( 1+\alpha \right)$ where $0<\alpha <1.$ For details see the
paper of Blahota and Tephnadze \cite{bt2}.

\begin{theorem}
\label{Theorem2nor} Let $f\in H_{p},$ where $0$\thinspace $<p\,<1/\left(
1+\alpha \right) ,$ $0<\alpha \leq 1$ and $\{q_{k}:k\in \mathbb{N}\},$ be a
sequence of non-increasing numbers satisfying the conditions (\ref{6a}) and (%
\ref{7a}). Then there exists an absolute constant $c_{\alpha ,p},$ depending
only on $\alpha $ and $p$ such that
\begin{equation*}
\overset{\infty }{\underset{k=1}{\sum }}\frac{\left\Vert t_{k}f\right\Vert
_{H_{p}}^{p}}{k^{2-\left( 1+\alpha \right) p}}\leq c_{\alpha ,p}\left\Vert
f\right\Vert _{H_{p}}^{p}.
\end{equation*}
\end{theorem}

{\bf Proof}:
By Lemma \ref{lemma2.1} the it suffices to show that%
\begin{equation*}
\overset{\infty }{\underset{m=1}{\sum }}\frac{\left\Vert t_{m}a\right\Vert
_{p}^{p}}{m^{2-\left( 1+\alpha \right) p}}\leq c_{\alpha ,p}<\infty
\end{equation*}%
for every $p$-atom $a.$ Analogously to the proofs of the previous theorems we
may assume that $a$ be an arbitrary $p$-atom with support $I,\ \mu \left(
I\right) =M_{N}^{-1}$ and $I=I_{N}$ and $m>M_{N}.$

Let $x\in I_{N}.$ Since $t_{m}$ is bounded from $L_{\infty }$ to $L_{\infty
} $ (the boundedness follows from Corollary \ref{corollary3n}) and $%
\left\Vert a\right\Vert _{\infty }\leq M_{N}^{1/p}$ we obtain that
\begin{equation*}
\int_{I_{N}}\left\vert t_{m}a\right\vert ^{p}d\mu \leq \left\Vert
a\right\Vert _{\infty }^{p}M_{N}^{-1}\leq 1.
\end{equation*}%
Hence
\begin{equation*}
\overset{\infty }{\underset{m=M_{N}}{\sum }}\frac{\int_{I_{N}}\left\vert
t_{m}a\right\vert ^{1/\left( 1+\alpha \right) }d\mu }{m^{2-\left( 1+\alpha
\right) p}}
\end{equation*}%
\begin{equation*}
\leq \overset{\infty }{\underset{m=1}{\sum }}\frac{1}{m^{2-\left( 1+\alpha
\right) p}}\leq c_{\alpha ,p}<\infty .
\end{equation*}

According to (\ref{1.1}) and (\ref{12nor})-(\ref{12anor}) we can conclude
that
\begin{equation*}
\overset{\infty }{\underset{m=M_{N}+1}{\sum }}\frac{\int_{\overline{I_{N}}%
}\left\vert t_{m}a\right\vert ^{p}d\mu }{m^{2-\left( 1+\alpha \right) p}}
\end{equation*}%
\begin{equation*}
=\overset{N-2}{\underset{k=0}{\sum }}\overset{N-1}{\underset{l=k+1}{\sum }}%
\sum\limits_{x_{j}=0,j\in \{l+1,\dots ,N-1\}}^{m_{j}-1}\frac{%
\int_{I_{N}^{k,l}}\left\vert t_{m}a\right\vert ^{p}d\mu }{m^{2-\left(
1+\alpha \right) p}}
\end{equation*}%
\begin{equation*}
+\overset{n}{\underset{m=M_{N}+1}{\sum }}\overset{N-1}{\underset{k=0}{\sum }}%
\frac{\int_{I_{N}^{k,N}}\left\vert t_{m}a\right\vert ^{p}d\mu }{m^{2-\left(
1+\alpha \right) p}}
\end{equation*}%
\begin{equation*}
\leq c_{\alpha ,p}\overset{\infty }{\underset{m=M_{N}+1}{\sum }}\left( \frac{%
M_{N}^{1-p}}{m^{2-p}}\overset{N-2}{\underset{k=0}{\sum }}\overset{N-1}{%
\underset{l=k+1}{\sum }}\frac{M_{l}^{p\alpha }M_{k}^{p}}{M_{l}}+\frac{%
M_{N}^{1-p}}{m^{2-\left( 1+\alpha \right) p}}\overset{N-1}{\underset{k=0}{%
\sum }}\frac{M_{k}^{p}}{M_{N}}\right) .
\end{equation*}

Since%
\begin{equation*}
\overset{N-2}{\underset{k=0}{\sum }}\overset{N-1}{\underset{l=k+1}{\sum }}%
\frac{M_{l}^{p\alpha }M_{k}^{p}}{M_{l}}\leq \overset{N-2}{\underset{k=0}{%
\sum }}M_{k}^{p}\overset{N-1}{\underset{l=k+1}{\sum }}\frac{1}{%
M_{l}^{1-p\alpha }}
\end{equation*}%
\begin{equation*}
\leq \overset{N-2}{\underset{k=0}{\sum }}\frac{1}{M_{k}^{1-p\alpha }}%
M_{k}^{p}\leq \overset{N-2}{\underset{k=0}{\sum }}\frac{1}{M_{k}^{1-p\left(
\alpha +1\right) }}<c<\infty
\end{equation*}%
and%
\begin{equation*}
\overset{N-1}{\underset{k=0}{\sum }}\frac{M_{k}^{p}}{M_{N}}\leq
M_{N}^{p-1}<\infty
\end{equation*}%
we obtain that%
\begin{equation*}
\overset{\infty }{\underset{m=M_{N}+1}{\sum }}\frac{\int_{\overline{I_{N}}%
}\left\vert t_{m}a\right\vert ^{p}d\mu }{m^{2-\left( 1+\alpha \right) p}}
\end{equation*}%
\begin{equation*}
<c_{\alpha ,p}M_{N}^{1-p}\overset{\infty }{\underset{m=M_{N}+1}{\sum }}\frac{%
1}{m^{2-p}}+c_{\alpha ,p}\overset{\infty }{\underset{m=M_{N}+1}{\sum }}\frac{%
1}{m^{2-\left( 1+\alpha \right) p}}\leq c_{\alpha ,p}<\infty .
\end{equation*}%
The proof is complete.
\QED

The final result in this section is due to Blahota, Persson and
Tephnadze \cite{bpt1}.

\begin{theorem}
\label{Theorem3nor}Let $f\in H_{1/(1+\alpha )}$ where $0<\alpha \leq 1$ and
$\{q_{k}:k\in \mathbb{N}\}$ be a sequence of non-increasing numbers
satisfying the conditions (\ref{6a}) and (\ref{7a}). Then there exists an
absolute constant $c_{\alpha }$ depending only on $\alpha$ such that
\begin{equation*}
\frac{1}{\log n}\overset{n}{\underset{m=1}{\sum }}\frac{\left\Vert
t_{m}f\right\Vert _{H_{1/\left( 1+\alpha \right) }}^{1/\left( 1+\alpha
\right) }}{m}\leq c_{\alpha }\left\Vert f\right\Vert _{H_{1/\left( 1+\alpha
\right) }}^{1/\left( 1+\alpha \right)}.
\end{equation*}
\end{theorem}

{\bf Proof}:
\textbf{\ }By Lemma \ref{lemma2.1} the proof is complete if we show
that%
\begin{equation*}
\frac{1}{\log n}\overset{n}{\underset{m=1}{\sum }}\frac{\left\Vert
t_{m}a\right\Vert _{1/\left( 1+\alpha \right) }^{1/\left( 1+\alpha \right) }%
}{m}\leq c_{\alpha }<\infty
\end{equation*}%
for every $1/\left( 1+\alpha \right) $-atom $a.$ Analogously to the proofs
of previous resalts we may assume that $a$ is an arbitrary $1/\left(
1+\alpha \right) $-atom with support$\ I,$ $\mu \left( I\right) =M_{N}^{-1}$
and $I=I_{N}$ and $m>M_{N}.$

Let $x\in I_{N}.$ Since $t_{m}$ is bounded from $L_{\infty }$ to $L_{\infty
} $ (the boundedness follows from Corollary \ref{corollary3n}) and $%
\left\Vert a\right\Vert _{\infty }\leq M_{N}^{1+\alpha }$ we obtain that
\begin{equation*}
\int_{I_{N}}\left\vert t_{m}a\right\vert ^{1/\left( 1+\alpha \right) }d\mu
\end{equation*}%
\begin{equation*}
\leq \left\Vert a\right\Vert _{\infty }^{1/\left( 1+\alpha \right)
}/M_{N}\leq 1.
\end{equation*}%
Hence
\begin{equation*}
\frac{1}{\log n}\overset{n}{\underset{m=M_{N}}{\sum }}\frac{%
\int_{I_{N}}\left\vert t_{m}a\left( x\right) \right\vert ^{1/\left( 1+\alpha
\right) }d\mu }{m}
\end{equation*}%
\begin{equation*}
\leq \frac{1}{\log n}\overset{n}{\underset{m=1}{\sum }}\frac{1}{m}\leq
c_{\alpha }<\infty .
\end{equation*}

According to (\ref{1.1}) together with (\ref{saa})-(\ref{saaa}) we can conclude that%
\begin{equation*}
\frac{1}{\log n}\overset{n}{\underset{m=M_{N}+1}{\sum }}\frac{\int_{%
\overline{I_{N}}}\left\vert t_{m}a\right\vert ^{1/\left( 1+\alpha \right)
}d\mu }{m}
\end{equation*}%
\begin{equation*}
=\frac{1}{\log n}\overset{n}{\underset{m=M_{N}+1}{\sum }}\overset{N-2}{%
\underset{k=0}{\sum }}\overset{N-1}{\underset{l=k+1}{\sum }}%
\sum\limits_{x_{j}=0,\text{ }j\in \{l+1,\dots ,N-1\}}^{m_{j-1}}\frac{%
\int_{I_{N}^{k,l}}\left\vert t_{m}a\right\vert ^{1/\left( 1+\alpha \right)
}d\mu }{m}
\end{equation*}%
\begin{equation*}
+\frac{1}{\log n}\overset{n}{\underset{m=M_{N}+1}{\sum }}\overset{N-1}{%
\underset{k=0}{\sum }}\frac{\int_{I_{N}^{k,N}}\left\vert t_{m}a\right\vert
^{1/\left( 1+\alpha \right) }d\mu }{m}
\end{equation*}%
\begin{equation*}
\leq \frac{c_{\alpha }}{\log n}\overset{n}{\underset{m=M_{N}+1}{\sum }}\frac{%
M_{N}^{\alpha /\left( 1+\alpha \right) }}{m^{\alpha /\left( 1+\alpha \right)
+1}}\overset{N-2}{\underset{k=0}{\sum }}\overset{N-1}{\underset{l=k+1}{\sum }%
}M_{l}^{\alpha /\left( 1+\alpha \right) }M_{k}^{1/\left( 1+\alpha \right) }%
\frac{m_{l+1}\dotsm m_{N-1}}{M_{N}}
\end{equation*}%
\begin{equation*}
+\frac{1}{\log n}\overset{n}{\underset{m=M_{N}+1}{\sum }}\frac{1}{mM_{N}}%
\overset{N-1}{\underset{k=0}{\sum }}M_{k}^{1/\left( 1+\alpha \right)
}M_{N}^{\alpha /\left( 1+\alpha \right) }.
\end{equation*}%
Since%
\begin{equation*}
\overset{N-2}{\underset{k=0}{\sum }}\overset{N-1}{\underset{l=k+1}{\sum }}%
M_{l}^{\alpha /\left( 1+\alpha \right) }M_{k}^{1/\left( 1+\alpha \right) }%
\frac{m_{l+1}\dotsm m_{N-1}}{M_{N}}
\end{equation*}%
\begin{equation*}
=\overset{N-2}{\underset{k=0}{\sum }}\overset{N-1}{\underset{l=k+1}{\sum }}%
M_{l}^{\alpha /\left( 1+\alpha \right) }M_{k}^{1/\left( 1+\alpha \right) }%
\frac{1}{M_{l}}
\end{equation*}%
\begin{equation*}
=\overset{N-2}{\underset{k=0}{\sum }}M_{k}^{1/\left( 1+\alpha \right) }%
\overset{N-1}{\underset{l=k+1}{\sum }}\frac{1}{M_{l}^{1/\left( 1+\alpha
\right) }}
\end{equation*}%
\begin{equation*}
\leq \overset{N-2}{\underset{k=0}{\sum }}M_{k}^{1/\left( 1+\alpha \right) }%
\frac{1}{M_{k}^{1/\left( 1+\alpha \right) }}\leq \overset{N-2}{\underset{k=0}%
{\sum }}1\leq N
\end{equation*}%
and%
\begin{equation*}
\overset{N-1}{\underset{k=0}{\sum }}M_{k}^{1/\left( 1+\alpha \right)
}M_{N}^{\alpha /\left( 1+\alpha \right) }
\end{equation*}%
\begin{equation*}
\leq M_{N}^{1/\left( 1+\alpha \right) }M_{N}^{\alpha /\left( 1+\alpha
\right) }\leq M_{N}
\end{equation*}%
we obtain that%
\begin{equation*}
\frac{1}{\log n}\overset{n}{\underset{m=M_{N}+1}{\sum }}\frac{\int_{%
\overline{I_{N}}}\left\vert t_{m}a\right\vert ^{1/\left( 1+\alpha \right)
}d\mu }{m}
\end{equation*}%
\begin{equation*}
\leq \frac{c_{\alpha }}{\log n}\left( \overset{n}{\underset{m=M_{N}+1}{\sum }%
}\frac{NM_{N}^{\alpha /\left( 1+\alpha \right) }}{m^{\alpha /\left( 1+\alpha
\right) +1}}+\overset{n}{\underset{m=M_{N}+1}{\sum }}\frac{1}{m}\right)
<c_{\alpha }<\infty .
\end{equation*}%
The proof is complete.
\QED

\subsection{Maximal operators of Riesz and N\"orlund logarithmic means on martingale Hardy spaces}

In our previous sections we investigated N\"orlund means with
non-increasing sequences $\left\{ q_{k}:k\in \mathbb{N}\right\} $, but the
case when $q_{k}=1/k$ was excluded, since this sequence does not satisfies the
condition (\ref{6a}) for any $0<\alpha \leq 1$. On the other hand, Riesz
logarithmic means are not examples of N\"orlund means. In this subsection we
fill up this gap simultaneously for both cases.

Both theorems in this
section are due to Tephnadze \cite{tep11}.

\begin{theorem}
\label{theorem1r}a)\bigskip\ The maximal operator$\ $of Riesz logarithmic
means $R^{\ast }$ \textit{is bounded from the Hardy space }$H_{1/2}$\textit{%
\ to the space }$weak-L_{1/2}.$

b) Let $0<p\leq 1/2$. Then there exists a martingale $f\in H_{p}$ such that
\begin{equation*}
\left\Vert R^{\ast }f\right\Vert _{p}=+\infty .
\end{equation*}
\end{theorem}

{\bf Proof}:
By using Abel transformation we obtain that
\begin{equation*}
R_{n}f=\frac{1}{l_{n}}\overset{n-1}{\underset{j=1}{\sum }}\frac{\sigma _{j}f%
}{j+1}+\frac{\sigma _{n}f}{l_{n}}.
\end{equation*}%
Consequently,%
\begin{equation}
R^{\ast }f\leq c\sigma ^{\ast }f.  \label{24a}
\end{equation}

Since $\sigma ^{\ast }$ is bounded from the martingale Hardy space $H_{1/2}$
to the space $weak-L_{1/2}$ by using (\ref{24a}) we can conclude that
\begin{equation*}
\left\Vert R^{\ast }f\right\Vert _{weak-L_{1/2}}\leq c\left\Vert
f\right\Vert _{H_{1/2}}
\end{equation*}%
and the proof of part a) is complete.

Let $f=\left( f^{\left( n\right) }:n\in
\mathbb{N}
\right) $ be martingale defined in Example \ref{example2.5.1}.

Set $q_{n}^{s}=M_{2n}+M_{2s}-1,$ $n>s.$ Then we can write that%
\begin{equation*}
R_{q_{\alpha _{k}}^{s}}f=\frac{1}{l_{q_{\alpha _{k}}^{s}}}\underset{j=1}{%
\overset{q_{\alpha _{k}}^{s}}{\sum }}\frac{S_{j}f}{j}
\end{equation*}%
\begin{equation*}
=\frac{1}{l_{q_{\alpha _{k}}^{s}}}\underset{j=1}{\overset{M_{2\alpha _{k}}}{%
\sum }}\frac{S_{j}f}{j}+\frac{1}{l_{q_{\alpha _{k}}^{s}}}\underset{%
j=M_{2\alpha _{k}}+1}{\overset{q_{\alpha _{k}}^{s}}{\sum }}\frac{S_{j}f}{j}%
:=I+II.
\end{equation*}

According to (\ref{sn101}) we have that%
\begin{equation*}
\left\vert I\right\vert \leq \frac{1}{l_{q_{\alpha _{k}}^{s}}}\underset{j=1}{%
\overset{M_{2\alpha _{k}}-1}{\sum }}\frac{\left\vert S_{j}f\left( x\right)
\right\vert }{j}
\end{equation*}%
\begin{equation*}
\leq \frac{1}{\alpha _{k}}\frac{2\lambda M_{2\alpha _{k-1}}^{1/p}}{\alpha
_{k-1}^{1/2}}\sum_{j=1}^{M_{2\alpha _{k}}-1}\frac{1}{j}
\end{equation*}%
\begin{equation*}
\leq \frac{2\lambda M_{2\alpha _{k-1}}^{1/p}}{\alpha _{k-1}^{1/2}}\leq \frac{%
2\lambda M_{\alpha _{k}}^{1/p}}{\alpha _{k}^{3/2}}.
\end{equation*}

Let $M_{2\alpha _{k}}\leq j\leq q_{\alpha _{k}}^{s}.$ According to the
second inequality of (\ref{1t6.1}) in the case when $l=k,$ we deduce that
\begin{equation}
S_{j}f=S_{M_{_{2\alpha _{k}}}}f+\frac{M_{2\alpha _{k}}^{1/p-1}\psi
_{_{M_{2\alpha _{k}}}}D_{j-M_{_{2\alpha _{k}}}}}{\alpha _{k}^{1/2}}.
\label{11y}
\end{equation}

Hence, we can rewrite $II$ as%
\begin{equation}
II=\frac{1}{l_{q_{\alpha _{k}}^{s}}}\underset{j=M_{2\alpha _{k}}}{\overset{%
q_{\alpha _{k}}^{s}}{\sum }}\ \frac{S_{M_{_{2\alpha _{k}}}}f}{j}  \label{12y}
\end{equation}%
\begin{equation}
+\frac{1}{l_{q_{\alpha _{k}}^{s}}}\frac{M_{2\alpha _{k}}^{1/p-1}\psi
_{_{M_{2\alpha _{k}}}}}{\alpha _{k}^{1/2}}\sum_{j=M_{2\alpha
_{k}}}^{q_{\alpha _{k}}^{s}}\frac{D_{j-M_{_{2\alpha _{k}}}}}{j}  \notag
\end{equation}%
\begin{equation*}
:=II_{1}+II_{2}.
\end{equation*}

In view of (\ref{sn101}) we find that%
\begin{equation}
\left\vert II_{1}\right\vert \leq \frac{1}{l_{q_{\alpha _{k}}^{s}}}\underset{%
j=M_{2\alpha _{k}}}{\overset{q_{\alpha _{k}}^{s}}{\sum }}\ \frac{1}{j}%
\left\vert S_{M_{_{2\alpha _{k}}}}f\right\vert \label{122y}
\end{equation}%
\begin{equation*}
\leq \left\vert S_{M_{_{2\alpha _{k}}}}f\right\vert \leq \frac{2\lambda
M_{2\alpha _{k-1}}^{1/p}}{\alpha _{k-1}^{1/2}}.
\end{equation*}

Let $0<p\leq 1/2,$ $x\in I_{2s}\backslash I_{2s+1}$ for $s=\left[ 2\alpha
_{k}/3\right] ,\ldots ,\alpha _{k}.$ Since%
\begin{equation*}
\sum_{j=0}^{M_{2s}-1}\frac{D_{_{j}}\left( x\right) }{_{j+M_{2\alpha _{k}}}}%
\geq \sum_{j=0}^{M_{2s}-1}\frac{j}{_{j+M_{2\alpha _{k}}}}
\end{equation*}%
\begin{equation*}
\geq \sum_{j=0}^{M_{2s}-1}\frac{j}{_{2M_{2\alpha _{k}}}}\geq \frac{%
cM_{2s}^{2}}{M_{2\alpha _{k}}}
\end{equation*}%
we obtain that%
\begin{equation}
\left\vert II_{2}\right\vert =\frac{1}{l_{q_{\alpha _{k}}^{s}}}\frac{%
M_{2\alpha _{k}}^{1/p-1}}{\alpha _{k}^{1/2}}\left\vert \psi _{_{M_{2\alpha
_{k}}}}\sum_{j=0}^{M_{2s}-1}\frac{D_{j}}{_{j+M_{2\alpha _{k}}}}\right\vert
\label{15y}
\end{equation}%
\begin{equation*}
\geq \frac{c}{\alpha _{k}}\frac{M_{2\alpha _{k}}^{1/p-1}}{\alpha _{k}^{1/2}}%
\frac{M_{2s}^{2}}{M_{2\alpha _{k}}}\geq \frac{cM_{2\alpha
_{k}}^{1/p-2}M_{2s}^{2}}{\alpha _{k}^{3/2}}.
\end{equation*}

By combining (\ref{1t2})-(\ref{1t4}) with (\ref{12y})-(\ref{15y}%
) we get that%
\begin{equation*}
\left\vert R_{q_{\alpha _{k}}^{s}}f\right\vert =\left\vert
II_{2}-I-II_{1}\right\vert
\end{equation*}%
\begin{equation*}
\geq \left\vert II_{2}\right\vert -\left\vert I\right\vert -\left\vert
II_{1}\right\vert \geq \left\vert II_{2}\right\vert -\frac{2\lambda
M_{\alpha _{k}}^{1/p}}{\alpha _{k}^{3/2}}
\end{equation*}%
\begin{equation*}
\geq \frac{cM_{2\alpha _{k}}^{1/p-2}M_{2s}^{2}}{\alpha _{k}^{3/2}}-\frac{%
cM_{\alpha _{k}}^{1/p}}{\alpha _{k}^{3/2}}
\end{equation*}%
\begin{equation*}
\geq \frac{cM_{2\alpha _{k}}^{1/p-2}M_{2s}^{2}}{\alpha _{k}^{3/2}}.
\end{equation*}

Hence,%
\begin{equation*}
\int_{G_{m}}\left\vert R^{\ast }f\left( x\right) \right\vert ^{p}d\mu \left(
x\right)
\end{equation*}%
\begin{equation*}
\geq c_{p}\sum_{s=\left[ 2\alpha _{k}/3\right] }^{\alpha _{k}}\underset{%
I_{2s}\backslash I_{2s+1}}{\int }\left\vert R_{q_{\alpha _{k}}^{s}}f\left(
x\right) \right\vert ^{p}d\mu \left( x\right)
\end{equation*}%
\begin{equation*}
\geq c_{p}\overset{\alpha _{k}}{\underset{s=\left[ 2\alpha _{k}/3\right] }{%
\sum }}\underset{I_{2s}\backslash I_{2s+1}}{\int }\frac{cM_{2\alpha
_{k}}^{1-2p}M_{2s}^{2p}}{\alpha _{k}^{3p/2}}d\mu \left( x\right)
\end{equation*}%
\begin{equation*}
\geq c_{p}\sum_{s=\left[ 2\alpha _{k}/3\right] }^{\alpha _{k}}\frac{%
M_{2\alpha _{k}}^{1-2p}M_{2s}^{2p-1}}{\alpha _{k}^{3p/2}}
\end{equation*}%
\begin{equation*}
\geq \left\{
\begin{array}{ll}
\frac{c_{p}2^{\alpha _{k}(1-2p)/3}}{\alpha _{k}^{3p/2}}, & \text{\thinspace
when \ }0<p<1/2, \\
c\alpha _{k}^{1/4}, & \text{when \ }p=1/2,%
\end{array}%
\right. \rightarrow \infty \text{, when \ }k\rightarrow \infty .
\end{equation*}

The proof is complete.
\QED

\begin{theorem}
\label{theorem2nl}Let $0<p\leq 1$. Then there exists a martingale $f\in
H_{p} $ such that%
\begin{equation*}
\left\Vert L^{\ast }f\right\Vert _{p}=+\infty .
\end{equation*}
\end{theorem}

{\bf Proof}:
We write that%
\begin{equation}
L_{q_{\alpha _{k}}^{s}}f=\frac{1}{l_{q_{\alpha _{k},s}}}\underset{j=1}{%
\overset{q_{\alpha _{k}}^{s}}{\sum }}\frac{S_{j}f}{q_{\alpha _{k}}^{s}-j}
\label{16y}
\end{equation}%
\begin{equation}
=\frac{1}{l_{q_{\alpha _{k}}^{s}}}\underset{j=1}{\overset{M_{2\alpha _{k}}-1}%
{\sum }}\frac{S_{j}f}{q_{\alpha _{k}}^{s}-j}+\frac{1}{q_{\alpha _{k}}^{s}}%
\underset{j=M_{2\alpha _{k}}}{\overset{q_{\alpha _{k}}^{s}}{\sum }}\frac{%
S_{j}f}{q_{\alpha _{k}}^{s}-j}:=III+IV.  \notag
\end{equation}

In the view of (\ref{sn101}) for $III$ we get the following estimate:%
\begin{equation}
\left\vert III\right\vert \leq \frac{1}{\alpha _{k}}\sum_{j=0}^{M_{2\alpha
_{k-1}}}\frac{1}{q_{\alpha _{k}}^{s}-j}\frac{M_{2\alpha _{k-1}}^{1/p}}{%
\alpha _{k-1}^{1/2}}\text{ }\leq \frac{M_{2\alpha _{k-1}}^{1/p}}{\alpha
_{k-1}^{1/2}}\leq \frac{M_{\alpha _{k}}^{1/p}}{\alpha _{k}^{3/2}}.
\label{17y}
\end{equation}

Moreover, according to the second inequality of (\ref{1t6.1}) in the case
when $l=k,$ (see also (\ref{11y})) we can rewrite $IV$ as
\begin{equation}
IV=\frac{1}{l_{q_{\alpha _{k}}^{s}}}\sum_{j=M_{2\alpha _{k}}}^{q_{\alpha
_{k}}^{s}}\frac{S_{M_{_{2\alpha _{k}}}}f}{q_{\alpha _{k},s}-j}  \label{18y}
\end{equation}%
\begin{equation*}
+\frac{1}{l_{q_{\alpha _{k}}^{s}}}\frac{M_{2\alpha _{k}}^{1/p-1}\psi
_{_{M_{2\alpha _{k}}}}}{\alpha _{k}^{1/2}}\sum_{j=M_{2\alpha
_{k}}}^{q_{\alpha _{k}}^{s}}\frac{D_{j-M_{_{2\alpha _{k}}}}}{q_{\alpha
_{k}}^{s}-j}:=IV_{1}+IV_{2}.
\end{equation*}

By now applying (\ref{sn101}) again we have that
\begin{equation}
\left\vert IV_{1}\right\vert \leq \frac{2\lambda M_{2\alpha _{k-1}}^{1/p}}{%
\alpha _{k-1}^{1/2}}\leq \frac{2\lambda M_{\alpha _{k}}^{1/p}}{\alpha
_{k}^{3/2}}.  \label{19y}
\end{equation}

Let $x\in I_{2s}\backslash I_{2s+1},$ $s=\left[ 2\alpha _{k}/3\right]
,...,\alpha _{k}.$ Since%
\begin{equation*}
\sum_{j=0}^{M_{2s}-1}\frac{D_{j}}{_{M_{2s}}-j}=\sum_{j=0}^{M_{2s}-1}\frac{j}{%
_{M_{2s}}-j}=\sum_{j=0}^{M_{2s}-1}\frac{M_{2s}}{_{M_{2s}}-j}%
-\sum_{j=0}^{M_{2s}-1}\frac{M_{2s}-j}{_{M_{2s}}-j}
\end{equation*}%
\begin{equation*}
=c_{p}sM_{2s}-M_{2s}\geq c_{p}sM_{2s},
\end{equation*}%
we obtain that%
\begin{equation}
\left\vert IV_{2}\right\vert =\frac{1}{l_{q_{\alpha _{k},s}}}\frac{%
M_{2\alpha _{k}}^{1/p-1}}{\alpha _{k}^{1/2}}\left\vert \sum_{j=0}^{M_{2s}-1}%
\frac{D_{j}}{_{M_{2s}}-j}\right\vert  \label{21y}
\end{equation}%
\begin{equation*}
\geq \frac{c_{p}M_{2\alpha _{k}}^{1/p-1}}{\alpha _{k}^{3/2}}sM_{2s},\qquad
x\in I_{2s}/I_{2s+1}.
\end{equation*}

By combining (\ref{1t4}) with (\ref{16y})-(\ref{21y}) for $x\in
I_{2s}\backslash I_{2s+1},s=\left[ 2\alpha _{k}/3\right] ,...,\alpha _{k}$
and $0<p\leq 1$ we get that%
\begin{equation*}
\left\vert L_{q_{\alpha _{k}}^{s}}f\left( x\right) \right\vert =\left\vert
III+IV_{1}+IV_{2}\right\vert \geq \left\vert IV_{2}\right\vert -\left\vert
III\right\vert -\left\vert IV_{1}\right\vert
\end{equation*}%
\begin{equation*}
\geq \frac{c_{p}M_{2\alpha _{k}}^{1/p-1}}{\alpha _{k}^{3/2}}sM_{2s}-\frac{%
2\lambda M_{\alpha _{k}}}{\alpha _{k}^{3/2}}\geq \frac{c_{p}M_{2\alpha
_{k}}^{1/p-1}}{\alpha _{k}^{3/2}}sM_{2s}.
\end{equation*}

Consequently,%
\begin{equation*}
\int_{G_{m}}\left\vert L^{\ast }f\left( x\right) \right\vert ^{p}d\mu \left(
x\right) \geq \sum_{s=\left[ 2\alpha _{k}/3\right] }^{\alpha
_{k}}\int_{I_{2s}\backslash I_{2s+1}}\left\vert L^{\ast }f\left( x\right)
\right\vert ^{p}d\mu \left( x\right)
\end{equation*}%
\begin{equation*}
\geq \sum_{s=\left[ 2\alpha _{k}/3\right] }^{\alpha
_{k}}\int_{I_{2s}\backslash I_{2s+1}}\left\vert L_{q_{\alpha
_{k}}^{s}}f\left( x\right) \right\vert ^{p}d\mu \left( x\right)
\end{equation*}%
\begin{equation*}
\geq c_{p}\sum_{s=\left[ 2\alpha _{k}/3\right] }^{\alpha
_{k}}\int_{I_{2s}\backslash I_{2s+1}}\frac{M_{2\alpha _{k-1}}^{1-p}}{\alpha
_{k}^{3p/2}}s^{p}M_{2s}^{p}d\mu \geq c_{p}\underset{s=\left[ 2\alpha _{k}/3%
\right] }{\overset{\alpha _{k}}{\sum }}\frac{M_{2\alpha _{k-1}}^{1-p}}{%
\alpha _{k}^{p/2}}M_{2s}^{p-1}
\end{equation*}%
\begin{equation*}
\geq \left\{
\begin{array}{ll}
\frac{c_{p}2^{\alpha _{k}(1-p)/3}}{\alpha _{k}^{p/2}}, & \text{when }0<p<1,
\\
c\alpha _{k}^{1/2}, & \text{when \ }p=1,%
\end{array}%
\right. \rightarrow \infty \text{, when \ }k\rightarrow \infty .
\end{equation*}

The proof is complete.
\QED

\subsection{Applications}

\bigskip

First we consider N\"orlund means $t_{n}$ with monotone and
bounded sequence $\left\{ q_{k}:k\in \mathbb{N}\right\} .$ The results in
our previous sections in particular imply the following results:

\begin{theorem}
a) Let $f\in H_{1/2}$ and $t_{n}$ be N\"orlund means with monotone and
bounded sequence $\left\{ q_{k}:k\in \mathbb{N}\right\} .$ Then there exists an
absolute constant $c$ such that%
\begin{equation*}
\left\Vert t^{\ast }f\right\Vert _{weak-L_{1/2}}\leq c\left\Vert
f\right\Vert _{H_{1/2}}.
\end{equation*}%
b) There exists a martingale $f\in H_{1/2},$ \textit{such that}
\begin{equation*}
\left\Vert t^{\ast }f\right\Vert _{1/2}=\infty .
\end{equation*}
\end{theorem}

\begin{theorem}
a) Let $p>1/2$, $f\in H_{p}$ and $t_{n}$ be N\"orlund means with monotone
and bounded sequence $\left\{ q_{k}:k\in \mathbb{N}\right\} .$ Then there
exists an absolute constant $c_{p}$ depending only on $p$ such that%
\begin{equation*}
\left\Vert \overset{\sim }{t_{p,1}}^{\ast }f\right\Vert _{p}\leq
c_{p}\left\Vert f\right\Vert _{H_{p}}.
\end{equation*}%
b) Let $\varphi :\mathbb{N}_{+}\rightarrow \lbrack 1,\infty )$ be a
non-decreasing function satisfying the condition (\ref{cond2}). Then the following
maximal operator%
\begin{equation*}
\sup_{n\in \mathbb{N}}\frac{\left\vert t_{n}f\right\vert }{\Phi _{n}}
\end{equation*}%
\textit{is not bounded from the Hardy space }$H_{p}$\textit{\ to the space }%
$weak-L_{p}.$
\end{theorem}

\begin{theorem}
a) Let $f\in H_{1/2}$ and $t_{n}$ be N\"orlund means with monotone and
bounded sequence $\left\{ q_{k}:k\in \mathbb{N}\right\} .$ Then there exists an
absolute constant $c$ such that
\begin{equation*}
\left\Vert \overset{\sim }{t_{1}}^{\ast }f\right\Vert _{1/2}\leq c\left\Vert
f\right\Vert _{H_{1/2}}.
\end{equation*}%
b) Let $\varphi :\mathbb{N}_{+}\rightarrow \lbrack 1,\infty )$ be a
non-decreasing function satisfying the condition (\ref{cond1}). Then the following
maximal operator%
\begin{equation*}
\sup_{n\in \mathbb{N}}\frac{\left\vert B_{n}f\right\vert }{\Phi _{n}}
\end{equation*}%
\textit{is not bounded from the Hardy space }$H_{1/2}$\textit{\ to the
Lebesgue space }$L_{1/2}.$
\end{theorem}

\begin{theorem}
Let $0<p<1/2$, $f\in H_{p}$ and $t_{n}$ be N\"orlund means with monotone and
bounded sequence $\left\{ q_{k}:k\in \mathbb{N}\right\} .$ Then there exists
an absolute constant $c_{p}$ depending only on $p$ such that%
\begin{equation*}
\overset{\infty }{\underset{m=1}{\sum }}\frac{\left\Vert t_{m}f\right\Vert
_{p}^{p}}{m^{2-2p}}\leq c_{p}\left\Vert f\right\Vert _{H_{p}}^{p}.
\end{equation*}
\end{theorem}

\begin{theorem}
Let $f\in H_{1/2}$ and $t_{n}$ be N\"orlund means with monotone and bounded
sequence $\left\{ q_{k}:k\in \mathbb{N}\right\} .$ Then there exists an
absolute constant $c$ such that%
\begin{equation*}
\frac{1}{\log n}\overset{n}{\underset{m=1}{\sum }}\frac{\left\Vert
t_{m}f\right\Vert _{1/2}^{1/2}}{m}\leq c\left\Vert f\right\Vert
_{H_{1/2}}^{1/2},\text{ \ }n=2,3,\ldots .
\end{equation*}
\end{theorem}

Next we remark that $\sigma _{n}$ are N\"orlund means with monotone and
bounded sequence $\left\{ q_{k}:k\in \mathbb{N}\right\} $ and $\kappa
_{n}^{\alpha ,\beta }$ are concrete examples of N\"orlund means with
non-decreasing and unbounded sequence $\left\{ q_{k}:k\in \mathbb{N}\right\}
$ but they have similar boundedness properties. The results in our previous sections in particular imply the following
results:

\begin{theorem}
a) Let $f\in H_{1/2}.$ Then there exists absolute constant $c$ such that%
\begin{equation*}
\left\Vert \sigma ^{\ast }f\right\Vert _{weak-L_{1/2}}\leq c\left\Vert
f\right\Vert _{H_{1/2}}
\end{equation*}%
and%
\begin{equation*}
\left\Vert \kappa ^{\alpha ,\beta ,\ast }f\right\Vert _{weak-L_{1/2}}\leq
c\left\Vert f\right\Vert _{H_{1/2}}.
\end{equation*}%
b) T\textit{here exists a martingale} $f\in H_{1/2}$ \textit{such that}
\begin{equation*}
\left\Vert \sigma ^{\ast }f\right\Vert _{1/2}=\infty
\end{equation*}%
and%
\begin{equation*}
\left\Vert \kappa ^{\alpha ,\beta ,\ast }f\right\Vert _{1/2}=\infty .
\end{equation*}
\end{theorem}

\begin{theorem}
a) Let $0<p<1/2$ and $f\in H_{p}.$ Then there exists an absolute constant $%
c_{p}$ depending only on $p$ such that%
\begin{equation*}
\left\Vert \overset{\sim }{\sigma _{p}}^{\ast }f\right\Vert _{p}\leq
c_{p}\left\Vert f\right\Vert _{H_{p}}
\end{equation*}%
and%
\begin{equation*}
\left\Vert \widetilde{\kappa }_{p}^{\alpha ,\beta ,\ast }f\right\Vert
_{p}\leq c_{p}\left\Vert f\right\Vert _{H_{p}}.
\end{equation*}%
b) Let $\varphi :\mathbb{N}_{+}\rightarrow \lbrack 1,\infty )$ be a
non-decreasing function satisfying the condition (\ref{cond2}). Then the following
maximal operators%
\begin{equation*}
\text{\ }\sup_{n\in \mathbb{N}}\frac{\left\vert \sigma _{n}f\right\vert }{%
\Phi _{n}}\text{ \ \ \ and \ \ \ }\sup_{n\in \mathbb{N}}\frac{\left\vert
\kappa _{n}^{\alpha ,\beta }\right\vert }{\Phi _{n}}
\end{equation*}%
\textit{are not bounded from the Hardy space }$H_{p}$\textit{\ to the space }%
$weak-L_{p}.$
\end{theorem}

\begin{theorem}
a) Let $f\in H_{1/2}.$ Then there exists an absolute constant $c$ such that
\begin{equation*}
\left\Vert \overset{\sim }{\sigma }^{\ast }f\right\Vert _{1/2}\leq
c\left\Vert f\right\Vert _{H_{1/2}}
\end{equation*}%
and%
\begin{equation*}
\left\Vert \widetilde{\kappa }^{\alpha ,\beta ,\ast }f\right\Vert _{1/2}\leq
c\left\Vert f\right\Vert _{H_{1/2}}.
\end{equation*}%
b) Let $\varphi :\mathbb{N}_{+}\rightarrow \lbrack 1,\infty )$ be a
nondecreasing function satisfying the condition (\ref{cond1}). Then the following
maximal operators%
\begin{equation*}
\text{\ }\sup_{n\in \mathbb{N}}\frac{\left\vert \sigma _{n}f\right\vert }{%
\Phi _{n}}\text{ \ \ \ and \ \ \ }\sup_{n\in \mathbb{N}}\frac{\left\vert
\kappa _{n}^{\alpha ,\beta }\right\vert }{\Phi _{n}}
\end{equation*}%
\textit{are not bounded from the Hardy space }$H_{1/2}$\textit{\ to the
Lebesgue space }$L_{1/2}.$
\end{theorem}

\begin{theorem}
Let $0<p<1/2$ and $f\in H_{p}.$ Then there exists an absolute constant $c_{p}$ depending only on $p$ such that%
\begin{equation*}
\overset{\infty }{\underset{m=1}{\sum }}\frac{\left\Vert \sigma
_{m}f\right\Vert _{p}^{p}}{m^{2-2p}}\leq c_{p}\left\Vert f\right\Vert
_{H_{p}}^{p}
\end{equation*}%
and%
\begin{equation*}
\overset{\infty }{\underset{m=1}{\sum }}\frac{\left\Vert \kappa _{m}^{\alpha
,\beta }f\right\Vert _{p}^{p}}{m^{2-2p}}\leq c_{p}\left\Vert f\right\Vert
_{H_{p}}^{p}.
\end{equation*}
\end{theorem}

\begin{theorem}
Let $f\in H_{1/2}.$ Then there exists an absolute constant $c$ such that%
\begin{equation*}
\frac{1}{\log n}\overset{n}{\underset{m=1}{\sum }}\frac{\left\Vert \sigma
_{m}f\right\Vert _{1/2}^{1/2}}{m}\leq c\left\Vert f\right\Vert
_{H_{1/2}}^{1/2},\text{ \ }n=2,3,\ldots ,
\end{equation*}%
and
\begin{equation*}
\frac{1}{\log n}\overset{n}{\underset{m=1}{\sum }}\frac{\left\Vert \kappa
_{m}^{\alpha ,\beta }f\right\Vert _{1/2}^{1/2}}{m}\leq c\left\Vert
f\right\Vert _{H_{1/2}}^{1/2},\text{ }n=2,3,\ldots .
\end{equation*}
\end{theorem}

Furthemore, we note that it is obvious that Ces\`{a}ro $\left( C,\alpha \right) $ and Riesz $\left(
R,\alpha \right) $ means are examples of N\"orlund means with non-increasing
sequence satisfying the conditions (\ref{6a}) and (\ref{7a}). It follows that
all results concerning such summability methods are true also for the Ces%
\`{a}ro $\left( C,\alpha \right) $ and Riesz $\left( R,\alpha \right) $
means. Hence, the results in our previous sections in particular imply the
following results:

\begin{theorem}
a) Let $0<p<1/2.$ Then there exists a martingale $f\in H_{p}$ such that%
\begin{equation*}
\sup_{n\in \mathbb{N}}\left\Vert \sigma _{n}^{\alpha }f\right\Vert
_{weak-L_{p}}=\infty
\end{equation*}%
and%
\begin{equation*}
\sup_{n\in \mathbb{N}}\left\Vert \beta _{n}^{\alpha }f\right\Vert
_{weak-L_{p}}=\infty .
\end{equation*}
\end{theorem}

\begin{theorem}
a) Let $0<\alpha \leq 1$ and $f\in H_{1/(1+\alpha )}.$ Then there exists
absolute constant $c_{\alpha }$ depending only on $\alpha$ such that
\begin{equation*}
\left\Vert \sigma ^{\alpha ,\ast }f\right\Vert _{weak-L_{1/\left( 1+\alpha
\right) }}\leq c_{\alpha }\left\Vert f\right\Vert _{H_{1/\left( 1+\alpha
\right) }}
\end{equation*}%
and%
\begin{equation*}
\left\Vert \beta ^{\alpha ,\ast }f\right\Vert _{weak-L_{1/\left( 1+\alpha
\right) }}\leq c_{\alpha }\left\Vert f\right\Vert _{H_{1/\left( 1+\alpha
\right) }}.
\end{equation*}%
b) Let $0<\alpha \leq 1.$ \textit{Then there exists a martingale} $f\in
H_{1/(1+\alpha )}$ \textit{such that}
\begin{equation*}
\left\Vert \sigma ^{\alpha ,\ast }f\right\Vert _{1/\left( 1+\alpha \right)
}=\infty
\end{equation*}%
and
\begin{equation*}
\left\Vert \beta ^{\alpha ,\ast }f\right\Vert _{1/\left( 1+\alpha \right)
}=\infty .
\end{equation*}
\end{theorem}

\begin{theorem}
a) Let $0<\alpha \leq 1$, $0<p<1/\left( 1+\alpha \right) $ and $f\in H_{p}.$
Then there exists an absolute constant $c_{\alpha ,p}$ depending only on $%
\alpha $ and $p$ such that
\begin{equation*}
\left\Vert \overset{\sim }{\sigma _{p}}^{\alpha ,\ast }f\right\Vert _{p}\leq
c_{\alpha ,p}\left\Vert f\right\Vert _{H_{p}}
\end{equation*}%
and%
\begin{equation*}
\left\Vert \overset{\sim }{\beta _{p}}^{\alpha ,\ast }f\right\Vert _{p}\leq
c_{\alpha ,p}\left\Vert f\right\Vert _{H_{p}}.
\end{equation*}%
b) Let $0<\alpha \leq 1$ and $\varphi :\mathbb{N}_{+}\rightarrow \lbrack
1,\infty )$ be a non-decreasing function satisfying the condition (\ref%
{6lbbb}). Then the following maximal operators%
\begin{equation*}
\sup_{n\in \mathbb{N}}\frac{\left\vert \sigma _{n}^{\alpha }f\right\vert }{%
\Phi _{n}}\text{ \ \ \ and \ \ \ }\sup_{n\in \mathbb{N}}\frac{\left\vert
\beta _{n}^{\alpha }f\right\vert }{\Phi _{n}}
\end{equation*}%
\textit{are not bounded from the Hardy space }$H_{p}$\textit{\ to the space }%
$weak-L_{p}.$
\end{theorem}

\begin{theorem}
a) Let $0<\alpha \leq 1$ and $f\in H_{1/(1+\alpha )}.$ Then there exists an
absolute constant $c_{\alpha }$ depending only on $\alpha $ such that
\begin{equation*}
\left\Vert \overset{\sim }{\sigma }^{\alpha ,\ast }f\right\Vert _{1/\left(
1+\alpha \right) }\leq c_{\alpha }\left\Vert f\right\Vert _{H_{1/\left(
1+\alpha \right) }}
\end{equation*}%
and
\begin{equation*}
\left\Vert \overset{\sim }{\beta }^{\alpha ,\ast }f\right\Vert _{1/\left(
1+\alpha \right) }\leq c_{\alpha }\left\Vert f\right\Vert _{H_{1/\left(
1+\alpha \right) }}.
\end{equation*}%
b) Let $0<\alpha \leq 1$ and $\varphi :\mathbb{N}_{+}\rightarrow \lbrack
1,\infty )$ be a nondecreasing function satisfying the condition (\ref{nom1}%
). Then the following maximal operators%
\begin{equation*}
\sup_{n\in \mathbb{N}}\frac{\left\vert \sigma _{n}^{\alpha }f\right\vert }{%
\Phi _{n}}\text{ \ \ and \ \ \ \ \ }\sup_{n\in \mathbb{N}}\frac{\left\vert
\beta _{n}^{\alpha }f\right\vert }{\Phi _{n}}
\end{equation*}%
\textit{are not bounded from the Hardy space }$H_{1/\left( 1+\alpha \right)
} $\textit{\ to the Lebesgue space }$L_{1/\left( 1+\alpha \right) }.$
\end{theorem}

\begin{theorem}
Let $0<\alpha \leq 1$, $0<p<1/\left( 1+\alpha \right) $ and $f\in H_{p}.$
Then there exists an absolute constant $c_{\alpha ,p}$ depending only on $%
\alpha $ and $p$ such that%
\begin{equation*}
\overset{\infty }{\underset{m=1}{\sum }}\frac{\left\Vert \sigma _{m}^{\alpha
}a\right\Vert _{p}^{p}}{m^{\left( 1+\alpha \right) \left( 1-p\right) }}\leq
c_{\alpha ,p}\left\Vert f\right\Vert _{H_{p}}^{p}
\end{equation*}%
and
\begin{equation*}
\overset{\infty }{\underset{m=1}{\sum }}\frac{\left\Vert \beta _{m}^{\alpha
}a\right\Vert _{p}^{p}}{m^{\left( 1+\alpha \right) \left( 1-p\right) }}\leq
c_{\alpha ,p}\left\Vert f\right\Vert _{H_{p}}^{p}.
\end{equation*}
\end{theorem}

\begin{theorem}
Let $0<\alpha \leq 1$ and $f\in H_{1/(1+\alpha )}.$ Then there exists an
absolute constant $c_{\alpha }$ depending only on $\alpha$ such that
\begin{equation*}
\frac{1}{\log n}\overset{n}{\underset{m=1}{\sum }}\frac{\left\Vert \sigma
_{m}^{\alpha }f\right\Vert _{1/\left( 1+\alpha \right) }^{1/\left( 1+\alpha
\right) }}{m}\leq c_{\alpha }\left\Vert f\right\Vert _{H_{1/\left( 1+\alpha
\right) }}^{1/\left( 1+\alpha \right) }
\end{equation*}%
and
\begin{equation*}
\frac{1}{\log n}\overset{n}{\underset{m=1}{\sum }}\frac{\left\Vert \beta
_{m}^{\alpha }f\right\Vert _{1/\left( 1+\alpha \right) }^{1/\left( 1+\alpha
\right) }}{m}\leq c_{\alpha }\left\Vert f\right\Vert _{H_{1/\left( 1+\alpha
\right) }}^{1/\left( 1+\alpha \right) }.
\end{equation*}
\end{theorem}

Finally, we present some applications concerning almost everywhere convergence of some summability methods. To study almost everywhere convergence of some summability methods is one of
the most difficult topics in Fourier analysis. They involve techniques from
function theory and Hardy spaces.

In most applications the a.e. convergence of $\left\{ T_{n}:n\in \mathbb{N}%
\right\} $ can be established for $f$ in some dense class of $L_{1}\left(
G_{m}\right) .$ In particular, the following result play an important role for studying this type of questions (see e.g. the books \cite{gol}, \cite{sws} and
\cite{13}).

\begin{lemma}
\label{lemmaae} Let $f\in L_{1}$ and $T_{n}:L_{1}\rightarrow L_{1}$ be some
sub-linear operators and
\begin{equation*}
T^{\ast }:=\sup_{n\in
\mathbb{N}
}\left\vert T_{n}\right\vert .
\end{equation*}%
If
\begin{equation*}
T_{n}f\rightarrow f\text{ a.e. \ \ for every }f\in S
\end{equation*}%
where the set $S$ is dense in the space $\ L_{1}$ and the maximal operator $%
T^{\ast }$ is bounded from the space $L_{1}$ to the space $weak-L_{1},$ that
is%
\begin{equation*}
\sup_{\lambda >0}\lambda \mu \left\{ x\in G_{m}:\text{ }\left\vert T^{\ast
}f\left( x\right) \right\vert >\lambda \right\} \leq \left\Vert f\right\Vert
_{1},
\end{equation*}%
then
\begin{equation*}
T_{n}f\rightarrow f,\text{ a.e. \ \ for every }f\in L_{1}\left( G_{m}\right)
.
\end{equation*}
\end{lemma}

\begin{remark}
\label{remarkae}Since the Vilenkin function $\psi _{m}$ is constant on $%
I_{n}(x)$ for every $x\in G_{m}$ and $0\leq m<M_{n},$ it is clear that each
Vilenkin function is a complex-valued step function, that is, it is a finite
linear combination of the characteristic functions%
\begin{equation*}
\chi \left( E\right) =\left\{
\begin{array}{ll}
1, & x\in E, \\
0, & x\notin E.%
\end{array}%
\right.
\end{equation*}%
On the other hand, notice that, by Corollary \ref{dn2.3} (Paley lemma), it
yields that%
\begin{equation*}
\chi \left( I_{n}(t)\right) \left( x\right) =\frac{1}{M_{n}}%
\sum_{j=0}^{M_{n}-1}\psi _{j}\left( x-t\right) ,\text{ \ \ }x\in I_{n}(t),
\end{equation*}%
for each $x,t\in $ $G_{m}$ and $n\in \mathbb{N}$. Thus each step function is
a Vilenkin polynomial. Consequently, we obtain that the collection of step functions
coincides with a collection of Vilenkin polynomials $\mathcal{P}$. Since the
Lebesgue measure is regular it follows from Lusin theorem that given $%
f\in L_{1}$ there exist Vilenkin polynomials $P_{1},P_{2}...,$ such that $%
P_{n}\rightarrow f$ \ a.e. when $n\rightarrow \infty .$ This means that the
Vilenkin polynomials are dense in the space $L_{1}.$
\end{remark}

By using Lemma \ref{lemmaae}, remarks \ref{remarkae} and the results in our previous sections we in particular obtain the
following a.e. convergence results.

\begin{theorem}
\label{theoremae1}Let $f\in L_{1}$ and $t_{n}$ be the regular N\"orlund means
with non-decreasing sequence $\{q_{k}:k\in
\mathbb{N}
\}$. Then%
\begin{equation*}
t_{n}f\rightarrow f\text{ \ \ a.e., \ \ when \ }n\rightarrow \infty .\text{\ }
\end{equation*}
\end{theorem}

{\bf Proof}:
Since
\begin{equation*}
S_{n}P=P,\text{ for every }P\in \mathcal{P}
\end{equation*}%
we obtain that
\begin{equation*}
t_{n}P\rightarrow P\ \text{a.e.,}\ \ \ \text{when \ \ \ }\ n\rightarrow
\infty ,
\end{equation*}%
for every regular N\"orlund mean with non-decreasing sequence $\{q_{k}:k\in
\mathbb{N}
\}$, where $P\in \mathcal{P}$ is dense in the space $L_{1}$ (see Remark \ref%
{remarkae})$.$

On the other hand, by combining Lemma \ref{lemma2.3} and Theorem \ref%
{theorem1norlund} we obtain that the maximal operator of N\"orlund means, with
non-decreasing sequence $\{q_{k}:k\in
\mathbb{N}
\}$ satisfying conditions (\ref{fn0}) is bounded from the space $L_{1}$ to
the space $weak-L_{1},$ that is,%
\begin{equation*}
\sup_{\lambda >0}\lambda \mu \left\{ x\in G_{m}:\text{ }\left\vert t^{\ast
}f\left( x\right) \right\vert >\lambda \right\} \leq \left\Vert f\right\Vert
_{1}.
\end{equation*}

According to Lemma \ref{lemmaae} we obtain that under condition (\ref{fn0})
we have a.e. convergence%
\begin{equation*}
t_{n}f\rightarrow f\text{ \ \ a.e., \ \ when \ }n\rightarrow \infty .
\end{equation*}

The proof is complete.
\QED

\begin{corollary}
\label{corollaryae1}Let $f\in L_{1}$. Then%
\begin{equation*}
\sigma _{n}f\rightarrow f\text{ \ \ \ a.e., \ \ \ when \ }n\rightarrow \infty
,\text{\ \ \ }
\end{equation*}%
and
\begin{equation*}
\kappa _{n}^{\alpha ,\beta }f\rightarrow f\text{ \ \ \ a.e., \ \ when \ }%
n\rightarrow \infty .
\end{equation*}
\end{corollary}

{\bf Proof}:
Since $\sigma _{n}$ and $\kappa _{n}$ are regular N\"orlund means with
non-decreasing sequence $\{q_{k}:k\in
\mathbb{N}
\}$ (see (\ref{monotone1}) and (\ref{node00})) the proof is complete by
just using Theorem \ref{theoremae1}.
\QED

\begin{theorem}
\label{theoremae2}Let $f\in L_{1}$ and $t_{n}$ be N\"orlund means with
non-increasing sequence $\{q_{k}:k\in
\mathbb{N}
\}$ satisfying conditions (\ref{6a}) and (\ref{7a}). Then
\begin{equation*}
t_{n}f\rightarrow f\text{ \ \ a.e., \ \ when\ }n\rightarrow \infty .
\end{equation*}
\end{theorem}

{\bf Proof}:
By using part b) of Theorem \ref{theorem1n} we get that every N\"orlund mean,
with non-increasing sequence $\{q_{k}:k\in
\mathbb{N}
\}$ is regular. Since
\begin{equation*}
S_{n}P=P,\text{ for every }P\in \mathcal{P}
\end{equation*}%
we obtain that $t_{n}P\rightarrow P,$ \ when \ $n\rightarrow \infty ,$ where $%
P\in \mathcal{P}$ is dense in the space $L_{1}$ (see Remark \ref{remarkae})$%
. $

On the other hand, by combining Lemma \ref{lemma2.3} and part a) of Theorem %
\ref{Theorem3normax} we obtain that the maximal operator of N\"orlund means with non-increasing sequence $\{q_{k}:k\in
\mathbb{N}
\}$ satisfying conditions (\ref{6a}) and (\ref{7a}) is bounded from the
space $L_{1}$ to the space $weak-L_{1},$ that is%
\begin{equation*}
\sup_{\lambda >0}\lambda \mu \left\{ x\in G_{m}:\text{ }\left\vert t^{\ast
}f\left( x\right) \right\vert >\lambda \right\} \leq \left\Vert f\right\Vert
_{1}.
\end{equation*}

According to Lemma \ref{lemmaae} we obtain that under conditions (\ref{6a})
and (\ref{7a}) we have a.e. convergence%
\begin{equation*}
t_{n}f\rightarrow f\text{ \ \ a.e., \ \ when \ }n\rightarrow \infty .
\end{equation*}

The proof is complete.
\QED

\begin{corollary}
\label{corollaryae2}Let $f\in L_{1}$. Then%
\begin{equation*}
\text{\ \ \ }\sigma _{n}^{\alpha }f\rightarrow f\text{ \ \ \ a.e., \ \ when \
}n\rightarrow \infty ,\text{ \ \ when \ }0<\alpha <1.
\end{equation*}%
and
\begin{equation*}
\beta _{n}^{\alpha }f\rightarrow f\text{ \ \ \ a.e., \ \ when \ }n\rightarrow
\infty ,\text{ \ \ \ \ when \ }0<\alpha <1.
\end{equation*}
\end{corollary}

{\bf Proof}:
By combining (\ref{node0})-(\ref{node01}) with (\ref{node1})-(\ref{node2})
we obtain that the N\"orlund means $\sigma _{n}^{\alpha }$ and $\beta
_{n}^{\alpha }$, with non-increasing sequence $\{q_{k}:k\in
\mathbb{N}
\}$ satisfy conditions (\ref{6a}) and (\ref{7a}).

Hence, the proof is complete by just using Theorem \ref{theoremae2}.
\QED

\begin{corollary}
\label{corollaryae22a}Let $f\in L_{1}$ and $t_{n}$ be N\"orlund means with
monotone and bounded sequence $\left\{ q_{k}:k\in \mathbb{N}\right\} $. Then%
\begin{equation*}
\text{\ \ \ }t_{n}f\rightarrow f \text{ \ \ \ a.e., \ \ when \ }n\rightarrow
\infty .
\end{equation*}
\end{corollary}

{\bf Proof}:
The proof is analogously to the proofs of Theorems \ref{theoremae1} and \ref%
{theoremae2} so we leave out the details.
\QED

\begin{theorem}
\label{theoremae3}Let $f\in L_{1}$. Then%
\begin{equation*}
\text{\ \ \ }R_{n}f\rightarrow f,\text{ \ \ \ a.e. \ \ \ \ when \ }%
n\rightarrow \infty .
\end{equation*}
\end{theorem}

{\bf Proof}:
The proof is analogously to the proofs of Theorems \ref{theoremae1} and \ref%
{theoremae2} so we leave out the details.
\QED

\end{document}